%% file: Main.tex
\documentclass[12pt, a4paper]{article}

\input{Framework.tex}

\input{Macros.tex}

\begin{document}

\pdfoutput=1

\begin{center}
    \begingroup
    \setstretch{1.05}
    \fontsize{16.2pt}{18.5pt}\selectfont
    \libertinusDisplay
    \fontdimen2\font=0.3em
    Burklund--Lin--Wang--Xu Methods in the 
    Cofiber-of-Tau Formalism \\ and Applications to 
    Equivariant Slice Differentials\par
    \endgroup
\end{center}

\vspace{0.3em}

\begin{center}
    {\fontsize{14.4pt}{18.5pt}\selectfont Yuchen Wu}
\end{center}

\vspace{1.5em}

\begin{abstract}
    We reinvestigate the theory of spectral sequences by studying the $(\infty,1)$-category of filtered spectra $\Fil\Sp$ through the cofiber-of-\(\tau\) formalism of Burklund--Isaksen--Pstragowski--Wang--Xu. 
    In this framework, we define and analyze hidden extensions along arbitrary maps of filtered spectra, 
    establishing computational principles that extend the generalized Leibniz rule and the generalized Mahowald trick of Lin--Wang--Xu \cite{Lin-Wang-Xu-kervaire}, 
    as well as Burklund's Leibniz rule for total differentials \cite{Burklund-cookware}, 
    from the Adams spectral sequence to this broader setup. 
    Our formulation uses a more refined, layered notion of extension, 
    which slightly sharpens these statements even for the Adams spectral sequence. 
    As an application, we study equivariant slice spectral sequences 
    and obtain new families of ``exotic transfer'' differentials in the \(C_4\)-slice spectral sequences for the Hill--Hopkins--Ravenel theories \(\BPCfour\langle m\rangle\) for every \(m \ge 1\).
\end{abstract}

\vskip 2em

\setcounter{tocdepth}{2}
\tableofcontents

\newpage 

\section{Introduction}

\input{Intro.tex}

\input{Conventions.tex}

\subsection{Acknowledgment}
\input{Acknowledgment.tex}

\newpage

\section{Filtered spectra and BIPWX's cofiber-of-tau formalism}

\input{Basics.tex}

\newpage

\section{Hidden extensions and the generalized Leibniz rule}

\input{GLR.tex}

\newpage

\section{Generalized Mahowald trick}

\input{GMT.tex}

\newpage

\section{Applications in the equivariant slice SS}

\input{Slice.tex}

\appendix

\newpage

\section{Appendix · Rigidity of the synthetic Adams SS} \label{app:A}

\input{SynASS.tex}

\newpage

\section{Appendix · Proof of the coherent Mahowald trick} \label{app:B}

\input{cMT-Proof.tex}

\newpage

\section{Appendix · Lax monoidality of the standard SS functors} \label{app:C}

\input{Picard.tex}

\newpage 

\printbibliography

\end{document}

%% file: Framework.tex
\usepackage{amsthm,amsmath,amssymb}
\usepackage{relsize,euscript}
\usepackage{enumitem}
\usepackage{multicol}

\usepackage{graphicx}
\usepackage{tikz, caption, float}
\usetikzlibrary{arrows.meta, cd, calc}

\usepackage[utf8]{inputenc}
\usepackage[T1]{fontenc}
\usepackage[expansion=false,nopatch=footnote]{microtype}
\usepackage[lining,proportional,scale=.875]{FiraSans}
\usepackage[tt=false,semibold]{libertinus}
\usepackage[libertine,smallerops]{newtxmath}
\usepackage[lining,scaled=.8]{FiraMono}
\usepackage[cal=boondox,frak=euler]{mathalpha}

\newcommand*{\sbf}{\firamedium}

\DeclareFontEncoding{MDA}{}{}
\DeclareSymbolFont{mathdesignA}{MDA}{mdput}{m}{n}
\SetSymbolFont{mathdesignA}{bold}{MDA}{mdput}{b}{n}
\DeclareSymbolFontAlphabet{\mathbb}{mathdesignA}

\tikzcdset{
row sep/normal = 5ex, column sep/normal = 5ex,
arrow style=tikz,
arrows={/tikz/line width=.7pt},
diagrams={>={Straight Barb[scale=0.8]}}
}

\DeclareFontFamily{OMX}{MnSymbolE}{}
\DeclareSymbolFont{MnLargeSymbols}{OMX}{MnSymbolE}{m}{n}
\SetSymbolFont{MnLargeSymbols}{bold}{OMX}{MnSymbolE}{b}{n}
\DeclareFontShape{OMX}{MnSymbolE}{m}{n}{
    <-6>  MnSymbolE5
   <6-7>  MnSymbolE6
   <7-8>  MnSymbolE7
   <8-9>  MnSymbolE8
   <9-10> MnSymbolE9
  <10-12> MnSymbolE10
  <12->   MnSymbolE12
}{}
\DeclareFontShape{OMX}{MnSymbolE}{b}{n}{
    <-6>  MnSymbolE-Bold5
   <6-7>  MnSymbolE6
   <7-8>  MnSymbolE7
   <8-9>  MnSymbolE8
   <9-10> MnSymbolE10
  <10-12> MnSymbolE10
  <12->   MnSymbolE12
}{}
\DeclareFontFamily{U}{MnSymbolC}{}
\DeclareFontShape{U}{MnSymbolC}{m}{n}{
    <-6>  MnSymbolC5
   <6-7>  MnSymbolC6
   <7-8>  MnSymbolC7
   <8-9>  MnSymbolC8
   <9-10> MnSymbolC9
  <10-12> MnSymbolC10
  <12->   MnSymbolC12
}{}
\DeclareFontShape{U}{MnSymbolC}{b}{n}{
    <-6>  MnSymbolC-Bold5
   <6-7>  MnSymbolC-Bold6
   <7-8>  MnSymbolC-Bold7
   <8-9>  MnSymbolC-Bold8
   <9-10> MnSymbolC-Bold9
  <10-12> MnSymbolC-Bold10
  <12->   MnSymbolC-Bold12
}{}

\DeclareSymbolFont{YWMnSymbolC}{U}{MnSymbolC}{m}{n}
\SetSymbolFont{YWMnSymbolC}{bold}{U}{MnSymbolC}{b}{n}

\let\filledstar\relax
\DeclareMathSymbol{\filledstar}{\mathord}{YWMnSymbolC}{"81}

\let\oldrightsquigarrow\rightsquigarrow

\makeatletter
\renewcommand{\rightsquigarrow}{%
  \mathrel{\mathpalette\rightsquigarrow@{}}%
}
\newcommand{\rightsquigarrow@}[2]{%
  \vcenter{\hbox{\scalebox{0.82}{$\m@th#1\oldrightsquigarrow$}}}%
}
\makeatother

\renewcommand{\geq}{\geqslant}
\renewcommand{\leq}{\leqslant}

\usepackage{accents}

\renewcommand{\overrightarrow}[1]{%
  \mathchoice
    {\mathord{\accentset{\smash{\raisebox{-.25ex}{$\scriptstyle\rightharpoonup$}}}{#1}}}%
    {\mathord{\accentset{\smash{\raisebox{-.25ex}{$\scriptstyle\rightharpoonup$}}}{#1}}}%
    {\mathord{\accentset{\smash{\raisebox{-.18ex}{$\scriptscriptstyle\rightharpoonup$}}}{#1}}}%
    {\mathord{\accentset{\smash{\raisebox{-.15ex}{$\scriptscriptstyle\rightharpoonup$}}}{#1}}}%
}

\usepackage{titling}
\pretitle {\vspace{-48pt}\setstretch{1.05}\begin{center}\LARGE\libertinusDisplay\fontdimen2\font=0.3em} % fontdimen2 is interword space
\posttitle{\end{center}\vspace{0pt}}

\usepackage[dvipsnames]{xcolor}
\usepackage{adjustbox}
\usepackage{multirow}
\usepackage{array}
\usepackage{placeins}
\usepackage{diagbox}
\usepackage[explicit]{titlesec}
\titleformat{\section}[hang]{\large\normalfont\sbf\filcenter}{\thesection}{1em}{\MakeUppercase{#1}}
\titleformat{\subsection}[hang]{\large\normalfont\sbf\filcenter}{\thesubsection}{1em}{%\MakeUppercase
{#1}}
\AtBeginDocument{
    \setlength{\belowdisplayshortskip}{\belowdisplayskip}
}

\definecolor{cite}{rgb}{0.45,0.55,0.6}
\definecolor{link}{rgb}{0.45,0.55,0.6}
\definecolor{url}{rgb}{0.4,0.4,0.4}
\usepackage[%
  hyperfootnotes=false,
  colorlinks,%
  linkcolor=link,%
  citecolor=cite,%
  urlcolor=url,%
]{hyperref}
\def\thesection{\arabic{section}\texorpdfstring{}{.}} % pdf bookmark numbering
\def\thesubsection{\arabic{section}.\arabic{subsection}\texorpdfstring{}{.}} % pdf bookmark numbering
\usepackage[style=alphabetic,sorting=nyt,useprefix=true,maxcitenames=7,
  mincitenames=7,
  maxbibnames=7,
  minbibnames=7]{biblatex}

\AtBeginBibliography{%
  \setlength{\biblabelsep}{0.4em}%
}

\usepackage{xpatch}

\DeclareFieldFormat[article]{volume}{\mkbibbold{#1}}
\DeclareFieldFormat[book,inbook]{number}{\mkbibbold{#1}}
\DeclareFieldFormat[article]{number}{(#1)}
\DeclareFieldFormat*{year}{(#1)}
\DeclareFieldFormat{pages}{#1}
\renewbibmacro{in:}{}
\renewbibmacro*{volume+number+eid}{%
    \printfield{volume}%
    \setunit*{\addnbspace}% originally: \setunit*{\adddot}
    \printfield{number}%
    \setunit{\addcomma\space}%
    \printfield{eid}%
}
\xapptobibmacro{author/editor+others/translator+others}{%
    \setunit{\space}%
    \printfield{year}%
    \clearfield{year}%
}{}{}
\xapptobibmacro{author/translator+others}{%
    \setunit{\space}%
    \printfield{year}%
    \clearfield{year}%
}{}{}
    \renewbibmacro*{issue+date}{%
    \printfield{issue}%
    \newunit%
}
\AtBeginBibliography{
    \DeclareFieldFormat{labelalpha}{#1}
    \DeclareFieldFormat{extraalpha}{\mknumalph{#1}}
}
\AtEveryBibitem{
    \clearfield{eprintclass}% suppresses e.g. [math.AT] after arXiv number
    \ifentrytype{online}{%
        \clearfield{year}%
    }{}
}

\usepackage{setspace}
\renewenvironment{abstract}{%
    \centering
    \begin{minipage}{.85\textwidth}%
    \setlength{\parindent}{1.5em}%
    \centerline{\fontsize{13pt}{15pt}\selectfont\sbf {Abstract}}%
    \par\vspace{6pt}%
}{%
    \end{minipage}%
    \par\vspace{12pt}%
}

\newcommand{\parr}{\vspace{\topsep}}

\theoremstyle{definition}

\newtheorem{definition}[equation]{Definition}

\newtheorem{recollection}[equation]{Recollection}
\newtheorem{summary}[equation]{Summary}
\newtheorem{construction}[equation]{Construction}
\newtheorem{setup}[equation]{Setup}
\newtheorem{example}[equation]{Example}
\newtheorem{proposition-definition}[equation]{Proposition-Definition}

\newtheorem{convention}{Convention}

\theoremstyle{plain}

\newtheorem{theorem}[equation]{Theorem}

\newtheorem{proposition}[equation]{Proposition}
\newtheorem{corollary}[equation]{Corollary}
\newtheorem{lemma}[equation]{Lemma}
\newtheorem{fact}[equation]{Fact}

\newtheorem{introthm}{Theorem}

\theoremstyle{remark}

\newtheorem{remark}[equation]{Remark}

\newenvironment{solution}{\begin{proof}}{\end{proof}}

\newcommand{\pf}{\begin{proof}}
\newcommand{\epf}{\end{proof}}
\newcommand{\sol}{\begin{solution}}
\newcommand{\esol}{\end{solution}}

\numberwithin{equation}{section}

\DeclareMathOperator*{\colim}{\mathrm{colim}}
\DeclareMathOperator*{\coker}{\mathrm{coker}}
\DeclareMathOperator*{\cofib}{\mathrm{cof}}
\DeclareMathOperator*{\fib}{\mathrm{fib}}

%% file: Macros.tex
\usepackage{mathdots}

\newcommand{\AAb}{\mathbb{A}}
\newcommand{\Cb}{\mathbb{C}}

\newcommand{\Eb}{\mathbb{E}}
\newcommand{\Fb}{\mathbb{F}}
\newcommand{\Gb}{\mathbb{G}}

\newcommand{\Nb}{\mathbb{N}}

\newcommand{\Qb}{\mathbb{Q}}
\newcommand{\Rb}{\mathbb{R}}
\newcommand{\Sb}{\mathbb{S}}

\newcommand{\Zb}{\mathbb{Z}}

\newcommand{\CA}{\mathcal{A}}
\newcommand{\CB}{\mathcal{B}}
\newcommand{\CC}{\mathcal{C}}
\newcommand{\CD}{\mathcal{D}}
\newcommand{\CE}{\mathcal{E}}
\newcommand{\CF}{\mathcal{F}}
\newcommand{\CI}{\mathcal{I}}

\newcommand{\CK}{\mathcal{K}}
\newcommand{\CL}{\mathcal{L}}
\newcommand{\CM}{\mathcal{M}}
\newcommand{\CN}{\mathcal{N}}
\newcommand{\CO}{\mathcal{O}}
\newcommand{\CP}{\mathcal{P}}
\newcommand{\CQ}{\mathcal{Q}}

\newcommand{\CT}{\mathcal{T}}

\newcommand{\CV}{\mathcal{V}}
\newcommand{\CW}{\mathcal{W}}

\newcommand{\cell}{\mathrm{cell}}

\newcommand{\id}{\mathrm{id}}

\newcommand{\infl}{\mathrm{Infl}}
\newcommand{\op}{\mathrm{op}}
\newcommand{\cp}{\mathrm{cp}}
\newcommand{\oneb}{\mathbb{1}}

\newcommand{\ev}{\mathrm{ev}}
\newcommand{\Ev}{\mathrm{Ev}}

\newcommand{\Aut}{\mathrm{Aut}}
\newcommand{\Map}{\mathrm{Map}}
\newcommand{\Msp}{\mathrm{map}}
\newcommand{\Res}{\mathrm{Res}}

\newcommand{\Free}{\mathrm{Free}}

\newcommand{\Indrm}{\mathrm{Ind}}
\newcommand{\Day}{\mathrm{Day}}
\newcommand{\Kos}{\mathrm{Kos}}

\newcommand{\N}{\mathrm{N}}

\newcommand{\res}{\mathrm{res}}
\newcommand{\tr}{\mathrm{tr}}
\newcommand{\tbar}{\bar{t}}
\newcommand{\vbar}{\bar{v}}
\newcommand{\dfbar}{\bar{\mathfrak{d}}}

\newcommand{\st}{\mathrm{st}}
\newcommand{\sgn}{\mathrm{sgn}}

\newcommand{\triv}{\mathrm{triv}}

\newcommand{\defopara}{\tau}
\newcommand{\quarterrep}{\lambda}

\newcommand{\Ab}{\mathsf{Ab}}
\newcommand{\Ac}{\mathrm{Ac}}

\newcommand{\Alg}{\mathsf{Alg}}

\newcommand{\Aux}{\mathsf{Aux}}

\newcommand{\BP}{\mathrm{BP}}

\newcommand{\MO}{\mathrm{MO}}
\newcommand{\MU}{\mathrm{MU}}

\newcommand{\CAlg}{\mathsf{CAlg}}
\newcommand{\Cat}{\mathsf{Cat}}

\newcommand{\Ch}{\mathsf{Ch}}

\newcommand{\D}{\mathsf{D}}
\newcommand{\Dec}{\textrm{D\'ec}}

\newcommand{\DT}{\mathsf{DT}}

\newcommand{\Fil}{\mathsf{Fil}}
\newcommand{\biFil}{\mathsf{biFil}}

\newcommand{\Gr}{\mathsf{Gr}}
\newcommand{\gr}{\mathrm{gr}}

\newcommand{\HH}{\mathrm{H}}

\newcommand{\lax}{\mathrm{lax}}
\newcommand{\oplax}{\mathrm{oplax}}

\newcommand{\Fun}{\mathsf{Fun}}

\newcommand{\Ind}{\mathsf{Ind}}

\newcommand{\Mack}{\mathsf{Mack}}
\newcommand{\Mod}{\mathsf{Mod}}

\newcommand{\Pic}{\mathsf{Pic}}
\newcommand{\pic}{\mathrm{pic}}

\newcommand{\PPr}{\mathsf{Pr}}

\newcommand{\QCoh}{\mathsf{QCoh}}

\newcommand{\Sl}{\mathrm{Sl}}
\newcommand{\Fin}{\mathsf{Fin}}

\newcommand{\RO}{\mathrm{RO}}

\newcommand{\Set}{\mathsf{Set}}

\newcommand{\Sh}{\mathsf{Sh}}
\newcommand{\SH}{\mathsf{SH}}
\newcommand{\Sp}{\mathsf{Sp}}
\newcommand{\Span}{\mathsf{Span}}
\newcommand{\SpSeq}{\mathsf{SS}}
\newcommand{\Syn}{\mathsf{Syn}}

\newcommand{\std}{\mathrm{std}}

\newcommand{\Pe}{\mathsf{P}_e}
\newcommand{\PG}{\mathsf{P}_G}

\newcommand{\Deltaone}{[1]}
\newcommand{\Deltatwo}{[2]}

\DeclareMathOperator{\IIm}{\mathrm{Im}}

\theoremstyle{plain}

\usepackage{enumitem}

\usepackage{cjhebrew}

\usepackage{relsize}

\newcommand{\HHRsup}[1]{^{(\mkern-4.22mu(#1)\mkern-3.5mu)}}

\newcommand{\BPCfour}{\BP\HHRsup{C_4}}
\newcommand{\BPCeight}{\BP\HHRsup{C_8}}
\newcommand{\BPCtwon}{\BP\HHRsup{C_{2^n}}}
\newcommand{\BPG}{\BP\HHRsup{G}}

\newcommand{\MUCfour}{\MU\HHRsup{C_4}}

\newcommand{\MUG}{\MU\HHRsup{G}}

%% file: Intro.tex
\label{sec:1}

\subsection{Background}

Spectral sequences (SS) are ubiquitous computational devices across mathematics, and they play a particularly central role in algebraic topology. A striking illustration is the Kervaire invariant problem: Geometrically, this is about whether there exists a smooth closed framed manifold of dimension $2^{j + 1} - 2$  that cannot be converted to a homotopy sphere via framed surgery, while the works of Kervaire--Milnor \cite{KM63} and Browder \cite{Bro69} translate this into deciding whether the Kervaire class $h_j^2$ in the $\mathrm{H}\Fb_2$-Adams spectral sequence of the sphere spectrum is a permanent cycle. For $j \leq 5$, these classes are permanent cycles due to the computation of Mahowald--Tangora \cite{MT67} and Barratt--Jones--Mahowald \cite{BJM84}, while the remaining cases were open for decades, until two major breakthroughs in the next century. \parr 

In \cite{HHR-paper}, Hill, Hopkins and Ravenel solve the Kervaire invariant problem for $j \geq 7$:

\begin{theorem}[Hill--Hopkins--Ravenel]
    For $j \geq 7$, $h_j^2$ is not a permanent cycle. Thus, every closed framed manifold of dimension $2^{j + 1} - 2$ can be converted to a homotopy sphere via framed surgery.
\end{theorem}

The resolution of Hill--Hopkins--Ravenel is based on equivariant homotopy theory. They carry out an extensive study of the \textbf{equivariant slice spectral sequences}, and they show that the Kervaire classes are detected in the slice spectral sequence of a certain genuine equivariant spectrum, namely the Hill--Hopkins--Ravenel theory $\BPG\langle m \rangle$ with $G = C_8, m = 1$. The slice spectral sequence analysis further shows that a suitable localization of this spectrum has a gap and a $256$-periodicity in its equivariant homotopy groups, from which they conclude that $h_j^2$ does not survive for each $j \geq 7$. \parr 

The last open case, i.e. the case $j = 6$, is solved by Lin, Wang and Xu in \cite{Lin-Wang-Xu-kervaire}:

\begin{theorem}[Lin--Wang--Xu]
    For $j = 6$, $h_6^2$ is a permanent cycle. Thus, there exist closed framed manifolds of dimension $126$ that cannot be converted to a homotopy sphere via framed surgery. 
\end{theorem}

The resolution of Lin--Wang--Xu stays in the realm of Adams SS. Their key technique lies in the systematic study of \textbf{(hidden) extensions} in the Adams SS, which record information about a map $f\colon X \to Y$ between finite spectra that is not detected by the induced map $f\colon E_{2,\mathrm{ASS}}^{*,*}(X) \to E_{2, \mathrm{ASS}}^{*,*}(Y)$. More precisely, for such a map $f\colon X \to Y$, they established 
\begin{itemize}
    \item The \textbf{Generalized Mahowald trick (GMT)}, which yields a correspondence between extensions along the map $f$ and the Adams differentials of $\cofib(f)$.  
    \item The \textbf{Generalized Leibniz rule (GLR)}, which combines short Adams differentials for $X$ and extensions along $f$ to produce long Adams differentials for $Y$.
\end{itemize}
Given as input the Adams $E_2$-pages and $d_2$ differentials of many finite spectra (computed by Lin's program using 
secondary Steenrod operations, cf. \cite{Lin-Wang-Xu-machine}), these two techniques generate lots of longer differentials, which clarify the pattern of Adams differentials of $\Sb^0$ around stem $126$. 

\subsection{Summary of the results}

This paper extends the study of hidden extensions from the Adams SS to the standard SS associated to arbitrary filtered spectra, and applies this framework to equivariant slice SS computations. \parr 

In \cite{Lin-Wang-Xu-kervaire}, hidden extensions in the Adams SS are characterized using $\mathrm{H}\Fb_2$-synthetic spectra. The underlying strategy goes back to the cofiber-of-$\tau$ formalism in the stable $\Cb$-motivic homotopy category by Gheorghe, Isaksen, Wang and Xu in \cite{Isaksen-stable-stems} and \cite{GWX}, while it has been further developed in a series of works by, among others, Burklund, Isaksen, Pstragowski, Wang and Xu in \cite{Pst23}, \cite{GIKR}, \cite{BHS1}, \cite{Burklund-Xu}, \cite{BIX} and \cite{Lin-Wang-Xu-kervaire}. Throughout this literature, the common formal input is a stable presentably (symmetric) monoidal $\infty$-category $\CC$ equipped with a (symmetric) monoidal functor from the poset $(\Zb, \leq)$ to $\CC$. In this paper, we work in the universal such setting, namely the $\infty$-category of filtered spectra $\Fil\Sp = \Fun((\Zb, \leq)^{\op}, \Sp)$ together with the spectral Yoneda embedding $\gamma\colon (\Zb, \leq) \to \Fil\Sp$. \parr 

Many key constructions in synthetic spectra admit a universal, and often simpler, formulation in $\Fil\Sp$. There is a bigraded family of invertible objects
\(\Sb^{n,w}=\Sigma^{n}\gamma(w),\)
which induces bigraded suspensions and bigraded homotopy groups. 
The $\gamma$-image of the map $(-1)\to 0$ defines a ``deformation parameter''\footnote{We choose to write $\defopara$ because of its origination from the map  $\tau$ in the stable motivic homotopy category over $\Cb$.} \(\defopara\colon \Sb^{0,-1}\to \Sb^{0,0},\) and hence natural transformations $\defopara^{k}\colon \Sigma^{0,-k}\to \id$ for $k\ge 1$. 
Moreover, from the composable pair 
\(\gamma((-r)\to (-k)\to 0) = \Sb^{0, -r}\to \Sb^{0,-k}\to \Sb^{0,0}\)
one extracts a fundamental distinguished triangle
\[
\Sigma^{0,-k}X/\defopara^{\,r-k}\xrightarrow{\defopara^{k}}
X/\defopara^{\,r}\xrightarrow{\rho_{k}^{r}}
X/\defopara^{\,k}\xrightarrow{\delta_{k}^{\,r-k}}
\Sigma^{1,-k}X/\defopara^{\,r-k},
\]
which defines the ``projection'' $\rho=\rho_{k}^{r}$ and the ``total differential'' $\delta_k=\delta_{k}^{\,r-k}$.
On the other hand, each $X\in\Fil\Sp$ has an associated ``standard spectral sequence'' $\{E_r^{*,*}(X)\}_{r\ge 2}$, and this construction encompasses most spectral sequences encountered in practice; see Example~\ref{SS-examples-1}, Example~\ref{SS-examples-2}, and Definition~\ref{slice-filtrations}. 
The \textbf{cofiber-of-$\tau$ formalism} in the sense of Burklund--Isaksen--Pstragowski--Wang--Xu (BIPWX) 
is the study of the standard SS through the bigraded homotopy groups of $X/\defopara^{\,r}$ for $1\le r\le\infty$, together with the operators $\defopara^{k}$, $\rho$, and $\delta_k$. In Theorem \ref{standard-SS} we construct the standard SS from these data. Conversely, the bigraded homotopy groups and these operators can be read off from the standard SS due to

\begin{introthm}[Bockstein dictionary; Precise formulation in Theorems \ref{Bockstein-dictionary-finite}, \ref{Bockstein-dictionary-infinite}, \ref{delta-as-total-diff}; See also Summary \ref{Bockstein-dictionary-summary}] \label{intro-thm-A}
    Suppose $1 \leq r \leq \infty, 1 \leq k < \infty$ and $X \in \Fil\Sp$ (subject to some mild assumptions if $r = \infty$). Then 
    \begin{enumerate}
        \item We have $E_2^{s, t}(X) \cong \pi_{t - s, t}(X / \defopara)$. 
        \item Each $d_r$-cycle $x \in Z_r^{s,t} \subset E_2^{s, t}$ admits a lift $[x] \in \pi_{t - s, t}(X / \defopara^r)$ along the projection $\rho\colon X / \defopara^r \to X / \defopara$. Conversely, each bigraded stem $\alpha \in \pi_{t - s, t}(X / \defopara^r)$ is of the form $\defopara^a[y]$ for some $y \in Z_{r - a}^{s + a, t + a}$.
        \item The operator $\defopara^k\colon X / \defopara^r \to X / \defopara^{r + k}$ kills the lifts of $(k + 1)$-boundaries $B_{k + 1}^{s, t} \subset E_2^{s, t}$ (up to some $\defopara$-divisible stems), while the projection $\rho = \rho^r_k\colon X / \defopara^r \to X / \defopara^k$ removes the $\defopara^{k}$-divisible stems. 
        \item The total differential $\delta_k  = \delta_k^r\colon X / \defopara^k \to \Sigma^{1, -k} X / \defopara^r$ repackages the data of $d_2, \ldots, d_{k + r}$, in the sense that $\delta_k^r \defopara^a [x] = \defopara^b [y]$ if and only if $d_{k - a + b + 1}(x) = y$ in the standard SS.
    \end{enumerate}
\end{introthm}

\begin{remark}
    The Bockstein dictionary results are closely related to \emph{rigidity} results for trigraded SS. For a more detailed account of the relevant history, including work of Hu--Kriz--Ormsby, Isaksen, Burklund--Hahn--Senger, Burklund--Xu, and Carrick--Davies--van Nigtevecht, see Remark \ref{credits-for-theorem-A}.
\end{remark}

The functoriality of these bigraded homotopy groups allows us to make sense of ``extensions of standard SS'' along any map $f\colon X \to Y$ in $\Fil\Sp$: there is a $j$-th layer extension along $f$ on the $E_{r + 1}$-page of filtration jump $k$ from $x \in E_2^{*,*}(X)$ to $y \in E_2^{*,*}(Y)$ if there exist lifts $[x] \in \pi_{**}(X / \defopara^r), [y] \in \pi_{**}(Y / \defopara^{r - k})$ so that $\defopara^j f[x] = \defopara^{k + j}[y]$. In terms of this, we have the following

\begin{introthm}[Generalized Leibniz rule; Precise formulation in Theorems \ref{GLR-for-SS} and \ref{GLR-for-SS-2}] \label{intro-thm-B}
    Suppose $f\colon X \to Y$ is a map in $\Fil\Sp$, $x, x_{\infty}\in E_2^{*,*}(X)$ and $y, y_{\infty} \in E_2^{*,*}(Y)$, such that 
    \begin{itemize}
        \item $d_{r + 1}(x) = x_{\infty}$ in the standard SS of $X$. 
        \item There is a finite page $f$-extension $x \rightsquigarrow y$ of filtration jump $m$. 
    \end{itemize}
    Then we have the following: 
    \begin{enumerate}
        \item If there is an $E_{\infty}$-page $f$-extension $x_{\infty} \rightsquigarrow y_{\infty}$ of filtration jump $l$, and certain no crossing conditions are satisfied, then $d_{r - m + l + 1}(y) = y_{\infty}$ in the standard SS of $Y$. 
        \item Conversely, if $d_{r - m + l + 1}(y) = y_{\infty}$ in the standard SS of $Y$, and some other no crossing conditions are satisfied, then there is an $E_{\infty}$-page $f$-extension $x_{\infty} \rightsquigarrow y_{\infty}$ of filtration jump $l$.
    \end{enumerate}
\end{introthm}

\begin{introthm}[Generalized Mahowald trick; Precise formulation in Theorems \ref{GMT-for-SS} and \ref{GMT-for-SS-2}] \label{intro-thm-C}
    Suppose $f\colon X \to Y$ is a map in $\Fil\Sp$ with cofiber $Z$, $x \in E_2^{*,*}(X)$, $y \in E_2^{*,*}(Y)$, $\overline{x}, \overline{y} \in E_2^{*,*}(Z)$, so that 
    \begin{itemize}
        \item There is a finite page extension along $g\colon Y \to Z$ from $y$ to $\overline{y}$ of filtration jump $m$. 
        \item There is a finite page extension along $h\colon Z \to \Sigma^{1, 0} X$ from $\overline{x}$ to $x$ of filtration jump $l$. 
    \end{itemize}
    Write $r = n + m + l$ for some $n \geq 0$. Then we have the following: 
    \begin{enumerate}
        \item If $d_{r + 1}(\overline{x}) = \overline{y}$ in the standard SS of $Z$, and certain no crossing conditions are satisfied, then there is a finite page $f$-extension $x \rightsquigarrow y$ of filtration jump $n$. 
        \item Conversely, if there is a finite page $f$-extension $x \rightsquigarrow y$ of filtration jump $n$, and some other no crossing conditions are satisfied, then $d_{r + 1}(\overline{x}) = \overline{y}$ in the standard SS of $Z$. 
    \end{enumerate}
\end{introthm}

Here ``no crossing conditions'' serve as necessary inputs to eliminate the potential discrepancy (nonuniqueness in the choice of lift) along the Bockstein dictionary (Theorem \ref{intro-thm-A}) translation. 

\begin{remark}
    These results generalize Lin--Wang--Xu's Adams SS versions \cite[Theorems 6.1 and 6.12]{Lin-Wang-Xu-kervaire}. For historical attributions concerning Theorem \ref{intro-thm-B}, including works of Chua, Hill--Hopkins--Ravenel, Burklund--Hahn--Senger, Burklund--Xu, Marek, Isaksen--Kong--Li--Ruan--Zou, and Carrick--Davies--van Nigtevecht, see Remark \ref{credits-for-theorem-B}. For the corresponding discussion of Theorem \ref{intro-thm-C}, including works of Mahowald, Hill--Hopkins--Ravenel, Behrens, and Ma, see Remark \ref{credits-for-theorem-C}.
\end{remark}

Even for the Adams SS, Theorems \ref{intro-thm-B} and \ref{intro-thm-C} are slightly sharper than the corresponding results of \cite{Lin-Wang-Xu-kervaire}, since they are formulated in terms of \emph{layered extensions}, which keep track of the relevant $\defopara$-power in a hidden extension. This refinement distinguishes some configurations that count as crossings in the sense of \cite{Lin-Wang-Xu-kervaire} but not in the layered framework; see the discussion on crossings at the end of Example \ref{example-extn-crossing}. Consequently, the same conclusions can hold under weaker hypotheses, as demonstrated in Example \ref{GLR-example}.

\begin{remark}
   The mechanism underlying this layered refinement already appears implicitly in concrete computations: multiplying by a suitable power of $\defopara$ can sometimes remove a crossing obstruction at the relevant bidegree. Instances of this idea occur in \cite{Burklund-hidden-extension}, \cite{BIX} and \cite[\S~7]{Lin-Wang-Xu-kervaire}. The role of layered extensions is to isolate this mechanism systematically.
\end{remark}

Before moving on to the remaining main results, we pause to highlight two technical inputs that underlie the theorems, which are developed in two appendices.

\begin{enumerate}
    \item[(A).] We compare our approach to hidden extensions with the approach in \cite{Lin-Wang-Xu-kervaire} in Remark \ref{comparison-with-extension-SS}. This relies on the identification of the synthetic Adams SS of $\nu X / \defopara^r$ ($1 \leq r \leq \infty$) with the corresponding trigraded $\defopara$-Bockstein SS that is compatible with the abutments, and we prove such a result in Theorem \ref{synthetic-Adams-vs-Bockstein}. In Appendix \ref{app:A} we also establish the ``rigidity'' of trigraded $\defopara$-Bockstein SS in general,
    and we apply this to reprove the fact that the synthetic Adams SS of $\nu X$ (and also $\nu X / \defopara^r$) determines, and is determined by the classical Adams SS of $X$. 
    \item[(B).] Theorem \ref{intro-thm-B} and Theorem \ref{intro-thm-C} both come from a careful Bockstein dictionary translation applied to the (more intuitive) blueprint results in terms of bigraded stems. The Blueprint \ref{blueprint-GLR} for Theorem \ref{intro-thm-B} follows from the naturality of the total differential operator $\delta$; on the other hand, the Blueprint \ref{blueprint-GMT} for Theorem \ref{intro-thm-C} is technically more subtle. 
    
    Indeed, to identify the bigraded homotopy groups lifting the extensions and/or differentials, one would need to construct a delicate homotopy between different internal isomorphisms in a $2$-dimensional coherent enhancement of distinguished triangles, which we refer to as the ``coherent Mahowald trick''. 
    
    The paper \cite{Lin-Wang-Xu-kervaire} resolves this challenge by restricting to the special situation 
    where the ``$2$d triangle'' comes from the smash product of two classical distinguished triangles, for which a model categorical argument is sketched in \cite[Proof of (TC3) in \S~6]{May-additivity}. 
    In this paper, we achieve further generality: we make sense of the concept ``$2$d triangle'' in an arbitrary stable $\infty$-category in \S~\hyperref[subsec:4.1]{4.1}, state the coherent Mahowald trick (cf. Theorem \ref{coherent-Mahowald-trick}) in this context, and provide a model-independent proof in Appendix~\ref{app:B}.
\end{enumerate}

Besides the generalized Leibniz rule and the generalized Mahowald trick, Burklund--Isaksen--Pstragowski--Wang--Xu's cofiber-of-$\tau$ formalism also allows us to generalize the classical Leibniz rule along another direction in $\Fil\Sp$, extending the \emph{Leibniz rule for total differentials} proposed by Burklund in \cite[Claim 3.3.3]{Burklund-cookware}:

\begin{introthm}[Leibniz rule for total differentials; Precise formulation in Theorem \ref{Burklund's-Leibniz-rule}] \label{intro-thm-D}
    Suppose $F\colon X \otimes Y \to T$ is a map in $\Fil\Sp$. For each $r \geq 1$,
    the induced pairing
    \[F = F / \defopara^r\colon \pi_{**}(X / \defopara^r) \times \pi_{**}(Y/\defopara^r) \to \pi_{**}(T/\defopara^r)\] 
    satisfies the Leibniz rule for total differentials $\delta_r = \delta_r^r$, in the sense that 
    \[\delta_r F(\alpha, \beta) = F(\delta_r(\alpha), \beta) + (-1)^{|\alpha|} F(\alpha, \delta_r(\beta)).\]
    In particular, the induced pairing on the $E_{r + 1}$-page satisfies the classical Leibniz rule.
\end{introthm}

\begin{remark} \label{credits-for-theorem-D}
    A special case of Theorem \ref{Burklund's-Leibniz-rule}, where $F$ is the multiplication for a ($\mathrm{H}\Fb_2$-)synthetic ring spectrum, is proved by Carrick--Davies--van Nigtevecht in \cite[Theorem 2.34]{CDvN2}.
\end{remark}

In Appendix \ref{app:C}, we recast the Leibniz rule for ordinary/total differentials as a categorical statement, expressed in terms of the lax $\Eb_n$-monoidality of certain functors. 

\begin{introthm}[Lax monoidality for the standard SS functor and its total-differential refinements; Precise formulation in Theorems \ref{categorified-total-Leibniz-rule-with-general-coefficient}, \ref{categorified-total-Leibniz-rule-with-picard-trigrading} and \ref{categorified-total-Leibniz-rule-with-picard-bigrading}] \label{intro-thm-E}
    Take $n \in \Nb_{\geq 1} \cup \{\infty\}$. If $\CE$ is a presentably $\Eb_{n}$-monoidal $t$-$\infty$-category, then the assignments
    \begin{align*}
        A_{r + 1}^{*,*}\colon \Fil(\CE) \to \Gr\Ch(\CE^\heartsuit)^{E_{r + 1}}, &\quad X \mapsto \{(A_{r + 1}^{s, t}(X) = \pi_{t - s, t}(X / \defopara^{r}), \delta_r^r)\}_{s, t \in \Zb} \\
        E_{r + 1}^{*,*}\colon \Fil(\CE) \to \Gr\Ch(\CE^\heartsuit)^{E_{r + 1}}, &\quad X \mapsto \{(E_{r + 1}^{s, t}(X), d_{r + 1})\}_{s, t \in \Zb}\\
        E_*^{*,*}\colon \Fil(\CE) \to \SpSeq(\CE^\heartsuit), &\quad X \mapsto \{(E_{r + 1}^{s, t}(X), d_{r + 1})\}_{r \geq 1, s, t \in \Zb}
    \end{align*}
    define lax $\Eb_{n}$-monoidal functors for each $r \in \Nb_{\geq 1}$. Here the target $\Gr\Ch(\CE^\heartsuit)^{E_{r + 1}}$ is the $1$-category of graded chain complexes $\{(M^{s, t}, d_{r + 1})\}$ valued in $\CE^\heartsuit$, 
    and $\SpSeq(\CE^\heartsuit)$ is the $1$-truncated $\infty$-operad over $\Eb_n$ of spectral sequences valued in $\CE^\heartsuit$. \parr

    Suppose moreover $\CE$ is $\Eb_{n + 1}$-monoidal and it is equipped with an $\Eb_{n + 1}$-monoidal retract $V \mapsto |V|$ of the canonical map $\Pic(\Sb) \to \Pic(\CE)$. Then the assignments 
    \begin{align*}
        A_{r + 1}^{*,\filledstar}\colon \Fil(\CE) \to \Gr_{\Pic(\CE)}^{E_{r + 1}}\Ch(\CE^\heartsuit), &\quad X \mapsto \{(A_{r + 1}^{s, V}(X) = \pi_{V - s, |V|}(X / \defopara^{r}), \delta_r^r)\}_{s \in \Zb, V \in \Pic(\CE)} \\
        E_{r + 1}^{*,\filledstar}\colon \Fil(\CE) \to \Gr_{\Pic(\CE)}^{E_{r + 1}}\Ch(\CE^\heartsuit), &\quad X \mapsto \{(E_{r + 1}^{s, V}(X) = E_{r + 1}^{s, 0}(\hom_{\Fil(\CE)}(\Sb^{V, |V|}, X)), d_{r + 1})\}_{s \in \Zb, V \in \Pic(\CE)} \\
        E_*^{*,\filledstar}\colon \Fil(\CE) \to \SpSeq_{\Pic(\CE)}(\CE^\heartsuit), &\quad X \mapsto \{(E_{r + 1}^{s, V}(X), d_{r + 1})\}_{r \geq 1, s \in \Zb, V \in \Pic(\CE)}
    \end{align*}
    define lax $\Eb_n$-monoidal functors for each $r \in \Nb_{\geq 1}$. Here the target $ \Gr_{\Pic(\CE)}^{E_{r + 1}}\Ch(\CE^\heartsuit)$ is the $1$-category of $\Pic(\CE)$-graded chain complexes $\{(M^{s, V}, d_{r + 1})\}$ valued in $\CE^\heartsuit$ for which the sign $C_2$-action on $V \in \Pic(\CE)$ induces the sign $C_2$-action on each $M^{s, V} \in \CE^\heartsuit$, and $\SpSeq_{\Pic(\CE)}(\CE^\heartsuit)$ is the $1$-truncated ${\infty}$-operad over $\Eb_n$ of $\Pic(\CE)$-graded spectral sequences valued in $\CE^\heartsuit$. 
\end{introthm}

\begin{remark}
    The lax symmetric monoidality of $E_*^{*,*}\colon \Fil\Sp \to \SpSeq$ was first proved by Hedenlund \cite[Theorem II.3.5]{Hed20}, and subsequently revisited by Antieau \cite[Theorem 8.9]{Ant24}. Their proofs are based on the streamlined décalage functor $\Dec\colon \Fil\Sp \to \Fil\Sp$, while our proof, based on total differentials, is different in nature. The décalage construction is well suited for extracting the standard SS itself; however, since the $r$-fold décalage of a filtered spectrum $X$ does not retain the earliest $r$ pages of the standard SS, hidden extensions and total differentials are not visible from décalage alone.
\end{remark}

From a different perspective, in \S~\hyperref[subsec:2.4]{2.4} we demonstrate that our results in $\Fil\Sp$ apply equally well to 
the study of an arbitrary presentable $\Fil\Sp$-module, namely any presentable stable $\infty$-category $\CC$ equipped with a cocontinuous monad $T \in \Alg(\Fun^L(\CC, \CC))$ such that $\fib(\eta\colon \id \to T)$ is an invertible functor (cf. Theorem \ref{GHMG} and Corollary \ref{GHMG-for-modules}). \parr

Having all these computational methods at hand, we test them in $\CC = \Fil(\Sp^G)$ with equivariant slice spectral sequences for the Hill--Hopkins--Ravenel theories  $\BPG \langle m \rangle$, which are genuine $G$-spectra (for any $2$-group $G$ equipped with a fixed inclusion $C_2 \subset G$) introduced in \cite{HHR-paper}. Such computations have received significant attention for two reasons:

\begin{itemize}
    \item Hill, Hopkins and Ravenel resolve the Kervaire invariant problem in the negative for $j \geq 7$ through a partial computation for the $C_8$-slice SS of $\BPCeight \langle 1 \rangle$.
    \item For each $m, n \geq 1$, $\BPCtwon \langle m \rangle$ is a decompleted model for the Lubin--Tate theory $E_h$ (here $h = 2^{n - 1} \cdot m$) with $C_{2^n}$-action, and the $C_{2^n}$-slice SS of $\BPCtwon \langle m \rangle$ contains the information about all differentials in the $C_{2^n}$-homotopy fixed point SS of $E_h$, cf. \cite{BHSZ}. 
\end{itemize}

The first equivariant slice SS computation for $\BPG\langle m\rangle$ appears in Dugger's thesis \cite{Dug05}, which corresponds to the case $G = C_2$ (so $\BPG=\BP_{\Rb}$) and $m = 1$; subsequently, Hu and Kriz \cite{Hu-Kriz} compute the $C_2$-slice SS of $\BP_{\Rb}\langle m\rangle$ for all $1\le m\le\infty$, Hill, Hopkins and Ravenel \cite{HHR-BPC4<1>} compute the $C_4$-slice SS for $\BPCfour\langle 1\rangle$, and Hill, Shi, Wang and Xu \cite{HSWX-BPC4<2>} compute the $C_4$-slice SS for $\BPCfour\langle 2\rangle$. Beyond these cases, we do not currently have a complete computation of the $G$-slice SS for $\BPG\langle m\rangle$ for any other $G$ or $m$. Nevertheless, a substantial body of structural results on the equivariant slice filtration and HHR theories---height stratifications, transchromatic shearing isomorphisms, vanishing lines, periodicities, and more---has been developed by Hill, Meier, Shi, Li, Duan, Zeng, Wang, Xu, et al. in works such as \cite{MSZ-localized-slice-SS, MSZ-stratification, MSZ-transchromatic, DLS-vanishing-lines, DHLLSWX-periodicity}, among others. BIPWX's cofiber-of-$\tau$ formalism provides a natural framework for organizing these results and extracting further computational information from them.\parr 

In this paper we focus on the computation for $G = C_4$. There are three known families of differentials in the $C_4$-slice SS of $\BPCfour\langle m \rangle$: 
\begin{itemize}
    \item the ``{norm differentials}'' $d_{5}, d_{13}, d_{29}, d_{61},\ldots, d_{2^{2m + 2} - 3}$ constructed from applying the HHR norm functor to the (completely known) $C_2$-slice SS, as in \cite[Theorem 4.7]{HHR-BPC4<1>}. 
    \item the ``{transfer differentials}'' $d_{3}, d_{7}, d_{15}, d_{31},\ldots, d_{2^{2m + 1} - 1}$ constructed from applying the transfer functor to the $C_2$-slice SS, as in \cite[Theorem 2.3]{DLS-vanishing-lines}. 
    \item the ``{sheared differentials}'' $d_{5}, d_{13}, \ldots, d_{2^{m + 2} - 3}$ constructed from the full knowledge of the $C_4$-slice SS above the line of slope $1$, as in \cite[Theorem 4.1]{MSZ-transchromatic}.
\end{itemize}
Apart from these families, the other differentials in existing computations (e.g. \cite{HSWX-BPC4<2>}) mostly come from ad-hoc case-by-case arguments, which admit no apparent generalization for $m \geq 3$. However, using Theorem \ref{intro-thm-B} and Theorem \ref{intro-thm-C}, we are able to produce new families of differentials in the $C_4$-slice SS of $\BPCfour\langle m \rangle$ for all $m$: 

\begin{introthm}[Exotic transfer differentials; Precise formulation in Theorem \ref{exotic-transfer-paradigm}] \label{intro-thm-F}
    For $1 \leq m < \infty$, consider the $C_4$-slice SS of $\BPCfour \langle m \rangle$. Suppose for $1 \leq h \leq 2m$ we have
    \begin{itemize}
        \item a homogeneous polynomial $P \in B_m =  \Zb[T_1, T_3, \ldots, T_{2^m - 1}]$ of degree $p$
        \item two homogeneous classes $x, c \in R_m = \Zb[\tbar_1, \gamma \tbar_1, \tbar_3, \gamma \tbar_3, \ldots, \tbar_{2^m - 1}, \gamma \tbar_{2^m - 1}]$ 
    \end{itemize}
    subject to the algebraic identity
    \[x \cdot \vbar_h = P(\tbar_1 \gamma \tbar_1, \tbar_3 \gamma \tbar_3, \ldots, \tbar_{2^m - 1} \gamma \tbar_{2^m - 1}) + c + \gamma(c) \mod I_h = (2, \vbar_1, \ldots, \vbar_{h - 1})\]
    for the Araki generators $\vbar_k \in R_m$ $(1 \leq k \leq 2m)$ in Theorem \ref{diffs-in-C2-sliceSS}. Then for each $q \in \{1, \ldots, m\}, j \in \Nb$, there is a differential in the $C_4$-slice SS of $\BPCfour \langle m \rangle$ 
    \[d_{2^{h + 1} + 2^{q + 2} - 5}(\tr(x a_{\sigma_2})u_{2^{h - 1}\quarterrep} u_{(2^{q + 1}\! j + 2^q - p - 1)\sigma}) = P(\dfbar_1, \dfbar_3, \ldots, \dfbar_{2^m - 1}) \dfbar_{2^q - 1} u_{2^{q + 1} \! j \sigma} a_{(2^h + 2^q - 1)\quarterrep}a_{(2^{q + 1} - 2)\sigma}\]
    constructed from ``splicing together'' the $h$-th transfer differential and the $q$-th sheared differential. 
\end{introthm}

With appropriately chosen inputs (cf. Examples \ref{exotic-transfer-diff-longest} and \ref{exotic-transfer-diff-second-longest}), this recovers 
\begin{itemize}
    \item all $d_{11}$ differentials in the $C_4$-slice SS of $\BPCfour \langle 1 \rangle$
    \item all $d_{43}$ differentials in the $C_4$-slice SS of $\BPCfour \langle 2 \rangle$ 
    \item all transchromatic\footnote{i.e. crossing the stratification line $\CL_1$ of slope $1$ in \cite[Figures 4 and 5]{MSZ-stratification}} $d_{19}$ differentials in the $C_4$-slice SS of $\BPCfour \langle 2 \rangle$ 
\end{itemize}
and extends them to families of differentials in the $C_4$-slice SS of $\BPCfour \langle m \rangle$ for $m \geq 3$. \parr 

Finally, in \S~\hyperref[subsec:5.6]{5.6}, we apply Theorems \ref{intro-thm-D} and \ref{intro-thm-F} to obtain a satisfactory interpretation for all differentials in the $C_4$-slice SS of $\BPCfour \langle 1 \rangle$: each differential in the SS comes from one of the seven generating differentials under the classical Leibniz rule, while each generating differential is established by a systematic (i.e. non-ad-hoc) construction. This is a promising first step towards a systematic description for all slice differentials of $\BPG \langle m \rangle$. 

\subsection{Organization and conventions}

In \S~\hyperref[sec:2]{2} we set up BIPWX's cofiber-of-$\tau$ formalism in the realm of filtered spectra. In particular, we prove the Bockstein dictionary (Theorem \ref{intro-thm-A}) and the Leibniz rule for total differentials (Theorem \ref{intro-thm-D}). In \S~\hyperref[sec:3]{3} we define hidden extensions and their no crossing conditions in terms of the deformation parameter $\defopara$, and prove the generalized Leibniz rule (Theorem \ref{intro-thm-B}). In \S~\hyperref[sec:4]{4} we prove the generalized Mahowald trick (Theorem \ref{intro-thm-C}) up to the delicate homotopy referred to as the ``coherent Mahowald trick''. In \S~\hyperref[sec:5]{5} we apply these results to equivariant slice SS to construct the exotic transfer differentials in Theorem \ref{intro-thm-F}. In Appendix \ref{app:A} we investigate the notion of trigraded $\defopara$-Bockstein SS and its rigidity (Theorem \ref{rigidity-for-trigraded-defopara-Bockstein-SS} and Remark \ref{rigidity-for-filtered-SS-by-definition}), and identify the synthetic Adams SS as a special case in Theorem \ref{synthetic-Adams-vs-Bockstein}. In Appendix \ref{app:B} we prove the coherent Mahowald trick stated in Theorem \ref{coherent-Mahowald-trick}. Lastly, in Appendix \ref{app:C} we prove the categorical refinements of Leibniz rules (Theorem \ref{intro-thm-E}). \parr 

We summarize our convention throughout this paper below:

%% file: Conventions.tex
\begin{convention}
    Conventions for categories in general.
\begin{itemize}
    \item Throughout the paper, we work with the theory of $(\infty,1)$-categories as developed in \cite{HTT, HA}. All definitions, constructions and arguments in this paper are supposed to be model-independent, and may equally be carried out in the axiomatic framework of \cite{CCNW26, Cno26}.
    \item We suppress most size considerations of categories throughout the paper; all such issues are standard and will be made explicit only when relevant. 
    \item We denote by $\overrightarrow{\Zb}$ or $(\Zb, \leq)$ the poset of integers, and we write $\Zb$ or $\Zb^\delta$ for its underlying set. Similar conventions apply to the set of nonnegative integers $\Nb$.
    \item We write categories, or $1$-categories, for $(1, 1)$-categories. 
    \item We treat posets as special instances of $1$-categories. For each $n \in \Nb$, we denote by $[n]$ or $\Delta^n$ the totally ordered set $\{0 \leq 1 \leq \cdots \leq n\}$. We write $\Delta$ for the $1$-category of all finite totally ordered sets. 
    \item We write $\infty$-categories for $(\infty, 1)$-categories, and we treat $1$-categories as special instances of $\infty$-categories. For any $\infty$-category $\CC$, we denote by $h\CC$ its homotopy $1$-category.  
    \item If $\CK, \CC$ are $\infty$-categories, we write $\CC^{\CK}$ or $\Fun(\CK, \CC)$ for the functor $\infty$-category. 
    \item For any $\infty$-category $\CC$, we write $\Fil(\CC) = \Fun({\overrightarrow{\Zb}^{\op}}, \CC)$ and $\Gr(\CC) = \Fun((\Zb^\delta)^{\op}, \CC)$. 
    \item We write $\mathsf{Spc}$ for the $\infty$-category of spaces (i.e. $\infty$-groupoids, or $(\infty, 0)$-categories). For any $\infty$-category $\CC$, we write $\CC^\simeq$ for the \emph{core groupoid} (i.e. the maximal sub-$\infty$-groupoid) of $\CC$. 
    \item We write $\mathsf{Spc}_*$ for the $\infty$-category of pointed spaces. 
    \item We write $\mathsf{Set}$ for the $1$-category of sets, and $\Set_*$ for its pointed variant. 
    \item We denote by $\Sp$ the $\infty$-category of spectra. We write $\Map_{\CC}$ for mapping spaces and $\Msp_{\CC}$ for mapping spectra. Furthermore, we write $[X, Y]$ or $[X, Y]_{\CC}$ for $\pi_0 \Map_{\CC}(X, Y)$. 
    \item We denote by $\Ab$ the $1$-category of abelian groups. For each $A \in \Ab$, we write $\mathrm{H} A \in \Sp$ for its Eilenberg--MacLane spectrum. 
    \item For any pointed $1$-category $\CA$, we write $\Ch(\CA)$ for the $1$-category of (cohomological) chain complexes, i.e. the category of pairs $(\{M^s \in \CA\}_{
    s \in \Zb}, d\colon M^s \to M^{s + 1})$ such that $dd = 0$.
    \item A \textbf{$t$-$\infty$-category} is a stable $\infty$-category with a $t$-structure. 
    \item A \textbf{homological functor} is a functor $F\colon \CC \to \CA$ whose source $\CC$ is a stable $\infty$-category, whose target $\CA$ is an abelian $1$-category, such that for any fiber sequence $A \to B \to C$ in $\CC$, the image $F(A) \to F(B) \to F(C)$ is exact in $\CA$. 
    \item We write $\Cat$ (resp. $\Cat_{\st}$) for the $\infty$-category of (stable) $\infty$-categories. 
    \item We denote by $\Alg(\Cat)^{\lax}$ the $\infty$-category of monoidal $\infty$-categories with lax monoidal functors in between, and by $\Alg(\Cat)$ the wide subcategory spanned by monoidal functors. Similar conventions apply to the $\Eb_n$-monoidal variants $\Alg_{\Eb_n}(\Cat) \subset \Alg_{\Eb_n}(\Cat)^{\lax}$, in particular (for $n = \infty$) the symmetric monoidal variants $\CAlg(\Cat) \subset \CAlg(\Cat)^{\lax}$.
    \item For $\CC, \CD \in \Alg(\Cat)$, we write $\Fun^{\otimes \lax}(\CC, \CD)$ for the $\infty$-category of lax monoidal functors, and we write $\Fun^{\otimes \oplax}(\CC, \CD) \coloneq \Fun^{\otimes \lax}(\CC^{\op}, \CD^{\op})^{\op}$ for the $\infty$-category of oplax monoidal functors. Similar conventions apply to the $\Eb_n$-monoidal variants.
    \item For each $\Eb_n$-monoidal $\infty$-category $\CC$, we write $\Pic(\CC)$ for the \textbf{picard $\infty$-groupoid} of $\CC$, i.e. the sub-$\infty$-groupoid of $\CC^\simeq$ spanned by invertible objects. In the case $n = \infty$, we also write $\pic(\CC)$ for the corresponding connective spectra. For $R \in \Alg_{\Eb_{n + 1}}(\Sp)$, we write $\Pic(R)$ for $\Pic(\Mod_R)$; if $n = \infty$ we also write $\pic(R)$ for $\pic(\Mod_R)$. 
    \item We write $\PPr$ (resp. $\PPr_{\st}$) for the $\infty$-category of presentable (resp. presentable stable) $\infty$-categories and left adjoints in between, equipped with the Lurie tensor product $\otimes$. 
    \item We write $X \times_{f, Z, g} Y$ or $X \times_Z Y$ for the pullback of $X \xrightarrow{f} Z \xleftarrow{g} Y$. Similarly for pushouts. 
    \item Suppose $\CK$ is a small $\infty$-category and $\CC$ admits $\CK$-indexed colimits. Write $\CK^\triangleright = \CK \filledstar [0] \coloneq (\CK \times [1]) \sqcup_{\CK \times \{1\}} [0]$. We can construct a colimit comparison functor
    \[\cp = \cp_{\CK}\colon \CC^{\CK^\triangleright} \to \CC^{[1]}, \qquad T \mapsto (\colim T|_{\CK} \to T(\text{cone point})).\]
    Concretely, $\cp$ is the composite $\CC^{\CK^\triangleright} \to \CC^{\CK \times [1]} = \Fun([1], \CC^{\CK}) \xrightarrow{\colim_{\CK}} \Fun([1], \CC)$. 
    \item Suppose $\CC$ is a stable $\infty$-category. Then $h\CC$ admits a triangulated structure, whose shift endomorphism is the suspension functor $\Sigma$, and in which $X \xrightarrow{f} Y \xrightarrow{g} Z \xrightarrow{h} \Sigma X$ is a \textbf{distinguished triangle (DT)} if it extends to a diagram $[2] \times [1] \to \CC$ of the form
    \[\begin{tikzcd}
        X \ar[r] \ar[d] & Y \ar[d] \ar[r] & 0 \ar[d] \\
        0 \ar[r] & Z \ar[r] & \Sigma X
    \end{tikzcd}\]
    so that both squares are cocartesian in $\CC$, and the outer rectangle corresponds to $\id\colon \Sigma X \to \Sigma X$ under the pushout comparison functor, cf. \cite[\S~1.1.2]{HA}. We refer to such an enhancement (if it exists) as an \textbf{extended cofiber sequence} lifting the DT.
    \item In a triangulated category, we write $(h, g, f)$ for a triangle 
    \[\begin{tikzcd}
        X \ar[r, "f"] & Y \ar[r, "g"] & Z \ar[r, "h"] & \Sigma X.
    \end{tikzcd}\]
    We say it is an \textbf{anti-distinguished triangle (ADT)} if $(-h, -g, -f)$ is a distinguished triangle. 
\end{itemize}
\end{convention}

\begin{convention}
    Conventions on equivariant homotopy theory.
\begin{itemize}
    \item We write $C_n$ for the cyclic group of order $n \in \Nb$. 
    \item Suppose $G$ is a finite group. We write $\Sp^{G}$ for the $\infty$-category of \textbf{genuine $G$-spectra}, and we denote by $\oneb \in \Sp^G$ its tensor unit. For $X \in \Sp^G$ and $K \subset G$, we write $X^K$ (resp. $X^{hK}, \Phi^K X$) for its \textbf{categorical} (resp. \textbf{homotopic}, \textbf{geometric}) \textbf{$K$-fixed point}. 
    \item We write $X \mapsto X^e$ for the forgetful functor $\Sp^G \to \Fun(\mathrm{B}G, \Sp)$, and refer to $X^e$ as the \textbf{underlying (Borel $G$-)spectrum} of $X$. This functor admits a fully faithful right adjoint, whose image consists of the \textbf{Borel complete $G$-spectra}. 
    \item Due to \cite{GM11} (see also \cite[Theorem A.1]{CMNN20}), $\Sp^G$ can be identified with the $\infty$-category of \emph{spectral $G$-Mackey functors} $\Mack_G(\Sp) = \Fun^{\times}(\Span(\Fin_G), \Sp)$, i.e. finite product preserving functors from the $\infty$-category of \emph{spans of finite $G$-sets} in \cite[Appendix C]{BH21} to $\Sp$. Therefore, we have a \emph{pointwise $t$-structure} on $\Sp^G$ whose heart is equivalent to $\Mack_G(\Ab) = \Fun^\times(\Span(\Fin_G), \Ab)$, the $1$-category of ordinary $G$-Mackey functors. For $M \in \Mack_G(\Ab)$, we refer to its image $\HH M \in \Sp^G$ as the \textbf{Eilenberg--MacLane $G$-spectrum} of $M$. 
    \item For any group homomorphism $\rho\colon K \to G$, the forgetful functor $\Fin_G \to \Fin_K$ induces a unique symmetric monoidal left adjoint $\rho^*\colon \Sp^G \to \Sp^K$. If $\rho$ is surjective with $G = K / A$, then $\rho^* = \infl^{K}_{K / A}$ is the \textbf{inflation} whose right adjoint is taking categorical $A$-fixed points. If $\rho$ is injective, then $\rho^* = \Res^G_K$ is the \textbf{restriction} whose right adjoint is the \textbf{induction} $\Indrm_K^G$. By the Wirthmüller isomorphism (cf. \cite{MS06, Cno23}), $\Indrm^G_K$ is also left adjoint to $\Res^G_K$. 
    \item We denote by $\Sigma_+^{\infty}\colon \Span(\Fin_G) \cong \Span(\Fin_G)^{\op} \to \Sp^G$ the spectral Yoneda embedding, and for each $X \in \Sp^G, P \in \Fin_G$ we write $X[P] = X \otimes \Sigma_+^{\infty} P$. Note that $X[G/H] \cong \Indrm_H^G \Res^G_H X$. 
    \item We write $\RO(G)$ for the real representation ring of $G$. For any (virtual) $G$-representation $V$, we write $|V| = \dim_{\Rb} V$,  and we denote by $\Sb^V \in \Sp^G$ the \textbf{representation sphere} associated to $V$. If $V$ is an actual $G$-representation, we denote by $a_V\colon \oneb = \Sb^0 \to \Sb^V$ the \textbf{Euler class} of $V$ (induced by the inclusion $0 \subset V$). 
    \item For $K \subset G, X \in \Sp^G$ and $V \in \RO(K)$, we write $\pi_V^{K}(X) = [\Sb^V, \Res^G_K X]_{\Sp^K}$ and refer to these as the \textbf{$\RO(K)$-graded homotopy groups}. 
    \item Suppose $X \in \Sp^G$. For any subgroup $H \subset G$, there is a \textbf{restriction map} 
    \[\res = \res^G_H\colon \pi^G_V(X) \to \pi^H_{\Res^G_H V}(X) = \pi^G_{V}(X[G / H])\] 
    induced by the unit $\res \colon X \to X[G / H]$ of the adjunction $\Indrm_H^G \dashv \Res^G_H$, and a \textbf{transfer map} 
    \[\tr = \tr_H^G\colon \pi^H_{\Res^G_H V}(X) = \pi_{V}^G(X[G / H]) \to \pi^G_V(X)\] 
    induced by the counit $\tr\colon X[G / H] \to X$ of the adjunction $\Res^G_H \dashv \Indrm_H^G$. Note that $\res_{X} \cong X \otimes \res_{\oneb}$, $\tr_X \cong X \otimes \tr_{\oneb}$, and under the Spanier--Whitehead self-duality of $\oneb, \oneb[G/H]$ the two maps $\res_{\oneb}, \tr_{\oneb}$ are dual to each other. As for the maps on homotopy groups, both of them are group homomorphisms, while for products (say $X = R$ is an $\mathbb{A}_2$-ring) we have $\res(a\cdot b) = \res(a) \cdot  \res(b)$ and the \textbf{Frobenius relations} $\tr(a) \cdot b = \tr(a \cdot \res(b)), b \cdot \tr(a) = \tr(\res(b) \cdot a)$. 
\end{itemize}
\end{convention}

\begin{convention}
    Grading conventions in this paper, to be detailed in \S~\hyperref[sec:2]{2}.
\begin{itemize}
    \item For bigraded spheres in $\Fil\Sp$, our convention is determined by three principles:
    \begin{itemize}
        \item The tensor unit $\oneb \in \Fil\Sp$ is $\Sb^{0, 0}$. 
        \item The homological suspension $\Sigma\Sb^{0, 0} = \cofib(0 \to \Sb^{0, 0})$ is $\Sb^{1, 0}$. 
        \item The ``deformation parameter'' (cf. Definition \ref{deformation-parameter}) is $\defopara\colon \Sb^{0, -1} \to \Sb^{0, 0}$.
    \end{itemize}
    \item For spectral sequences, our convention is determined by two principles: 
    \begin{itemize}
        \item Every SS starts from the $E_2$-page. 
        \item We use cohomological Adams convention for differentials, namely $d_r\colon E_r^{s, t} \to E_r^{s + r, t + r - 1}$. 
    \end{itemize}
    \item From any filtered spectrum $X$ we can extract a ``standard spectral sequence'' (cf. Theorem \ref{standard-SS}) whose grading convention is determined by $E_2^{s, t}(X) \cong \pi_{t - s, t}(X / \defopara)$.
\end{itemize}

To compare this with some other conventions in the literature, we use the table below: its first row lists the objects (bigraded homotopy groups and spectral sequence components) in this paper and the following rows describe the same objects in the language of the cited works.

\begin{table}[htbp]
        \centering
        \renewcommand{\arraystretch}{1.35}
        \begin{tabular}{|c||c|c|}
            \hline 
            & Bigraded stems & SS \\
            \hline 
            \hline 
            {\cite{Lin-Wang-Xu-kervaire} and this paper} & $\pi_{n, w}(X)$ & $E_r^{s, t}$ \\
            \hline
            {\cite[\S~7]{Burklund-Xu} and \cite{CDvN2}} & $\pi_{n, w - n}(X)$ & $E_r^{t - s, s}$ \\
            \hline 
            {\cite[\S~4]{HZ25}} & $\pi_{n, w}(X)$ & $E_{r - 1}^{t - s, t}$ \\
            \hline
            {\cite[Part II]{Hed20}} & $\pi_{n}(X(w))$ & $E^{r - 1}_{t - s, t}$ \\
            \hline
            {\cite[\S~3 and \S~4]{Ant24}%\footnote{In the first page of \cite{Ant24} there is also a translation table for different indexing conventions.}
            } & $\pi_{n}(F^w X)$ & $E^{r - 1}_{-t, 2t - s}$  \\
            \hline
            {\cite{HA}} & $\pi_{n}(X(-w))$ & $E_{r - 1}^{-t, 2t - s}$ \\
            \hline
            {\cite[\S~8]{Boa99}} & the $n$-th component of $A^w = \pi_*(F^w X)$ & $E_{r - 1}^{t, s - 2t}$ \\
            \hline
        \end{tabular}
        % \captionsetup{labelformat=empty}
        \caption{Grading conventions in the literature.}
        \label{table-grading-conventions}
    \end{table}
    \FloatBarrier 

    As an example, in \cite[\S~4]{HZ25} the authors adopt the same convention for bigraded spheres, while their spectral sequences start from the $E_1$-page and follow the Atiyah--Hirzebruch convention for differentials (namely $d_r\colon E_r^{n, s} \to E_r^{n - 1, s + r}$). According to the table, their grading convention for the standard spectral sequences is determined by $E_1^{n, s}(X) \cong \pi_{n, s}(X / \defopara)$, and in their convention Theorem \ref{Bockstein-dictionary-infinite} becomes the statement that each $y \in B_{k + 1}^{n, s}(X) \setminus B_{k}^{n, s}(X)$ admits a lift $[y] \in \pi_{n, s}(X)$ so that $\defopara^k  [y] \neq 0$ and $\defopara^{k + 1}[y] = 0$, etc. 
\end{convention}

%% file: Acknowledgment.tex
The author thanks Weinan Lin, Guozhen Wang, and Zhouli Xu for discussions related to their work on hidden extensions in the Adams spectral sequence; Zhipeng Duan, Guchuan Li, and Xiaolin Danny Shi for discussions related to equivariant slice spectral sequence computations; and Eva Belmont, Yueshi Hou, Hana Jia Kong, Cheng Li, Yuxuan Li, Tongtong Liang, Yunze Lu, and Shangjie Zhang for helpful conversations during the preparation of this paper. The author also thanks Dan Isaksen, Xiaolin Danny Shi, Zhouli Xu and Shangjie Zhang for comments on earlier drafts.

%% file: Basics.tex
\label{sec:2}

In this section we set up BIPWX's cofiber-of-$\tau$ formalism in the realm of filtered spectra. More precisely, in \S~\hyperref[subsec:2.1]{2.1} we introduce the $\infty$-category of filtered spectra $\Fil\Sp$, the bigraded spheres, the fundamental operators $\defopara^k, \rho, \delta_k$ and their basic relations; in \S~\hyperref[subsec:2.2]{2.2} we construct the standard spectral sequence of a filtered spectrum and prove the Bockstein dictionary (Theorems \ref{Bockstein-dictionary-finite}, \ref{Bockstein-dictionary-infinite} and \ref{delta-as-total-diff}) that translates between standard SS and the bigraded homotopy groups; in \S~\hyperref[subsec:2.3]{2.3} we deal with the multiplicativity of standard SS and prove the Leibniz rule for total differentials in Theorem \ref{Burklund's-Leibniz-rule}; finally, in \S~\hyperref[subsec:2.4]{2.4} we demonstrate how our results facilitate the study of any presentable $\infty$-category with a cocontinuous $\Fil\Sp$-action, and construct numerous examples. 

\begin{remark} \label{standard-SS-with-general-coefficient-cats}
    Although we work solely with $\Fil\Sp$ in the first three subsections, most of the discussions there hold true for $\Fil(\CE)$ where $\CE$ is a presentably $\Eb_n$ ($n \geq 1$) monoidal $t$-$\infty$-category satisfying some mild assumptions, cf. the construction in Corollary \ref{standard-SS-with-general-coefficients}, the Bockstein dictionary results in Corollary \ref{Bockstein-dictionary-finite-page-with-general-coefficients} and the multiplicativity results in Appendix \ref{app:C}.
\end{remark}

\subsection{Filtered spectra} 
\label{subsec:2.1}

\begin{definition}
    The $\infty$-category of \textbf{filtered spectra} is $\Fil\Sp = \Fun(\overrightarrow{\Zb}^{\op}, \Sp)$. 
\end{definition}

Concretely, an object in $\Fil\Sp$ is a sequence of maps between spectra 
\[X \colon \begin{tikzcd}
    \cdots \ar[r] & X(2) \ar[r] & X(1) \ar[r] & X(0) \ar[r] & X(-1) \ar[r] & \cdots
\end{tikzcd}\]
while a morphism $f\colon X \to Y$ amounts to a family of squares
\[\begin{tikzcd}
    \cdots \ar[r] & X(2) \ar[r] \ar[d, "{f_2}"] & X(1) \ar[r]  \ar[d, "{f_1}"] & X(0) \ar[r]  \ar[d, "{f_0}"] & X(-1) \ar[r]  \ar[d, "{f_{-1}}"] & \cdots \\
    \cdots \ar[r] & Y(2) \ar[r] & Y(1) \ar[r] & Y(0) \ar[r] & Y(-1) \ar[r] & \cdots
\end{tikzcd}\]

\begin{remark}
    For $X \in \Fil\Sp$ treated as a diagram $\overrightarrow{\Zb}^{\op} \to \Sp$, we write $X(-{\infty})$ (resp. $X(\infty)$) for its colimit (resp. limit) in $\Sp$. We say $X \in \Fil\Sp$ is \textbf{complete} if $X(\infty) = 0$. 
\end{remark}

\begin{fact}
    The $\infty$-category $\Fil\Sp$ is presentable and stable. Furthermore, it acquires a presentably symmetric monoidal structure via Day convolution (cf. \cite[Example 2.2.6.17]{HA}). 
\end{fact}

Concretely, the Day convolution tensor unit $\oneb \in \Fil\Sp$ is given as follows: 
\[\begin{tikzcd}
    \cdots \ar[r] & \oneb(2) = 0 \ar[r] & \oneb(1) = 0 \ar[r] & \oneb(0) = \Sb^0 \ar[r, "{\id}"] & \oneb(-1) = \Sb^0 \ar[r, "{\id}"] & \cdots
\end{tikzcd}\]
Moreover, the tensor product $X \otimes Y$ is computed by $(X \otimes Y)(n) = \colim_{a + b \geq n} X(a) \otimes Y(b)$. 

\begin{remark}
    The above discussions also work for the $\infty$-category of \textbf{graded spectra} $\Gr\Sp = \Fun((\Zb^\delta)^{\op}, \Sp)$. Concretely, the Day convolution tensor unit of $\Gr\Sp$ is a single $\Sb^0$ sitting in stage zero, and for $X, Y \in \Gr\Sp$ their tensor product $X \otimes Y$ is given as $(X \otimes Y)(n) = \bigoplus_{a + b = n} X(a) \otimes Y(b)$. 
\end{remark}

\begin{definition}[Bigraded spheres] \label{bigraded-spheres}
    There are two families of invertible objects in $\Fil\Sp$. 
    \begin{itemize}
        \item By construction, the spectral Yoneda embedding functor 
        \[\gamma\colon \overrightarrow{\Zb} \to \Fil\Sp, \quad w \mapsto \gamma_w = \Sigma^{\infty}_+ \Map(-, w)\] 
        is symmetric monoidal. Thus, each $\gamma_w$ is invertible in $\Fil\Sp$ with inverse $\gamma_{-w}$. We write $\Sb^{0, w} = \gamma_w$. Concretely, the image of $k \in \overrightarrow{\Zb}^{\op}$ under $\Sb^{0, w}$ is $\Sb^0$ if $k \leq w$ or $0$ otherwise, and the transition maps are either $0$ or $\id$. In particular, $\Sb^{0, 0} = \oneb$. 
        \item On the other hand, each suspension $\Sigma^{n} \oneb$ is invertible with inverse $\Sigma^{-n} \oneb$. We write $\Sb^{n, 0} = \Sigma^{n} \oneb$. 
    \end{itemize}
    In general, we set $\Sb^{n, w} = \Sb^{n, 0} \otimes \Sb^{0, w}$, and we write $\Sigma^{n, w} X = X \otimes \Sb^{n, w}$ for each $X \in \Fil\Sp$.
\end{definition}

For $X \in \Fil\Sp$ we write $\pi_{n, w}(X)$ for $\pi_0(\Map(\Sb^{n, w}, X))$. Concretely, $\pi_{n, w}(X) = \pi_n(X(w))$. 

\begin{definition}[Deformation parameter] \label{deformation-parameter}
    We write $\defopara\colon \Sb^{0, -1} \to \Sb^{0, 0}$ for the image of $-1 \to 0$ under the spectral Yoneda embedding $\gamma\colon \overrightarrow{\Zb} \to \Fil\Sp$. Concretely, this corresponds to the ladder diagram 
    \[\begin{tikzcd}
        \Sb^{0, -1}  = \ar[d, "\defopara"] \\
        \Sb^{0, 0} \ \ =
    \end{tikzcd} \begin{tikzcd}
        \cdots \ar[r] & 0 \ar[r] \ar[d] & 0 \ar[r] \ar[d] & \Sb^0 \ar[r, "{\id}"] \ar[d, "{\id}"] & \Sb^{0} \ar[r, "{\id}"] \ar[d, "{\id}"] & \cdots \\
        \cdots \ar[r] & 0 \ar[r] & \Sb^{0} \ar[r, "{\id}"]  & \Sb^{0} \ar[r, "{\id}"]  & \Sb^{0} \ar[r, "{\id}"]  & \cdots
    \end{tikzcd}\]
    For $X \in \Fil\Sp$, we write $\defopara$ for $X \otimes \defopara\colon \Sigma^{0, -1} X \to X$. We also write $\defopara^k\colon \Sigma^{0, -k} X \to X$ for the $k$-fold composite of $\defopara$, i.e. the image of $-k \to 0$ under $X \otimes \gamma\colon \overrightarrow{\Zb} \to \Fil\Sp$. 
\end{definition}

For technical convenience, we will adopt the following conventions: We write $\defopara^0$ for the map $\id_X\colon X \to X$, and we write $\defopara^{\infty}$ for the map $0\colon 0 \to X$. We can therefore make sense of $X / \defopara^{r}$ for each $r \in \Nb \cup \{\infty\}$, with $X / \defopara^0 = 0$ and $X / \defopara^{\infty} = X$.

\begin{proposition}[Generic fiber and special fiber]
    \begin{itemize}
        \item The ``$\defopara$-generic fiber'' of $\Fil\Sp$ is $\Sp$. More precisely, $\oneb[\defopara^{-1}] = \colim_{w \to \infty} \Sb^{0, w}$ is an idempotent $\Eb_{\infty}$ algebra in $\Fil\Sp$, and the functor $\Fil\Sp \to \Sp, X \mapsto X(-\infty)$ induces a symmetric monoidal equivalence $\Mod_{\oneb[\defopara^{-1}]} \Fil\Sp \cong \Sp$. 
        \item The ``$\defopara$-special fiber'' of $\Fil\Sp$ is $\Gr\Sp$. More precisely, the cofiber $\oneb / \defopara$ admits an $\Eb_{\infty}$ algebra structure in $\Fil\Sp$, and the functor $\gr_\bullet\colon \Fil\Sp \to \Gr\Sp , X \mapsto \{\gr_w(X) = X(w) / X(w + 1)\}_{w \in \Zb}$ induces a symmetric monoidal equivalence $\Mod_{\oneb / \defopara} \Fil\Sp \cong \Gr\Sp$. 
    \end{itemize}
\end{proposition}

\begin{proof}
    We first construct the two symmetric monoidal left adjoints $\Fil\Sp \to \Sp$ and $\Fil\Sp \to \Gr\Sp$. For the first one, as $\Fil\Sp$ is the free stable cocompletion of $\overrightarrow{\Zb}$, it suffices to construct a symmetric monoidal functor $\overrightarrow{\Zb} \to \Sp$ sending each $w$ to $\Sb^0$. This comes from the composite $\overrightarrow{\Zb}  \to * \to \Sp$ where the second map hits the tensor unit $\Sb^0 \in \Sp$. For the second one, it suffices to construct a symmetric monoidal functor $\overrightarrow{\Zb} \to \Gr\Sp$ sending each $w$ to $\Sb^{0, w}$. This comes from pre-composing $\Gr(\Sigma^{\infty})\colon \Gr(\mathsf{Spc}_*) \to \Gr\Sp$ with the symmetric monoidal functor $\overrightarrow{\Zb} \to \Gr(\mathsf{Set}_*) \subset \Gr(\mathsf{Spc}_*)$, where the latter functor is the unique enhancement of the pointed Yoneda embedding $\Zb^\delta \rightarrow \Gr(\mathsf{Set}) \to \Gr(\mathsf{Set}_*)$ that assigns each non-identity map $z < w$ in $\overrightarrow{\Zb}$ to the zero map in the codomain. \parr 
    
    From the construction it is immediate that the right adjoint $\Sp \to \Fil\Sp$ sends $\Sb^0$ to the constant filtration with value $\Sb^0$, which is the same as $\oneb[\defopara^{-1}]$. Similarly, the right adjoint $\Gr\Sp \to \Fil\Sp$ sends $\Sb^{0, 0}$ to the filtration $\cdots \to 0 \to \Sb^0 \to 0 \to \cdots$ with $\Sb^0$ concentrated in degree $0$, which is the same as $\oneb / \defopara$. These right adjoints are lax symmetric monoidal (since their left adjoints are symmetric monoidal), so they equip $\oneb[\defopara^{-1}]$ and $\oneb / \defopara$ with $\Eb_{\infty}$-algebra structures in $\Fil\Sp$. Furthermore, $\oneb[\defopara^{-1}] \otimes \oneb[\defopara^{-1}] \cong \oneb[\defopara^{-1}]$ since $\oneb[\defopara^{-1}] = \colim_{w \to \infty} \Sb^{0, w}$ where the index category $\overrightarrow{\Zb}$ is filtered, i.e. $\oneb[\defopara^{-1}]$ is idempotent, so the $\Eb_{\infty}$-algebra structure on $\oneb[\defopara^{-1}]$ is unique (\cite[Proposition 4.8.2.9]{HA}). It remains to prove the two induced functors $\Sp \to \Mod_{\oneb[\defopara^{-1}]} \Fil\Sp$ and $\Gr\Sp \to \Mod_{\oneb / \defopara} \Fil\Sp$ are equivalences. These equivalences are recorded in \cite[Examples A.5 and A.6]{BHS2}, as straightforward applications of the criterion \cite[Proposition A.4]{BHS2}.
\end{proof}

\begin{remark}
    See Example \ref{SS-examples-2} item 4 for a justification of the names  ``generic fiber'' and ``special fiber'' from the perspective of algebraic geometry.
\end{remark}

\begin{remark}
    Since $\oneb[\defopara^{-1}]$ is idempotent, the forgetful functor $\Mod_{\oneb[\defopara^{-1}]} \Fil\Sp \to \Fil\Sp$ is fully faithful. For $X \in \Fil\Sp$, we will not distinguish the constant filtration $X[\defopara^{-1}] = X \otimes \oneb[\defopara^{-1}]$ and its image $X[\defopara^{-1}](-\infty) = X(-\infty) \in \Sp$. In particular, we write $\pi_{n}(X[\defopara^{-1}])$ for $\pi_{n}(X(-\infty))$, which is also isomorphic to $\pi_{n, w}(X[\defopara^{-1}])$ for each $w \in \Zb$. 
\end{remark}

\begin{remark}
    For $X \in \Fil\Sp$, the limit $\lim_{w \to \infty} \Sigma^{0, -w} X$ along multiplication by $\defopara$ is the constant filtration on $X({\infty})$, so $X$ is complete iff $\lim_{w \to \infty} \Sigma^{0, -w} X = 0$. The inclusion of complete filtered spectra $\Fil\Sp_{\hat{\defopara}} \hookrightarrow \Fil\Sp$ admits a left adjoint that sends $X$ to its \textbf{completion} $X_{\hat{\defopara}}$, here
    \[X_{\hat{\defopara}} = \cofib(\lim_{w \to \infty} \Sigma^{0, -w} X \to X) = \lim_{w \to \infty} X / \defopara^w.\]
    It follows that $X \to X_{\hat{\defopara}}$ induces an isomorphism $X / \defopara^w \cong X_{\hat{\defopara}} / \defopara^w$ for each $0 \leq w < \infty$. 
\end{remark}

\begin{remark} \label{recollements-and-Nakayama-lemma}
    We record some categorical facts for $\Fil\Sp$: 
    \begin{itemize}
        \item The two inclusions $\Sp = \Fil\Sp[\defopara^{-1}] \hookrightarrow \Fil\Sp \hookleftarrow \Fil\Sp_{\hat{\defopara}}$ exhibit $\Fil\Sp$ as a \emph{recollement} of $\Sp$ and $\Fil\Sp_{\hat{\defopara}}$ in the sense of \cite[Definition A.8.1]{HA} due to \cite[Theorem A.8.20]{HA}. Furthermore, we can equip $\Fil\Sp_{\hat{\defopara}}$ with the complete tensor product $(X, Y) \to (X \otimes Y)_{\hat{\defopara}}$, which upgrades this to a \emph{symmetric monoidal recollement} in the sense of \cite[Definition 2.20]{Sha21}. 
        \item If a map $f\colon X \to Y$ in $\Fil\Sp$ induces an isomorphism $X / \defopara \cong Y / \defopara$, then it also induces an isomorphism $X / \defopara^w \cong Y / \defopara^w$ for each $w \geq 1$ and an isomorphism of completions $X_{\hat{\defopara}} \cong Y_{\hat{\defopara}}$. This is because each $X / \defopara^w$ is an iterated extension of $X / \defopara$, and $X_{\hat{\defopara}} = \lim_{w \to \infty} X / \defopara^w$.
    \end{itemize}
\end{remark}

\begin{construction}[Fundamental distinguished triangles] \label{fundamental-distinguished-triangle}
    For any pair of arrows $A \to B \to C$ in a stable $\infty$-category, there is an ``octahedral'' distinguished triangle $B / A \to C / A \to C / B \to \Sigma (B / A)$, cf. \cite[Proof of (TR4) in Theorem 1.1.2.14]{HA}. Applying this to the pair 
    \[\Sigma^{0, -r} \oneb \xrightarrow{\defopara^{r - k}} \Sigma^{0, - k} \oneb \xrightarrow{\defopara^{k}} \oneb\] 
    for $\oneb \in \Fil\Sp$, $0 \leq k < r \leq \infty$ and tensoring this with $X \in \Fil\Sp$ leads to the distinguished triangle 
    \[\Sigma^{0, -k} X / \defopara^{r - k} \xrightarrow{\defopara^{k}} X / \defopara^{r} \xrightarrow{\rho_{k}^{r}} X / \defopara^{k} \xrightarrow{\delta_{k}^{r - k}} \Sigma^{1, -k} X / \defopara^{r - k}. \]
    In practice, we usually suppress some of the indices by writing $\rho = \rho_k^r$ and $\delta_k = \delta_k^{r - k}$.  
\end{construction}

In some sense, $\defopara$, $\rho$ and $\delta$ are the ``gist'' of $\Fil\Sp$. In the next section, we will study their effect in terms of spectral sequences. We end this section by collecting some relations on their composites. 

\begin{proposition} \label{defopara-n-refines-defopara-n}
    For $a, b \geq 1$, the composite $\Sigma^{0, -b} X / \defopara^a \xrightarrow{\defopara^b} X / \defopara^{a + b} \xrightarrow{\rho} X / \defopara^a$ is homotopic to the image of $\defopara^b\colon \Sb^{0, -b} \to \Sb^{0, 0}$ under $(X / \defopara^a) \otimes - \colon \Fil\Sp \to \Fil\Sp$. 
\end{proposition}

\begin{remark}
    It is for this reason that we reserve the name $\defopara^b$ for the map $\Sigma^{0, -b} X / \defopara^a \to X / \defopara^{a + b}$ despite potential ambiguity. When we have to distinguish this with $(X / \defopara^a) \otimes (\defopara^b\colon \Sb^{0, -b} \to \Sb^{0, 0})$, we choose to write $\cdot \defopara^b$ or $\defopara^b \cdot$ for the latter map as it is the actual ``multiplication by $\defopara^b$''. 
\end{remark}

\begin{proof}
    It suffices to prove this for $X = \oneb$. The map $(\oneb / \defopara^a) \otimes \defopara^b$ is the vertical cofiber of the square 
    \[\defopara^a \otimes \defopara^b \colon \begin{tikzcd}
        \Sb^{0, - a} \otimes \Sb^{0, -b} \ar[r, "{\id \otimes \defopara^b}"] \ar[d, "{\defopara^a \otimes \id}"] & \Sb^{0, - a} \otimes \Sb^{0, 0} \ar[d, "{\defopara^a \otimes \id}"] \\
         \Sb^{0, 0} \otimes \Sb^{0, -b} \ar[r, "{\id \otimes \defopara^b}"] & \Sb^{0, 0} \otimes \Sb^{0, 0}
    \end{tikzcd}\]
    while the composite $\rho \defopara^b$ is the vertical cofiber of the outer rectangle in 
    \[\begin{tikzcd}
        \Sb^{0, - a - b} \ar[r, "{\id}"] \ar[d, "{\defopara^a}"] & \Sb^{0, - a - b} \ar[r, "{\defopara^b}"] \ar[d, "{\defopara^{a + b}}"] &  \Sb^{0, -a} \ar[d, "{\defopara^a}"] \\
        \Sb^{0, - b} \ar[r, "{\defopara^b}"]  & \Sb^{0, 0} \ar[r, "{\id}"]  &  \Sb^{0, 0} 
    \end{tikzcd}\]
    so it reduces to producing an isomorphism between the two squares. In general, it would be difficult as one has to compare the two filling homotopies. However, in this case both squares lie in the image of the monoidal functor $\gamma\colon \overrightarrow{\Zb} \to \Fil\Sp$, and their preimages can both be identified with the unique square 
    \[\begin{tikzcd}
        -a - b \ar[r] \ar[d] & -a \ar[d] \\
        -b \ar[r] & 0
    \end{tikzcd}\]
    since the source is a $1$-category in which every commutative square has a unique filler.
\end{proof}

\begin{remark}
    This result is not as innocent as it seems. For instance, taking $a = b$, it implies $\cdot \defopara^a\colon \Sigma^{0, -a} X / \defopara^a \to X / \defopara^a$ is nullhomotopic, which is a bit surprising compared to the fact in $\Sp$ that $\cdot 2 \colon \Sb^0 / 2 \to \Sb^0 / 2$ is not nullhomotopic. The above argument does not work in the Moore spectrum situation because there is no \emph{monoidal} functor $\overrightarrow{\Zb} \to \Sp$ sending $-1 \to 0$ to $\cdot 2 \colon \Sb^0 \to \Sb^0$, or equivalently (cf. Theorem \ref{GHMG}) $\Sb^0 / 2$ does not admit the structure of an $\Eb_1$-algebra in $\Sp$.  
\end{remark}

\begin{fact} \label{basic-commutations}
    It is clear from construction that $\defopara^a \defopara^b = \defopara^{a + b}$, $\rho_{a}^{a + b} \rho_{a + b}^{a + b + k} = \rho_a^{a + b + k}$ and $\rho_{a}^{a + b}\delta_k^{a + b} = \delta_k^a$.     
\end{fact}

\begin{lemma} \label{defopara-rho-cartesian}
   For any $X \in \Fil\Sp$, $1 \leq a \leq \infty, 0 \leq b < \infty , 0 \leq k \leq \infty$, we have a cartesian square:
    \[\begin{tikzcd}
        \Sigma^{0, -b} X / \defopara^{a + k} \ar[r, "{\defopara^b}"] \ar[d, "{\rho}"] & X / \defopara^{a + k + b} \ar[d, "{\rho}"] \\
        \Sigma^{0, -b} X / \defopara^{a} \ar[r, "{\defopara^b}"] & X / \defopara^{a + b}
    \end{tikzcd}\]
\end{lemma}

\begin{proof}
    Such a square is constructed by taking horizontal fibers of the square
    \[\begin{tikzcd}
        X / \defopara^{a + k + b} \ar[r, "{\rho}"] \ar[d, "{\rho}"] & X / \defopara^{b}  \ar[d, "{\id}"] \\
        X / \defopara^{a + b} \ar[r, "{\rho}"] &  X / \defopara^{b}
    \end{tikzcd}\]    
    and it is cartesian since the vertical right arrow of the latter square is an isomorphism.  
\end{proof}

\begin{corollary} \label{defopara-totaldiff}
    Under the same setup, if $a, k < \infty$, there is also a commutative square in $\Fil\Sp$ 
    \[\begin{tikzcd}
        \Sigma^{0, -b} X / \defopara^{a} \ar[r, "{\defopara^b}"] \ar[d, "{\delta_a}"] & X / \defopara^{a + b} \ar[d, "{\delta_{a + b}}"]\\
        \Sigma^{1, - a -b} X / \defopara^{k} \ar[r, "{\id}"] & \Sigma^{1, - a -b} X / \defopara^{k}
    \end{tikzcd}\]
    deduced from the original square by taking vertical cofibers. 
\end{corollary}

\begin{lemma} \label{projection-totaldiff-cartesian}
    For $X \in \Fil\Sp$ and $0 \leq a \leq \infty, 1 \leq b < \infty, 0 \leq k < \infty$, we have a cartesian square: 
    \[\begin{tikzcd}
        X / \defopara^{b + k} \ar[r, "{\delta_{b + k}}"] \ar[d, "{\rho}"] &\Sigma^{1, -b - k} X / \defopara^a \ar[d, "{\defopara^k}"] \\
        X / \defopara^{b} \ar[r,"{\delta_b}"] & \Sigma^{1, -b} X / \defopara^{a + k}
    \end{tikzcd}\]
\end{lemma}

\begin{proof}
    Such a square is constructed by taking horizontal cofibers of the square
    \[\begin{tikzcd}
        X / \defopara^{a + k + b} \ar[r, "{\rho}"] \ar[d, "{\id}"] & X / \defopara^{b + k}  \ar[d, "{\rho}"] \\
        X / \defopara^{a + k + b} \ar[r, "{\rho}"] &  X / \defopara^{b}
    \end{tikzcd}\]    
    and it is cartesian since the vertical left arrow of the latter square is an isomorphism.  
\end{proof}

\subsection{SS and the Bockstein dictionary}
\label{subsec:2.2}

\begin{definition} \label{SS}
    A \textbf{spectral sequence (SS)} of abelian groups $E = \{E^{s, t}_r\}$ consists of the following data: 
    \begin{itemize}
        \item For each $r \geq 2$, $s, t \in \Zb$, an abelian group $E_r^{s, t}$. 
        \item For each $r \geq 2$, $s, t \in \Zb$, a map of abelian groups $d_r\colon E_r^{s, t} \to E_r^{s + r, t + r - 1}$, so that $d_r d_r = 0$. 
        \item For each $r \geq 2$, $s, t \in \Zb$, an isomorphism $\frac{\ker(d_r \colon E_r^{s, t} \to E_r^{s + r, t + r - 1})}{\IIm(d_r\colon E_r^{s - r, t - r + 1} \to E_r^{s, t})} \cong E_{r + 1}^{s, t}$. 
    \end{itemize}
    A map between SS $E \to E'$ is a family of maps of abelian groups $E_r^{s, t} \to (E')_{r}^{s, t}$ that commutes with the $d_r$ maps and the structural isomorphisms. These assemble into a $1$-category $\SpSeq$. 
\end{definition}

\begin{remark}
    When drawing a SS in the plane, we take $n = t - s$ (resp. $s$) as the horizontal (resp. vertical) coordinate. Thus a $d_r$-differential in the picture moves one step to the left and $r$ steps up.
\end{remark}

\begin{remark} \label{cycles-and-boundaries}
    Given a SS $E = \{E_r^{s, t}\}$, we refer to the maps $d_r$ as the \textbf{differentials} on the $E_r$-page. Furthermore, by lifting subquotients, we can produce a canonical family of subgroup inclusions for each $E_2^{s, t}$: 
    \[0 = B_1^{s, t} \subset B_2^{s, t} \subset B_3^{s, t} \subset \cdots \subset B_r^{s, t} \subset \cdots \subset Z_r^{s, t} \subset \cdots \subset Z_3^{s, t} \subset Z_2^{s, t} \subset Z_1^{s, t} = E_2^{s, t}\]
    so that $E_{r + 1}^{s, t} \cong Z_{r}^{s, t} / B_{r}^{s, t}$ for each $r \geq 1$. We refer to $Z_r^{s, t}$ (resp. $B_r^{s, t}$) as the group of \textbf{$d_r$-cycles} (resp. \textbf{$d_r$-boundaries}) at bidegree $(s, t)$. By construction, any map of SS $f\colon E \to E'$ induces a family of maps between $d_r$-cycles (resp. $d_r$-boundaries). 
\end{remark}

To each filtered spectrum $X$ we can associate a ``standard SS'' $E_*^{*,*}(X)$. We start with a lemma:

\begin{lemma} \label{cartesian-square-composite-image}
    Suppose $\CC$ is a stable $\infty$-category and there is a cartesian square in $\CC$ 
    \[\begin{tikzcd}
        A \ar[r] \ar[d] & C \ar[d] \\
        B \ar[r] & D
    \end{tikzcd}\]
    then for each homological functor $F\colon \CC \to \CA$, we have canonical identifications in $\CA$:
    \begin{align*}
        \IIm(F(A)\to F(D)) &\cong \frac{F(A)}{\ker(F(A) \to F(B)) \oplus \ker(F(A) \to F(C))} \\
        &\cong \frac{\IIm(F(A) \to F(B))}{\IIm(\ker(F(A) \to F(C)) \to F(A) \to  F(B))} \\
        &\cong \frac{\IIm(F(A) \to F(C))}{\IIm(\ker(F(A) \to F(B)) \to F(A) \to F(C))}.
    \end{align*}
\end{lemma}

\begin{proof}
    Write $X = \fib(A \to B) = \fib(C \to D)$ and $Y = \fib(A \to C) = \fib(B \to D)$. We will prove that $X \to A$, $Y \to A$ and $A \to D$ fit into a cofiber sequence $X \oplus Y \to A \to D$. Consider the diagrams
    \[\begin{tikzcd}
        A \ar[r, "{\id}"] \ar[d] & A \ar[r] \ar[d] & B \ar[d] \\
        B \ar[r] & D \ar[r, "{\id}"] & D
    \end{tikzcd} \qquad \text{and} \qquad
    \begin{tikzcd}
        A \ar[r, "{\id}"] \ar[d] & A \ar[r] \ar[d] & C \ar[d] \\
        C \ar[r] & D \ar[r, "{\id}"] & D
    \end{tikzcd}\]
    Taking vertical fibers for these, we obtain fiber sequences $X \to \fib(A \to D) \to Y$ and $Y \to \fib(A \to D) \to X$ so that the composites $X \to \fib(A \to D) \to X$ and $Y \to \fib(A \to D) \to Y$ are $\id$. From these we obtain an isomorphism $\fib(A \to D) \cong X \oplus Y$ compatible with their maps into $A$. Taking the $F(-)$ long exact sequence for $X \oplus Y \to A \to D$ yields the first isomorphism we want, while the other two follow by a straightforward image computation in $\CA$.
\end{proof}

\begin{remark} \label{cartesian-square-composite-image-dual}
    Under the above setup, the same argument also proves a dual statement:
    \begin{align*}
        \IIm(F(A) \to F(D)) &\cong \ker(F(D) \to F(D) / F(B) \oplus F(D) / F(C)) \\
        &\cong \ker(\IIm(F(B) \to F(D)) \to \IIm(F(B) \to F(D) \to F(D) / F(C))) \\
        &\cong \ker(\IIm(F(C) \to F(D)) \to \IIm(F(C) \to F(D) \to F(D) / F(B))).
    \end{align*}
\end{remark}

\begin{theorem} \label{standard-SS}
    For each $X \in \Fil\Sp$ there is a \textbf{standard SS} $E_*^{*,*}(X) = \{E_r^{s, t}(X)\}_{r \geq 2}$ as follows: 
    \begin{enumerate}
        \item $E_2^{s, t}(X) = \pi_{t - s, t}(X / \defopara)$. 
        \item For $r \geq 2$, each $E_r^{s, t}(X)$ is a quotient of $\pi_{t - s, t}(X / \defopara^{r - 1})$. Actually, $E_r^{s, t} = Z_{r - 1}^{s, t} / B_{r - 1}^{s, t}$, where 
        \begin{itemize}
            \item $Z_{r - 1}^{s, t} \subset E_2^{s, t} = \pi_{t - s, t}(X / \defopara)$ is the image of $\rho\colon \pi_{t - s, t}(X / \defopara^{r - 1}) \to \pi_{t - s, t}(X / \defopara)$. 
            \item $B_{r - 1}^{s, t} \subset E_2^{s, t} = \pi_{t - s, t}(X / \defopara)$ is the image of $\delta_{r - 2}\colon \pi_{t - s + 1, t - r + 2}(X / \defopara^{r - 2}) \to \pi_{t - s, t}(X / \defopara)$. 
        \end{itemize}
        Here we adopt the convention $Z_1^{s, t} = E_2^{s, t}$ and $B_1^{s, t} = 0$. 
        \item $|d_r| = (r, r - 1)$ for each $r \geq 2$. For $x \in Z_{r - 1}^{s, t}$, if $[x] \in \pi_{t - s, t}(X / \defopara^{r - 1})$ is a lift of $x$, then $\delta_{r - 1}[x] \in \pi_{t - s - 1, t + r - 1}(X / \defopara) = E_2^{s + r, t + r - 1}(X)$ is a representative of $d_r(x)$. 
    \end{enumerate}
    Furthermore, these assemble into a functor $E_*^{*,*}\colon \Fil\Sp \to \SpSeq, X \mapsto E_*^{*,*}(X)$. 
\end{theorem} 

\begin{proof}
    Consider the following square in $\Fil\Sp$ 
    \[\begin{tikzcd}
        X / \defopara^{r - 1} \ar[r, "{\defopara^{r - 2}}"] \ar[d, "{\rho^{r - 1}_1}"] & \Sigma^{0, r - 2} X / \defopara^{2r - 3} \ar[d, "{\rho^{2 r  - 3}_{r - 1}}"] \\
        X / \defopara \ar[r, "{\defopara^{r - 2}}"] &  \Sigma^{0, r - 2} X / \defopara^{r - 1}
    \end{tikzcd}\]
    which is cartesian by Lemma \ref{defopara-rho-cartesian}. Here the composite $\defopara^{r - 2} \rho^{r - 1}_1\colon X / \defopara^{r - 1} \to \Sigma^{0, r - 2} X / \defopara^{r - 1}$ is homotopic to $\cdot \defopara^{r - 2}$ by Proposition \ref{defopara-n-refines-defopara-n}. Applying Lemma \ref{cartesian-square-composite-image} to this square, we conclude that the image of $\cdot \defopara^{r - 2}\colon \pi_{t - s, t}(X / \defopara^{r - 1}) \to \pi_{t - s, t - r + 2}(X / \defopara^{r - 1})$ is $\IIm(\rho^{r - 1}_1) / \IIm(\rho^{r - 1}_{r - 2} \delta_{r - 2}^{r - 1}) = Z_{r - 1}^{s, t} / B_{r - 1}^{s, t} = E_r^{s, t}$. Furthermore, we have a commutative square 
    \[\begin{tikzcd}
        X / \defopara^{r - 1} \ar[r, "{\cdot\defopara^{r - 2}}"] \ar[d, "{\delta_{r - 1}^{r - 1}}"] & \Sigma^{0, r - 2} X / \defopara^{r - 1} \ar[d, "{\delta_{r - 1}^{r - 1}}"] \\
        \Sigma^{1, 1 - r} X / \defopara^{r - 1} \ar[r, "{\cdot\defopara^{r - 2}}"] & \Sigma^{1, -1} X / \defopara^{r - 1}
    \end{tikzcd}\]
    since $- \otimes \defopara^{r - 2}$ defines a natural transformation. Taking $\pi_{t - s, t}$, we then obtain the expected $d_r$
    \[\begin{tikzcd}
        \pi_{t - s, t}(X / \defopara^{r - 1}) \ar[r] \ar[d, "{\delta_{r - 1}^{r - 1}}"] &  E_r^{s, t} \ar[d, dashed, "{d_r}"]  \ar[r] &  \pi_{t - s, t - r + 2}(X / \defopara^{r - 1}) \ar[d, "{\delta_{r - 1}^{r - 1}}"] \\
        \pi_{t - s - 1, t + r - 1}(X / \defopara^{r - 1}) \ar[r] & E_r^{s + r, t + r - 1} \ar[r] & \pi_{t - s - 1, t + 1}(X / \defopara^{r - 1})
    \end{tikzcd}\]
    from functoriality of the epi-mono factorization applied to the horizontal arrows. It remains to show $d_r d_r = 0$ and $\ker(d_r) / \IIm(d_r) \cong E_{r + 1}^{s, t}$. Actually, $d_rd_r = 0$ because $\delta_{r - 1}^{r - 1} \delta_{r - 1}^{r - 1} \simeq \delta_{r - 1}^{r - 1} \rho^{2r - 2}_{r - 1} \delta_{r - 1}^{2r - 2} \simeq 0$ due to Fact \ref{basic-commutations}. Furthermore, the composite map 
    \[\pi_{t - s, t}(X / \defopara^r) \xrightarrow{\rho_{r - 1}^r} \pi_{t - s, t}(X / \defopara^{r - 1}) \to E_r^{s, t}\] 
    factors through an epimorphism onto $\ker(d_r) \subset E_r^{s, t}$, since $\rho^{r - 1}_1\delta_{r - 1}^{r - 1}\rho^r_{r - 1} \simeq \delta_{r - 1}^1\rho^r_{r - 1} \simeq 0$ and the two maps $\rho_{r - 1}^r\colon \pi_{t - s, t}(X / \defopara^r) \to \pi_{t - s, t}(X / \defopara^{r - 1})$, $\defopara\colon \pi_{t - s, t + 1}(X / \defopara^{r - 2}) \to \pi_{t - s, t}(X / \defopara^{r - 1})$ collectively surject onto the kernel of $\defopara^{r - 2} \rho_1^{r - 1} \delta_{r - 1}^{r - 1} \simeq \defopara^{r - 2} \delta_{r - 1}^1 \simeq \delta_1^{r - 1} \rho_1^{r - 1}\colon \pi_{t - s - 1, t}(X / \defopara^{r - 1}) \to \pi_{t - s - 1, t + 1}(X / \defopara^{r - 1})$ due to an application of Lemma \ref{cartesian-square-composite-image} on the cartesian square witnessing $ \defopara^{r - 2} \delta_{r - 1}^1 \simeq \delta_1^{r - 1} \rho_1^{r - 1}$ from Lemma \ref{projection-totaldiff-cartesian}. A dual argument (replacing Lemma \ref{cartesian-square-composite-image} by Remark \ref{cartesian-square-composite-image-dual}) shows the composite
    \[E_{r}^{s, t} \to \pi_{t - s, t - r + 2}(X / \defopara^{r - 1}) \xrightarrow{\defopara} \pi_{t - s, t - r + 1}(X / \defopara^r)\]
    factors through a monomorphism out of $E^{s, t}_r / \IIm(d_r)$. Splicing these together, we obtain 
    \[\pi_{t - s, t}(X / \defopara^r) \to \ker(d_r) \to \ker(d_r) / \IIm(d_r) \to  E^{s, t}_r / \IIm(d_r) \to \pi_{t - s, t - r + 1}(X / \defopara^r)\] 
    which serves as an epi-mono factorization of the composite $\defopara \defopara^{r - 2} \rho^{r - 1}_1 \rho^r_{r - 1} \simeq \cdot \defopara^{r - 1} \colon \pi_{t - s, t}(X / \defopara^r) \to \pi_{t - s, t - r + 1}(X / \defopara^r)$. Thus, $\ker(d_r) / \IIm(d_r) \cong E_{r + 1}^{s, t}$ by uniqueness of the epi-mono factorization. \parr 
    
    Finally, the functoriality of $E_*^{*,*}\colon \Fil\Sp \to \SpSeq$ is straightforward from the construction. 
\end{proof}

In fact, the definition of SS makes sense for every abelian $1$-category. On the other hand, if $\CE$ is a presentable stable $\infty$-category, then the Lurie tensor product $\Fil(\CE) = \CE \otimes \Fil\Sp$ admits a $\Fil\Sp$-action, so the notions of $\Sigma^{n, w} X = X \otimes \Sb^{n, w}$, $X / \defopara^r$, and the fundamental distinguished triangle also make sense in $\Fil(\CE)$; see Definition \ref{FilSp-module-def} and Construction \ref{fundamental-stuff-for-FilSp-modules} for details. The argument in Theorem \ref{standard-SS} thus yields the following general result.  

\begin{corollary} \label{standard-SS-with-general-coefficients}
    Let $\CE$ be a presentable stable $\infty$-category, $\CE^\heartsuit$ be an abelian $1$-category and $\pi_0\colon \CE \to \CE^\heartsuit$ be a homological functor.\footnote{In practice, usually $\CE$ is a $t$-$\infty$-category, and $\pi_0\colon \CE \to \CE^\heartsuit$ takes the $0$-th homotopy object in its heart.} We denote by $\pi_n\colon \CE \to \CE^\heartsuit$ ($n \in \Zb$) the composite $\pi_0\circ \Sigma^{-n}$. Then for each $X \in \Fil(\CE)$, there is an \textbf{$\CE^\heartsuit \!$-valued standard SS} $E_*^{*,*}(X) = \{E_r^{s, t}(X)\}_{r \geq 2}$ as follows:  
    \begin{enumerate}
        \item $E_2^{s, t}(X) = \pi_{t - s, t}(X / \defopara)$. Here by $\pi_{n, w} Y$ we mean $\pi_n(Y(w)) \in \CE^\heartsuit$. 
        \item For $r \geq 2$, each $E_r^{s, t}(X)$ is a quotient of $\pi_{t - s, t}(X / \defopara^{r - 1})$. Actually, $E_r^{s, t} = Z_{r - 1}^{s, t} / B_{r - 1}^{s, t}$, where 
        \begin{itemize}
            \item $Z_{r - 1}^{s, t} \subset E_2^{s, t} = \pi_{t - s, t}(X / \defopara)$ is the image of $\rho\colon \pi_{t - s, t}(X / \defopara^{r - 1}) \to \pi_{t - s, t}(X / \defopara)$. 
            \item $B_{r - 1}^{s, t} \subset E_2^{s, t} = \pi_{t - s, t}(X / \defopara)$ is the image of $\delta_{r - 2}\colon \pi_{t - s + 1, t - r + 2}(X / \defopara^{r - 2}) \to \pi_{t - s, t}(X / \defopara)$. 
        \end{itemize}
        Alternatively, $E_r^{s, t}$ is the middle term in the epi-mono factorization of $\cdot \defopara^{r - 2}\colon \pi_{t - s, t}(X / \defopara^{r - 1}) \to \pi_{t - s, t - r + 2}(X / \defopara^{r - 1})$ in the abelian $1$-category $\CE^\heartsuit$.
        \item The differential $d_r\colon E_r^{s, t} \to E_{r}^{s + r, t + r - 1}$ comes from the commutative square
        \[\begin{tikzcd}
            \pi_{t - s, t}(X / \defopara^{r - 1}) \ar[r, "{\cdot \defopara^{r - 2}}"] \ar[d, "{\delta_{r - 1}^{r - 1}}"] &  \pi_{t - s, t - r + 2}(X / \defopara^{r - 1}) \ar[d, "{\delta_{r - 1}^{r - 1}}"] \\
            \pi_{t - s - 1, t + r - 1}(X / \defopara^{r - 1}) \ar[r, "{\cdot \defopara^{r - 2}}"] & \pi_{t - s - 1, t + 1}(X / \defopara^{r - 1})
        \end{tikzcd}\]
        through functoriality of the epi-mono factorization.
    \end{enumerate}
    Furthermore, these assemble into a functor $E_*^{*,*}\colon \Fil(\CE) \to \SpSeq(\CE^\heartsuit), X \mapsto E_*^{*,*}(X)$. 
\end{corollary}

\begin{remark} \label{Lurie-comparison}
    Up to reindexing, the standard SS construction in the proof of Theorem \ref{standard-SS} coincides with \cite[Construction 1.2.2.6]{HA}, and the subsequent argument amounts to a reformulation of \cite[Theorem 1.2.2.7]{HA} in terms of BIPWX's cofiber-of-$\tau$ formalism. 
\end{remark}

\begin{remark} \label{Boardman-comparison}
    Alternatively, for $X \in \Fil\Sp$, $\pi_{**}(X)$ and $\pi_{**}(X / \defopara)$ form an unrolled exact couple, and the SS it induces (cf. \cite{Boa99}) is precisely the standard SS.
\end{remark}

\begin{remark} \label{Boardman-completion}
    The standard SS of $X \in \Fil\Sp$ depends only on its completion $X_{\hat{\defopara}}$ since taking completion does not change any of the mod-$\defopara^r$ cofibers. On the other hand, if $X$ is complete, then by the Milnor sequence and the fact $X(\infty) = 0$, we see the unrolled exact couple $(\pi_{**}(X), \pi_{**}(X / \defopara))$ imposes \textbf{conditional convergence} to $\pi_{*}(X[\defopara^{-1}])$ on the standard SS in the sense of \cite[Definition 5.10]{Boa99}. We return to convergence later in Theorem \ref{Bockstein-dictionary-infinite}.
\end{remark}

\begin{remark} \label{all-SS-are-Bockstein-SS}
    For $X \in \Fil\Sp$ one can consider the $\overrightarrow{\Zb}^{\op}$-indexed $\defopara$-Bockstein filtration in $\Fil\Sp$:
    \[\cdots \xrightarrow{\defopara} \Sigma^{0, -2} X \xrightarrow{\defopara} \Sigma^{0, -1} X \xrightarrow{\defopara} X \xrightarrow{\defopara} \Sigma^{0, 1} X\xrightarrow{\defopara} \Sigma^{0, 2} X \xrightarrow{\defopara} \cdots\]
    its image under $\Msp(\oneb, -)\colon \Fil(\Fil\Sp) \to \Fil\Sp$ is the filtration $X$ itself. Together with the standard SS construction, this justifies the slogan ``Every SS (from filtered spectra) is a Bockstein SS''.
\end{remark}

We can translate back and forth between $\pi_{**}(X / \defopara^r)$ and the first $(r + 1)$ pages of the standard SS. The forward direction is clear: the standard SS is defined using the bigraded stems. The other direction is given by the \textbf{Bockstein dictionaries}\footnote{We choose this name since the results amount to a $\defopara$-Bockstein SS computation, although we will prove them in a slightly different way. Similar results in \cite{BHS1} are referred to as the ``omnibus theorem''.} below.

\begin{theorem}[Bockstein dictionary, finite page]\label{Bockstein-dictionary-finite}
    Take $X \in \Fil\Sp$ and $1 \leq r < \infty$. For each $(s, t) \in \Zb \times \Zb$, there is a finite filtration on $\pi_{t - s, t}(X / \defopara^r)$ induced by all $\defopara^a \colon \pi_{t - s, t + a} (X / \defopara^{r - a}) \to \pi_{t - s, t}(X / \defopara^r)$
    \[0 = \IIm(\defopara^r) \subset \IIm(\defopara^{r - 1}) \subset \cdots \subset \IIm(\defopara^2) \subset \IIm(\defopara^1) \subset \IIm(\defopara^0) = \pi_{t - s, t}(X / \defopara^r)\]
    for which $\IIm(\defopara^a) / \IIm(\defopara^{a + 1}) \cong Z_{r - a}^{s + a, t + a} / B_{a + 1}^{s + a, t + a}$. In particular: 
    \begin{itemize}
        \item Every $x \in Z_{r}^{s, t}(X)$ admits a lift $[x] \in \pi_{t - s, t}(X / \defopara^r)$. Conversely, any nonzero bigraded stem $\beta \in \pi_{t - s, t} (X / \defopara^r)$ lies in some $\IIm(\defopara^a) \setminus \IIm(\defopara^{a + 1})$, and there is a unique nonzero $x \in Z_{r - a}^{s + a, t + a} / B_{a + 1}^{s + a, t + a}$ so that $\beta = \defopara^{a}[x]$ for a certain lift $[x] \in \pi_{t - s, t + a}(X / \defopara^{r - a})$. 
        \item Suppose $\{x_{a, j}\}_{j \in J_a}$ is a family of generators of $Z_{r - a}^{s + a, t + a} / B_{a + 1}^{s + a, t + a}$ for each $0 \leq a \leq r - 1$. Choose a lift $[x_{a, j}]_{\std} \in \pi_{t - s, t + a}(X / \defopara^{r - a})$ for each $x_{a, j}$. Then every $\beta \in \pi_{t - s, t}(X / \defopara^r)$ can be presented as $\sum_{0 \leq a < r} \sum_{j \in J_a} \defopara^a n_{a, j}[x_{a, j}]_{\std}$ for a certain family of integers $\{n_{a, j}\}_{0 \leq a < r, j \in J_a}$. 
        \item Suppose $y \in Z_c^{s, t}$ for $c \geq r$ and $y \neq 0$ until the $E_{k + 1}$-page for some $k \leq c$. Then for every lift $[y] \in \pi_{t - s, t}(X / \defopara^r)$ and any $0 \leq a < k$ the class $\defopara^a [y] \in \pi_{t - s, t + a}(X / \defopara^{r + a})$ is nonzero. If furthermore we have $d_{k + 1}(x) = y$, then there is a specific lift $[y]_0 = \delta_{k}^{r}[x] \in \pi_{t - s, t}(X / \defopara^r)$ for which $\defopara^k [y]_0 = 0$, while for a generic lift $[y]$ we only know that $\defopara^k[y] \in \IIm(\defopara^{k + 1}) \subset \pi_{t - s, t + k}(X / \defopara^{r + k})$. 
    \end{itemize}
\end{theorem}

\begin{proof}
    The items after ``in particular'' follow directly from the main body. To establish the main body statement, for $0 \leq a < r < \infty$ we take the cartesian square from Lemma \ref{defopara-rho-cartesian}
    \[\begin{tikzcd}
        \Sigma^{0, -a} X / \defopara^{r - a} \ar[r, "{\defopara^a}"] \ar[d, "{\rho}"] & X / \defopara^r \ar[d, "{\rho}"] \\
        \Sigma^{0, -a} X / \defopara \ar[r, "{\defopara^a}"] & X / \defopara^{a + 1}
    \end{tikzcd}\]
    and apply Lemma \ref{cartesian-square-composite-image}. 
    It follows that on $\pi_{t - s, t}$ we get $\IIm(\defopara^a) / \IIm(\defopara^{a + 1}) \cong \IIm(\rho) / \IIm(\rho \delta_a) = Z_{r - a}^{s + a, t + a} / B_{a + 1}^{s + a, t + a}$ as expected.
\end{proof}

Similarly, we can translate back and forth between $\pi_{**}(X)$ and the \emph{permanent cycles} in the standard SS. To phrase this we shall introduce some notions about $E_{\infty}$-page and (strong) convergence.

\begin{definition}
    \begin{itemize}
        \item Suppose $E = \{E_r^{s, t}\}$ is a SS. We denote by $Z_{\infty}^{s, t} = \lim_r Z_r^{s, t} \subset E_2^{s, t}$ the group of \textbf{permanent cycles}, and $B_{\infty}^{s, t} = \colim_r B_r^{s, t} \subset E_2^{s, t}$ the group of all boundaries. We refer to the collection of all quotients $E_{\infty}^{s, t} = Z_{\infty}^{s, t} / B_{\infty}^{s, t}$ as the \textbf{$E_{\infty}$-page} of this SS. 
        \item Suppose $F^* A = \{F^w \! A\}_{w \in \Zb}$ is a $\Zb$-indexed strict\footnote{This means each $F^w A$ is a subgroup of $A$.} decreasing filtration on $A \in \Ab$. 
        \begin{itemize}
            \item We say $F^* A$ is \textbf{exhaustive} if $\bigcup_w F^w \!  A = A$. 
            \item We say $F^* A$ is \textbf{Hausdorff} if the comparison map $A \to \lim_w A / F^w \! A$ is injective, where the limit is computed in $\Ab$. In other words, $F^* A$ is Hausdorff if $\bigcap_w F^w \!  A = 0$. 
            \item Similarly, we say $F^* A$ is \textbf{precomplete}\footnote{We can treat strictly filtered abelian groups as $\defopara$-torsionfree objects in $\Mod_{\Zb[\defopara]} \Gr\Ab$. One may wonder what is the relationship between derived $\defopara$-completeness there and the notion of precompleteness above. Actually, the $\defopara$-torsionfree object $F^* A$ is derived $\defopara$-complete iff it is classically $\defopara$-complete, iff $F^* A$ is both precomplete and Hausdorff. For instance, if $F^* A$ is the constant filtration on $A = \Zb$, then it is precomplete but not derived $\defopara$-complete.} if the comparison map $A \to \lim_w A / F^w \! A$ is surjective. In other words, $F^* A$ is precomplete if $\lim^1_w F^w \! A = 0$. We say $F^* A$ is \textbf{complete} if it is both precomplete and Hausdorff.
        \end{itemize}
        \item Suppose $E = \{E_r^{s, t}\}$ is a SS and $A = \{A_n\} \in \Gr\Ab$. We say $E$ \textbf{converges (strongly)} to $A$, or $E_r^{s, t} \Rightarrow A_{t - s}$, if for each $n \in \Zb$ there is a complete exhaustive filtration $F^* A_n$ together with a family of isomorphisms $F^s \! A_n / F^{s + 1} \! A_n \cong E_{\infty}^{s, n + s}$. We also refer to $A$ as the \textbf{abutment} of $E$. 
    \end{itemize}
\end{definition}

\begin{theorem}[Bockstein dictionary, infinite page]\label{Bockstein-dictionary-infinite}
    Take $X \in \Fil\Sp$. 
    \begin{enumerate}
        \item The following strict filtration on $\pi_{t - s, t}(X)$ induced by all $\defopara^a \colon \pi_{t - s, t + a}(X) \to \pi_{t - s, t}(X)$
        \[\cdots \subset \IIm(\defopara^a) \subset \cdots \subset \IIm(\defopara^2) \subset \IIm(\defopara^1) \subset \IIm(\defopara^0) = \pi_{t - s, t}(X)\]
        is exhaustive with associated graded $\IIm(\defopara^a) / \IIm(\defopara^{a + 1}) \cong \IIm(\rho^{\infty}_1)/ B_{a + 1}^{s + a, t + a}$, where 
        \[\IIm(\rho^{\infty}_1 \colon \pi_{t - s, t + a}(X) \to \pi_{t - s, t + a}(X / \defopara)) \subset Z_{\infty}^{s + a, t + a}.\] 
        Suppose $X$ is moreover complete and for each bidegree $(s, t)$ the \textbf{weak obstruction} $RE_{\infty}^{s, t}(X) = \lim^1_{r} Z_r^{s, t}$ vanishes. Then the image of $\rho\colon \pi_{t - s, t}(X) \to \pi_{t - s, t}(X / \defopara) = E_2^{s, t}$ is precisely $Z_{\infty}^{s, t}$, and this filtration is also complete, with $\IIm(\defopara^a) / \IIm(\defopara^{a + 1}) \cong Z_{\infty}^{s + a, t + a} / B_{a + 1}^{s + a, t + a}$. In particular:
        \begin{itemize}
            \item Every $x \in Z_{\infty}^{s, t}(X)$ admits a lift $[x] \in \pi_{t - s, t}(X)$. Conversely, any nonzero bigraded stem $\beta \in \pi_{t - s, t} (X)$ lies in some $\IIm(\defopara^a) \setminus \IIm(\defopara^{a + 1})$, and there is a unique nonzero $x \in Z_{\infty}^{s + a, t + a} / B_{a + 1}^{s + a, t + a}$ so that $\beta = \defopara^{a}[x]$ for a certain lift $[x] \in \pi_{t - s, t + a}(X)$. 
            \item Suppose $\{x_{a, j}\}_{j \in J_a}$ is a family of generators of $Z_{\infty}^{s + a, t + a} / B_{a + 1}^{s + a, t + a}$ for each $a \geq 0$. Choose a lift $[x_{a, j}]_{\std} \in \pi_{t - s, t + a}(X)$ for each $x_{a, j}$. Then for any levelwise finitely supported family of integers $\{n_{a, j}\}_{a \geq 0, j \in J_a}$, the expression $\sum_{a \geq 0} \sum_{j \in J_a} \defopara^a n_{a, j}[x_{a, j}]_{\std}$ determines a unique element in $\pi_{t - s, t}(X)$. Conversely, every $\beta \in \pi_{t - s, t}(X)$ can be presented in this form.
            \item If $y \in E_2^{s, t}$ survives to the $E_{k + 1}$-page and is killed by a $d_{k + 1}$-differential, then for every lift $[y] \in \pi_{t - s, t}(X)$ and any $0 \leq a < k$ the class $\defopara^a [y] \in \pi_{t - s, t + a}(X)$ is nonzero. Furthermore, if $d_{k + 1}(x) = y$, then there is a specific lift $[y]_0 = \delta_{k}^{\infty}[x] \in \pi_{t - s, t}(X)$ for which $\defopara^k [y]_0 = 0$, while for a generic lift $[y]$ we only know that $\defopara^k[y] \in \IIm(\defopara^{k + 1}) \subset \pi_{t - s, t + k}(X)$. 
        \end{itemize}
        \item Under the assumption of item 1, for each $n \in \Zb$, we further consider the filtration on $\pi_n(X[\defopara^{-1}])$
        \[\cdots \subset \IIm(\iota_2) \subset \IIm(\iota_1) \subset \IIm(\iota_{0}) \subset \IIm(\iota_{-1}) \subset \IIm(\iota_s) \subset \cdots\]
        with $\IIm(\iota_s)$ being the image of $\iota_s\colon \pi_{n, n + s} (X) \to \pi_{n}(X[\defopara^{-1}])$. This is precomplete exhaustive, so that $\IIm(\iota_s) / \IIm(\iota_{s + 1})\cong E_{\infty}^{s, n + s}$. If furthermore for each $n \in \Zb$ the \textbf{whole plane obstruction}
        \[W_n(X) = \ker\left(\colim_{a \to -\infty} \lim_{c \to \infty} \pi_{n, a}(X / \defopara^{c - a}) \to \lim_{c \to \infty} \colim_{a \to -\infty} \pi_{n, a}(X / \defopara^{c - a})\right)\]
        vanishes, then these filtrations are also Hausdorff, thus $E_r^{s, t}(X) \Rightarrow \pi_{t - s} (X[\defopara^{-1}])$. In particular: 
        \begin{itemize}
            \item Any nonzero stem $\beta \in \pi_{n} (X[\defopara^{-1}])$ lies in some $\IIm(\iota_s) \setminus \IIm(\iota_{s + 1})$, and there is a unique nonzero $x \in E_{\infty}^{s, n + s}$ so that a certain lift $[x] \in \pi_{n, n + s}(X)$ recovers $\beta$ after $\defopara$-inversion. In this case we say $\beta$ is of filtration $s$ and $x \in E_{\infty}^{s, n + s}$ \textbf{detects} $\beta  \in \pi_{n} (X[\defopara^{-1}])$.
            \item Suppose $\{x_{s, j}\}_{j \in J_s}$ is a family of generators of $E_{\infty}^{s, n + s}$ for each $s \in \Zb$. Choose a lift $[x_{s, j}]_{\std} \in \pi_{n, n + s}(X)$ for each $x_{s, j}$. Then for any $w \in \Zb$ and any levelwise finitely supported family of integers $\{k_{s, j}\}_{s \geq w, j \in J_s}$, the expression $\sum_{s \geq w} \sum_{j \in J_s} \iota_s (k_{s, j}[x_{s, j}]_{\std})$ determines a unique element in $\pi_{n} (X[\defopara^{-1}])$. Conversely, every $\beta \in \pi_{n} (X[\defopara^{-1}])$ can be presented in this form.
        \end{itemize}
    \end{enumerate}
\end{theorem}

\begin{proof}
    The items after ``in particular'' follow directly from the main body. For the main body of part 1, we take the cartesian square from Lemma \ref{defopara-rho-cartesian}
    \[\begin{tikzcd}
        \Sigma^{0, -a} X \ar[r, "{\defopara^a}"] \ar[d, "{\rho}"] & X \ar[d, "{\rho}"] \\
        \Sigma^{0, -a} X / \defopara \ar[r, "{\defopara^a}"] & X / \defopara^{a + 1}
    \end{tikzcd}\]
    and apply Lemma \ref{cartesian-square-composite-image} to this square. On $\pi_{t - s, t}$ we get $\IIm(\defopara^a) / \IIm(\defopara^{a + 1}) \cong \IIm(\rho) / \IIm(\rho \delta_a) = \IIm(\rho) / B_{a + 1}^{s + a, t + a}$, here $\IIm(\rho\colon \pi_{t - s, t + a}(X) \to \pi_{t - s, t + a}(X / \defopara))$ is a subgroup of $Z_{\infty}^{s + a, t + a}$. It remains to prove that if $X$ is complete and all $RE_{\infty}^{s, t}$ vanish, then the map $\IIm(\rho) \to Z_{\infty}^{s, t}$ is surjective, $\lim_{a} \IIm(\defopara^a) = 0$ and $\lim\nolimits^1_a \IIm(\defopara^a) = 0$ for each bidegree $(s, t)$. In fact, the second assertion implies the first: Suppose there is $x \in Z_{\infty}^{s, t} \subset \pi_{t - s, t}(X / \defopara)$ that does not lie in $\IIm(\rho)$, then $\delta_1^{\infty} (x)$ is nonzero in $\pi_{t - s + 1, t + 1}(X)$, so $\lim_{a} \IIm(\defopara^a) = 0$ implies $\delta_1^{\infty} (x) \not\in \IIm(\defopara^a)$ for $a \gg 1$, but we can take a lift $[x] \in \pi_{t - s, t}(X / \defopara^a)$ and obtain the contradictory equality $\delta_1^{\infty} (x) = \defopara^a \delta_a^{\infty}[x]$ via Corollary \ref{defopara-totaldiff}. \parr 

    To proceed, we translate everything to the setup in \cite{Boa99} as pointed out in Remark \ref{Boardman-comparison}. In this language, $\pi_{t - s, t}(X)$ is a component of the graded abelian group $A^s$ and $\lim_a \IIm(\defopara^a)_{s, t}$ (resp. $\lim^1_a \IIm(\defopara^a)_{s, t}$) is a component of the graded abelian group $Q^s$ (resp. $R Q^s$). Since $X$ is complete, a Milnor sequence argument shows that $\lim_s A^s = 0$ and $\lim^1_s A^s = 0$. If further $RE_{\infty}^{s, t} = 0$ for each $(s, t) \in \Zb^2$, then $Q^s = 0$ and $RQ^s = 0$ for each $s \in \Zb$ due to \cite[Lemma 5.7 and Lemma 5.11]{Boa99}. Thus, we have established the claims in part 1. \parr 

    For the main body of part 2, we consider the cartesian square 
    \[\begin{tikzcd}
        X \ar[r] \ar[d, "{\rho}"] & X[\defopara^{-1}] \ar[d] \\
        X / \defopara \ar[r] & \colim\limits_{w \to \infty} \Sigma^{0, w} X / \defopara^{w + 1}
    \end{tikzcd} = \quad \colim_{w \to \infty}\left[
    \begin{tikzcd}
        X \ar[r, "{\defopara^w}"] \ar[d, "{\rho}"] & \Sigma^{0, w} X \ar[d, "{\rho}"] \\
        X / \defopara \ar[r, "{\defopara^w}"] & \Sigma^{0, w} X / \defopara^{w + 1}
    \end{tikzcd}\right]\]
    and apply Lemma \ref{cartesian-square-composite-image} to this square on $\pi_n$. It follows that $\IIm(\iota_s) / \IIm(\iota_{s + 1})$ is isomorphic to the quotient of $\IIm(\rho) = Z_{\infty}^{s, n + s}$ by the image of $\pi_{n, n + s} (\fib(X \to X[\defopara^{-1}])) \to \pi_{n, n + s}(X) \to \pi_{n, n + s}(X / \defopara)$. Since $\fib(X \to X[\defopara^{-1}]) \cong \colim_{w \to {\infty}} \fib(\defopara^w \colon X \to \Sigma^{0, w} X) \cong \colim_{w \to \infty} \Sigma^{-1, w} X / \defopara^w$ and each composite $\Sigma^{-1, w} X / \defopara^w \to X \to X / \defopara$ is homotopic to $\delta_w$ up to sign, this image can be identified with $\colim_{w \to \infty} B_w^{s, n + s} = B_{\infty}^{s, n 
    + s}$. In other words, $\IIm(\iota_s) / \IIm(\iota_{s + 1}) \cong Z_{\infty}^{s, n + s} / B_{\infty}^{s, n 
    + s} = E_{\infty}^{s, n + s}$. The filtration is exhaustive as $\pi_{n}(X[\defopara^{-1}]) \cong \colim_{s \to \infty} \pi_{n, n + s}(X)$, which follows from compactness of $\Sb^{n, n}$. It remains to show $\lim^1_s \IIm(\iota_s) = 0$, and $\lim_s \IIm(\iota_s) = 0$ if we further have $W_n = 0$ for each $n \in \Zb$. For these we switch to the language in \cite{Boa99}. Here each $\IIm(\iota_s)_n$ is a component of the graded abelian group $F^s$. As we already have $Q^s = 0 = RQ^s$ for each $s$, the exact sequence in \cite[Lemma 8.11]{Boa99} shows that $\lim^1_{s} F^s = 0$ and $\lim_s F^s \cong W$ for the graded abelian group $W$ in \cite[Lemma 8.5]{Boa99}. The $n$-th component of $W$ takes the stated form due to \cite[Theorem 7.5]{HR19}, which concludes part 2.
\end{proof}

\begin{remark} \label{vanishing-of-the-obstructions-to-strong-convergence}
    \begin{itemize}
        \item For $(s, t) \in \Zb^2$, the weak obstruction $RE_{\infty}^{s, t} = 0$ if the family $\{Z_{r}^{s, t}\}_{r \geq 1}$ is pro-constant, in particular if there exists some $K \gg 1$ so that $Z_r^{s, t} = Z_{K}^{s, t}$ for all $r \geq K$ (i.e. there are only finitely many pages for which there exists a nontrivial differential out of bidegree $(s, t)$). 
        \item For $n \in \Zb$, the whole plane obstruction $W_n = 0$ if there exists $K \gg 1$ so that for any $s_1, s_2 > K$, $d_{s_1 + s_2}\colon E_{s_1 + s_2}^{-s_2, n + 1 - s_2} \to E_{s_1 + s_2}^{s_1, n + s_1}$ is the zero map (i.e. there does not exist an infinite family of nontrivial pairwise crossing differentials with increasing length hitting stem $n$), cf. \cite[Lemma 8.1]{Boa99}, \cite[Proposition 5.3 and Remark 5.4]{HR19}.
    \end{itemize}
\end{remark}

\begin{summary} \label{Bockstein-dictionary-summary}
    We pause to summarize the behavior of \(\defopara^k\), \(\rho\), and \(\delta_k\) in terms of the standard SS.

    \begin{itemize}
        \item The bigraded homotopy group \(\pi_{n,n+s}(X/\defopara^r)\) is built out of \(r\) layers:
        \begin{itemize}
            \item the \(0\)th layer consists of classes \([x]\) coming from \(Z_r\)-cycles in bidegree \((s,n+s)\);
            \item the \(1\)st layer consists of classes \(\defopara[x]\) coming from \(Z_{r-1}\)-cycles in bidegree \((s+1,n+s+1)\), modulo the \(d_2\)-boundaries;
            \item the \(2\)nd layer consists of classes \(\defopara^2[x]\) coming from \(Z_{r-2}\)-cycles in bidegree \((s+2,n+s+2)\), modulo the \(d_3\)-boundaries;
            \item \(\cdots\)
            \item the \((r-1)\)st layer consists of classes \(\defopara^{r-1}[x]\) coming from \(Z_1\)-cycles \(=E_2\)-classes in bidegree \((s+r-1,n+s+r-1)\), modulo the \(d_r\)-boundaries.
        \end{itemize}
    
        \item Multiplication by \(\defopara^k\) raises the layer by \(k\) and kills \(k\) additional pages of boundaries.
    
        \item The projection operator \(\rho^r_k\) removes all layers of order at least \(k\).
    
        \item The operator \(\delta_k^r\colon \pi_{**}(X/\defopara^k)\to \pi_{**}(X/\defopara^r)\) behaves like a \emph{total differential} combining information from \(d_2\) through \(d_{k+r}\), in the sense that $\delta_k^r \defopara^a [x] = \defopara^b [y]$ if and only if $d_{k - a + b + 1}(x) = y$. This is to be made precise below in Theorem \ref{delta-as-total-diff}.
    \end{itemize}
    
\end{summary}

\begin{lemma} \label{pi*-pullback-to-pullback-pi*}
   Suppose $\CC$ is a stable $\infty$-category. For any cospan $A \to C \gets B$ in $\CC$ and any homological functor $F\colon \CC \to \CA$, the limit comparison map $F(A \times_C B) \to F(A) \times_{F(C)} F(B)$ is surjective, and its kernel is given by the cokernel of $F(\Sigma^{-1}A) \oplus F(\Sigma^{-1}B) \to F(\Sigma^{-1}C)$. 
\end{lemma}

\begin{proof}
    The fiber sequence $A \times_C B \to A \oplus B \to C$ leads to a long exact sequence 
    \[\cdots \to F(\Sigma^{-1}A) \oplus F(\Sigma^{-1}B) \to F(\Sigma^{-1}C) \to F(A \times_C B) \to F(A) \oplus F(B) \to F(C) \to \cdots\]
    The desired result follows by breaking this into short exact sequences. 
\end{proof}

\begin{theorem} \label{delta-as-total-diff}
    For $1 \leq k < \infty$ and $1 \leq r \leq \infty$, consider $\delta_{k}^r\colon \pi_{n + 1, n + s + 1}(X / \defopara^k) \to \pi_{n, n + s + k}(X / \defopara^r)$. 
    \begin{itemize}
        \item Suppose $0 \leq a < k$, $0 \leq b < r$, $x \in E_2^{s, n + s + 1}(X)$ and $y \in E_2^{s + k - a + b + 1, n + s + k - a + b + 1}(X)$. Then there exist lifts $[x] \in \pi_{n + 1, n + s + 1}(X / \defopara^{k - a}), [y] \in \pi_{n, n + s + k - a + b + 1}(X / \defopara^{r - b})$ with $\delta_k^r \defopara^a [x] = \defopara^b[y]$ if and only if $x, y$ are both $(k - a + b)$-cycles and $d_{k - a + b + 1}(x) = y$ in the standard SS. 
        \item Suppose $0 \leq a < k$ and $x \in E_2^{s, n + s + 1}(X)$. For $r < \infty$, there exists a lift $[x] \in \pi_{n + 1, n + s + 1}(X / \defopara^{k - a})$ such that $\delta_k^r \defopara^a[x] = 0$ if and only if $x$ is an $(k - a + r)$-cycle. For $r = \infty$, the ``only if'' direction also holds true, while the ``if'' direction is true provided that $X$ is complete and all weak obstructions $RE_{\infty}^{s, t}$ (cf. Theorem \ref{Bockstein-dictionary-infinite}) vanish.
    \end{itemize}
\end{theorem}

\begin{proof} 
    For the ``only if'' direction in the first item, observe that $\delta_k^r\defopara^a[x] = \delta_{k - a}^r[x]$ by Corollary \ref{defopara-totaldiff}, and $\delta_{k - a}^r [x] = \defopara^b[y]$ implies by Lemma \ref{projection-totaldiff-cartesian} and Lemma \ref{pi*-pullback-to-pullback-pi*} that there is a lift $[x]' \in \pi_{n + 1, n + s + 1}(X / \defopara^{k - a + b})$ refining $[x]$ with $\delta_{k - a + b}^1[x]' = y$, thus $d_{k - a + b + 1}(x) = y$ in the standard spectral sequence. For the ``if'' direction, $d_{k - a + b + 1}(x) = y$ means there is a lift $[x]_1 \in \pi_{n + 1, n + s + 1}(X / \defopara^{k - a + b})$ so that $\delta_{k - a + b}[x]_1 = y_1$, with $y - y_1 \in B_{k - a + b}^{s + k -a + b + 1, n + s + k -a + b + 1}$. By construction, 
    \[\delta_{k - a + b - 1}^1\colon \pi_{n + 1, n + s + 2}(X / \defopara^{k - a + b - 1}) \to B_{k - a + b}^{s + k -a + b + 1, n + s + k -a + b + 1}\]
    is surjective, so we can find $[z] \in \pi_{n + 1, n + s + 2}(X / \defopara^{k - a + b - 1})$ with $\delta_{k - a + b - 1}^1[z] = y - y_1$. Thus, taking $[x] = [x]_1 + \defopara[z]$, we have $\delta^1_{k - a + b}[x] = y$ on the $E_2$-page, therefore $[y] = \delta_{k - a + b}^{r - b}[x]$ is a lift of $y$ for which $\defopara^b[y] = \defopara^b\delta_{k - a + b}^{r - b}[x] = \delta_{k - a}^r[x] = \delta_k^r \defopara^a [x]$ by reversing the argument for the ``only if'' direction. The second item with $r < \infty$ follows from the first item applied to $a, b = r - 1, x, y = 0 \in E_2^{s + k - a + r, n + s + k - a + r}$. For $r = \infty$, the ``only if'' direction follows from the $r < \infty$ cases, while the ``if'' direction holds true as $x \in Z_{\infty}^{s, n + s + 1}$ together with completeness and vanishing weak obstructions yield the existence of a lift $[x]' \in \pi_{n + 1, n + s + 1}(X)$ for which $[x] = \rho_{k - a}^{\infty}[x]'$ satisfies $\delta_k^{\infty} \defopara^a[x] = \delta_{k - a}^{\infty} [x] = 0$.
\end{proof}

\begin{remark} \label{essential-diff}
    A differential $d_r(x) = y$ is \textbf{essential} if $y \neq 0$ on the $E_r$-page. In the language of total differential $\delta$ above, $d_{k - a + b + 1}(x) = y$ is essential iff the lift $\delta_k \defopara^a [x] = \defopara^b[y] \not\in \IIm(\defopara^{b + 1})$. 
\end{remark}

It is also worthwhile to single out the ``limit case'' behavior of the fundamental distinguished triangle, which appears in the proof of Theorem \ref{Bockstein-dictionary-infinite}. Heuristically, this serves as the unification of all possible differentials in the standard SS of a filtered spectrum $X$.

\begin{construction}  \label{limit-distinguished-triangle}
    For each $X \in \Fil\Sp$, consider the tensor product of $X$ with the octahedral distinguished triangle coming from the colimit of all consecutive maps 
    \[0 \to \oneb \xrightarrow{\defopara^k} \Sigma^{0, k} \oneb\] 
    as $k \in \Nb$ increases monotonically to $\infty$. The resulting triangle takes the form
    \[X \xrightarrow{\iota \, = \, \langle \defopara^{\infty} \rangle} X[\defopara^{-1}]\xrightarrow{\rho_{\infty}^{\infty}} \Sigma^{1, 0} X^{\defopara\textup{-tors}} = \colim_{k \to \infty} \Sigma^{0, k} X / \defopara^{k} \xrightarrow{\delta_{\infty}^{\infty}} \Sigma^{1, 0} X\]
    where $X[\defopara^{-1}] = \colim_{k \to {\infty}} \Sigma^{0, k} X $, and $X^{\defopara\textup{-tors}} = \colim_{k \to \infty} \Sigma^{-1, k} X / \defopara^{k}$ is the \textbf{$\defopara$-torsionization} of $X$ abstractly constructed as the value of $X$ under the left adjoint of the completion functor $\Fil\Sp \to \Fil\Sp, X \mapsto X_{\hat{\defopara}}$. We may also describe the three maps in explicit terms:
    \begin{itemize}
        \item The first map $\iota = \langle \defopara^{\infty}\rangle$ is the $\defopara$-inversion map, i.e. the colimit of all maps $\defopara^k\colon X \to \Sigma^{0, k} X$. Concretely, for any permanent cycle $x$, it takes $[x] \in \pi_{**}(X)$ to the stem it detects in the abutment $\pi_*(X[\defopara^{-1}])$. 
        \item The second map $\rho_{\infty}^{\infty}$ is the colimit of all projections $\rho^{\infty}_k\colon  \Sigma^{0, k} X \to \Sigma^{0, k} X / \defopara^k$. Concretely, in good cases each $\alpha \in \pi_*(X[\defopara^{-1}])$ comes from a surviving permanent cycle $x$ of filtration $\geq -k$ for some $k \geq 0$ in that a certain lift $[x]$ detects $\alpha$, while $\rho^{\infty}_{\infty} (\alpha)$ amounts to forgetting $x$ to a $k$-cycle.  
        \item The third map $\delta^{\infty}_{\infty}$ is the $\defopara$-torsionization map, i.e. the colimit of all maps $\delta_k^{\infty}\colon \Sigma^{0, k} X / \defopara^k \to \Sigma^{1, 0} X$. Concretely, it takes every class $[x] \in \pi_{**}(X / \defopara^k) \hookrightarrow \pi_{**}(X^{\defopara\textup{-tors}})$ to the total differential $\delta_k^{\infty}[x]$, which is independent of the choice of $k$ due to Corollary \ref{defopara-totaldiff}.
    \end{itemize}
\end{construction}

We wrap up this subsection by generalizing some of the results above to the situation where $\Sp$ is replaced by an arbitrary $t$-$\infty$-category $\CE$, in terms of the $\CE^\heartsuit\!$-valued standard spectral sequence from Corollary \ref{standard-SS-with-general-coefficients}.

\begin{corollary} \label{Bockstein-dictionary-finite-page-with-general-coefficients}
    Let $\CE$ be a presentable stable $\infty$-category, $\CE^\heartsuit$ be an abelian $1$-category and $\pi_0\colon \CE \to \CE^\heartsuit$ be a homological functor. Take $X \in \Fil(\CE)$. 
    \begin{itemize}
        \item For each $r \in \Nb, s, t \in \Zb$ with $r \geq 1$, there is a finite filtration on $\pi_{t - s, t}(X / \defopara^r) \in \CE^\heartsuit$ induced by all maps in $\CE^\heartsuit$ of the form $\defopara^a \colon \pi_{t - s, t + a} (X / \defopara^{r - a}) \to \pi_{t - s, t}(X / \defopara^r)$
        \[0 = \IIm(\defopara^r) \subset \IIm(\defopara^{r - 1}) \subset \cdots \subset \IIm(\defopara^2) \subset \IIm(\defopara^1) \subset \IIm(\defopara^0) = \pi_{t - s, t}(X / \defopara^r)\]
        with $\IIm(\defopara^a) / \IIm(\defopara^{a + 1}) \cong Z_{r - a}^{s + a, t + a} / B_{a + 1}^{s + a, t + a}$. This follows by the same proof as in Theorem \ref{Bockstein-dictionary-finite}. 
        \item For each $r, k, a, b \in \Nb, n, s \in \Zb$ with $a < k, b < r$, suppose there is an object $T \in \CE^\heartsuit$ with two maps $x\colon T \to E_2^{s, n + s + 1}, y\colon T \to E_2^{s + k - a + b + 1, n + s + k - a + b + 1}$. Then the following are equivalent: 
        \begin{itemize}
            \item There exist $[x]\colon T \to \pi_{n + 1, n + s + 1}(X / \defopara^{k - a})$ and $[y] \colon T \to \pi_{n, n + s + k - a + b + 1}(X / \defopara^{r - b})$ such that $\rho^{k - a}_1[x] = x, \rho^{r - b}_1[y] = y$ and $\delta_k^r \defopara^a [x] = \defopara^b [y]$ as maps in $\CE^\heartsuit$. 
            \item The maps $x, y$ lift (uniquely) to $x\colon T \to Z_{k - a + b}^{s, n + s + 1}$, $y\colon T \to Z_{k - a + b}^{s + k - a + b + 1, n + s + k - a + b + 1}$, so that
            \[\begin{tikzcd}
                T \ar[r, "{x}"] \ar[d, "{y}"] & Z_{k - a + b}^{s, n + s + 1} \ar[r] & E_{k - a + b + 1}^{s, n + s + 1} \ar[d, "{d_{k - a + b + 1}}"] \\
                Z_{k - a + b}^{s + k - a + b + 1, n + s + k - a + b + 1} \ar[rr]& & E_{k - a + b + 1}^{s + k - a + b + 1, n + s + k - a + b + 1}
            \end{tikzcd}\]
            is a commutative square in $\CE^\heartsuit$.
        \end{itemize}
        This follows by the same proof as in Theorem \ref{delta-as-total-diff}. 
        \item On the other hand, the $r = \infty$ case of our Bockstein dictionary does not have a direct generalization. Indeed, suppose $\CE^\heartsuit$ admits sequential limits, exact sequential colimits, and the functor $\pi_0$ preserves sequential colimits. Then for each $s, t \in \Zb$, we can still construct the exhaustive strict filtration
        on $\pi_{t - s, t}(X)$
        \[\cdots \subset \IIm(\defopara^a) \subset \cdots \subset \IIm(\defopara^2) \subset \IIm(\defopara^1) \subset \IIm(\defopara^0) = \pi_{t - s, t}(X)\]
        induced by all maps in $\CE^\heartsuit$ of the form $\defopara^a \colon \pi_{t - s, t + a}(X) \to \pi_{t - s, t}(X)$, while the identification $\IIm(\defopara^a) / \IIm(\defopara^{a + 1}) \cong \IIm(\rho^{\infty}_1)/ B_{a + 1}^{s + a, t + a}$ with
        \[\IIm(\rho^{\infty}_1 \colon \pi_{t - s, t + a}(X) \to \pi_{t - s, t + a}(X / \defopara)) \subset Z_{\infty}^{s + a, t + a}\]
        still holds true due to the same proof as in Theorem \ref{Bockstein-dictionary-infinite}. However, to further study the completeness of this filtration, and subsequently the strong convergence property of the $\CE^\heartsuit \!$-valued standard SS, would require a detailed analysis of certain obstruction classes ($RE^{s, t}_{\infty}, W_n$, etc.) coming from the higher derived functors of the sequential limit functor in $\CE^\heartsuit$. In the case $\CE^\heartsuit = \Ab$, all higher derived limits vanish except for $\lim^1$, and the required analysis is achieved in \cite{Boa99}. For general $\CE^\heartsuit$, however, the contributions from $\lim^2$, etc. introduce substantial additional obstructions, and we do not pursue the corresponding analysis here.
    \end{itemize}
\end{corollary}

\begin{remark} \label{credits-for-theorem-A}
    We end this subsection with some historical comments and attributions concerning the Bockstein dictionary results proved above.
    \begin{itemize} 
        \item Historically, people usually define the standard SS of a filtered spectrum via approaches that do not involve bigraded stems, e.g. the exact couple formalism in \cite{Boa99}. Consequently, they tend to understand the bigraded stems of $X$ (or $X / \defopara^r$) as the abutment of the trigraded $\defopara$-Bockstein SS associated to the $\defopara$-filtration on the object $X$ (or $X / \defopara^r$) in $\Fil\Sp$. In this way, 
        \begin{enumerate}
            \item[(a).] The content of Theorem \ref{Bockstein-dictionary-infinite} becomes that the trigraded $\defopara$-Bockstein SS of $X$ is \textbf{rigid}, i.e. it admits the action of a formal parameter $\defopara$ and is completely determined by $E_*^{*,*}(X) = E_*^{*,*,*}(X)[\defopara^{-1}]$. See Theorem \ref{rigidity-for-trigraded-defopara-Bockstein-SS} item 1. 
            \item[(b).] The content of Theorem \ref{Bockstein-dictionary-finite} becomes that the trigraded $\defopara$-Bockstein SS of $X / \defopara^r$ is \textbf{rigid}, or more precisely the translation rule between this SS and the standard SS $E_*^{*,*}(X)$. See Theorem \ref{rigidity-for-trigraded-defopara-Bockstein-SS} item 3. 
        \end{enumerate}
        Actually, as we will discuss in Appendix \ref{app:A}, all examples of rigid trigraded SS come from $\defopara$-Bockstein filtrations, so the history of the rigidity results for SS is the same as that of the Bockstein dictionary type results for filtered objects. 
        \item The first rigidity result of type (a) is due to Hu--Kriz--Ormsby \cite[Lemma 15]{Hu-Kriz-Ormsby}, which proves the rigidity of motivic Adams--Novikov SS (mANSS) of $\oneb \in \SH(\Cb)$, i.e. the complete translation between data in this mANSS and the data in the classical ANSS of $\Sb^0$. 
        \item After this, Isaksen \cite[Proposition 6.2.5]{Isaksen-stable-stems} proves the rigidity of the mANSS of $\oneb / \tau \in \SH(\Cb)$, i.e. this SS has no nontrivial differential and the abutment $\pi_{**} (\oneb / \tau) $ is the classical ANSS $E_2$-page of $\Sb^0$. This is the first rigidity result of type (b).
        \item Burklund--Hahn--Senger \cite[Theorem A.8]{BHS1} proves the rigidity of $E$-synthetic ASS of $\nu X$. This generalizes the $\Cb$-motivic rigidity results since Pstragowski \cite[Theorem 7.34]{Pst23} identifies the $\infty$-category of ($p$-complete cellular) $\Cb$-motivic spectra with the full sub-$\infty$-category of ($p$-complete) $\MU$-synthetic spectra $\Syn_{\MU}$ generated under colimits by the even bigraded spheres
        $\{\Sb^{n, 2w}\}_{n, w \in \Zb}$.
        
        Furthermore, they convert (for the first time) the rigidity result into a type (a) Bockstein dictionary statement for bigraded stems \cite[Theorem A.1]{BHS1}, under the name of the ``omnibus theorem''.  They also briefly mention the truncated case (i.e. $E$-synthetic ASS of $\nu X / \defopara^r$) in \cite[Corollary A.11]{BHS1}. 
        \item Burklund--Xu \cite[Theorem 2.8]{Burklund-Xu} proves the rigidity of the motivic Cartan--Eilenberg SS. 
        \item Carrick--Davies--van Nigtevecht \cite[\S~2]{CDvN2}, \cite[\S~3]{vN25} establish the ``omnibus theorem'' and its truncated variant in filtered spectra by adapting the $\defopara$-Bockstein SS argument in \cite{BHS1} to $\Fil\Sp$. The results they obtain, namely \cite[Theorem 3.62, Theorem 3.67 and Theorem 3.70]{vN25}, are essentially equivalent to our Theorem \ref{Bockstein-dictionary-finite} and Theorem \ref{Bockstein-dictionary-infinite}. 
        
        Note that the present paper approaches these results from a different perspective. Rather than treating the bigraded stems of $X$ and $X / \defopara^r$ as auxiliary tools for understanding the standard SS of $X$, we regard them as part of the foundational data of the standard SS itself. Accordingly, our proof of the Bockstein dictionary results does not use spectral sequence arguments, but reduces instead to a simple lemma on homological functors, Lemma \ref{cartesian-square-composite-image}.
    \end{itemize}
    In the rigidity story, one remaining missing piece is the identification of the synthetic Adams SS of $\nu X / \defopara^r$ with a truncated variant of the trigraded $\defopara$-Bockstein SS. In Appendix \ref{app:A}, we give a more extensive review of the rigidity results and fill this gap in Theorem \ref{synthetic-Adams-vs-Bockstein}.
\end{remark}

\subsection{Multiplicativity}
\label{subsec:2.3}

Next we discuss the multiplicativity of the standard SS. Recall that  $\oneb / \defopara$ is an $\Eb_{\infty}$-algebra in $\Fil\Sp$. 

\begin{construction}[Accelerations] \label{accelerations}
    Take $k \geq 1$. 
    \begin{itemize}
        \item The map of abelian groups $\Zb \to \Zb, t \mapsto kt$ is compatible with the poset structures, so it left Kan extends to a symmetric monoidal left adjoint $\Sl_k\colon \Fil\Sp \to \Fil\Sp$ that sends $\Sb^{n, w}$ to $\Sb^{n, kw}$. For general $X$, $\Sl_k(X)(w) = X(\lceil w/k \rceil)$, so we refer to it as ``\textbf{$k$-fold slowdown}'' of the filtration $X$. Concretely, $\pi_{n, w} (\Sl_k X) = \pi_{n, \lfloor w / k \rfloor} X$.  Note that $\Sl_k$ sends $\defopara\colon \Sb^{0, -1} \to \Sb^{0, 0}$ to $\defopara^k\colon \Sb^{0, -k} \to \Sb^{0, 0}$, so $\Sl_k(\oneb / \defopara) = \oneb / \defopara^k$ also admits an $\Eb_{\infty}$-algebra structure.  
        \item The corresponding right adjoint $\Ac_k$ (which is automatically lax symmetric monoidal) sends $X \in \Fil\Sp$ to $\Ac_k X$ with $(\Ac_k X)(w) = \Msp(\Sb^{0, w}, \Ac_k X) = \Msp(\Sl_k(\Sb^{0, w}), X) = X(kw)$, so it amounts to a ``\textbf{$k$-fold acceleration}'' of the filtration $X$. Concretely, $\pi_{n, w} (\Ac_k X) = \pi_{n, kw} X$. It follows that $\Ac_k (\oneb / \defopara^k) \cong \oneb / \defopara$ and $\Ac_k(X / \defopara^k) \cong (\Ac_k X) / \defopara$ as $\oneb / \defopara$-modules for each $X \in \Fil\Sp$. Furthermore, the image of the total differential $\delta_k^k\colon X / \defopara^k \to \Sigma^{1, -k} X / \defopara^k$ under $\Ac_k$ turns out to be $\delta_1^1$ on $(\Ac_k X) / \defopara$. 
    \end{itemize}
\end{construction}

\begin{remark}
    Both $\oneb$ and $\oneb / \defopara$ lie in the full subcategory $\Fil_{\leq 0} \Sp \subset \Fil\Sp$ generated under colimits by $\{\Sb^{n, w}\}_{n \in \Zb, w \leq 0}$, for which we have an identification $\Fil_{\leq 0} \Sp \cong \Fun(\overrightarrow{\Zb}^{\op}_{\leq 0}, \Sp)$ where $\overrightarrow{\Zb}_{\leq 0}$ is the poset of nonpositive integers. We can reproduce the construction of $\Sl_k$ (and $\Ac_k$) in $\Fil_{\leq 0} \Sp$, which leads to a commutative square of symmetric monoidal left adjoints
    \[\begin{tikzcd}
        \Fil_{\leq 0} \Sp \ar[r, "{\Sl_k}"] \ar[d, "{\mathrm{incl}}"] & \Fil_{\leq 0} \Sp \ar[d, "{\mathrm{incl}}"] \\
        \Fil \Sp \ar[r, "{\Sl_k}"] & \Fil \Sp 
    \end{tikzcd}\]
    while we have an extra advantage in $\Fil_{\leq 0} \Sp$: among the symmetric monoidal functors (i.e. maps of posets) $\cdot k \colon \overrightarrow{\Zb}_{\leq 0} \to \overrightarrow{\Zb}_{\leq 0}$ there are symmetric monoidal natural transformations $\cdot (k + 1) \to \cdot {k}$ for all $k \geq 1$, and they extend to symmetric monoidal transformations $\Sl_{k + 1} \to \Sl_k$ by the universal property of spectral presheaves. This refines the diagram
    \[\cdots \xrightarrow{\rho^5_4} \Sl_4(\oneb / \defopara) =  \oneb / \defopara^4 \xrightarrow{\rho^4_3} \Sl_3(\oneb / \defopara) =  \oneb / \defopara^3 \xrightarrow{\rho^3_2} \Sl_2(\oneb / \defopara) =  \oneb / \defopara^2 \xrightarrow{\rho^2_1} \Sl_1(\oneb / \defopara) = \oneb / \defopara\]
    to a sequence of $\Eb_{\infty}$-ring maps, whose limit recovers $\oneb = \oneb_{\hat{\defopara}}$ as an $\Eb_{\infty}$-algebra.
\end{remark} 

\begin{recollection}[Koszul sign rule]
    Take $n \in \Nb_{\geq 1} \cup \{\infty\}$. 
    \begin{itemize}
        \item Suppose $\CA$ is an additively $\Eb_n$-monoidal $1$-category. Then there is a \textbf{Koszul $\Eb_n$-monoidal structure} on $\Gr(\CA)$. For $n = 1$ this coincides with the Day convolution monoidal structure. In general, the difference between $\Gr(\CA)^{\Day}$ and $\Gr(\CA)^{\Kos}$ lies solely in the swap map. As every object $X = \{X^w\}_{w \in \Zb} \in \Gr(\CA)$ is the direct sum of each $X^w[w]$, the object $X^w \in \CA$ concentrated in degree $w$, it suffices to treat the swap map intertwining $X[m_1]$ and $Y[m_2]$ for $m_1, m_2 \in \Zb^\delta$ and $X, Y \in \CA$. In fact, if $\epsilon = \epsilon_{X, Y}\colon X \otimes Y \to Y \otimes X$ is the swap map in $\CA$, then 
        \[\epsilon[m_1 + m_2]\colon X[m_1] \otimes Y[m_2] = (X \otimes Y)[m_1 + m_2] \to (Y \otimes X)[m_1 + m_2] = Y[m_2] \otimes X[m_1]\]
        is the swap map intertwining $X[m_1]$ and $Y[m_2]$ in $\Gr(\CA)^{\Day}$, while 
        \[(-1)^{m_1m_2}\epsilon[m_1 + m_2]\colon X[m_1] \otimes Y[m_2] \to Y[m_2] \otimes X[m_1]\] 
        is the swap map intertwining $X[m_1]$ and $Y[m_2]$ in $\Gr(\CA)^{\Kos}$. 
        \item There is a unique \textbf{Koszul $\Eb_n$-monoidal structure} on $\Ch(\CA)$ so that the forgetful functor $U\colon \Ch(\CA) \to \Gr(\CA), (\{M^s\}_{s \in \Zb}, d\colon M^s \to M^{s + 1}) \mapsto \{M^s\}_{s \in \Zb}$ is $\Eb_n$-monoidal with Koszul $\Eb_n$-monoidal structure on the target. Concretely, 
        \begin{itemize}
            \item For $M, N \in \Ch(\CA)$, the tensor product $M \otimes N$ is given on the underlying object by $(M \otimes N)^s = \bigoplus_{p + q = s} M^p \otimes N^q$, while the differential $d\colon (M \otimes N)^s \to (M \otimes N)^{s + 1}$ is the sum of maps $d \otimes \id + (-1)^p \id \otimes d\colon M^p \otimes N^q \to M^{p + 1} \otimes N^q \oplus  M^{p} \otimes N^{q + 1}$. 
            \item In the case $n \geq 2$, for each pair $M, N \in \Ch(\CA)$, the swap map on the underlying object $\epsilon\colon U(M) \otimes U(N) \to U(N) \otimes U(M)$ is automatically compatible with differentials on both sides, so these assemble into an $\Eb_n$-monoidal structure $\Ch(\CA)^{\Kos}$. 
        \end{itemize}  
        Also, taking cycles/homology yields lax $\Eb_n$-monoidal functors $Z, H\colon \Ch(\CA)^{\Kos} \to \Gr(\CA)^{\Kos}$.
    \end{itemize}
\end{recollection}

\begin{recollection}
     There is a \textbf{Beilinson $t$-structure} on the $\infty$-category $\Fil\Sp$, in which
     \begin{itemize}
         \item $\Fil\Sp^{B}_{\geq 0}$ consists of filtrations $X$ so that each $(X / \defopara) (-w)$ is $w$-connective. 
         \item $\Fil\Sp^{B}_{\leq 0}$ consists of complete filtrations $X$ so that each $(X / \defopara) (-w)$ is $w$-coconnective. 
     \end{itemize}
     Its heart $\Fil\Sp^{B, \heartsuit}$ can be identified with $\Ch(\Ab)$, the $1$-category of chain complexes of abelian groups, 
     and the Day convolution symmetric monoidal structure on $\Fil\Sp$ induces the Koszul symmetric monoidal structure on $\Ch(\Ab)$, cf. \cite[\S~3 and \S~8]{Ant24}. 
     For each $X \in \Fil\Sp, n \in \Zb$, the $n$-th homotopy object $\pi^B_n X \in \Ch(\Ab)$ is the chain complex $(\pi_{n - a, a}(X / \defopara), \delta_1^1\colon \pi_{n - a, a}(X / \defopara) \to \pi_{n - a -1, a + 1}(X / \defopara))$. Thus, the collection $\{\pi_n^B(X)\}_{n \in \Zb}$ amounts to the $E_2$-page of the standard SS together with its $d_2$ differentials.
\end{recollection}

\begin{proposition} \label{Antieau-lemma}
     Take $n \in \Nb_{\geq 1} \cup \{\infty\}$. Suppose $\CE$ is a stably $\Eb_n$-monoidal $t$-$\infty$-category such that 
     \begin{itemize}
         \item The $\infty$-category $\CE$ admits sequential colimits. 
         \item The tensor product on $\CE$ is exact and preserves sequential colimits in each variable.
         \item The $t$-structure on $\CE$ is compatible with the tensor product (cf. \cite[Example 2.2.1.3]{HA}).
     \end{itemize}
     Then the functor $\pi_*\colon \CE \to \Gr(\CE^\heartsuit)$ upgrades to a lax $\Eb_n$-monoidal functor $\pi_*\colon \CE \to \Gr(\CE^\heartsuit)^{\Kos}$. 
\end{proposition}

\begin{proof}
    This is \cite[Theorem 8.7]{Ant24}. Roughly speaking, there is a \textbf{diagonal $t$-structure} on $\Fil(\CE)$ whose connective part $\Fil(\CE)^{D}_{\geq 0}$ (resp. whose coconnective part $\Fil(\CE)^{D}_{\leq 0}$) consists of filtrations $X$ so that each $X(w)$ is $w$-connective (resp. $w$-coconnective) in $\CE$, and its heart $\Fil(\CE)^{D, \heartsuit}$ is equivalent to $\Gr(\CE^{\heartsuit})^{\Kos}$ as an $\Eb_n$-monoidal $1$-category. The functor we desire is the composite 
    \[\CE \xrightarrow{\mathrm{const}} \Fil(\CE) \xrightarrow{\pi_0^{D}} \Fil(\CE)^{D, \heartsuit} \cong \Gr(\CE^{\heartsuit})^{\Kos}\]
    whose lax $\Eb_n$-monoidality follows from \cite[Example 2.2.1.10]{HA}.
\end{proof}

\begin{example} \label{lax-symmetric-monoidality-of-pi-star}
    \begin{itemize}
        \item Taking $\CE = \Sp$, we obtain a lax symmetric monoidal functor $\pi_*\colon \Sp \to \Gr(\Ab)^{\Kos}$. As the functor $i^*\colon \Fil\Sp \to \Gr\Sp$ induced by precomposing $i\colon \Zb^\delta \to \overrightarrow{\Zb}$ is also lax symmetric monoidal, we obtain a lax symmetric monoidal composite 
        \[\pi_{**}\colon \Fil\Sp \to \Gr\Sp \to \Gr(\Gr(\Ab)^{\Kos})^{\Day} \cong \Gr(\Gr(\Ab)^{\Day})^{\Kos}, \quad X \mapsto \{\pi_{n, w}(X)\}_{n, w \in \Zb}.\]
        Here for bigraded objects such as $\{\pi_{n, w}(X)\}_{n, w \in \Zb} \in \Gr(\Gr(\Ab)^{\Day})^{\Kos}$, we always use the first index (in this example: $n$) to denote the external (in this example: Koszul) grading and use the second index (in this example: $w$) to denote the internal (in this example: Day) grading. 
        \item Taking $\CE = \Fil\Sp$ with its Beilinson $t$-structure, we deduce that
        \[E_2 = \pi^B_{*}\colon \Fil\Sp \to \Gr(\Ch(\Ab)^{\Kos})^{\Kos}, \quad  X \mapsto (\{\pi_{n - w, w}(X / \defopara) = E_2^{2w - n, w}(X)\}_{n, w \in \Zb}, d_2)\] 
        is lax symmetric monoidal. Also, the functor $E_2 \Ac_k\colon \Fil\Sp \to \Gr(\Ch(\Ab)^{\Kos})^{\Kos}$ is lax symmetric monoidal for each $k \in \Nb$, $k \geq 1$. 
        \end{itemize}
\end{example}

\begin{lemma} \label{products-on-the-E2-page}
    There is a commutative square in $\CAlg(\Cat)^{\lax}$
    \[\begin{tikzcd}
        \Fil\Sp \ar[r, "{-\otimes \oneb / \defopara}"] \ar[d, "{E_2 = \pi^B_*}"]& \Fil\Sp \ar[r, "{\pi_{**}}"] &  \Gr(\Gr(\Ab)^{\Day})^{\Kos} \ar[d, "{v^*}"]\\
        \Gr(\Ch(\Ab)^{\Kos})^{\Kos} \ar[rr, "{\mathrm{forget}}"] && \Gr(\Gr(\Ab)^{\Kos})^{\Kos}
    \end{tikzcd}\]
    in which the symmetric monoidal equivalence $v^*\colon \Gr(\Gr(\Ab)^{\Day})^{\Kos} \cong \Gr(\Gr(\Ab)^{\Kos})^{\Kos}$ is given by precomposing the isomorphism 
    $v\colon \overrightarrow{\Zb} \times \overrightarrow{\Zb} \to \overrightarrow{\Zb} \times \overrightarrow{\Zb}, (n, a) \mapsto (n - a, a)$.
\end{lemma}
\begin{proof}
    To prove this we introduce two $t$-structures on $\Gr\Sp$: 
    \begin{itemize}
        \item The \textbf{pointwise $t$-structure}, whose connective part $\Gr\Sp^{P}_{\geq 0}$ (resp. coconnective part $\Gr\Sp^{P}_{\leq 0}$) consists of $X \in \Gr\Sp$ such that each $X(w)$ is $0$-connective (resp. $0$-coconnective). Its heart $\Gr\Sp^{P,\heartsuit}$ is equivalent to $\Gr(\Ab)^{\Day}$ as a symmetric monoidal $1$-category. 
        \item The \textbf{Beilinson $t$-structure}, whose connective part $\Gr\Sp^{B}_{\geq 0}$ (resp. coconnective part $\Gr\Sp^{B}_{\leq 0}$) consists of $X \in \Gr\Sp$ such that each $X(-w)$ is $w$-connective (resp. $w$-coconnective). Its heart $\Gr\Sp^{B,\heartsuit}$ is equivalent to $\Gr(\Ab)^{\Kos}$ as a symmetric monoidal $1$-category. 
    \end{itemize}
    The symmetric monoidal functor $- / \defopara\colon \Fil\Sp \to \Gr\Sp$ is $t$-exact for Beilinson $t$-structures on both sides, so it induces a commutative square in $\CAlg(\Cat)^{\lax}$
    \[\begin{tikzcd}
        \Fil\Sp \ar[d, "{E_2 = \pi_*^{B}}"] \ar[r, "{- / \defopara}"] & \Gr\Sp \ar[d, "{\pi_*^{B}}"] \\
        \Gr(\Ch(\Ab)^{\Kos})^{\Kos} \ar[r, "{\mathrm{forget}}"] & \Gr(\Gr(\Ab)^{\Kos})^{\Kos}
    \end{tikzcd}\]
    On the other hand, the composite $\Fil\Sp \xrightarrow{- \otimes \oneb / \defopara} \Fil\Sp \xrightarrow{i^*} \Gr\Sp$ can be identified with $\Fil\Sp \xrightarrow{-/\defopara} \Gr\Sp$, and $\Gr\Sp \to \Gr(\Gr(\Ab)^{\Day})^{\Kos}$ can be identified with $\pi_*^P\colon \Gr\Sp \to \Gr(\Fil\Sp^{P, \heartsuit})^{\Kos} \cong \Gr(\Gr(\Ab)^{\Day})^{\Kos}$ by examining the construction in the proof of Proposition \ref{Antieau-lemma}. What remains is to establish the commutative square
    \[\begin{tikzcd}
        \Gr\Sp \ar[d, "{\id}"] \ar[r, "{\pi_*^P}"] &  \Gr(\Gr(\Ab)^{\Day})^{\Kos}  \ar[d, "{v^*}"]\\
        \Gr\Sp \ar[r, "{\pi_*^B}"] & \Gr(\Gr(\Ab)^{\Kos})^{\Kos}
    \end{tikzcd}\]
    in $\CAlg(\Cat)^{\lax}$. Actually, as each arrow of this tentative square is already symmetric monoidal and the lower-right corner is a $1$-category, it suffices to fill the square in $\Alg(\Cat)^{\lax}$, where $(-)^{\Kos} = (-)^{\Day}$. 
    We achieve this by pasting together two squares
    \[\begin{tikzcd}
        \Gr\Sp \ar[d, "{V}"] \ar[r, "{\pi_*^P}"] &  \Gr(\Gr(\Ab))  \ar[d, "{\id}"]\\
        \Gr\Sp \ar[r, "{\pi_*^B}"] & \Gr(\Gr(\Ab))
    \end{tikzcd} \qquad \text{and} \qquad 
    \begin{tikzcd}
        \Gr\Sp \ar[d, "{V^{-1}}"] \ar[rr, "{\pi_*^P = \pi_{**}}"] &&  \Gr(\Gr(\Ab))  \ar[d, "{v^*}"]\\
        \Gr\Sp \ar[rr, "{\pi_*^P = \pi_{**}}"] && \Gr(\Gr(\Ab))
    \end{tikzcd}\]
    Here, the LHS square is induced by the monoidal equivalence $V\colon \Gr\Sp \to \Gr\Sp, \{X(w)\}_{w \in \Zb} \mapsto \{\Sigma^{-w} X(w)\}_{w \in \Zb}$, which is $t$-exact for the pointwise $t$-structure on the source and the Beilinson $t$-structure on the target. For the RHS square, note that the lax monoidal functor $\pi^P_{*} = \pi_{**} \colon \Gr\Sp \to \Gr(\Gr(\Ab))$ admits an alternative description as the composite
    \[\Gr\Sp \xrightarrow{\kappa} \Gr(\Gr\Sp) \xrightarrow{\Gr(\Gr(\pi_0))} \Gr(\Gr(\Ab))\]
    in which $\kappa$ is the monoidal left adjoint sending $\{X(w)\}_{w \in \Zb}$ to $\{\Sigma^n X(w)\}_{n, w \in \Zb}$. The RHS square is then deduced from the square of monoidal left adjoints
    \[\begin{tikzcd}
        \Gr\Sp \ar[d, "{V^{-1}}"] \ar[r, "{\kappa}"] &  \Gr(\Gr\Sp)  \ar[d, "{v^*}"]\\
        \Gr\Sp \ar[r, "{\kappa}"] & \Gr(\Gr\Sp)
    \end{tikzcd}\]
    whose existence follows from the monoidal universal property of $\Gr\Sp$.
\end{proof}

\begin{theorem}\label{Burklund's-Leibniz-rule}
    Take a map $F\colon X \otimes Y \to T$ in $\Fil\Sp$. For any $k \geq 1$, the bilinear pairing 
    \[F = F / \defopara^k\colon \pi_{n_1, w_1}(X / \defopara^k) \times \pi_{n_2, w_2}(Y/\defopara^k) \to \pi_{n_1 + n_2, w_1 + w_2}(T/\defopara^k)\] 
    extracted from the lax monoidality of $\Fil\Sp \xrightarrow{- \otimes \oneb / \defopara^k} \Fil\Sp \xrightarrow{\pi_{**}} \Gr(\Gr(\Ab)^{\Kos})^{\Day}$ satisfies the Leibniz rule for total differentials $\delta_k = \delta_k^k$, in the sense that 
    \[\delta_k F(\alpha, \beta) = F(\delta_k(\alpha), \beta) + (-1)^{n_1} F(\alpha, \delta_k(\beta)).\] 
    In particular, this leads to a bilinear map $F\colon E_{k + 1}^{s_1, t_1}(X) \times E_{k + 1}^{s_2, t_2}(Y) \to E_{k + 1}^{s_1 + s_2, t_1 + t_2}(T)$ satisfying the classical Leibniz rule, in the sense that $d_{k + 1}F(x, y) = F(d_{k + 1} (x), y) + (-1)^{t_1 - s_1} F(x, d_{k + 1}(y))$. 
\end{theorem}

\begin{proof}
    For $k = 1$ both claims follow from the lax monoidality of $E_2\colon  \Fil\Sp \to \Gr(\Ch(\Ab)^{\Kos})^{\Kos}$ together with Lemma \ref{products-on-the-E2-page}. In general, we work with the composite $E_2 \Ac_k$. By replacing $F\colon X \otimes Y \to T$ with $F\colon \Sigma^{0, w_1} X \otimes \Sigma^{0, w_2} Y \to \Sigma^{0, w_1 + w_2} T$ (this does not cause any sign issue as the swap map on $\Sb^{0, 1} \otimes \Sb^{0, 1}$ is homotopic to the identity, which is already true in $\overrightarrow{\Zb}$), we can assume that $w_1 = w_2 = 0$. Note that $\pi_{n, 0}(X / \defopara^k) = \pi_{n, 0}(\Ac_k(X / \defopara^k)) = \pi_{n, 0}(\Ac_k(X) / \defopara)$ and similarly for $Y$, while $\delta_k$ becomes $\delta_1$ after taking $\Ac_k$. Therefore, the image of $F\colon X \otimes Y \to T$ under the lax monoidal functor $E_2 \Ac_k$ gives rise to a refinement of the bilinear pairing on mod $\defopara^k$ stems in this case, with Leibniz rule for total differentials as expected. To construct the pairing on the $E_{k + 1}$-page, recall in the proof of Theorem~\ref{standard-SS} we identify $E_{k + 1}^{s, t}(X)$ with the image of $\cdot \defopara^{k - 1}\colon \pi_{t - s, t}(X / \defopara^{k}) \to \pi_{t - s, t - k + 1}(X / \defopara^{k})$. From this, we get the commutative square 
    \[\begin{tikzcd}
        X / \defopara^k \otimes Y / \defopara^k \ar[r, "{F}"] \ar[d, "{\cdot \defopara^{k - 1}}"] & T / \defopara^k \ar[d, "{\cdot \defopara^{k - 1}}"] \\
        \Sigma^{0, k - 1} (X / \defopara^{k} \otimes Y/ \defopara^{k}) \ar[r, "{F}"] &  \Sigma^{0, k - 1} T / \defopara^k
    \end{tikzcd}\]
    which further induces the diagram below after taking $\pi_{t - s, t}$ (with $s = s_1 + s_2, t = t_1 + t_2$)
    \[\begin{tikzcd}[column sep = small]
        \pi_{t_1 - s_1, t_1}(X / \defopara^k) \otimes_{\Ab} \pi_{t_2 - s_2, t_2}(Y / \defopara^k) \ar[r] &  \pi_{t - s, t}(X / \defopara^k \otimes Y / \defopara^k) \ar[r, "{F}"] \ar[d, "{\cdot \defopara^{k - 1}}"] & \pi_{t - s, t}(T / \defopara^k) \ar[d, "{\cdot \defopara^{k - 1}}"] \\
        & \pi_{t - s, t - k + 1}(X / \defopara^k \otimes Y / \defopara^k) \ar[r, "{F}"] & \pi_{t - s, t - k + 1}(T / \defopara^k)
    \end{tikzcd}\]
    The composite for the arrows on the left $\pi_{t_1 - s_1, t_1}(X / \defopara^k) \otimes_{\Ab} \pi_{t_2 - s_2, t_2}(Y / \defopara^k) \to  \pi_{t - s, t}(X / \defopara^k \otimes Y / \defopara^k) \to \pi_{t - s, t - k + 1}(X / \defopara^k \otimes Y / \defopara^k)$ factors through the epimorphism to $E_{k + 1}^{s_1, t_1}(X) \otimes_{\Ab} E_{k + 1}^{s_2, t_2}(Y)$ since its kernel contains both $\ker(\cdot \defopara^{k - 1}) \otimes \id$ and $\id \otimes \ker(\cdot \defopara^{k - 1})$. On the other hand, the epi-mono factorization for the vertical right arrow is precisely $\pi_{t - s, t}(T / \defopara^k) \to E_{k + 1}^{s, t}(T) \to \pi_{t - s, t - k + 1}(T / \defopara^k)$. Thus, we obtain an induced map $F\colon E_{k + 1}^{s_1, t_1}(X) \otimes_{\Ab} E_{k + 1}^{s_2, t_2}(Y) \to E_{k + 1}^{s, t}(T)$ in $\Ab$, and the Leibniz rule for total differentials implies immediately the Leibniz rule for ordinary differentials. 
\end{proof}

\begin{remark} \label{Burklund's-Leibniz-rule-multilinear}
    The same proof as in Theorem \ref{Burklund's-Leibniz-rule} also shows that, for every map $X_1 \otimes \cdots \otimes X_n \to T$, the induced multilinear pairings on $\pi_{**}(-/\defopara^r)$ and on $E_{r + 1}^{*,*}$ satisfy the Leibniz rule for (total) differentials.
\end{remark}

\begin{remark}
    The family of pairings $\{F = F / \defopara^k\}_{k \geq 1}$ also acquires compatibility with $\defopara^m$ and $\rho$. 
    \begin{itemize}
        \item For $1 \leq m \leq k \leq \infty, \alpha \in \pi_{**}(X / \defopara^k), \beta \in \pi_{**}(Y / \defopara^k)$, we have $\rho^k_mF(\alpha, \beta) = F(\rho^k_m(\alpha), \rho^k_m(\beta))$. Actually, the $\Eb_{\infty}$-ring map $\rho^k_m \colon \oneb / \defopara^k \to \oneb / \defopara^m$ induces a symmetric monoidal transformation between lax symmetric monoidal endofunctors $ - \otimes \oneb / \defopara^k \to - \otimes  \oneb / \defopara^m$, and the relation follows by postcomposing this with $\pi_{**}\colon \Fil\Sp \to \Gr(\Gr(\Ab)^{\Kos})^{\Day}$. 
        \item For $0 \leq m < k \leq \infty, \alpha \in \pi_{**}(X / \defopara^{k - m}), \beta \in \pi_{**}(Y / \defopara^k)$, we have $F(\defopara^{m}\alpha, \beta) = \defopara^m F(\alpha, \rho^k_{k - m}(\beta))$. To prove this, it suffices to construct a commutative square
        \[\begin{tikzcd}
            X / \defopara^{k - m} \otimes Y / \defopara^k \ar[rr, "{\id \otimes \rho^k_{k - m}}"] \ar[d, "{\defopara^m \otimes \id}"] && X / \defopara^{k - m} \otimes Y / \defopara^{k - m} \ar[r, "{F}"] & T / \defopara^{k - m} \ar[d, "{\defopara^m}"] \\
            X / \defopara^{k} \otimes Y / \defopara^k \ar[rrr, "{F}"] &&& T / \defopara^k
        \end{tikzcd}\]
        in $\Fil\Sp$. As the upper-right composite of the tentative square can be identified with $X / \defopara^{k - m} \otimes Y / \defopara^k \to X / \defopara^{k - m} \otimes_{\oneb / \defopara^{k}} Y / \defopara^k \cong (X \otimes Y) / \defopara^{k - m} \xrightarrow{F} T / \defopara^{k - m} \xrightarrow{\defopara^m} T / \defopara^k$, and the right-lower composite can be identified with $X / \defopara^{k - m} \otimes Y / \defopara^k \to  X / \defopara^{k - m} \otimes_{\oneb / \defopara^{k}} Y / \defopara^k \cong (X \otimes Y) / \defopara^{k - m} \xrightarrow{\defopara^m} (X \otimes Y) / \defopara^{k} \xrightarrow{F} T / \defopara^k$, one can fill this square with the homotopy $\defopara^m \circ F \simeq F \circ \defopara^m$. 
    \end{itemize}
\end{remark}

\begin{remark}\label{decalage}
    Using the Beilinson $t$-structure, one can construct a \textbf{décalage}\footnote{which stands for ``page-turning'' in SS} functor $\Dec \colon \Fil\Sp \to \Fil\Sp$ which sends $X$ to the filtration $\Dec(X)$ with $\Dec(X)(w) = (\tau^B_{\geq w}X)(-\infty)$. By construction, $\Dec(X)(-\infty) = X(-\infty)$ and $\Dec(X)$ is complete if $X$ is complete. Moreover, $\Dec$ admits a lax symmetric monoidal structure according to \cite[Lemma 8.7 and Theorem 8.9]{Ant24}. One can then produce a different description of the standard SS: $E_r^{s, t}(X) = E_2^{s + (r - 2)(t - s), s + (r - 1)(t - s)}(\Dec^{(r - 2)}(X))$ and\footnote{here $\Dec^{(k)}$ stands for $k$-fold iteration of $\Dec$, in particular $\Dec^{(0)} = \id$} $d_r$ is given by reindexing $\delta_1 = d_2$ for $\Dec^{(r - 2)}(X)$, cf. Remark \ref{Lurie-comparison} and \cite[Theorem 4.11]{Ant24}. Consequently, the lax monoidality of $\Dec$ leads to another proof for the classical Leibniz rule. 
\end{remark}

In Appendix \ref{app:C}, we will recast Theorem \ref{Burklund's-Leibniz-rule} as a categorical statement, expressed in terms of the lax symmetric monoidality of certain functors.

\subsection{The stages beyond filtered spectra}
\label{subsec:2.4}

In this subsection we record a vastly general setup for BIPWX's cofiber-of-$\tau$ formalism. We start by showing there are plenty of methods to construct filtrations that recover known spectral sequences.

\begin{example} \label{SS-examples-1}
    \begin{enumerate}
        \item Suppose $X$ is a pointed CW complex. Then the standard SS of the filtration $\Sigma^{\infty} \mathrm{sk}_{*} X  \in \Fil\Sp$ is a reindexed version\footnote{In the standard convention, this is the Atiyah--Hirzebruch SS of $X$ that starts from its $E_1$ page} of the Atiyah--Hirzebruch SS of $X$. 
        \item Suppose $R \to S$ is a map of $\Eb_n$-ring spectra and $M$ is an $R$-module, such that $M$ is \textbf{$S$-nilpotent complete} in the sense that $M \cong \lim_{[n] \in \Delta} M \otimes_R S^{\otimes_R (n + 1)}$. Then the totalization of Whitehead towers $\lim_{[n] \in \Delta} \tau_{\geq *}(M \otimes_R S^{\otimes_R (n + 1)})$ is a filtration whose standard SS is the \textbf{descent SS} of $M$. In particular, if $R = \Sb^0$ and $S = E$ is an Adams type ring spectrum (cf. \cite[Definition 3.14]{Pst23}), then this recovers the $E$-Adams SS of $M \in \Sp$ starting from its $E_2$-page. 
        \item Suppose $\CC$ is a $t$-$\infty$-category and $F\colon \CC \to \Sp$ is an exact functor. For $X \in \CC$, the standard SS of $F(\tau_{\geq *} X)$ recovers many familiar constructions: 
        \begin{itemize}
            \item If $\CC = \Fun(\mathrm{B}G, \Sp)$ for some compact Lie group $G$, then for $F = (-)_{hG}$ (resp. $F = (-)^{hG}$, or $F = (-)^{tG}$), this recovers the homotopy orbit SS (resp. homotopy fixed point SS, or Tate SS) of the Borel $G$-spectrum $X$. 
            \item If $\CC = \Sh(\CE; \Sp)$ is the $\infty$-category of spectral sheaves over some $\infty$-topos $\CE$ and $F = \Gamma$ is the global section functor, then this recovers (a generalization of) the Leray SS computing the hypercohomology of  $X$.
        \end{itemize}
        \item Suppose $\CA \xrightarrow{F} \CB \xrightarrow{G} \CC$ is a composable pair of left exact functors between Grothendieck abelian $1$-categories, which induces a pair of derived functors $\D^+(\CA) \xrightarrow{\mathrm{R} F} \D^+(\CB) \xrightarrow{\mathrm{R} G} \D(\CC)$, cf. \cite[\S~1.3.3]{HA}. For $X \in \D^+(\CA)$, $\mathrm{R} G(\tau_{\geq * } \mathrm{R} F(X)) \in \Fil(\D(\CC))$ yields a $\CC$-valued standard SS (cf. Corollary \ref{standard-SS-with-general-coefficients}), which recovers the Grothendieck SS computing $\pi_*\mathrm{R}(G \circ F)(X)$.
    \end{enumerate}
\end{example}

To fully make use of BIPWX's cofiber-of-$\tau$ formalism, it is not enough to produce filtrations only. In particular, the computational methods based on hidden extensions (in the next two sections) require more structured inputs such as cofiber sequences of filtrations. We shall introduce the ``minimal'' coherent structure that generates such functoriality data. Recall that $\Fil\Sp$ is a stable presentably symmetric monoidal $\infty$-category, in other words $\Fil\Sp \in \CAlg(\PPr_{\st})$. 

\begin{definition} \label{FilSp-module-def}
    A \textbf{presentable $\Fil\Sp$-module} is a $\Fil\Sp$-module in\footnote{here the definition will not change if we replace $\PPr_{\st}$ by $\PPr$ due to idempotency of $\Sp \in \PPr$} $\PPr_{\st}$. Concretely, this is a presentable stable $\infty$-category $\CC$ equipped with a functor $\CC \times \Fil\Sp \to \CC$ that commutes with colimits on each side, subject to associativity and unitality coherences. For $\CC \in \Mod_{\Fil\Sp}(\PPr_{\st})$, we write $\CC[\defopara^{-1}] = \Mod_{\oneb[\defopara^{-1}]}(\CC) = \CC \otimes_{\Fil\Sp} \Sp \in \PPr_{\st}$ for the \textbf{$\defopara$-generic fiber} of $\CC$, and we write $\CC/\defopara = \Mod_{\oneb / \defopara}(\CC) = \CC \otimes_{\Fil\Sp} \Gr\Sp \in \Mod_{\Gr\Sp}(\PPr_{\st})$ for the \textbf{$\defopara$-special fiber} of $\CC$. We denote by $\defopara^{-1}\colon \CC \to \CC[\defopara^{-1}]$, $-/\defopara\colon \CC \to \CC / \defopara$ the two basechange comparison functors. 
\end{definition}

\begin{remark}
    Here we invoke the equivalence of $\infty$-categories $\CC\otimes_{\CV\!}\Mod_A(\CV) \cong \Mod_{A}(\CC)$ following \cite[Theorem 4.8.4.6]{HA}. 
\end{remark}

\begin{remark}
    Alternatively, a presentable $\Fil\Sp$-module is a \emph{locally filtered stable $\infty$-category} in the sense of \cite[Definition 3.1.10]{Rotational-Invariance} whose underlying $\infty$-category is presentable. This follows from comparing the characterization in  \cite[Remark 3.1.12]{Rotational-Invariance} and in Corollary \ref{GHMG-for-modules}.
\end{remark}

\begin{construction} \label{fundamental-stuff-for-FilSp-modules}
    For $\CC \in \Mod_{\Fil\Sp}(\PPr_{\st})$ and $X \in \CC$, we set $\Sigma^{n, w} X = X \otimes \Sb^{n, w}$ and set $\defopara^k = X \otimes \defopara^k_{\oneb}\colon \Sigma^{0, -k} X \to X$. We also construct the fundamental distinguished triangles (cf. Construction \ref{fundamental-distinguished-triangle}) for $X$ by tensoring $X$ with the fundamental distinguished triangles of $\oneb \in \Fil\Sp$. Clearly, all results we proved in \S~\hyperref[subsec:2.1]{2.1} also hold true in this generality.
\end{construction}

This construction imports the basic language of $\Fil\Sp$ into $\CC$. On the other hand, from $\CC$ we can also extract information into $\Fil\Sp$, using the standard fact that a $\Fil\Sp$-action on a presentable $\infty$-category amounts to a $\Fil\Sp$-enrichment in the sense of \cite{GH15}. 

\begin{theorem}
    Suppose $\CV \in \CAlg(\PPr_{\st})$, and $\CM \in \Mod_{\CV}(\PPr_{\st})$. Then for each $m \in \CM$, the functor $\CV \to \CM, v \mapsto m \otimes v$ admits a right adjoint $\CM \to \CV, m' \mapsto \hom_{\CM}^{\CV}(m, m')$, and these assemble to a $\CV\!$-enriched $\infty$-category structure on $\CM$. Furthermore, this renders an equivalence between $\Mod_{\CV\!}(\PPr_{\st})$ and the subcategory of large $\CV\!$-enriched $\infty$-categories spanned by those $\CM$ admitting $\CV\!$-tensors whose underlying $\infty$-categories are presentable stable, and left adjoint $\CV\!$-enriched functors in between.  
\end{theorem}

\begin{proof}
    This follows from \cite[Theorem 1.2]{Hei23}.
\end{proof}

\begin{remark}
    For $\CV = \Fil\Sp$ we shall write $\Msp^{\Fil}$ for $\hom^{\Fil\Sp}$. Similarly for $\CV = \Gr\Sp$. 
\end{remark}

By construction, $\Msp^{\Fil}$ satisfies the following standard properties:

\begin{fact} \label{filtered-mapping-spectra-facts}
    Suppose $\CC$ is a presentable $\Fil\Sp$-module, $X, Y \in \CC$.
    \begin{itemize} 
        \item The filtered mapping spectrum $\Msp^{\Fil}_{\CC}(X, Y) \in \Fil\Sp$ is concretely given by $\Msp^{\Fil}_{\CC}(X, Y)(w) \cong \Msp_{\CC}(\Sigma^{0, w} X, Y) \cong \Msp_{\CC}(X, \Sigma^{0, -w}Y)$. 
        \item The functor $\Msp^{\Fil}(-, -)\colon \CC^{\op} \times \CC \to \Fil\Sp$ preserves limits and finite colimits for each variable.  
        \item The functor $\Msp^{\Fil}(-, -)$ is compatible with bigraded suspensions, in the sense that 
        \[\Sigma^{n, w}\Msp^{\Fil}_{\CC}(X, Y) \cong \Msp^{\Fil}_{\CC}(\Sigma^{-n, -w}X, Y) \cong \Msp^{\Fil}_{\CC}(X, \Sigma^{n, w} Y).\] 
        \item The $\Msp^{\Fil}(X, -)$ image of $\defopara^k \colon \Sigma^{0, -k} Y \to Y$ is $\defopara^k\colon \Sigma^{0, -k} \Msp^{\Fil}(X, Y) \to \Msp^{\Fil}(X, Y)$. Combining this with its limit preservation, it follows that the $\Msp^{\Fil}(X, -)$ image of any fundamental distinguished triangle (cf. Construction \ref{fundamental-stuff-for-FilSp-modules}) of $Y$ is the corresponding fundamental distinguished triangle for $\Msp^{\Fil}(X, Y)$. Similar results hold true for $\Msp^{\Fil}(-, Y)$ as well. % \footnote{In the contravariant variable, the analogous statements may involve shifts; for instance, $\Msp^{\Fil}(X / \defopara^a, Y) \cong \Sigma^{-1, a}(\Msp^{\Fil}(X, Y)/ \defopara^a)$.}.
    \end{itemize}
\end{fact}

We further relate the generic/special fiber constructions with filtered mapping spectra.

\begin{lemma} \label{enriched-hom-comparison}
    Suppose $F\colon \CV \to \CW$ is a map in $\CAlg(\PPr_{\st})$, $\CM \in \Mod_{\CV}(\PPr_{\st}), \CN \in \Mod_{\CW}(\PPr_{\st})$ and $F_1 \colon \CM \to \CN$ is a $\CV\!$-linear left adjoint. Write $F^R, F_1^R$ for the right adjoints of $F$ and $F_1$. 
    \begin{itemize}
        \item For $x \in \CM, z \in \CN$, there is a natural isomorphism $F^R(\hom^{\CW}\!\!(F_1(x), z)) \cong \hom^{\CV}\!(x, F_1^R(z))$. 
        \item For $x, y \in \CM$, there is a natural map $F(\hom^{\CV}\!\!(x, y)) \to \hom^{\CW}\!\!(F_1(x), F_1(y))$.
        \item For $x, y \in \CM$, there is a natural map $F^R F(\hom^{\CV}\!\!(x, y)) \to \hom^{\CV}\!\!(x, F_1^R F_1(y))$. If $F^R$ is conservative, then this is invertible iff $F(\hom^{\CV}\!\!(x, y)) \to \hom^{\CW}\!\!(F_1(x), F_1(y))$ is an isomorphism.
    \end{itemize} 
\end{lemma}

\begin{proof}
    For the first item, take any $v \in \CV$, we have 
    \begin{align*}
        \Map_{\CV}(v, F^R(\hom^{\CW}\!\!(F_1(x), z))) &\cong \Map_{\CW}(F(v), \hom^{\CW}\!\!(F_1(x), z)) \\
        &\cong \Map_{\CN}(F_1(x) \otimes F(v), z) \\
        &\cong \Map_{\CN}(F_1(x \otimes v), z) \\
        &\cong \Map_{\CM}(x \otimes v, F_1^R(z)) \\
        &\cong \Map_{\CV}(v, \hom^{\CV}\!\!(x, F_1^R(z))).
    \end{align*}
    This induces the isomorphism $F^R(\hom^{\CW}\!\!(F_1(x), z)) \cong \hom^{\CV}\!\!(x, F_1^R(z))$. The second map is constructed as the transpose of the map $F_1(x) \otimes F(\hom^{\CV}\!\!(x, y)) \cong F_1(x \otimes \hom^{\CV}\!\!(x, y))  \to F_1(y)$, and the third map is constructed as the composite 
    \[F^R F(\hom^{\CV}\!\!(x, y)) \to F^R(\hom^{\CW}\!\!(F_1(x), F_1(y))) \cong \hom^{\CV}\!\!(x, F_1^R F_1(y))\]
    of the second map and the first isomorphism. Therefore, the third map is an isomorphism iff the second map becomes an isomorphism after taking $F^R$ image, which is further equivalent to the second map itself being invertible if $F^R$ is conservative. 
\end{proof}

\begin{theorem} \label{filtered-mapping-spectra-interactions}
    Suppose $\CC$ is a presentable $\Fil\Sp$-module, and $X, Y \in \CC$.  
    \begin{itemize}
        \item For the $\defopara$-special fiber, we have $\Msp_{\CC}^{\Fil}(X, Y) / \defopara \cong \Msp_{\CC / \defopara}^{\Gr}(X / \defopara, Y / \defopara)$. 
        \item For the $\defopara$-generic fiber, if $X$ is sequentially compact in $\CC$, then 
        \[\Msp^{\Fil}_{\CC}(X, Y) [\defopara^{-1}] \cong \Msp_{\CC[\defopara^{-1}]}(X[\defopara^{-1}], Y[\defopara^{-1}]).\]
        \item Also, $\Msp^{\Fil}_{\CC}(X, Y)$ is complete if $Y$ is $\defopara$-complete (i.e. $\lim_{w \to \infty} \Sigma^{0, -w} Y = 0$), and for completions $\Msp^{\Fil}_{\CC}(X, Y)_{\hat{\defopara}} \cong \Msp^{\Fil}_{\CC}(X, Y_{\hat{\defopara}}) \cong \lim_{w \to \infty} \Msp^{\Fil}_{\CC}(X, Y/\defopara^w)$. 
    \end{itemize} 
\end{theorem}

\begin{proof}
    To establish the first identification, by Lemma \ref{enriched-hom-comparison} it suffices to show $\Msp^{\Fil}_{\CC}(X, Y)/\defopara \to \Msp^{\Fil}_{\CC}(X, Y / \defopara)$ is invertible. This is true as both of them are the cofiber of $\defopara\colon \Sigma^{0, -1}\Msp^{\Fil}_{\CC}(X, Y) \to \Msp^{\Fil}_{\CC}(X, Y)$. For the second part, one can work with $\Msp^{\Fil}_{\CC}(X, Y)[\defopara^{-1}] \to \Msp^{\Fil}_{\CC}(X, Y[\defopara^{-1}])$. As inverting $\defopara$ is essentially taking a sequential colimit, this is an isomorphism if $\Msp_{\CC}(X, -)$ preserves all sequential colimits, i.e. if $X$ is sequentially compact in $\CC$. Lastly, the results on completeness follow from the fact that $\Msp^{\Fil}_{\CC}(X, -)$ preserves limits and cofiber sequences. 
\end{proof}

\begin{remark} \label{recollements-in-general}
    Write $\CC_{\hat{\defopara}}$ for the full subcategory of $\CC$ spanned by $\defopara$-complete objects. It follows from \cite[Proposition 5.2.3]{CSY21} that $\CC_{\hat{\defopara}} = \CC \otimes_{\Fil\Sp} \Fil\Sp_{\hat{\defopara}}$, and 
    the inclusions $\CC[\defopara^{-1}] \hookrightarrow \CC \hookleftarrow \CC_{\hat{\defopara}}$ exhibit $\CC$ as a recollement of $\CC[\defopara^{-1}]$ and $\CC_{\hat{\defopara}}$. Note that the proof of Theorem \ref{filtered-mapping-spectra-interactions} also implies that $\Msp^{\Fil}_{\CC}(X, Y)_{\hat{\defopara}} \cong \Msp^{\Fil}_{\CC_{\hat{\defopara}}}(X_{\hat{\defopara}}, Y_{\hat{\defopara}})$.
\end{remark}

Theorem \ref{filtered-mapping-spectra-interactions}, together with Fact \ref{filtered-mapping-spectra-facts}, builds the bridge between $\Fil\Sp$ and an arbitrary presentable $\infty$-category $\CC$ with cocontinuous $\Fil\Sp$-action. It remains to find examples for such $\CC$. 

\begin{example} \label{SS-examples-2}
    If $\CC$ is a stable presentably symmetric monoidal $\infty$-category, and $\Fil\Sp \to \CC$ is a symmetric monoidal left adjoint, then $\CC$ admits a cocontinuous $\Fil\Sp$-action. For instance:
    \begin{enumerate}
        \item If $\CE$ is a stable presentably symmetric monoidal $\infty$-category, then the unique symmetric monoidal left adjoint $\Sp \to \CE$ yields a symmetric monoidal left adjoint $\Fil\Sp \to \CC = \Fil(\CE)$ by postcomposition. In this case, $\CC[\defopara^{-1}] \cong \CE$ and $\CC/\defopara \cong \Gr(\CE)$.
        \item If $R$ is an $\Eb_{\infty}$ algebra in $\Fil\Sp$, then there is a symmetric monoidal left adjoint $- \otimes R \colon \Fil\Sp \to \CC = \Mod_R (\Fil\Sp)$. In this case, $\CC[\defopara^{-1}] = \Mod_{R[\defopara^{-1}]}(\Sp)$ and $\CC / \defopara \cong \Mod_{R / \defopara}(\Gr\Sp)$. By \cite[Proposition C.22]{BHS2} this covers the case of \textbf{(cellular) $E$-synthetic spectra} in \cite{Pst23} since $\Syn_E^{\cell} \cong \Mod_{R_E} (\Fil\Sp)$, where $R_E = \lim_{[n] \in \Delta} \tau_{\geq *}(E^{\otimes(n + 1)})$. 
        In this case, write $\nu\colon \Sp \to \Syn_E$ for the \textbf{synthetic analogue} functor in \cite[Definition 4.3]{Pst23}, then for each $E$-nilpotent complete spectrum $X$ this equivalence of categories sends $\nu X$ to the $E$-Adams filtration of $X$ in Example \ref{SS-examples-1} item 2, cf. Theorem \ref{synthetic-Adams-vs-Bockstein}.
        \item More generally, if $(\CC = \CE^{\mathrm{def}}, \CE)$ is a \emph{deformation pair} in the sense of \cite[Definition C.14]{BHS2}, then there is a symmetric monoidal left adjoint $\Fil\Sp \to \CE^{\mathrm{def}}$ due to \cite[Construction C.18]{BHS2}. In this case, $\CE^{\mathrm{def}}[\defopara^{-1}] \cong \CE$. 
        \item On the other hand, in terms of spectral algebraic geometry, we have $\Fil(\Mod_R(\Sp)) \cong \QCoh(\AAb^1_R / \Gb_{m, R})$ for each $R \in \CAlg(\Sp)$, cf. \cite[Theorem 1.1]{Mou21}. Thus, every spectral\footnote{It is easy to construct maps of classical stacks $f\colon X \to \AAb^1 / \Gb_{m}$ over a classical commutative ring $R$ (or maps of \emph{derived} stacks $f\colon X \to \AAb^1 / \Gb_{m}$ over an \emph{animated} commutative ring $R$): such a map is determined by the data of a line bundle $\CL$ on $X$  (which corresponds to $f^*(\Sigma^{0, 1} R)$) together with a section $\CO_X \to \CL$ (which corresponds to the map $f^* \defopara$), cf. \cite[Proposition 3.2.6]{Khan-Rydh}. In general, the spectral stack $\AAb^1 / \Gb_{m}$ classifies \emph{strict} line bundles with \emph{strict} sections, whose construction would require more homotopy coherence data.} stack $X$ over $\AAb^1 / \Gb_{m}$ gives rise to a symmetric monoidal left adjoint $f^*\colon \Fil\Sp \to \QCoh(X)$. In this case, $\QCoh(X)[\defopara^{-1}] \cong \QCoh(X \times_{\AAb^1 / \Gb_{m}} *)$ and  $\QCoh(X) / \defopara \cong \QCoh(X \times_{\AAb^1 / \Gb_{m}} \mathrm{B} \Gb_m)$. As a concrete example, if $R = \Cb$ and $X = Y^{\mathrm{dR}, +}$ is Simpson's filtered de Rham stack (as in \cite[Definition 2.3.5]{FGauge}) of a smooth complex variety $Y$, then the right adjoint $f_*\colon \QCoh(X) \to \Fil\Sp$ sends the structure sheaf $\CO_X$ to the filtration whose standard SS is the Hodge-to-de Rham SS of $Y$, cf. \cite[Theorem 2.3.6]{FGauge}. 
    \end{enumerate}
\end{example}

\begin{remark}[Comparison with synthetic spectra] \label{synthetic-spectra-as-special-case}
    Consider the $\infty$-category $\Syn_E$ of $E$-synthetic spectra for an Adams type ring spectrum $E$. Write $U = \Msp^{\Fil}(\Sb^{0, 0}_E, -)\colon \Syn_E \to \Fil\Sp$ for the forgetful functor, which is right adjoint to the functor $\Fil\Sp \to \Syn^{\cell}_E \subset \Syn_E$ above. By Fact \ref{filtered-mapping-spectra-facts}, the operators $\tau^k$ (or $\lambda^k$), $\rho$, and $\delta_k$ in \cite{BHS1}, \cite[\S~7]{Burklund-Xu}, and \cite{Lin-Wang-Xu-kervaire} map under $U$ to the operators $\defopara^k$, $\rho$, and $\delta_k$ in $\Fil\Sp$ discussed above. Consequently, our Bockstein dictionary (i.e. Theorem \ref{Bockstein-dictionary-infinite}) recovers the ``omnibus theorem'' in $\Syn_E$, namely \cite[Theorem 9.19]{BHS1}. 
\end{remark}

In general it requires a lot of coherence data to construct \emph{symmetric monoidal} left adjoints from $\Fil\Sp$ to a certain abstract $\infty$-category $\CC$. On the other hand, to produce cocontinuous $\Fil\Sp$-actions it suffices to construct \emph{monoidal} left adjoints $\Fil\Sp \to \CC$, which is much simpler. 

\begin{theorem} \label{GHMG}
    Suppose $\CC$ is a stable presentably monoidal $\infty$-category. Then the assignment $F \mapsto F(\oneb / \defopara)$ defines an equivalence from the $\infty$-category of monoidal left adjoint functors $\Fil\Sp \to \CC$ to the full subcategory of $\Alg(\CC)$ spanned by those $\Eb_1$-algebras $R$ for which $V = \fib(\oneb \to R)$ is invertible.
\end{theorem}

\begin{proof}
    This follows from \cite[Proposition 1.5]{GHMG}. 
\end{proof}

\begin{corollary} \label{GHMG-for-modules}
    According to \cite[Corollary 4.7.1.41]{HA}, a presentable $\Fil\Sp$-module structure on a presentable stable $\infty$-category $\CC$ amounts to a monoidal left adjoint from $\Fil\Sp$ to the (stable presentably monoidal) $\infty$-category of left adjoint endofunctors $\Fun^L(\CC, \CC)$. Thus, by Theorem \ref{GHMG}, the data of a presentable $\Fil\Sp$-module structure on $\CC$ is the same as the data of a cocontinuous monad $T \in \Alg(\Fun^L(\CC, \CC))$ whose fiber of the unit transformation $V = \fib (\eta\colon \id_{\CC} \to T)$ is an invertible functor. There is an abundant source of such monads from \textbf{spherical adjunctions}, cf. \cite{DKSS, GHMG}. 
\end{corollary}

%% file: GLR.tex
\label{sec:3}

\newdimen\widthA
\newcommand{\resizelabel}[2]{%
  \settowidth{\widthA}{#2}%
  \pgfmathsetmacro{\ratio}{min((0.8cm)/max(\widthA, 0.8cm), #1)}%
  \scalebox{\ratio}{#2}%
}

In this section we begin our study of the hidden extensions using BIPWX's cofiber-of-$\tau$ formalism. More precisely, in \S~\hyperref[subsec:3.1]{3.1} we define extensions on the $E_r$-page for $2\leq r \leq \infty$ and the corresponding no crossing conditions, demonstrate how to translate different no-crossing conditions in terms of bigraded homotopy groups, and establish compatibility of extensions across different pages; in \S~\hyperref[subsec:3.2]{3.2} we prove the generalized Leibniz rule (Theorems \ref{GLR-for-SS} and \ref{GLR-for-SS-2}). 

\subsection{Hidden extensions and crossings}
\label{subsec:3.1}

We first introduce the notion of \emph{(hidden) extensions}, which captures the data of a map $f\colon X \to Y$ in $\Fil\Sp$---especially the part invisible after passing to the induced map on graded pieces---in terms of the standard spectral sequences of $X$ and $Y$. We also introduce \emph{crossings}, which quantify a certain subtlety in the study of extensions; their role will become precise after Lemma \ref{free-choice-lemma-for-refined-extensions}. 

\begin{definition}[Extensions and crossings] \label{extensions-def}
    Let $f\colon X \to Y$ be a map in $\Fil\Sp$. Take $r \in \Nb \cup \{\infty\}$, $j, k\in \Nb$ such that $k < r$. 
    \begin{itemize}
        \item There is a \textbf{$j$-th layer extension along $f$} on the $E_{r + 1}$-page of filtration jump $k$ from $x \in Z_r^{s, t}(X)$ to $y \in Z_{r - k}^{s + k, t + k}(Y) / B_{1 + k + j}^{s + k, t + k}(Y)$ if there exist $[x] \in \pi_{t - s, t}(X / \defopara^r), [y] \in \pi_{t - s, t + k}(Y / \defopara^{r - k})$ lifting $x, y$ so that
        \[\defopara^{j}f[x] = \defopara^{k + j} [y].\]
        In this case we write $d_k^{f, E_{r + 1}, j}(x) = y$, and we refer to the pair $([x], [y])$ as a \textbf{witness} of the $j$-th layer extension. If the index $j$ is clear from the context, we just refer to this as an \textbf{extension}, and if furthermore $k \geq 1$ we say the extension is \textbf{hidden}. In the case $j = 0$, we write $d_k^{f, E_{r + 1}}(x) = y$ as an abbreviation for $d_k^{f, E_{r + 1}, 0}(x) = y$. 
        \item Such a $j$-th layer extension is \textbf{trivial} if $y = 0$ in $Z_{r - k}^{s + k, t + k}(Y) / B_{1 + k + j}^{s + k, t + k}(Y)$, in other words the witness $\defopara^{k + j} [y] \in \IIm(\defopara^{k + j + 1})$. It is \textbf{nontrivial} if it is not trivial.
        \item Such a $j$-th layer extension is \textbf{inessential} if it is trivial or there exists $x' \in Z^{s + a, t + a}_{r - a}(X)$ with $0 < a \leq k$ so that $d_{k - a}^{f, E_{r + 1 - a}, j + a}(x') = y$. It is \textbf{essential} if it is not inessential.
        \item A \textbf{precrossing} for such a $j$-th layer extension is a pair 
        \[x' \in Z_{r - a}^{s + a, t + a}(X), y' \neq 0 \in Z_{r - k + b}^{s + k - b, t + k - b}(Y) / B^{s + k - b, t + k - b}_{1 + k + j - b}(Y)\] 
        with $0 < a \leq k, 0 \leq b \leq k - a$ so that $d_{k - a - b}^{f, E_{r + 1 - a}, j + a}(x') = y'$. Such a precrossing is a \textbf{crossing} if it is furthermore an essential extension.
        \item The extension has \textbf{(full) no crossing} if no such crossing extension exists. More precisely, for $S \subset \{1, 2, \ldots, k\}, T \subset \{0, 1,\ldots, k - 1\}$, the extension has \textbf{(partial) no crossing with source range $S$, target range $T$} if there is no such crossing extension with $a \in S, b \in T$.
    \end{itemize}
\end{definition}

\begin{remark}[Classical hidden extensions]
    There is also a notion of ``extensions on the abutment'' which matches the usage of the term ``(hidden) extensions'' in common literature such as \cite{IWX}: 
    \begin{itemize}
        \item Let $f\colon X \to Y$ be a map in $\Fil\Sp$, fix $s, t \in \Zb$, $k \in \Nb$ and let $x \in Z_{\infty}^{s, t}(X)$, $y \in Z_{\infty}^{s + k, t + k}(Y)$ be two surviving permanent cycles. There is an \textbf{extension on the abutment} along $f$ of filtration jump $k$ from $x$ to $y$ if there exists $[x] \in \pi_{t - s, t}(X)$ and $[y] \in \pi_{t - s, t + k}(Y)$ lifting $x, y$, so that $f \iota[x] = \iota[y]$ in $\pi_{t - s}(Y[\defopara^{-1}])$. Here $\iota$ is the $\defopara$-inversion map in Construction \ref{limit-distinguished-triangle}.
    \end{itemize}
    As $\iota\colon Y \to Y[\defopara^{-1}]$ is the colimit of all maps $\defopara^w\colon Y \to \Sigma^{0, w} Y$, the information of an $f$-extension on the abutment is equivalent to that of an $f$-extension on the $E_{\infty}$-page in a sufficiently high layer. 
\end{remark}

\begin{remark} \label{comparison-with-extension-SS}
    The paper \cite{Lin-Wang-Xu-kervaire} develops the theory of hidden extensions from a slightly different perspective, using the language of \emph{extension SS}. We compare the two approaches below:
    \begin{itemize}
        \item For each map between strictly filtered abelian groups $f\colon F^* A \to F^* B$, one can consider the length $1$ chain complex $[A \to B]$ equipped with the induced filtration (i.e. $x$ is of filtration $\geq s$ if $x \in F^{s} B$ or $x \in F^{s + 1} A$). The \textbf{extension SS} is the standard SS for this filtration, suitably reindexed so that it starts from the $E_0$-page (which is the direct sum of $\gr^*(FA)$ and $\gr^*(FB)$). 
        \item Suppose $X \to Y$ is a map in $\Fil\Sp$. For $1 \leq r < \infty$, take $f\colon F^*_{\defopara} \pi_{**}(X / \defopara^r) \to F^*_{\defopara} \pi_{**}(Y / \defopara^r)$ to be the map between $\IIm(\defopara^a)$-filtrations in the sense of Theorem \ref{Bockstein-dictionary-finite}. Then a $j$-th layer extension along $f$ on the $E_{r + 1}$-page from $x$ to $y$ of filtration jump $k$ is a differential $d_k(\defopara^j x) = \defopara^{k + j} y$ in the extension SS. 
        \item Similarly, for $r = \infty$, take $f\colon F^*_{\defopara} \pi_{**}(X) \to F^*_{\defopara} \pi_{**}(Y)$ to be the map between $\IIm(\defopara^a)$-filtrations in the sense of Theorem \ref{Bockstein-dictionary-infinite}. Then a $j$-th layer extension along $f$ on the $E_{\infty}$-page from $x$ to $y$ of filtration jump $k$ is a differential $d_k(\defopara^{j}x) = \defopara^{k + j} y$ in the extension SS.  
        \item This characterization of $f$-extensions is equivalent to Definition \ref{extensions-def} above. To prove the identification, note that we can construct from $f\colon F^* A \to F^* B$ a distinguished triangle in $\Fil\Sp$
        \[ \Sigma^{0, -1} F^* A\xrightarrow{\defopara f} F^* B \xrightarrow{g} F^*[A \to B] \xrightarrow{h} \Sigma^{1, -1} F^* A.\]
        From this, the translation between Definition \ref{extensions-def} and the extension differential approach follows from the generalized Mahowald trick (GMT), i.e. Theorem \ref{GMT-for-SS} and Theorem \ref{GMT-for-SS-2}. Here 
        \begin{itemize}
            \item Assumption item 4 in each case holds true for degree reasons.
            \item Assumption item 5 in Theorem \ref{GMT-for-SS-2} is also automatic since there is no nontrivial boundary term for the strict filtration $F^* B$.
        \end{itemize}
        \item To finish the comparison with \cite{Lin-Wang-Xu-kervaire} there is one remaining point: the $\IIm(\defopara^a)$-filtrations on $\pi_{**}(X)$ and $\pi_{**}(X / \defopara^r)$ come from the abutment of the trigraded $\defopara$-Bockstein SS (cf. Construction \ref{trigraded-defopara-Bockstein-SS} and Theorem \ref{rigidity-for-trigraded-defopara-Bockstein-SS}), while in \cite[Definition 5.4]{Lin-Wang-Xu-kervaire} the extension SS are constructed from the abutment of the synthetic $\mathrm{H}\Fb_2$-Adams SS. In fact, these are also equivalent due to Theorem \ref{synthetic-Adams-vs-Bockstein}.  
        \item The extension SS approach explains the terminology $d_k^{f, E_{r + 1}, j}(x) = y$ we use. Also, an extension is essential iff it is an essential (cf. Remark \ref{essential-diff}) extension differential, and crossings for extensions are realized as diagrammatic crossings of extension differentials. 
    \end{itemize}
\end{remark}

\begin{remark} \label{comparison-with-extension-SS-contd}
    Modulo the translation between bigraded stem functoriality and extension differentials as discussed in Remark \ref{comparison-with-extension-SS}, there are some further differences in describing hidden extensions and crossings between \cite[\S~5]{Lin-Wang-Xu-kervaire} and this section. 
    \begin{itemize}
        \item The \textit{extensions} $d_k^{f, E_{r + 1}}$ in \cite[\S~5]{Lin-Wang-Xu-kervaire} correspond to the \textit{$0$-th extensions} $d_k^{f, E_{r + 1}} = d_k^{f, E_{r + 1}, 0}$ in Definition \ref{extensions-def}. The $j$-th layer extensions for $j \geq 1$ in the current paper briefly appear in \cite[\S~4]{Lin-Wang-Xu-kervaire} as certain extension differentials between the $E_{\infty}$ pages of synthetic Adams SS, which are not the focus of their study.  
        \item A \textit{crossing} for $d_{k}^{f, E_{r + 1}, 0}(x) = y$ in the sense of \cite[\S~5]{Lin-Wang-Xu-kervaire} is a pair 
        \[x' \in Z_{r - a}^{s + a, t + a}(X), y' \neq 0 \in Z_{r - k + b}^{s + k - b, t + k - b}(Y) / B^{s + k - b, t + k - b}_{1 + k - a - b}(Y)\] 
        with $0 < a \leq k, 0 \leq b \leq k - a$ so that we have an essential\footnote{here \textit{essential} in \cite[\S~5]{Lin-Wang-Xu-kervaire} means there is no $z \in Z_{r - a_1}^{s + a_1, t + a_1}(X)$ with $a < a_1 \leq k - b$ so that $d_{k - a_1 - b}^{f, E_{r + 1 - a_1}, 0}(z) = y'$} extension $d_{k - a - b}^{f, E_{r + 1 - a}, 0}(x') = y'$, and the original extension $d_{k}^{f, E_{r + 1}, 0}(x) = y$ has no crossing if no such crossing extension exists. It follows from the same discussion as in Fact \ref{crossing-facts} below that if $d_{k}^{f, E_{r + 1}, 0}(x) = y$ has no crossing in the sense of \cite[\S~5]{Lin-Wang-Xu-kervaire}, then it has no crossing in the sense of Definition \ref{extensions-def}. The converse does not hold true in general.
    \end{itemize}
\end{remark}

\begin{example} \label{example-extn-crossing}

    We now present an example illustrating most of the possible behaviors of extensions and crossings. Its conclusions may serve as a first guide to these phenomena, while its arguments are best appreciated after the general results developed in the remainder of this section and the next have been established. To formulate this example, we import data from \(\mathrm{H}\mathbb{F}_2\)-synthetic spectra via the functor \(U\colon \Syn_{\mathrm{H}\mathbb{F}_2}\to \Fil\Sp\), according to Remark \ref{synthetic-spectra-as-special-case} (here everything is $2$-local). 
    \begin{itemize}
        \item  Write \(\mathbb{S}^{0,0}_{\mathrm{H}\mathbb{F}_2}=U(\nu_{\mathrm{H}\mathbb{F}_2}\mathbb{S}^0)\), and let \(\mathbb{S}^{n,w}_{\mathrm{H}\mathbb{F}_2}=\Sigma^{n,w}\mathbb{S}^{0,0}_{\mathrm{H}\mathbb{F}_2}\). 
        \item Since the stem \(\nu\in\pi_3(\mathbb{S})\) has Adams filtration \(1\), by \cite[Lemma 9.15]{BHS1} it lifts to a map \([h_2]\colon \mathbb{S}^{3,4}_{\mathrm{H}\mathbb{F}_2}\to \mathbb{S}^{0,0}_{\mathrm{H}\mathbb{F}_2}\) whose induced map on graded pieces, equivalently whose extension of filtration jump \(0\), recovers multiplication by \(h_2\) on the $\mathrm{H}\Fb_2$-Adams \(E_2\)-page.
        \item By \cite[Proposition 3.20]{Lin-Wang-Xu-kervaire}, this map fits into a distinguished triangle
        \begin{equation*} \label{delta-dt}\tag{$\mkern0.5mu\blacktriangle\mkern-0.5mu$}
        \mathbb{S}^{3,4}_{\mathrm{H}\mathbb{F}_2}\xrightarrow{[h_2]}\mathbb{S}^{0,0}_{\mathrm{H}\mathbb{F}_2}\xrightarrow{i_\blacktriangle}\cofib(\nu)_{\mathrm{H}\mathbb{F}_2}\xrightarrow{p_\blacktriangle}\mathbb{S}^{4,4}_{\mathrm{H}\mathbb{F}_2}
        \end{equation*}
        in which $\cofib(\nu)_{\mathrm{H}\mathbb{F}_2} = U(\nu_{\mathrm{H}\mathbb{F}_2}\cofib(\nu\colon \Sb^3 \to \Sb^0))$.
    \end{itemize} 
    We will compute all extensions along the (shifted) map
    \[
    f=[h_2]\colon X=\mathbb{S}^{0,0}_{\mathrm{H}\mathbb{F}_2}\to Y=\mathbb{S}^{-3,-4}_{\mathrm{H}\mathbb{F}_2}
    \]
    in stem \(45\), and record the result in Figure~\ref{figure-crossings-extn}. 
    \begin{itemize}
        \item Here \(\{E_r^{s,t}(X)\}_{r \geq 2}\) is the \(\mathrm{H}\mathbb{F}_2\)-Adams spectral sequence of \(\mathbb{S}^0\), while \(E_r^{s,t}(Y) = E_{r, {\mathrm{ASS}}}^{s + 1, t + 4}(\Sb^0)\). 
        \item We draw our figures with horizontal coordinate \(t-s\) and vertical coordinate \(s\), and refer to \((t-s,s)\) as the bidegree in the picture. We also label both copies according to the indexing on the Adams spectral sequence of \(\mathbb{S}^0\). Thus, for example, a class marked as having bidegree \((48,7)\) in the right-hand column actually lies in \(E_2^{6,45+6}(Y) = E_{2, {\mathrm{ASS}}}^{7,48+7}(\Sb^0)\). 
    \end{itemize}
    In the discussion of this example, we will freely use the generalized Mahowald trick (more precisely, Theorem~\ref{GMT-for-SS}) to construct extensions, together with Theorem~\ref{stretching-extensions-across-pages} to adjust their pages and/or layers. We will also appeal to the \(\mathrm{H}\mathbb{F}_2\)-Adams spectral sequence computations in \cite{Lin-Wang-Xu-machine}, which are available on Lin's \href{https://waynelin92.github.io/ss/kervaire-49.html}{website}. 

    \begin{figure}[htbp]
        \centering
        \scalebox{1}{
        \begin{tikzpicture}[line width=0.1pt]
            \tikzset{
                diff/.style={-{Stealth},line width=0.3pt,shorten <=3pt,shorten >=3pt},
                diff1/.style={-{Stealth},dashed,line width=0.6pt,shorten <=3pt,shorten >=3pt},
                extarr/.style={-{Stealth},line width=0.3pt,shorten <=3pt,shorten >=3pt},
                extarr1/.style={-{Stealth},dashed,line width=0.6pt,shorten <=3pt,shorten >=3pt},
                mapbelow/.style={-{Straight Barb[scale=0.8]},line width=0.7pt,shorten <=.7cm,shorten >=.7cm},
            }
            
            %======================
            % left column
            %======================
            \draw (44.5,2.5) rectangle (45.5,17.5);
            \node at (45,2) {$45$};
            \node at (44,3) {$3$};
            \node at (44,4) {$4$};
            \node at (44,5) {$5$};
            \node at (44,6) {$6$};
            \node at (44,7) {$7$};
            \node at (44,8) {$8$};
            \node at (44,9) {$9$};
            \node at (44,10) {$10$};
            \node at (44,11) {$11$};
            \node at (44,12) {$12$};
            \node at (44,13) {$13$};
            \node at (44,14) {$14$};
            \node at (44,15) {$15$};
            \node at (44,16) {$16$};
            \node at (44,17) {$17$};
            
            \begin{scope}
            \clip (44.5,2.5) rectangle (45.5,17.5);
            \draw[black!10] (44.5,3.5) -- (45.5,3.5);
            \draw[black!10] (44.5,4.5) -- (45.5,4.5);
            \draw[black!10] (44.5,5.5) -- (45.5,5.5);
            \draw[black!10] (44.5,6.5) -- (45.5,6.5);
            \draw[black!10] (44.5,7.5) -- (45.5,7.5);
            \draw[black!10] (44.5,8.5) -- (45.5,8.5);
            \draw[black!10] (44.5,9.5) -- (45.5,9.5);
            \draw[black!10] (44.5,10.5) -- (45.5,10.5);
            \draw[black!10] (44.5,11.5) -- (45.5,11.5);
            \draw[black!10] (44.5,12.5) -- (45.5,12.5);
            \draw[black!10] (44.5,13.5) -- (45.5,13.5);
            \draw[black!10] (44.5,14.5) -- (45.5,14.5);
            \draw[black!10] (44.5,15.5) -- (45.5,15.5);
            \draw[black!10] (44.5,16.5) -- (45.5,16.5);
            
            % lower filtration part: unchanged from previous example
            \draw (45,3) -- (45,4);
            \draw (45.2,5.041) -- (45,6);
            \draw (45,6) -- (45,7);
            
            % upper filtration part: from extended column
            \draw (45,15) -- (45,16);
            \draw (45,16) -- (45,17);
            \draw[ForestGreen, diff] (45,12) -- (44,16);
            \draw[red, diff] (46,14) -- (45,17);
            \draw[cyan, diff] (45,15) -- (44,17);
            \draw[cyan, diff] (45,16) -- (44,18);
            
            \coordinate (ext1_L_45_3) at (45,3);
            \coordinate (ext1_L_45_4) at (45,4);
            \coordinate (ext1_L_45_5ur) at (44.8,4.959);
            \coordinate (ext1_L_45_5ul) at (45.2,5.041);
            \coordinate (ext1_L_45_6) at (45,6);
            \coordinate (ext1_L_45_7) at (45,7);
            \coordinate (ext1_L_45_9) at (45,9);
            \coordinate (ext1_L_45_12) at (45,12);
            \coordinate (ext1_L_45_15) at (45,15);
            \coordinate (ext1_L_45_16) at (45,16);
            \coordinate (ext1_L_45_17) at (45,17);
            
            \fill (ext1_L_45_3) circle (0.064) node[below] {\resizelabel{0.667}{$h_3^2h_5$}};
            \fill (ext1_L_45_4) circle (0.064);
            \fill (ext1_L_45_5ur) circle (0.064) node[above] {\resizelabel{0.6}{$h_1 g_2$}};
            \fill (ext1_L_45_5ul) circle (0.064) node[below] {\resizelabel{0.6}{$d_0h_5$}};
            \fill (ext1_L_45_6) circle (0.064);
            \fill (ext1_L_45_7) circle (0.064);
            \fill (ext1_L_45_9) circle (0.064) node[below] {\resizelabel{0.667}{$\Delta h_1g$}};
            \fill (ext1_L_45_12) circle (0.064) node[below] {\resizelabel{0.667}{$d_0^2e_0$}};
            \fill (ext1_L_45_15) circle (0.064) node[below] {\resizelabel{0.667}{$Pd_0i$}};
            \fill (ext1_L_45_16) circle (0.064);
            \fill (ext1_L_45_17) circle (0.064);
            \end{scope}
            
            %======================
            % right column
            %======================
            \begin{scope}[xshift=3cm,yshift=-1cm]
            
            \draw (47.5,3.5) rectangle (48.5,18.5);
            \node at (48,3) {$48$};
            \node at (49,4) {$4$};
            \node at (49,5) {$5$};
            \node at (49,6) {$6$};
            \node at (49,7) {$7$};
            \node at (49,8) {$8$};
            \node at (49,9) {$9$};
            \node at (49,10) {$10$};
            \node at (49,11) {$11$};
            \node at (49,12) {$12$};
            \node at (49,13) {$13$};
            \node at (49,14) {$14$};
            \node at (49,15) {$15$};
            \node at (49,16) {$16$};
            \node at (49,17) {$17$};
            \node at (49,18) {$18$};
            
            \begin{scope}
            \clip (47.5,3.5) rectangle (48.5,18.5);
            \draw[black!10] (47.5,4.5) -- (48.5,4.5);
            \draw[black!10] (47.5,5.5) -- (48.5,5.5);
            \draw[black!10] (47.5,6.5) -- (48.5,6.5);
            \draw[black!10] (47.5,7.5) -- (48.5,7.5);
            \draw[black!10] (47.5,8.5) -- (48.5,8.5);
            \draw[black!10] (47.5,9.5) -- (48.5,9.5);
            \draw[black!10] (47.5,10.5) -- (48.5,10.5);
            \draw[black!10] (47.5,11.5) -- (48.5,11.5);
            \draw[black!10] (47.5,12.5) -- (48.5,12.5);
            \draw[black!10] (47.5,13.5) -- (48.5,13.5);
            \draw[black!10] (47.5,14.5) -- (48.5,14.5);
            \draw[black!10] (47.5,15.5) -- (48.5,15.5);
            \draw[black!10] (47.5,16.5) -- (48.5,16.5);
            \draw[black!10] (47.5,17.5) -- (48.5,17.5);
            
            % lower filtration part: unchanged from previous example
            \draw (48,5) -- (48,6);
            \draw (48,6) -- (48.113,6.9);
            \draw (47.887,7.1) -- (48,8);
            \draw (48,8) -- (48.113,8.9);
            \draw[red,diff] (49,6) -- (48.113,8.9);
            \draw[cyan,diff] (49,5) -- (48.113,6.9);
            \draw[cyan, diff] (48,4) -- (47,6);
            \draw[cyan, diff] (48,5) -- (47,7);
            
            % upper filtration part: from extended column
            % \draw (47.113,12.959) -- (48,14);
            % \draw (46.887,13.041) -- (48,14);
            \draw (48,15) -- (48,16);
            % \draw (47.113,15.959) -- (48,17);
            \draw (48,16) -- (48,17);
            \draw[ForestGreen, diff] (48,17) -- (47.113,20.959);
            \draw[red, diff] (49,11) -- (48,14);
            \draw[cyan, diff] (48,15) -- (47.113,16.959);
            \draw[cyan, diff] (48,16) -- (47.113,17.959);
            
            \coordinate (ext1_R_48_4) at (48,4);
            \coordinate (ext1_R_48_5) at (48,5);
            \coordinate (ext1_R_48_6) at (48,6);
            \coordinate (ext1_R_48_7ur) at (48.113,6.9);
            \coordinate (ext1_R_48_7ul) at (47.887,7.1);
            \coordinate (ext1_R_48_8) at (48,8);
            \coordinate (ext1_R_48_9ur) at (48.113,8.9);
            \coordinate (ext1_R_48_9ul) at (47.887,9.1);
            \coordinate (ext1_R_48_12) at (48,12);
            \coordinate (ext1_R_48_14) at (48,14);
            \coordinate (ext1_R_48_15) at (48,15);
            \coordinate (ext1_R_48_16) at (48,16);
            \coordinate (ext1_R_48_17) at (48,17);
            
            \fill (ext1_R_48_4) circle (0.064) node[below] {\resizelabel{0.667}{$h_3c_2$}};
            \fill (ext1_R_48_5) circle (0.064) node[below] {\resizelabel{0.667}{$e_0h_5$}};
            \fill (ext1_R_48_6) circle (0.064);
            \fill (ext1_R_48_7ur) circle (0.064);
            \fill (ext1_R_48_8) circle (0.064);
            \fill (ext1_R_48_9ur) circle (0.064);
            \fill (ext1_R_48_9ul) circle (0.064) node[above] {\resizelabel{0.667}{$h_1 Pc_0 h_5$}};
            \fill (ext1_R_48_12) circle (0.064) node[below] {\resizelabel{0.667}{$d_0^2g$}};
            \fill (ext1_R_48_14) circle (0.064) node[above] {\resizelabel{0.667}{$y_{48,14}$}};
            \fill (ext1_R_48_15) circle (0.064) node[below] {\resizelabel{0.667}{$d_0Pj$}};
            \fill (ext1_R_48_16) circle (0.064);
            \fill (ext1_R_48_17) circle (0.064);
            \fill (ext1_R_48_7ul) circle (0.064) node[above] {\resizelabel{0.667}{$Mh_2$}};
            \end{scope}
            
            \end{scope}
            
            %======================
            % extension arrow
            %======================
            \draw[ForestGreen, extarr] (ext1_L_45_3) -- (ext1_R_48_7ul);
            \draw[ForestGreen, extarr] (ext1_L_45_4) -- (ext1_R_48_8);
            \draw[gray, extarr] (ext1_L_45_5ul) -- (ext1_R_48_6);
            \draw[gray, extarr] (ext1_L_45_6) -- (ext1_R_48_7ur);
            \draw[cyan, extarr] (ext1_L_45_7) -- (ext1_R_48_9ur);
            \draw[red, extarr] (ext1_L_45_9) -- (ext1_R_48_12);
            \draw[cyan, extarr] (ext1_L_45_12) -- (ext1_R_48_14);
            \draw[red, extarr] (ext1_L_45_12) -- (ext1_R_48_15);
            \draw[gray, extarr] (ext1_L_45_15) -- (ext1_R_48_16);
            \draw[gray, extarr] (ext1_L_45_16) -- (ext1_R_48_17);
            
            %======================
            % map X -> Y
            %======================

            \coordinate (label-X) at (45,0.8);
            \coordinate (label-Y) at (51,0.8);
            
            \node at (label-X) {$X$};
            \node at (label-Y) {$Y$};
            \draw[mapbelow] (label-X) -- (label-Y)
              node[midway,above] {\footnotesize $f = [h_2]$};
        \end{tikzpicture}
        }
        \caption{Extensions along $[h_2]\colon X = \Sb^{0, 0}_{\mathrm{H} \Fb_2} \to Y = \Sb^{-3, -4}_{\mathrm{H} \Fb_2}$ at stem $45$, and their crossings.}
        \label{figure-crossings-extn}
    \end{figure}
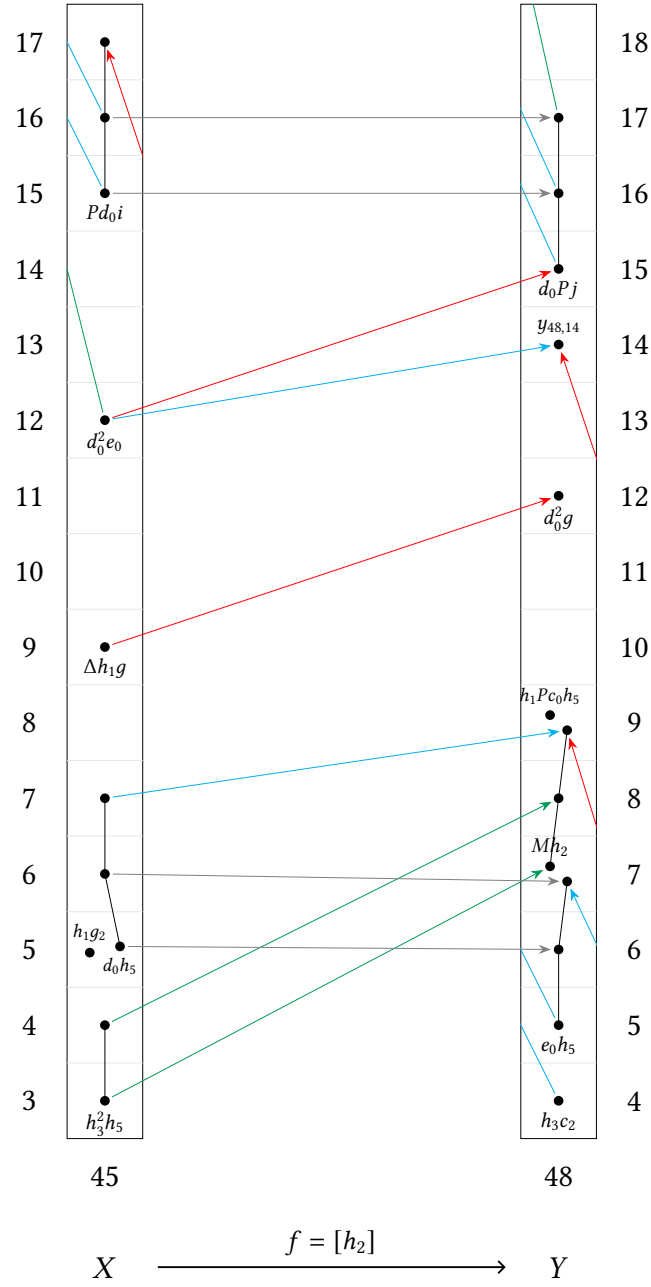

    \parr

    We start from the extensions with filtration jump $= 0$. 

    \begin{itemize}
        \item The class $P d_0 i$ of bidegree $(t - s, s) = (45, 15)$ supports an essential $d_2$-differential. It also supports an essential extension $\color{gray} d_0^{f, E_{2}, 0}(P d_0 i) = h_0 d_0 Pj$ due to the $E_2$-page relation $h_2 i = h_0 j$. In fact, $d_0^{f, E_{2}, u}(P d_0 i) = h_0 d_0 Pj$ is an essential extension for any $u \in \Nb$. 
        \item Similarly, the class $h_0 P d_0 i$ of bidegree $(t - s, s) = (45, 16)$ supports an essential $d_2$-differential, and it also supports an essential extension $\color{gray} d_0^{f, E_{2}, 0}(h_0 P d_0 i) = h_0^2 d_0 Pj$. Also, $d_0^{f, E_{2}, u}(h_0 P d_0 i) = h_0^2 d_0 Pj$ is an essential extension for any $u \in \Nb$. 
        \item The class $d_0 h_5$ of bidegree $(t - s, s) = (45, 5)$ is a permanent cycle, and it supports an essential extension $\color{gray} d_0^{f, E_{\infty}, 0}(d_0 h_5) = h_0 e_0 h_5$. Furthermore, $d_0^{f, E_{r + 1}, u}(d_0 h_5) = h_0 e_0 h_5$ is an essential extension for any $r \in \Nb \cup \{\infty\}, u \in \Nb$. To prove this, note that we have $d_0^{f, E_{2}, 0}(d_0 h_5) = h_0 e_0 h_5$ due to the $E_2$-page relation $h_2 d_0 = h_0 e_0$, and the rest follows from Theorem \ref{stretching-extensions-across-pages}.
        \item Similarly, the class $h_0 d_0 h_5$ of bidegree $(t - s, s) = (45, 6)$ is a permanent cycle, and it supports an essential extension $\color{gray} d_0^{f, E_{\infty}, 0}(h_0 d_0 h_5) = h_0^2 e_0 h_5$. However, since $h_0^2 e_0 h_5 = d_2([f_0] h_5)$ is a boundary class, $d_0^{f, E_{\infty}, 1}(h_0 d_0 h_5) = h_0^2 e_0 h_5$ is no longer essential. Actually, $d_k^{f, E_{r + 1}, u}(h_0 d_0 h_5) = 0$ makes sense for all $r \in \Nb \cup \{\infty\}, k, u \in \Nb$ with $u \geq 1$. 
        
        To prove this, by Theorem \ref{stretching-extensions-across-pages} it suffices to establish $d_k^{f, E_{\infty}, 1}(h_0 d_0 h_5) = 0$ for all $k$. Indeed, consider the distinguished triangle (\ref{delta-dt}). As $[f_0] h_5[0] = i_\blacktriangle([f_0] h_5)$ is a surviving permanent cycle in $E_{*, \mathrm{ASS}}^{*, *}(\cofib(\nu))$, it follows that $i_\blacktriangle(\delta_1^{\infty}([f_0] h_5)) = \delta_1^{\infty}([f_0] h_5[0]) = 0$, so $\delta_1^{\infty}([f_0] h_5)$ is in $\IIm(f)$. Since $\delta_1^{\infty}([f_0] h_5)$, being a lift of $h_0^2 e_0 h_5$, is not $\defopara$-divisible, the only possible situation is $\delta_1^{\infty}([f_0] h_5) = f[h_0 d_0 h_5]$ for some $[h_0 d_0 h_5] \in \pi_{45, 45 + 6}(X)$ lifting $h_0 d_0 h_5$. It follows that $\defopara f[h_0 d_0 h_5] = 0$, which is what we want.
    \end{itemize}

    Next we treat the classes that do not support any essential extension. 

    \begin{itemize}
        \item The class $h_0^2 P d_0 i$ of bidegree $(t - s, s) = (45, 17)$ is a permanent cycle, and it does not support any essential extension along $f$, in other words $d^{f, E_{r + 1}, u}_k(h_0^2 P d_0 i) = 0$ makes sense for all $r \in \Nb \cup \{\infty\}, k, u \in \Nb$. To prove this, by Theorem \ref{stretching-extensions-across-pages} it suffices to establish $d^{f, E_{\infty}, 0}_k(h_0^2 P d_0 i) = 0$ for all $k$. As the class $h_0^2 P d_0 i[4] = h_1 d_0P^2 d_0[4] \in E_{2, \mathrm{ASS}}^{17, 49+17}(\cofib(\nu))$ is a permanent cycle (indeed, a $d_4$-boundary), by Theorem \ref{Bockstein-dictionary-infinite} there exists a lift $[h_0^2 P d_0 i[4]]$ so that $p_\blacktriangle([h_0^2 P d_0 i[4]]) = [h_0^2 P d_0 i] \in \pi_{17, 45 + 17}(X)$ is a lift of $h_0^2 P d_0 i$. It follows that $f[h_0^2 P d_0 i] = 0$, which is what we want.
        \item The class \(h_1g_2\) of bidegree \((t-s,s)=(45,5)\) is a permanent cycle, and it supports no essential extension along \(f\), in other words \(d_k^{f,E_{r+1},j}(h_1g_2)=0\) makes sense for all \(r\in \mathbb{N}\cup\{\infty\}\) and \(k,j\in \mathbb{N}\). This follows from the same argument: since \(h_1g_2[4]\in E_{2,\mathrm{ASS}}^{5,49+5}(\cofib(\nu))\) is a permanent cycle, one obtains a lift \([h_1g_2]\in \pi_{45,45+5}(X)\) satisfying \(f[h_1g_2]=0\).
        \item It also makes sense to write, for instance, \(d_{10}^{f,E_{12},0}(h_1g_2)=h_0d_0Pj\). Indeed, one may choose a lift \([h_1g_2]\in \pi_{45,45+5}(X/\defopara^{11})\) such that \(f[h_1g_2]=0\), then \([h_1g_2]_1=[h_1g_2]+\defopara^{10}Pd_0i\) is another lift of \(h_1g_2\) such that \(f[h_1g_2]_1=\defopara^{10}h_0d_0Pj\). This is a typical inessential extension.
    \end{itemize}

    We then derive the essential extensions with nontrivial filtration jumps by applying the generalized Mahowald trick (Theorem \ref{GMT-for-SS}) to the distinguished triangle (\ref{delta-dt}).

    \begin{itemize}
        \item The class $h_3^2 h_5$ of bidegree $(t - s, s) = (45, 3)$ is a permanent cycle, and it supports an essential extension $\color{ForestGreen}{d_3^{f, E_{\infty}, 0}(h_3^2 h_5) = Mh_2}$. In fact, $d_3^{f, E_{r + 1}, u}(h_3^2 h_5) = Mh_2$ is an essential extension for any $4 \leq r \leq \infty, u \in \Nb$. 
        
        As for the proof, note that in the distinguished triangle (\ref{delta-dt}), there exist classes $Mh_2[0] \in E_{2,\mathrm{ASS}}^{7, 48 + 7}(\cofib(\nu))$, $h_3^2 h_5[4] \in E_{2,\mathrm{ASS}}^{3, 49 + 3}(\cofib(\nu))$ such that $i_\blacktriangle(Mh_2) = Mh_2[0]$ and $p_\blacktriangle(h_3^2 h_5[4]) = h_3^2 h_5$. Also, $\color{ForestGreen} d_4(h_3^2 h_5[4]) = Mh_2[0]$ is an essential differential in $E_{*,\mathrm{ASS}}^{*,*}(\cofib(\nu))$. Thus, the extension $\color{ForestGreen}{d_3^{f, E_{5}, 0}(h_3^2 h_5) = Mh_2}$ can be derived from this differential via Theorem \ref{GMT-for-SS}. The other claimed extensions follow from Theorem \ref{stretching-extensions-across-pages}, whose assumption is satisfied since all potential crossing classes for this $E_5$-page extension are permanent cycles.
        \item The class $h_0 h_3^2 h_5$ of bidegree $(t - s, s) = (45, 4)$ is a permanent cycle, and it supports an essential extension $\color{ForestGreen}{d_3^{f, E_{\infty}, 0}(h_0 h_3^2 h_5) = h_0 Mh_2}$. In fact, $d_3^{f, E_{r + 1}, u}(h_0 h_3^2 h_5) = h_0 Mh_2$ is an essential extension for any $4 \leq r \leq \infty, u \in \Nb$. This is proved by the same GMT argument, applied to the differential \(\color{ForestGreen} d_4(h_0h_3^2h_5[4])=h_0Mh_2[0]\) in \(E_{*,\mathrm{ASS}}^{*,*}(\cofib(\nu))\).
        \item The class $\Delta h_1 g$ of bidegree $(t - s, s) = (45, 9)$ is a permanent cycle, and it supports an essential extension $\color{red} d_2^{f, E_{\infty}, 0}(\Delta h_1 g) = d_0^2 g$. In fact, $d_2^{f, E_{r + 1}, u}(\Delta h_1 g) = d_0^2 g$ is an essential extension for any $3 \leq r \leq \infty, u \in \Nb$. This is also proved by the same GMT argument, applied to the differential \(\color{red} d_3(\Delta h_1 g[4]) = d_0^2 g[0]\) in \(E_{*,\mathrm{ASS}}^{*,*}(\cofib(\nu))\).
        \item The class $h_0^2 d_0 h_5$ of bidegree $(t - s, s) = (45, 7)$ is a permanent cycle, and it supports an essential extension $\color{cyan} d_1^{f, E_{\infty}, 0}(h_0^2 d_0 h_5) = h_0^2 M h_2$. In fact, $d_1^{f, E_{r + 1}, 0}(h_0^2 d_0 h_5) = h_0^2 M h_2$ is an essential extension for any $2 \leq r \leq \infty$. For $u \geq 1$, as $h_0^2 M h_2 = d_3(h_0[f_0]h_5)$ is a boundary term, the extension $d_1^{f, E_{\infty}, u}(h_0^2 d_0 h_5) = h_0^2 M h_2$ is no longer essential, and indeed $d_k^{f, E_{r + 1}, u}(h_0^2 d_0 h_5) = 0$ makes sense for all $r \in \Nb \cup\{\infty\}, k, u \in \Nb$ with $u \geq 1$. 

        To prove this, note that the $d_3$-differential implies $\delta_1^{\infty}(h_0[f_0]h_5) = \defopara [h_0^2 M h_2]$ due to Theorem \ref{delta-as-total-diff}. Thus, by Theorem \ref{stretching-extensions-across-pages} it suffices to find a lift $[h_0^2 d_0 h_5] \in \pi_{45, 45 + 7}(X)$ so that $f[h_0^2 d_0 h_5] = \delta_1^{\infty}(h_0[f_0]h_5)$. Indeed, as $h_0[f_0]h_5[0]$ is a permanent cycle in $E_{*,\mathrm{ASS}}^{*,*}(\cofib(\nu))$, it follows that $i_\blacktriangle(\delta_1^{\infty}(h_0[f_0]h_5)) = \delta_1^{\infty}(h_0[f_0]h_5[0]) = 0$, so $\delta_1^{\infty}(h_0[f_0]h_5)$ is in $\IIm(f)$. Since $\delta_1^{\infty}(h_0[f_0]h_5) = \defopara [h_0^2 M h_2]$ is not $\defopara^2$-divisible and $E_2^{8, 45 + 8}(X) = 0$, the only possible situation is $\delta_1^{\infty}(h_0[f_0]h_5) = f[h_0^2 d_0 h_5]$ for some $[h_0^2 d_0 h_5] \in \pi_{45, 45 + 7}(X)$ lifting $h_0^2 d_0 h_5$, as expected.
    \end{itemize}

    It remains to examine the \(f\)-extensions of the class \(d_0^2e_0\) of bidegree \((t-s,s)=(45,12)\), which supports an essential \(d_4\). This special case calls for the introduction of \(u\)-th extensions for \(u\geq 1\).

    \begin{itemize}
        \item There is an essential extension $\color{cyan} d_1^{f, E_{4}, 0}(d_0^2 e_0) = y_{48,14} = h_1(\Delta h_0^2 i + P\Delta h_1 d_0)$, which also gives an essential extension $d_1^{f, E_{3}, 0}(d_0^2 e_0) = y_{48,14}$. As $y_{48,14} = d_3(d_0 m)$, the extension $d_1^{f, E_{4}, 1}(d_0^2 e_0) = y_{48,14}$ is no longer essential; instead, we obtain a longer essential extension $\color{red} d_2^{f, E_{4}, 1}(d_0^2 e_0) = d_0 Pj$. Furthermore, $d_2^{f, E_{4}, u}(d_0^2 e_0) = d_0 Pj$ is an essential extension for all $u \geq 1$.

        To prove these claims, note that there is an essential differential $d_4(d_0 m[0]) = - d_0 P j[0]$ in $E_{*,\mathrm{ASS}}^{*,*}(\cofib(\nu))$. 
        It follows that $\delta_1^3(d_0 m[0]) = - \defopara^2 d_0 P j[0]$ due to Theorem \ref{delta-as-total-diff}.
        Thus, $i_\blacktriangle(\delta_1^3(d_0 m) + \defopara^2 d_0 P j) = \delta_1^3(d_0 m[0]) + \defopara^2 d_0 P j[0] = 0$, 
        in other words $\delta_1^3(d_0 m) + \defopara^2 d_0 P j$ is in $\IIm(f)$. 
        Note that the essential differential $d_3(d_0 m) = y_{48, 14}$ in $E_{*,\mathrm{ASS}}^{*,*}(\Sb^0)$ implies that $\delta_1^3(d_0 m) = \defopara [y_{48, 14}]$, 
        so the class $\delta_1^3(d_0 m) + \defopara^2 d_0 P j$ is not $\defopara^2$-divisible. 
        Together with the fact $E_2^{13, 45 + 13}(X) = 0$, we deduce that the only possible situation is $\delta_1^3(d_0 m) + \defopara^2 d_0 P j = f[d_0^2 e_0]$ for some $[d_0^2 e_0] \in \pi_{45, 45 + 12}(X / \defopara^3)$ lifting $d_0^2 e_0$. 
        This relation witnesses both essential extensions $d_1^{f, E_{4}, 0}(d_0^2 e_0) = y_{48,14}$ and $d_2^{f, E_{4}, 1}(d_0^2 e_0) = d_0 Pj$, while the rest follows from Theorem \ref{stretching-extensions-across-pages}.
        
    \end{itemize}

    These exhaust all essential extensions along the map $f = [h_2]$ at stem $45$. As for crossings: 

    \begin{itemize}
        \item A crossing of the extension \(\color{ForestGreen}{d_3^{f,E_{\infty},0}(h_3^2h_5)=Mh_2}\) is an essential extension \(d_{3-a-b}^{f,E_{\infty},a}(z)=w\), where {\(1\leq a\leq 3\), \(0\leq b\leq 3-a\)}, \(z\in Z_{\infty}^{3+a,45+3+a}(X)\), and \(w\in Z_{\infty}^{6-b,45+6-b}(Y)=Z_{{\infty},\mathrm{ASS}}^{7-b,48+7-b}(\mathbb{S}^0)\). From the figure, there are five candidate classes in the possible bidegrees of \(z\), namely \(h_0h_3^2h_5\), \(h_1g_2\), \(d_0h_5\), \(h_0d_0h_5\), and \(h_0^2d_0h_5\). Among these, only \(h_0d_0h_5\) and \(h_0^2d_0h_5\) support essential extensions whose targets lie in filtration at most \(7\) of \(E_{*,\mathrm{ASS}}^{*,*}(\mathbb{S}^0)\).
        
        The essential extension $\color{gray} d_0^{f, E_{\infty}, 2}(d_0 h_5) = h_0 e_0 h_5$ is a crossing for the original extension. On the other hand, the essential extension $\color{gray} d_0^{f, E_{\infty}, 0}(h_0 d_0 h_5) = h_0^2 e_0 h_5$ seems\footnote{It is indeed classified as a crossing in the sense of \cite[Definition 5.9]{Lin-Wang-Xu-kervaire}.} to contribute to a crossing as well. However, $d_0^{f, E_{\infty}, 3}(h_0 d_0 h_5) = h_0^2 e_0 h_5$ is not essential anymore, so the class $h_0 d_0 h_5$ cannot support any crossing for the original extension in the sense of Definition \ref{extensions-def}. 

        Therefore, $\color{gray} d_0^{f, E_{\infty}, 2}(h_0 d_0 h_5) = h_0^2 e_0 h_5$ is the only crossing extension for $\color{ForestGreen}{d_3^{f, E_{\infty}, 0}(h_3^2 h_5) = Mh_2}$. In general, for each $4 \leq r \leq \infty, u \in \Nb$, the extension $d_3^{f, E_{r + 1}, u}(h_3^2 h_5) = Mh_2$ has only one crossing $d_0^{f, E_{r - 1}, u + 2}(h_0 d_0 h_5) = h_0^2 e_0 h_5$.
        
        \item Similarly, $\color{gray} d_0^{f, E_{\infty}, 1}(h_0 d_0 h_5) = h_0^2 e_0 h_5$ is the only crossing extension for $\color{ForestGreen}{d_3^{f, E_{\infty}, 0}(h_0h_3^2 h_5) = h_0Mh_2}$, and in general for each $4 \leq r \leq \infty, u \in \Nb$, the extension $d_3^{f, E_{r + 1}, u}(h_0 h_3^2 h_5) = h_0 Mh_2$ has only one crossing $d_0^{f, E_{r}, u + 1}(h_0 d_0 h_5) = h_0^2 e_0 h_5$. Here the plausible interference $\color{gray} d_0^{f, E_{\infty}, 0}(h_0 d_0 h_5) = h_0^2 e_0 h_5$ again acts as a ``fake crossing'' since $d_0^{f, E_{\infty}, 2}(h_0 d_0 h_5) = h_0^2 e_0 h_5$ is not essential anymore.
        
        \item All the other essential extensions have full no crossing for degree reasons.
    \end{itemize}

    \FloatBarrier
    
\end{example}

\begin{fact} \label{extension-abundance}
    For any map $f\colon X \to Y$ in $\Fil\Sp$, there is an abundant supply of extensions.
    \begin{itemize}
        \item In the case $r < \infty$, it follows from Theorem~\ref{standard-SS} that every class $x \in Z_r^{s,t}(X)$ admits a lift $[x] \in \pi_{t-s,t}(X / \defopara^r)$. Furthermore, for any such lift $[x] \in \pi_{t-s,t}(X / \defopara^r)$ and any $j \in \Nb$, if $\defopara^j f[x] \in \pi_{t-s,t-j}(Y / \defopara^{r+j})$ is nonzero, then it lies in $\IIm(\defopara^{k+j}) \setminus \IIm(\defopara^{k+j+1})$ for a unique integer $0 \leq k < r$. Thus, by Theorem~\ref{Bockstein-dictionary-finite}, there exists a unique class $y \in Z_{r-k}^{s+k,t+k}(Y) / B_{1+k+j}^{s+k,t+k}(Y)$ together with a lift $[y] \in \pi_{t-s,t+k}(Y / \defopara^{r-k})$ so that $\defopara^j f[x] = \defopara^{k+j}[y]$, and this witnesses the extension $d_k^{f, E_{r+1}, j}(x) = y$.
    
        \item In the case $r = \infty$, if $x \in Z_{\infty}^{s,t}(X)$ lies in the image of $\rho_1^{\infty}\colon \pi_{t-s,t}(X) \to Z_{\infty}^{s,t}(X)$ (which is always the case if $X$ is complete and all weak obstructions in $E_*^{*,*}(X)$ vanish; cf.~Theorem~\ref{Bockstein-dictionary-infinite}), then tautologically there is some lift $[x] \in \pi_{t-s,t}(X)$. Furthermore, for any such lift $[x] \in \pi_{t-s,t}(X)$ and any $j \in \Nb$, either $\defopara^j f[x] \in \pi_{t-s,t-j}(Y)$ lies in $\bigcap_{a \in \Nb}\IIm(\defopara^a)$, or it lies in $\IIm(\defopara^{k+j}) \setminus \IIm(\defopara^{k+j+1})$ for a unique integer $k \in \Nb$. In the latter case (which is always the case if $Y$ is complete and all weak obstructions in $E_*^{*,*}(Y)$ vanish, provided that $\defopara^j f[x] \neq 0$), Theorem~\ref{Bockstein-dictionary-infinite} further yields a unique class $y \in Z_{\infty}^{s+k,t+k}(Y) / B_{1+k+j}^{s+k,t+k}(Y)$ together with a lift $[y] \in \pi_{t-s,t+k}(Y)$ such that $\defopara^j f[x] = \defopara^{k+j}[y]$, and this witnesses the extension $d_k^{f, E_{\infty}, j}(x) = y$.
    \end{itemize}
\end{fact}

\begin{fact} \label{crossing-facts}
    Take $f\colon X \to Y$ to be a map in $\Fil\Sp$. Take $r \in \Nb \cup \{\infty\}, j \in \Nb$. Suppose we have classes $x \in Z_r^{s, t}(X), y \in Z_{r - k}^{s + k, t + k}(Y)$ that support an extension $d_k^{f, E_{r + 1}, j}(x) = y$.
    \begin{enumerate}
        \item If $d_{k - a - b}^{f, E_{r + 1 - a}, j + a}(z) = w$ is a precrossing of this extension, then by definition we have $w \neq 0 \in Z_{r - k + b}^{s + k - b, t + k - b}(Y) / B_{1 + k + j - b}^{s + k - b, t + k - b}(Y)$, so it can be inessential only if there exists another extension $d_{k - a_1 - b}^{f, E_{r + 1 - a_1}, j + a_1}(z_1) = w$ with $a_1 > a$, which is also a precrossing. After sufficiently many iterations, we would end up with an essential extension $d_{k - a_n - b}^{f, E_{r + 1 - a_n}, j + a_n}(z_n) = w$, which turns out to be a crossing. In particular, the original extension has no crossing if and only if it has no precrossing. 
        \item As $\defopara^j f[x] = \defopara^{k + j} [y]$ implies $\defopara^{j + c} f[x] = \defopara^{k + j + c}[y]$ for $c \in \Nb$, an extension $d_k^{f, E_{r + 1}, j}(x) = y$ leads to extensions $d_k^{f, E_{r + 1}, j + c}(x) = y$ for all $c \in \Nb$. Furthermore, if $d_{k - a - b}^{f, E_{r + 1 - a}, j + a + c}(z) = w$ is a precrossing extension for $d_k^{f, E_{r + 1}, j + c}(x) = y$ witnessed by $\defopara^{j + a + c}f[z] = \defopara^{j + k + c - b}[w] \not\in \IIm(\defopara^{j + k + c - b + 1})$, then $\defopara^{j + a}f[z] \not\in \IIm(\defopara^{j + k - b + 1})$, so it would also witness a precrossing extension for $d_k^{f, E_{r + 1}, j}(x) = y$. In particular, if the original extension $d_k^{f, E_{r + 1}, j}(x) = y$ has full no crossing, then the induced extensions $d_k^{f, E_{r + 1}, j + c}(x) = y$ have full no crossing for all $c \in \Nb$. 
    \end{enumerate}
\end{fact}

To compute the extension of a certain class in the standard SS, it is necessary to choose one of its lifts in the bigraded homotopy group, while different choices of lifts could witness different extensions in principle. We interpret no crossing conditions as formulating the independence of such choices.

\begin{lemma} \label{free-choice-lemma-for-refined-extensions} 
    Suppose we have a $j$-th layer extension $d_k^{f, E_{r + 1}, j}(x) = y$ witnessed by $\defopara^j f[x]' = \defopara^{k + j} [y]'$. Then for each $[x] \in \pi_{t - s, t}(X / \defopara^r)$ lifting $x$, either $\defopara^j f[x] = \defopara^{k + j}[y]$ in which $[y] = [y]' + \defopara \varepsilon$ for some $\varepsilon \in \pi_{t - s, t + k + 1}(Y / \defopara^{r - k - 1})$, or $[x] = [x]' + \defopara^a [z]$ for some precrossing extension $d_{k - a - b}^{f, E_{r + 1 - a}, j + a}(z) = w$ witnessed by $\defopara^{j + a} f[z] = \defopara^{k + j - b} [w]$.
\end{lemma}

\begin{proof}
    Consider the difference $[x] - [x]'$. If it lies in $\IIm(\defopara^{k + 1})$ then the former case holds true. Otherwise, $[x] - [x]' \in \IIm(\defopara^a) \setminus \IIm(\defopara^{a + 1})$ for some $0 < a \leq k$, so Theorem \ref{Bockstein-dictionary-finite} (for $r < \infty$) or Theorem \ref{Bockstein-dictionary-infinite} (for $r = \infty$) implies that $[x] - [x]' = \defopara^a [z]$ for $[z] \in \pi_{t - s, t + a}(X / \defopara^{r - a})$ lifting some nonzero $z \in Z^{s + a, t + a}_{r - a} / B^{s + a, t + a}_{1 + a}$. If the image $\defopara^{j} f[z]$ lies in $\IIm(\defopara^{k + j + 1})$ then the former case also holds true. Otherwise, $\defopara^{j}f[z]  \in \IIm(\defopara^{k + j - b}) \setminus \IIm(\defopara^{k + j - b + 1})$ for some $0 \leq b \leq k - a$, so $\defopara^{j}f[z] = \defopara^{j + k - b}[w]$ for $[w]$ lifting some nonzero $w \in Z^{s + k - b, t + k - b}_{r - k + b} / B^{s + k - b, t + k - b}_{1 + k + j - b}$, therefore we extract a precrossing extension $d_{k - a - b}^{f, E_{r + 1 - a}, j + a}(z) = w$ as described in the latter case.
\end{proof}

\begin{corollary} \label{free-choice-lemma-for-refined-extensions-with-full-no-crossings} 
    The extension $d_k^{f, E_{r + 1}, j}(x) = y$ has no crossing if and only if for each $[x] \in \pi_{t - s, t}(X / \defopara^r)$ lifting $x$, there exists some $[y] \in \pi_{t - s, t + k}(Y / \defopara^{r - k})$ lifting $y$ so that $\defopara^j f[x] = \defopara^{k + j}[y]$.
\end{corollary}

\begin{lemma} \label{inessential-crossing-extensions}
    Suppose $f\colon X \to Y$ is a map in $\Fil\Sp$. Take $s, t \in \Zb$, $r, k, j \in \Nb$ so that $k < r$. Suppose we have classes $x \in Z_r^{s, t}(X)$, $y \in Z^{s + k, t + k}_{r - k}$  so that $d_k^{f, E_{r + 1}, j}(x) = y$. Then for each $0 < a < r$ and each $z \in Z^{s + a, t + a}_{r - a}(X)$, if $z$ is not the source of a crossing extension for $d_k^{f, E_{r + 1}, j}(x) = y$, then there is a lift $[z]_{\std} \in \pi_{t - s, t + a}(X / \defopara^{r - a})$ such that $\defopara^{j + a} f[z]_{\std} \in \IIm(\defopara^{k + j + 1})$.
\end{lemma}
\begin{proof}
    According to Lemma \ref{free-choice-lemma-for-refined-extensions}, for any lift $[z] \in \pi_{t - s, t + a}(X / \defopara^{r - a})$, either $\defopara^{j + a} f[z] \in \IIm(\defopara^{k + j + 1})$, or $a \leq k$ and $\defopara^{j + a} f[z] = \defopara^{k + j - b} [w] \in \IIm(\defopara^{k + j - b}) \setminus \IIm(\defopara^{k + j - b + 1})$ witnesses some precrossing extension $d_{k - a - b}^{f, E_{r + 1 - a}, j + a}(z) = w$. In the latter case, the precrossing extension has to be inessential, so there exists another precrossing extension $d_{k - a_1 - b}^{f, E_{r + 1 - a_1}, j + a_1}(z_1) = w$ with $a < a_1 \leq k - b$ witnessed by $\defopara^{j + a_1} f[z_1] = \defopara^{k + j - b} [w]_1$. Replacing $[z]$ with $[z] - \defopara^{a_1 - a}[z_1]$, which is still a lift of $z$, it follows that $\defopara^{j + a} f[z]_{\std} \in \IIm(\defopara^{k + j - b + 1})$. As this strictly increases the $\defopara$-divisibility of the target, after sufficiently many iterations we obtain a lift $[z] \in \pi_{t - s, t + a}(X / \defopara^{r - a})$ with $\defopara^{j + a} f[z] \in \IIm(\defopara^{k + j + 1})$ as expected. 
\end{proof}

We can translate $j$-th layer extensions for different $j$ and/or stretch them across pages.

\begin{lemma} \label{reducing-refined-extensions-to-extensions}
    Suppose $X \in \Fil\Sp$, $s, t \in \Zb$, $c, k, j \in \Nb$ and $k + c \leq r \leq \infty$. The following are equivalent: 
    \begin{enumerate}
        \item Every $\alpha \in \IIm(\defopara^k) \subset \pi_{t - s, t}(X / \defopara^r)$ that is $\defopara^j$-torsion also lies in $\IIm(\defopara^{k + c})$. 
        \item In the standard SS of $X$, $B^{s + a, t + a}_{1 + a + j} / B^{s + a, t + a}_{1 + a} = 0$ for each $k \leq a < k + c$.
    \end{enumerate}
\end{lemma}

\begin{proof}
    Suppose there exists $y \in B^{s + a, t + a}_{1 + a + j} \setminus B^{s + a, t + a}_{1 + a}$ for some $k \leq a < k + c$, then by Theorem \ref{delta-as-total-diff} we can find $[x] \in \pi_{t - s + 1, t - a - j}(X / \defopara^{a + j})$ so that $\delta_{a + j}^1 [x] = y$. It then follows that $\defopara^a \delta_{a + j}^{r - a}[x] \in \pi_{t - s, t} (X / \defopara^{r})$ is $\defopara^j$-torsion, while it lies in $\IIm(\defopara^a) \setminus \IIm(\defopara^{a + 1})$ due to Theorem \ref{Bockstein-dictionary-finite} (if $r < \infty$) or Theorem \ref{Bockstein-dictionary-infinite} (if $r = \infty$), so it is not in $\IIm(\defopara^{k + c})$. Thus, item 1 implies item 2. Conversely, if there exists $\alpha \in \IIm(\defopara^a) \setminus \IIm(\defopara^{a + 1}) \subset \pi_{t - s, t}(X / \defopara^r)$ for some $k \leq a < k + c$, then  Theorem \ref{Bockstein-dictionary-finite} (if $r < \infty$) or Theorem \ref{Bockstein-dictionary-infinite} (if $r = \infty$) implies that $\alpha = \defopara^a [x]$ where $[x]$ is a lift of some nonzero $x \in Z_{r - a}^{s + a, t + a} / B_{1 + a}^{s + a, t + a}$, and if furthermore $\defopara^j \alpha = 0$ then $x \in B^{s + a, t + a}_{1 + a + j} / B^{s + a, t + a}_{1 + a}$. Therefore, item 2 also implies item 1.
\end{proof}

\begin{theorem} \label{stretching-extensions-across-pages}
    Take $f\colon X \to Y$ a map in $\Fil\Sp$. Fix $s, t \in \Zb$, $r, r' \in \Nb \cup\{\infty\}$, $j, j', k \in \Nb$ subject to the relations $k < r' \leq r$ and $j' \leq j$.   
    \begin{enumerate}
        \item If for $x \in Z_{r}^{s, t}(X)$, $y \in Z_{r - k}^{s + k, t + k}(Y)$ we have $d_k^{f, E_{r + 1}, j'}(x) = y$, then $d_k^{f, E_{{r'} + 1}, j}(x) = y$ as well. 
        \item If for $x \in Z_{r'}^{s, t}(X)$, $y \in Z_{r' - k}^{s + k, t + k}(Y)$ we have $d_k^{f, E_{r' + 1}, j}(x) = y$ and 
        \begin{itemize}
            \item The class $x$ actually lies in $Z_r^{s, t}(X)$. 
            \item This extension has no crossing $d_{k - a - b}^{f, E_{{r'} + 1 - a}, j + a}(x') = y'$ with $x' \in Z_{r' - a}^{s + a, t + a}(X) \setminus Z_{r - a}^{s + a, t + a}(X)$.
        \end{itemize}
        Then $y$ actually lies in $Z_{r - k}^{s + k, t + k}(Y)$, and $d_k^{f, E_{r + 1}, j}(x) = y$ as well. 
        \item If for $x \in Z_{r}^{s, t}(X)$, $y \in Z_{r - k}^{s + k, t + k}(Y)$ we have $d_k^{f, E_{r + 1}, j}(x) = y$ such that
        \begin{itemize}
            \item The extension has a witness $\defopara^j f[x] = \defopara^{k + j} [y]$ in which $f[x] \in \IIm(\defopara^c)$ for some $0 \leq c \leq k$. 
            \item In the standard SS of $Y$, $B^{s + a, t + a}_{1 + a + j}(Y) / B^{s + a, t + a}_{1 + a + j'}(Y) = 0$ for each $c\leq a \leq k$.
        \end{itemize}
        Then 
        $d_k^{f, E_{r + 1}, j'}(x) = y$ as well. 
    \end{enumerate}
\end{theorem}

\begin{proof}
    For part 1, suppose $\defopara^{j'} \! f[x] = \defopara^{k + j'} [y]$ is a witness of the extension $d_k^{f, E_{r + 1}, j'}(x) = y$. Then by Lemma \ref{defopara-rho-cartesian} we deduce $\defopara^{j} f \rho^r_{r'}[x] = \defopara^{j - j'} \rho^r_{r'}f[x] = \defopara^{k + j} \rho^{r - k}_{r' - k}[y]$, which witnesses $d_k^{f, E_{{r'} + 1}, j}(x) = y$. \parr

    For part 2, we assume $r' < \infty$ with no loss of generality. Take a pair of lifts $[x]' \in \pi_{t - s, t}(X / \defopara^{r'})$ and $[y]' \in \pi_{t - s, t + k}(Y / \defopara^{{r'} - k})$ so that $\defopara^{j} f[x]' = \defopara^{k + j}[y]'$. As $x \in Z_r^{s, t}(X)$, we can choose another lift $[x] \in \pi_{t - s, t}(X / \defopara^r)$. We then induct on the $\defopara$-divisibility of $\rho^r_{r'}[x] - [x]'$. Suppose $\rho^r_{r'}[x] - [x]' \in \IIm(\defopara^a) \setminus \IIm(\defopara^{a + 1})$, then it is of the form $\defopara^a [z]'$ for some nonzero $z \in Z_{r' - a}^{s + a, t + a}(X)$, cf. Theorem \ref{Bockstein-dictionary-finite} (for $r' < \infty$) or Theorem \ref{Bockstein-dictionary-infinite} (for $r' = \infty$). By Lemma \ref{inessential-crossing-extensions}, either $z \in Z_{r - a}^{s + a, t + a}(X)$, or there exists a lift $[z]_{\std} \in \pi_{t - s, t + a}(X / \defopara^{r' - a})$ so that $\defopara^{j + a} f[z]_{\std} \in \IIm(\defopara^{k + j + 1})$. In the former case we take a certain $[z]_1 \in \pi_{t - s, t + a}(X / \defopara^{r - a})$ lifting $z$ and replace $[x]$ by $[x] - \defopara^a[z]_1$, while in the latter case we replace $[x]'$ by $[x]' + \defopara^a[z]_{\std}$. After this modification, $\rho^r_{r'}[x] - [x]' \in \IIm(\defopara^{a + 1})$, and we still have $\defopara^j f[x]' = \defopara^{k + j}[y]''$ for $[y]'' = [y] + \defopara \varepsilon$ lifting $y$. Continuing this process until $a \geq k$, we obtain $[x] \in \pi_{t - s, t}(X / \defopara^r)$ and $[y]'' \in \pi_{t - s, t + k}(Y / \defopara^{r' - k})$ lifting $x$ and $y$, so that 
    \[\rho^{r + j}_{r' + j}(\defopara^{j} f[x]) = \defopara^{j} f\rho^r_{r'}[x] =  \defopara^{k + j} [y]''.\] 
    Feeding this into the pullback square in Lemma \ref{defopara-rho-cartesian} (via Lemma \ref{pi*-pullback-to-pullback-pi*}), we obtain a stem $[y] \in \pi_{t - s, t + k}(Y / \defopara^{r - k})$ so that $\defopara^{k + j} [y] = \defopara^j f[x]$ and $\rho^{r - k}_{r' - k}[y] = [y]''$. Thus, $[y]$ is also a lift of $y$, and we conclude that $d_k^{f, E_{r + 1}, j}(x) = y$ as expected. \parr
    
    For part 3, the difference $\defopara^{j'} \! f[x] - \defopara^{k + j'} [y]$ lies in $\IIm(\defopara^{j' + c})$ and it is $\defopara^{j - j'}$-torsion, so Lemma \ref{reducing-refined-extensions-to-extensions} suggests that  $\defopara^{j'} \! f[x] - \defopara^{k + j'} [y] = \defopara^{j' + k + 1} \varepsilon$ for some $\varepsilon \in \pi_{t - s, t + k + 1}(Y / \defopara^{r - k - 1})$. Hence, $[y]' = [y] + \defopara \varepsilon$ is another lift of $y$ so that $\defopara^{j'} \! f[x] = \defopara^{j' + k}[y]'$, which witnesses the extension we expect.
\end{proof}

One especially important family of arrows in $\Fil\Sp$ consists of the ``total differential'' operators
\(
\delta_r^{\infty}\colon X / \defopara^r \to \Sigma^{1,-r}X,
\)
whose mapping data determine the differentials in ${E_r^{s, t}(X)}$ by Theorem \ref{delta-as-total-diff}. These maps, however, do not directly fit into the framework of (hidden) extensions on specified pages, since that framework is tailored to the intuition of maps between spectral sequences rather than to the internal structure of a single spectral sequence. We will instead analyze their mapping data in a parallel, though not identical, manner via the following notion of crossings.

\begin{definition}[crossing for differentials] \label{crossing-for-diff}
    Suppose $X \in \Fil\Sp$. Fix $r, k \in \Nb$ with $1 \leq k \leq r$. Take $x \in Z_{r}^{s, t}$ and $y \in Z_{\infty}^{s + r + 1, t + r}$ so that $d_{r + 1}(x) = y$ in the standard SS. 
    \begin{itemize}
        \item A \textbf{crossing on the $E_{k + 1}$-page} for this $d_{r + 1}$ differential is a pair $x' \in Z_{r - a - b}^{s + a, t + a}$, $y' \neq 0 \in E_{r - a - b + 1}^{s + r - b + 1, t + r - b}$ with $0 < a < k$ and $0 \leq b \leq r - k$, so that $d_{r - a - b + 1}(x') = y'$. 
    \end{itemize}
\end{definition}

\begin{remark}
    In the case of Adams SS, this recovers \cite[Definition 4.14]{Lin-Wang-Xu-kervaire}. 
\end{remark}

\begin{example} \label{example-diff-crossing}

Take again the $\mathrm{H} \Fb_2$-Adams SS of $\Sb^0$, which is the standard SS of $X = \Sb^{0, 0}_{\mathrm{H} \Fb_2}$. Here the class $h_5 P e_0$ of bidegree $(t - s, s) = (56, 9)$ supports a differential $\color{blue}{d_5(h_5 P e_0) = d_0 \Delta h_0^2 e_0}$. Its potential crossings are depicted in Figure~\ref{figure-crossings-diff} below.

\newcommand{\crossingsdiffwindow}[1]{%
  % \draw (54.5,8.5) rectangle (56.5,14.5);
  \node at (55,8) {$55$};
  \node at (56,8) {$56$};
  \node at (54,9) {$9$};
  \node at (54,10) {$10$};
  \node at (54,11) {$11$};
  \node at (54,12) {$12$};
  \node at (54,13) {$13$};
  \node at (54,14) {$14$};

  \begin{scope}
    \clip (54.5,8.5) rectangle (56.5,14.5);

    % shaded cells, passed in as an argument
    #1

    \draw[black!22] (55.5,8.5) -- (55.5,14.5);
    \draw[black!22] (54.5,9.5) -- (56.5,9.5);
    \draw[black!22] (54.5,10.5) -- (56.5,10.5);
    \draw[black!22] (54.5,11.5) -- (56.5,11.5);
    \draw[black!22] (54.5,12.5) -- (56.5,12.5);
    \draw[black!22] (54.5,13.5) -- (56.5,13.5);

    % \draw (53,10) -- (56,11);
    % \draw (53,11) -- (56,12);
    \draw (56,9) -- (57,10);
    \draw (54.113,10.959) -- (55,12);
    % \draw (53,12) -- (56.113,12.959);
    \draw (55.887,10.041) -- (56,11);
    \draw (55.887,10.041) -- (57,11);
    % \draw (55.887,10.041) -- (59,11);
    \draw (56,11) -- (56,12);
    % \draw (56,11) -- (59,12);
    \draw (55,12) -- (56.113,12.959);
    \draw (56,12) -- (56.113,12.959);

    \draw[blue,diff]  (56,9) -- (55,14);
    \draw[ForestGreen,diff] (55,11) -- (54.113,14.959);
    \draw[ForestGreen,diff] (55.887,13.041) -- (55,17);
    \draw[red,diff]   (57,7) -- (56.113,9.959);
    \draw[red,diff]   (56.887,8.041) -- (56,11);
    \draw[red,diff]   (56.887,9.041) -- (56,12);
    \draw[cyan,diff]  (55.887,10.041) -- (55,12);
    \draw[cyan,diff]  (57,11) -- (56.113,12.959);

    \fill (56,9) circle (0.064) node[below] {\resizelabel{0.75}{$h_5Pe_0$}};
    \fill (56.113,9.959) circle (0.064) node[right] {\resizelabel{0.6}{$gt$}};
    \fill (55,11) circle (0.064) node[below] {\resizelabel{0.75}{$gm$}};
    \fill (56,11) circle (0.064);
    \fill (55,12) circle (0.064) node[above] {\resizelabel{0.7}{$h_1^2 M\!P \ \ $}};
    \fill (56,12) circle (0.064);
    \fill (56.113,12.959) circle (0.064);
    \fill (55.887,13.041) circle (0.064) node[above] {\resizelabel{0.85}{$d_0\Delta h_1e_0$}};
    \fill (55,14) circle (0.064) node[above] {\resizelabel{0.9}{$d_0\Delta h_0^2e_0$}};
    \fill (55.887,10.041) circle (0.064) node[left] {\resizelabel{0.8}{$\Delta^2 h_1 h_3$}};
  \end{scope}
}

\begin{figure}[htbp]
\centering
\scalebox{1}{%
\begin{tikzpicture}[line width=0.1pt]
    \tikzset{
        diff/.style={-{Stealth},line width=0.3pt,shorten <=3pt,shorten >=3pt},
        diff1/.style={-{Stealth},dashed,line width=0.6pt,shorten <=3pt,shorten >=3pt},
        extarr/.style={-{Stealth},line width=0.3pt,shorten <=3pt,shorten >=3pt},
        extarr1/.style={-{Stealth},dashed,line width=0.6pt,shorten <=3pt,shorten >=3pt},
        mapbelow/.style={-{Straight Barb[scale=0.8]},line width=0.7pt,shorten <=.7cm,shorten >=.7cm},
    }

  % left copy:
  % centers (55,14), (55,13), (55,12), (56,10)
  \begin{scope}
    \crossingsdiffwindow{
      \fill[black!12] (54.5,13.5) rectangle (55.5,14.5);
      \fill[black!12] (54.5,12.5) rectangle (55.5,13.5);
      \fill[black!12] (54.5,11.5) rectangle (55.5,12.5);
      \fill[black!12] (55.5, 9.5) rectangle (56.5,10.5);
    }
    \draw (54.5,8.5) rectangle (56.5,14.5);
  \end{scope}

  % middle copy:
  % centers (55,14), (55,13), (56,11), (56,10)
  \begin{scope}[xshift=4.2cm]
    \crossingsdiffwindow{
      \fill[black!12] (54.5,13.5) rectangle (55.5,14.5);
      \fill[black!12] (54.5,12.5) rectangle (55.5,13.5);
      \fill[black!12] (55.5,10.5) rectangle (56.5,11.5);
      \fill[black!12] (55.5, 9.5) rectangle (56.5,10.5);
    }
    \draw (54.5,8.5) rectangle (56.5,14.5);
  \end{scope}

  % right copy:
  % centers (55,14), (56,12), (56,11), (56,10)
  \begin{scope}[xshift=8.4cm]
    \crossingsdiffwindow{
      \fill[black!12] (54.5,13.5) rectangle (55.5,14.5);
      \fill[black!12] (55.5,11.5) rectangle (56.5,12.5);
      \fill[black!12] (55.5,10.5) rectangle (56.5,11.5);
      \fill[black!12] (55.5, 9.5) rectangle (56.5,10.5);
    }
    \draw (54.5,8.5) rectangle (56.5,14.5);
  \end{scope}

\end{tikzpicture}%
}
\caption{Potential crossings for $\color{blue}{d_5(h_5 P e_0) = d_0 \Delta h_0^2 e_0}$ on the $E_3, E_4$ and $E_5$-page.}
\label{figure-crossings-diff}
\end{figure}
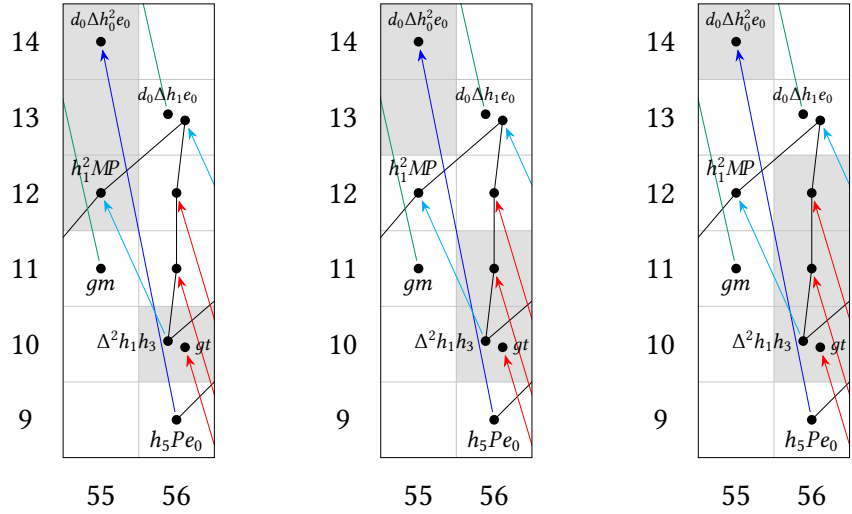
\begin{itemize}
    \item Every differential has no crossing on the $E_2$-page. 
    \item A crossing for this $d_5$ on the $E_3$-page is a pair $x' \in Z_{4 - a - b}^{9 + a, 56 + 9 + a}, y' \neq 0 \in E_{5 - a - b}^{14 - b, 55 + 14 - b}$ with \textbf{$1 \leq a \leq 1$, $0 \leq b \leq 2$} so that $d_{5 - a - b}(x') = y'$. In other words, this corresponds to an essential differential whose source and target lie in the shaded region for the figure \textbf{on the left}. As we can see, there is only one such crossing differential, namely $\color{cyan}{d_2(\Delta^2 h_1 h_3) = h_1^2 MP}$.
    \item A crossing for this $d_5$ on the $E_4$-page is a pair $x' \in Z_{4 - a - b}^{9 + a, 56 + 9 + a}, y' \neq 0 \in E_{5 - a - b}^{14 - b, 55 + 14 - b}$ with \textbf{$1 \leq a \leq 2$, $0 \leq b \leq 1$} so that $d_{5 - a - b}(x') = y'$. In other words, this corresponds to an essential differential whose source and target lie in the shaded region for the figure \textbf{in the middle}. From the figure we see this $d_5$ has no crossing on the $E_4$-page. 
    \item A crossing for this $d_5$ on the $E_5$-page is a pair $x' \in Z_{4 - a - b}^{9 + a, 56 + 9 + a}, y' \neq 0 \in E_{5 - a - b}^{14 - b, 55 + 14 - b}$ with \textbf{$1 \leq a \leq 3$, $0 \leq b \leq 0$} so that $d_{5 - a - b}(x') = y'$. In other words, this corresponds to an essential differential whose source and target lie in the shaded region for the figure \textbf{on the right}. From the figure we see this $d_5$ has no crossing on the $E_5$-page. 
\end{itemize}
\FloatBarrier
\end{example}

As in Lemma \ref{free-choice-lemma-for-refined-extensions}, we also interpret no crossing conditions for differentials as formulating the independence of bigraded stem lifts that witness a certain differential.

\begin{lemma} \label{free-choice-lemma-for-diff}
    Suppose for $r, k \in \Nb$ with $1 \leq k \leq r$, $X \in \Fil\Sp$ we have classes $x \in Z^{s, t}_r, y \in Z^{s + r + 1, t + r}_{\infty}$ in the standard SS of $X$ such that $d_{r + 1}(x) = y$, together with a pair of lifts  $[x]' \in \pi_{t - s, t}(X / \defopara^{k}), [y]' \in \pi_{t - s - 1, t + r}(X)$ so that $\delta_k^{\infty}[x]' = \defopara^{r - k}[y]'$ (which always exists due to Theorem \ref{delta-as-total-diff}). Then for each $[x] \in \pi_{t - s, t}(X / \defopara^{k})$ lifting $x$, either we have $\delta_k^{\infty}[x] = \defopara^{r - k} [y]$ in which $[y] = [y]' + \defopara \varepsilon$ for some $\varepsilon \in \pi_{t - s - 1, t + r + 1}(X)$, or we have $\delta_k^{\infty}([x] - [x]') = \defopara^{r - k - b} [w]$ in which $[w] \in \pi_{t - s - 1, t + r - b}(X)$ lifts the target of some $E_{k + 1}$ page crossing differential $d_{r - a - b + 1}(z) = w$. 
\end{lemma}

\begin{proof}
    The difference $[x] - [x]'$ lies in $\IIm(\defopara^a) \setminus \IIm(\defopara^{a + 1})$ for some $0 < a < k$, so Theorem \ref{Bockstein-dictionary-finite} implies that $[x] - [x]' = \defopara^a [z]$ for $[z] \in \pi_{t - s, t + a}(X / \defopara^{k - a})$ lifting some nonzero $z \in Z^{s + a, t + a}_{k - a} / B^{s + a, t + a}_{1 + a}$. If the image $\delta_k^{\infty}(\defopara^a [z]) = \delta_{k - a}^{\infty} [z]$ lies in $\IIm(\defopara^{r - k + 1})$ then the former case in the statement holds true. Otherwise, $\delta_{k - a}^{\infty} [z] \in \IIm(\defopara^{r - k - b}) \setminus \IIm(\defopara^{r - k - b + 1})$ for some $0 \leq b \leq r - k$, so Theorem \ref{Bockstein-dictionary-infinite} implies that $\delta_{k - a}^{\infty} [z] = \defopara^{r - k - b}[w]$ for $[w]$ lifting some nonzero $w \in Z^{s + r - b + 1, t + r - b}_{\infty} / B^{s + r - b + 1, t + r - b}_{r - k - b + 1}$. If the corresponding differential $d_{r - a - b + 1}(z) = w$ (under Theorem \ref{delta-as-total-diff}) is essential, then it is an $E_{k + 1}$ page crossing differential, so the latter case in the statement holds true. Otherwise, there is another lift $[z]_{\std}$ of $z$ so that $\delta_{k - a}^{\infty}[z]_{\std} \in \IIm(\defopara^{r - k - b + 1})$. Thus, the difference $\defopara^a([z] - [z]_{\std}) = \defopara^{a_1}[z_1] \in \IIm(\defopara^{a_1}) \setminus \IIm(\defopara^{a_1 + 1})$ for certain $a < a_1 < k$, and $\delta_{k}^{\infty}(\defopara^{a_1}[z_1]) = \defopara^{r - k - b}[w]_1$ where $[w]_1 = [w] + \defopara \varepsilon$ is another lift of $w$. If the corresponding differential $d_{r - a_1 - b + 1}(z_1) = w$ is essential then it is an $E_{k + 1}$ page crossing differential. Otherwise, we reiterate the process, and in the end we would obtain some $[z_n] \in \pi_{t - s, t + a_n}(X / \defopara^{k - a_n})$ with $0 < a_1 < \cdots < a_n < k$ whose corresponding differential $d_{r - a_n - b + 1}(z_n) = w$ is essential, thus an $E_{k + 1}$ page crossing differential as expected. 
\end{proof}

\begin{corollary} \label{free-choice-lemma-for-diff-with-full-no-crossing}
    The differential $d_{r + 1}(x) = y$ has no crossing on the $E_{k + 1}$ page if and only if for each $[x] \in \pi_{t - s, t}(X / \defopara^k)$ lifting $x$, there exists some $[y] \in \pi_{t - s - 1, t + r}(X)$ lifting $y$ so that $\delta_k^{\infty}[x] = \defopara^{r - k}[y]$. 
\end{corollary}

\begin{remark} \label{free-choice-lemma-for-diff-refined}
    Under the setup of Lemma \ref{free-choice-lemma-for-diff}, take $i \in \Nb$ with $1 \leq i \leq k$. The same argument there also shows that for each $[x] \in \pi_{t - s, t}(X / \defopara^{k})$ lifting $x$, either we have $\defopara^{k - i}\delta_{k}^{\infty}[x] = \delta_{i}^{\infty} \rho^{k}_{i}[x] = \defopara^{r - i} [y]$ in which $[y] = [y]' + \defopara \varepsilon$ for some $\varepsilon \in \pi_{t - s - 1, t + r + 1}(X)$, or we have $\defopara^{k - i}\delta_{k}^{\infty}([x] - [x]') = \defopara^{r - i - b} [w]$ in which $[w] \in \pi_{t - s - 1, t + r - b}(X)$ lifts the target of an essential differential $d_{r - a - b + 1}(z) = w$ with $0 < a < i, 0 \leq b \leq r - k$ (namely $d_{r - a - b + 1}(z) = w$ is an $E_{c + 1}$-page crossing of the original $d_{r + 1}$ for each $i \leq c \leq k$). 
\end{remark}

\begin{lemma} \label{inessential-crossing-differentials}
    Take $X \in \Fil\Sp$. Suppose for $s, t \in \Zb$, $r, k \in \Nb$ with $1 \leq k \leq r$ we have classes $x \in Z^{s, t}_r, y \in Z_{\infty}^{s + r + 1, t + r}$ in the standard SS of $X$ so that $d_{r + 1}(x) = y$. Then for each $0 < a < r$ and each $z \in Z_{k - a}^{s + a, t + a}$, if $z$ does not support an $E_{k + 1}$-page crossing differential for $d_{r + 1}(x) = y$, then there is a lift $[z]_{\std} \in \pi_{t - s, t + a}(X / \defopara^{k - a})$ so that $\delta_{k - a}^{\infty}[z]_{\std} \in \IIm(\defopara^{r - k + 1})$. 
\end{lemma}

\begin{proof}
    According to Lemma \ref{free-choice-lemma-for-diff}, for any lift $[z] \in \pi_{t - s, t + a}(X / \defopara^{k - a})$, either $\delta_{k - a}^{\infty}[z] \in \IIm(\defopara^{r - k + 1})$, or $\delta_{k - a}^{\infty}[z] = \defopara^{r - k - b} [w] \in \IIm(\defopara^{r - k - b}) \setminus \IIm(\defopara^{r - k - b + 1})$ so that $[w] \in \pi_{t - s - 1, t + r - b}(X)$ lifts the target of some $E_{k + 1}$-page crossing differential $d_{r - a_1 - b + 1}(z_1) = w$. 
    In the latter case, we have $a < a_1 < k$ as $z$ does not support any such crossing differential, so there exists $[z_1] \in \pi_{t - s, t + a_1}(X / \defopara^{k - a_1})$ for which $\delta_{k - a_1}^{\infty}[z_1] = [w]_1$ is another lift of $w$. Replacing $[z]$ with $[z] - \defopara^{a_1 - a}[z_1]$, which is still a lift of $z$, it follows that $\delta_{k - a}^{\infty}[z] \in \IIm(\defopara^{r - k - b + 1})$. 
    As this strictly increases the $\defopara$-divisibility of the target, after sufficiently many iterations we obtain a lift $[z] \in \pi_{t - s, t + a}(X / \defopara^{k - a})$ with $\delta_{k - a}^{\infty}[z] \in \IIm(\defopara^{r - k + 1})$ as expected.  
\end{proof}

\subsection{Generalized Leibniz rule for SS}
\label{subsec:3.2}

The ``prototype'' of the GLR is a simple consequence of the naturality for the total differentials. 
    
\begin{theorem}[Blueprint for GLR] \label{blueprint-GLR}
    Let $f\colon X \to Y$ be a map in $\Fil\Sp$. Fix $s, t \in \Zb$, $r, m, l, i, j, n \in \Nb$, $U \in \Nb \cup \{\infty\}$ so that $m < n \leq r$, $i \geq r - n$ and $U_1 = U - r + n > l$. Suppose there are bigraded stems
    \[[x] \in \pi_{t - s, t}(X / \defopara^n), [x_{\infty}] \in \pi_{t - s - 1, t + r}(X / \defopara^{U_1}), [y] \in \pi_{t - s, t + m}(Y / \defopara^{n - m}), [y_{\infty}] \in \pi_{t - s - 1, t + r + l}(Y / \defopara^{U_1 - l})\] 
    subject to the relations $\defopara^j f[x] = \defopara^{m + j}[y]$ and $\delta_n^{U}[x] = \defopara^{r - n}[x_{\infty}]$. Then $\defopara^{i} f[x_{\infty}] = \defopara^{i + l} [y_{\infty}]$ if and only if $\defopara^{i - r + n}\delta_{n - m}^{U}[y] = \defopara^{i + l} [y_{\infty}]$.
\end{theorem}

\begin{proof}
    Naturality of $\delta_{n}^{U}$ provides a commutative square 
    \[\begin{tikzcd}
        X / \defopara^n \ar[r, "{f}"] \ar[d, "{\delta_{n}^{U}}"] & Y / \defopara^n \ar[d, "{\delta_{n}^{U}}"] \\
        \Sigma^{1, -n} X / \defopara^{U} \ar[r, "{f}"] & \Sigma^{1, -n} Y / \defopara^{U}
    \end{tikzcd}\]
    thus $\delta_{n - m}^{U} [y] = \delta_{n + j}^{U} (\defopara^{m + j} [y]) = \delta_{n + j}^{U}(\defopara^j f[x])  = \delta_n^{U}f[x] =  f(\delta_n^{U}[x]) = \defopara^{r - n}  f[x_{\infty}]$ due to Corollary \ref{defopara-totaldiff}. The result then follows by multiplying $\defopara^{i - r + n}$ on both sides and comparing with $\defopara^{i + l}[y_{\infty}]$. 
\end{proof}

We then state the necessary no crossing conditions to translate Theorem \ref{blueprint-GLR} into a statement of standard spectral sequences. The following Definition \ref{complementary-crossings} aims for maximal generality, while we list some sufficient conditions easy to verify in practice in Example \ref{complementary-crossing-examples}.

\begin{definition}[Complementary crossings] \label{complementary-crossings}
    Let $f\colon X \to Y$ be a map in $\Fil\Sp$.
    \begin{itemize}
        \item Suppose for $s, t \in \Zb$, $r, n, m, j \in \Nb$ with $m < n \leq r $, we have classes $x \in Z_r^{s, t}(X)$, $x_{\infty} \in Z_{\infty}^{s + r + 1, t + r}(X)$ and $y \in Z_{n - m}^{s + m, t + m}(Y)$ so that $d_{r + 1}(x) = x_{\infty}$ and $d_m^{f, E_{n + 1}, j}(x) = y$. Then the differential and the $j$-th layer extension are \textbf{complementary on the source} % on the $E_{n + 1}$ page} 
        if for each $z \in Z^{s + a, t + a}_{n - a}(X)$ with $0 < a \leq m$, either $z$ does not support an $E_{n + 1}$-page crossing differential for $d_{r + 1}(x) = x_{\infty}$, or $z$ is not the source of a crossing extension for $d_m^{f, E_{n + 1}, j}(x) = y$.
        \item Suppose for $s, t \in \Zb$, $r, l, k, i \in \Nb$, $U \in \Nb \cup \{\infty\}$ with %$j + l < r$, $r - j - l\leq k$ 
        $k > 0$, $l < r$ and $l < U$, we have classes $x_{\infty} \in Z_{U}^{s - l, t - l}(X)$, $y \in Z_k^{s - r - 1, t - r}(Y)$ and $y_{\infty} \in Z_{U - l}^{s, t}(Y)$ so that $d_{l}^{f, E_{U + 1}, j}(x_{\infty}) = y_{\infty}$. Then the $i$-th layer extension and the class $y$ are \textbf{complementary on the target} if each $w \in Z_{U - l}^{s, t}(Y)$ is not the target of a crossing extension for $d_{l}^{f, E_{U + 1}, i}(x_{\infty}) = y_{\infty}$, and for each $w \in Z_{U - l + b}^{s - b, t - b}(Y)$ with $0 < b < l$, either $w$ is not the target of a crossing extension for $d_{l}^{f, E_{U + 1}, i}(x_{\infty}) = y_{\infty}$, or in the standard SS of $Y$ there is no essential differential of the form $d_{r + 1 - b}(y) = w$. 
        \item Suppose for $s, t \in \Zb$, $r, l, i \in \Nb$ with 
        $r > l + i$, we have $x_{\infty} \in Z_{\infty}^{s - l, t - l}(X)$, $y \in Z_r^{s - r - 1, t - r}(Y)$ and $y_{\infty} \in Z_{\infty}^{s, t}(Y)$ so that $d_{r + 1}(y) = y_{\infty}$. Then the differential and the permanent cycle $x_{\infty}$ are \textbf{$i$-complementary on the target} if each class $w \in Z_{\infty}^{s, t}(Y)$ is not the target of an essential differential $d_{r - a + 1}(z) = w$ with $0 < a < r - l - i$, and for each $w \in Z_{\infty}^{s - b, t - b}(Y)$ with $0 < b \leq l$, either $w$ is not the target of an essential differential $d_{r - a - b + 1}(z) = w$ with $0 < a < r - l - i$, or there is no nontrivial extension of the form $d_{l - b}^{f, E_{\infty}, i}(x_{\infty}) = w$. 
    \end{itemize}
\end{definition}

\begin{example} \label{complementary-crossing-examples}
    \begin{itemize}
        \item The differential $d_{r + 1}(x) = x_{\infty}$ and the extension $d_m^{f, E_{n + 1}, j}(x) = y$ are complementary on the source if there exist subsets $S_1,S_2 \subset \{1,\ldots,m\}$ such that the differential has no crossing on the $E_{n + 1}$-page with source range $S_1$, the extension has no crossing with source range $S_2$, and $\{1, \ldots, m\} \subset S_1 \cup S_2$. In particular, this holds if either the differential $d_{r + 1}(x) = x_{\infty}$ has no crossing on the $E_{n + 1}$-page, or the extension $d_m^{f, E_{n + 1}, j}(x) = y$ has no crossing.

        \item The extension $d_{l}^{f, E_{U + 1}, j}(x_{\infty}) = y_{\infty}$ and the class $y$ are complementary on the target if there exist subsets $T_1,T_2 \subset \{0,1,\ldots,l-1\}$ such that the extension has no crossing with target range $T_1$, the class $y$ does not support an essential $d_{r - m + l + 1 - b}$-differential for each $b \in T_2$, and $0 \in T_1$, $\{1, \ldots, l - 1\} \subset T_1 \cup T_2$. In particular, this holds if the extension $d_{l}^{f, E_{U + 1}, j}(x_{\infty}) = y_{\infty}$ has no crossing.

        \item The differential $d_{r + 1}(y) = y_{\infty}$ and the class $x_{\infty}$ are $i$-complementary on the target if there exist subsets $T_1,T_2 \subset \{0,1,\ldots,l\}$ such that the differential has no crossing of the form $d_{r - a - b + 1}(z) = w$ with $0 < a < r - l - i$ and $b \in T_1$, the class $x_{\infty}$ does not support a nontrivial $i$-th layer extension on the $E_{\infty}$-page with filtration jump $l - b$ for each $b \in T_2$, and $0 \in T_1$, $\{1, \ldots, l\} \subset T_1 \cup T_2$. In particular, this holds if the differential $d_{r + 1}(y) = y_{\infty}$ has no crossing on the $E_{c + 1}$-page for some $c \in \{r - l - i, \ldots, r - l\}$.
    \end{itemize}
\end{example}

The intuition behind the complementary crossing condition is as follows. We use no-crossing hypotheses to choose bigraded stems that simultaneously realize two tasks, such as witnessing an extension and a differential, or two extensions. Corollary \ref{free-choice-lemma-for-refined-extensions-with-full-no-crossings} and Corollary \ref{free-choice-lemma-for-diff-with-full-no-crossing} guarantee this under full no crossing from either side, but these assumptions are stronger than necessary: complementary crossing conditions already suffice to produce the desired common lift.

\begin{lemma} \label{combining-partial-nocrossings-in-source-refined}
    Suppose $f\colon X \to Y$ is a map in $\Fil\Sp$. Fix $s, t \in \Zb$, $r, n, m, j \in \Nb$ with $m < n \leq r $. Suppose we have classes $x \in Z_r^{s, t}(X)$, $x_{\infty} \in Z_{\infty}^{s + r + 1, t + r}(X)$ and $y \in Z_{n - m}^{s + m, t + m}(Y)$ which support a differential $d_{r + 1}(x) = x_{\infty}$ and an extension $d_m^{f, E_{n + 1}, j}(x) = y$ complementary on the source.
    Then there exist lifts $[x] \in \pi_{t - s, t}(X / \defopara^{n}), [x_{\infty}] \in \pi_{t - s - 1, t + r}(X), [y] \in \pi_{t - s, t + m}(Y / \defopara^{n - m})$ of $x, x_{\infty}$ and $y$, such that $\delta_n^{\infty}[x] = \defopara^{r - n}[x_{\infty}]$ and $\defopara^j f[x] = \defopara^{m + j} [y]$. 
\end{lemma}
\begin{proof}
    According to Lemma \ref{inessential-crossing-extensions}, Lemma \ref{inessential-crossing-differentials} and the ``complementary on the source'' assumption, any $z \in Z_{n - a}^{s + a, t + a}(X) / B_{a + 1}^{s + a, t + a}(X)$ with $0 < a < n$ admits a lift $[z]_{\std} \in \pi_{t - s, t + a}(X / \defopara^{n - a})$ for which either $\defopara^{j + a} f[z]_{\std} \in \IIm(\defopara^{m + j + 1})$ or $\delta_n^{\infty} \defopara^a[z]_{\std} = \delta_{n - a}^{\infty}[z]_{\std} \in \IIm(\defopara^{r - n + 1})$. The differential and the extension in the assumption imply the existence of $[x], [x]' \in \pi_{t - s, t}(X / \defopara^n), [x_{\infty}]' \in \pi_{t - s - 1, t + r}(X), [y]' \in \pi_{t - s, t + m}(Y / \defopara^{n - m})$ lifting $x, x_{\infty}, y$, so that $\delta_n^{\infty}[x] = \defopara^{r - n}[x_{\infty}]'$ and $\defopara^j f[x]' = \defopara^{m + j} [y]'$. We induct on the $\defopara$-divisibility of $[x]' - [x]$. Suppose $[x]' - [x] \in \IIm(\defopara^a) \setminus \IIm(\defopara^{a + 1})$, then it is detected by certain nonzero $z \in Z_{n - a}^{s + a, t + a} / B_{a + 1}^{s + a, t + a}$, cf. Theorem \ref{Bockstein-dictionary-finite}. If $\defopara^{j + a} f[z]_{\std} \in \IIm(\defopara^{m + j + 1})$ we replace $[x]$ by $[x] + \defopara^a [z]_{\std}$, otherwise we replace $[x]'$ by $[x]' - \defopara^a [z]_{\std}$. It follows that $\delta_n^{\infty}[x] = \defopara^{r - n}[x_{\infty}]''$ and $\defopara^j f[x]' = \defopara^{m + j} [y]''$ where $[x_{\infty}]'', [y]''$ are also lifts of $x_{\infty}, y$, while $[x]' - [x] \in \IIm(\defopara^{a + 1})$. Continuing this process until $a = n$, we obtain $[x] = [x]'$, for which we have both $\delta_n^{\infty}[x] = \defopara^{r - n}[x_{\infty}]$ and $\defopara^j f[x] = \defopara^{m + j} [y]$ with $[x_{\infty}], [y]$ lifting $x_{\infty}, y$.
\end{proof}

The Blueprint Theorem \ref{blueprint-GLR} organizes the information in two directions. On the one hand, if there is a differential $d(x) = x_{\infty}$ in the domain of $f\colon X \to Y$ in $\Fil\Sp$, together with two $f$-extensions $x \rightsquigarrow y$ and $x_{\infty} \rightsquigarrow y_{\infty}$, then it suggests the corresponding differential $d(y) = y_{\infty}$ in the codomain.

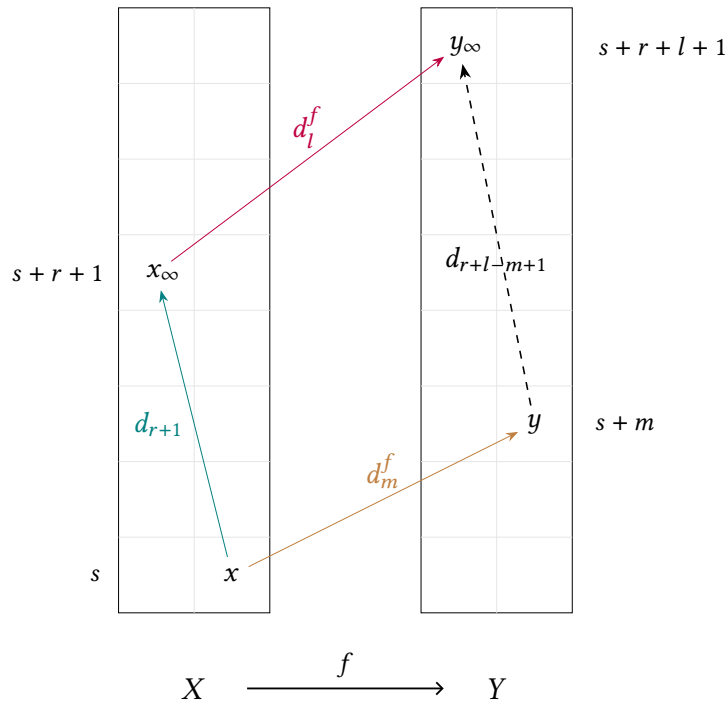
\begin{figure}[htbp]
    \centering
    \scalebox{1}{%
        \begin{tikzpicture}[line width=0.1pt]
            \tikzset{
                diff/.style={-{Stealth},line width=0.3pt,shorten <=7pt,shorten >=7pt},
                diff1/.style={-{Stealth},dashed,line width=0.6pt,shorten <=7pt,shorten >=7pt},
                extarr/.style={-{Stealth},line width=0.3pt,shorten <=7pt,shorten >=7pt},
                extarr1/.style={-{Stealth},dashed,line width=0.6pt,shorten <=7pt,shorten >=7pt},
                mapbelow/.style={-{Straight Barb[scale=0.8]},line width=0.7pt,shorten <=.7cm,shorten >=.7cm},
            }
            
            % left copy: 
            
            \begin{scope}
                \draw (49.5,11.5) rectangle (51.5,19.5);
                \node at (49.2,12) {\footnotesize $s$};
                \node at (48.7,16) {\footnotesize $s + r + 1$};
                \begin{scope}
                    % \clip (49.5,11.5) rectangle (51.5,19.5);
                    \draw[black!10] (50.5,11.5) -- (50.5,19.5);
                    \draw[black!10] (49.5,12.5) -- (51.5,12.5);
                    \draw[black!10] (49.5,13.5) -- (51.5,13.5);
                    \draw[black!10] (49.5,14.5) -- (51.5,14.5);
                    \draw[black!10] (49.5,15.5) -- (51.5,15.5);
                    \draw[black!10] (49.5,16.5) -- (51.5,16.5);
                    \draw[black!10] (49.5,17.5) -- (51.5,17.5);
                    \draw[black!10] (49.5,18.5) -- (51.5,18.5);

                    \coordinate (X-source) at (51,12);
                    \coordinate (X-target) at (50,16);
                    
                    \node at (X-source) {\small $x$};
                    \node at (X-target) {\small $\ \ x_{\infty}$};
                    \color{teal} \draw[diff] (X-source) -- (X-target) node[midway, left] {\small $d_{r + 1}$};
                \end{scope}
            \end{scope}
            
            % right copy: 
            
            \begin{scope}[xshift=4cm,yshift=0cm]
                \draw (49.5,11.5) rectangle (51.5,19.5);
                \node at (52.2,14) {\footnotesize $s + m$};
                \node at (52.7,19) {\footnotesize $s + r + l + 1$};
                \begin{scope}
                    % \clip (49.5,11.5) rectangle (51.5,19.5);
                    \draw[black!10] (50.5,11.5) -- (50.5,19.5);
                    \draw[black!10] (49.5,12.5) -- (51.5,12.5);
                    \draw[black!10] (49.5,13.5) -- (51.5,13.5);
                    \draw[black!10] (49.5,14.5) -- (51.5,14.5);
                    \draw[black!10] (49.5,15.5) -- (51.5,15.5);
                    \draw[black!10] (49.5,16.5) -- (51.5,16.5);
                    \draw[black!10] (49.5,17.5) -- (51.5,17.5);
                    \draw[black!10] (49.5,18.5) -- (51.5,18.5);

                    \coordinate (Y-source) at (51,14);
                    \coordinate (Y-target) at (50,19);
                    
                    \node at (Y-source) {\small $y$};
                    \node at (Y-target) {\small $\ \ y_{\infty}$};
                    \draw[diff1] (Y-source) -- (Y-target) node[midway, below] {\small $d_{r + l - m + 1}$};
                \end{scope}
            \end{scope}

            % extension arrows 

            \color{brown} \draw[extarr] (X-source) -- (Y-source) node[midway, above] {\small $d_{m}^f$};
            \color{purple} \draw[extarr] (X-target) -- (Y-target) node[midway, above] {\small $d_{l}^f$};
            
            % bottom map X -> Y labeled by f
            
            \coordinate (label-X) at (50.5, 10.5);
            \coordinate (label-Y) at (54.5, 10.5);

            \color{black}
            \node at (label-X) {$X$};
            \node at (label-Y) {$Y$};
            \draw[mapbelow] (label-X) -- (label-Y) node[midway,above] {\footnotesize $f$};
        \end{tikzpicture}
    }
    \caption{Generalized Leibniz rule.}
    \label{figure-GLR}
\end{figure}

\begin{theorem}[Generalized Leibniz rule, part 1] \label{GLR-for-SS}
     Let $f\colon X \to Y$ be a map in $\Fil\Sp$. Fix $s, t \in \Zb$, $r, m, l, j, n \in \Nb$, $U \in \Nb \cup \{\infty\}$ so that $m < n \leq r$ and $U_1 = U - r + n > l$. Suppose there are classes $x \in Z_r^{s, t}(X), x_{\infty} \in Z_{\infty}^{s + r + 1, t + r}(X), y \in Z_{n - m}^{s + m, t + m}(Y), y_{\infty} \in Z_{U_1 - l}^{s + r + l + 1, t + r + l}(Y)$ subject to the assumptions:
    \raggedcolumns
    \begin{multicols}{2}
    \begin{enumerate}
        \item ${\color{teal} d_{r + 1}(x) = x_{\infty}}$. 
        \item ${\color{brown} d^{f, E_{n + 1}, j}_{m}(x) = y}$. 
        \item ${\color{purple} d^{f, E_{U_1 + 1}, i}_{l}(x_{\infty}) = y_{\infty}}$, here $i = r - m - 1$.
        \item The {\color{teal} differential} in item 1 and the {\color{brown} extension} in item 2 are complementary on the source.
        \item The {\color{purple}  extension} in item 3 and the class $y$ are complementary on the target.
    \end{enumerate}
    \end{multicols}    
    \noindent Then $y \in Z_{r - m + l}^{s + m, t + m}(Y)$, and $d_{r - m + l + 1}(y) = y_{\infty}$ in the standard SS of $Y$ (see Figure \ref{figure-GLR}).
\end{theorem}

\begin{proof}
    Item 1,2 and 4 suggest the existence of $[x] \in \pi_{t - s, t}(X / \defopara^n), [x_{\infty}] \in \pi_{t - s - 1, t + r}(X / \defopara^{U_1})$ and $[y] \in \pi_{t - s, t + m}(Y / \defopara^{n - m})$ lifting $x, x_{\infty}$ and $y$, such that 
    \[\defopara^j f[x] = \defopara^{m + j} [y] \quad \text{ and } \quad \delta_n^{U}[x] = \defopara^{r - n}[x_{\infty}]\]
    due to Lemma \ref{combining-partial-nocrossings-in-source-refined}. Also, item 3 implies we can find $[x_{\infty}]' \in \pi_{t - s - 1, t + r}(X / \defopara^{U_1})$ and $[y_{\infty}]' \in \pi_{t - s - 1, t + r + l}(Y / \defopara^{U_1 - l})$ lifting $x_{\infty}$ and $y_{\infty}$, so that $\defopara^{i}f[x_{\infty}]' = \defopara^{i + l}  [y_{\infty}]'$. Thus by Lemma \ref{free-choice-lemma-for-refined-extensions}, either $\defopara^{i}f[x_{\infty}] = \defopara^{i + l} [y_{\infty}]$ for some $[y_{\infty}] = [y_{\infty}]' + \defopara \varepsilon$, 
    or $[x_{\infty}] = [x_{\infty}]' + \defopara^a [z]$ in which $[z]$ witnesses some precrossing extension $d_{l - a - b}^{f, E_{U_1 - a + 1}, i + a}(z) = w$. In the former case, the Blueprint Theorem \ref{blueprint-GLR} implies
    \[\delta_1^{U + n - m - 1} y = \delta_{1}^{U + n - m - 1} \rho_1^{n - m} [y] = \defopara^{n - m - 1} \delta_{n - m}^{U}[y] = \defopara^{i - r + n} \delta_{n - m}^{U}[y] = \defopara^{i + l}  [y_{\infty}].\] 
    Thus, $y \in Z_{r - m + l}^{s + m, t + m}(Y)$ and $d_{r - m + l + 1}(y) = y_{\infty}$ in the standard SS of $Y$ due to Theorem \ref{delta-as-total-diff}. In the latter case, the class $w$ is also the target of some crossing extension due to Fact \ref{crossing-facts}, so item 5 implies that $0 < b \leq l - a$ and there is no essential differential of the form $d_{r - m + l - b + 1}(y) = w$. On the other hand, according to the Blueprint Theorem \ref{blueprint-GLR},
    \[\delta_1^{U + n - m - 1} y = \defopara^{i - r + n} \delta_{n - m}^{U}[y] = \defopara^{i} f[x_{\infty}] = \defopara^{i} f(\defopara^a [z] + [x_{\infty}]') = \defopara^{i + l - b}([w] + \defopara^{b}[y_{\infty}]')\]
    which leads to a shorter differential $d_{r - m + l - b + 1}(y) = w$. As $w$ is the target of a crossing extension, we have $w \not\in B^{s + r + l - b + 1, t + r + l - b}_{r - m + l - b}(Y)$, therefore $d_{r - m + l - b + 1}(y) = w$ is essential. This contradicts our previous conclusion, so the latter case can never happen.
\end{proof}

\begin{remark} \label{explanation-of-the-parameters-GLR}
    We briefly discuss how to choose parameters when we apply Theorem \ref{GLR-for-SS}.
    \begin{itemize}
        \item The parameters $r, m, l$ are fixed by the desideratum. 
        \item The parameter $n$ is fine-tuned towards two goals: the existence of the extension from $x$ to $y$, and the ``complementary on the source'' condition. The smallest possible choice of $n$ yields $n + 1 = m + 2$, the earliest page on which extensions of filtration jump $m$ make sense. The largest possible choice is $n = r$, which is also the default choice in practice.
        \item The parameter $j$ can be as large as possible to guarantee the extension from $x$ to $y$ is unobstructed by early boundaries (i.e. $\defopara$-power torsion lifts). In practice, often $j = 0$ suffices.
        \item The parameter $U$ is fine-tuned towards two goals: the existence of the extension from $x_{\infty}$ to $y_{\infty}$, and the ``complementary on the target'' condition. The smallest possible choice of $U$ yields $U - r + n + 1 = l + 2$, the earliest page on which extensions of filtration jump $l$ make sense. The largest possible choice is $U = \infty$, which is also the default choice in practice.
    \end{itemize}
\end{remark}

\begin{remark} \label{reduction-to-GLR-in-LWX}
    The complementary crossing conditions can be simplified via Example \ref{complementary-crossing-examples}. 
    \begin{itemize}
        \item Assumption item 4 in Theorem \ref{GLR-for-SS} holds true if the differential in item 1 has no crossing on the $E_{n + 1}$ page or the extension in item 2 has no crossing. 
        \item Assumption item 5 holds true if the extension in item 4 has no crossing. 
    \end{itemize}
    Under the translation in Remark \ref{synthetic-spectra-as-special-case}, Remark \ref{comparison-with-extension-SS} and Remark \ref{comparison-with-extension-SS-contd}, we thus recover \cite[Theorem 6.1]{Lin-Wang-Xu-kervaire} for Adams SS. Note that our formulation of complementary crossing conditions in item 4 and item 5 match the expected generalization of GLR in \cite[Remark 6.5]{Lin-Wang-Xu-kervaire}.
\end{remark}

The Blueprint Theorem \ref{blueprint-GLR} can also be read from a complementary perspective: if $d(x) = x_{\infty}$ and $d(y) = y_{\infty}$ are two differentials lying respectively in the domain and the codomain of a map $f\colon X \to Y$ in $\Fil\Sp$, together with an $f$-extension $x \rightsquigarrow y$ between the sources, then it suggests the corresponding $f$-extension $x_{\infty} \rightsquigarrow y_{\infty}$ between the targets.

\begin{theorem}[Generalized Leibniz rule, part 2] \label{GLR-for-SS-2}
     Let $f\colon X \to Y$ be a map in $\Fil\Sp$. Fix $s, t \in \Zb$, $r, m, l, i, j, n \in \Nb$ so that $m < r - i \leq n \leq r$. Suppose there are classes $x \in Z_r^{s, t}(X), x_{\infty} \in Z_{\infty}^{s + r + 1, t + r}(X), y \in Z_{r - m + l}^{s + m, t + m}(Y), y_{\infty} \in Z_{\infty}^{s + r + l + 1, t + r + l}(Y)$ subject to the following assumptions:
    \raggedcolumns
    \begin{multicols}{2}
    \begin{enumerate}
        \item ${\color{black} d_{r + 1}(x) = x_{\infty}}$. 
        \item ${\color{black} d^{f, E_{n + 1}, j}_{m}(x) = y}$. 
        \item $d_{r - m + l + 1}(y) = y_{\infty}$.
        \item The {\color{black}  differential} in item 1 and the {\color{black} extension} in item 2 are complementary on the source.
        \item The differential in item 3 and the class $x_{\infty}$ are $i$-complementary on the target.
    \end{enumerate}
    \end{multicols}    
    \noindent Then there is an $f$-extension ${\color{black} d^{f, E_{\infty}, i}_{l}(x_{\infty}) = y_{\infty}}$.
\end{theorem}

\begin{proof}
    Item 1,2 and 4 suggest the existence of $[x] \in \pi_{t - s, t}(X / \defopara^n), [x_{\infty}] \in \pi_{t - s - 1, t + r}(X)$ and $[y] \in \pi_{t - s, t + m}(Y / \defopara^{n - m})$ lifting $x, x_{\infty}$ and $y$, such that 
    \[\defopara^j f[x] = \defopara^{m + j} [y] \quad \text{ and } \quad \delta_n^{\infty}[x] = \defopara^{r - n}[x_{\infty}]\]
    due to Lemma \ref{combining-partial-nocrossings-in-source-refined}. Also, item 3 implies via  Theorem \ref{delta-as-total-diff} that we can find $[y]' \in \pi_{t - s, t + m}(Y / \defopara^{n - m})$ and $[y_{\infty}]' \in \pi_{t - s - 1, t + r + l}(Y)$ lifting $y$ and $y_{\infty}$, so that $\delta_{n - m}^{\infty} [y]' = \defopara^{r - n + l} [y_{\infty}]'$. Therefore, by Remark \ref{free-choice-lemma-for-diff-refined}, either $\defopara^{i - r + n} \delta_{n - m}^{\infty} [y] = \defopara^{i + l} [y_{\infty}]$ for $[y_{\infty}] = [y_{\infty}]' + \defopara \epsilon$, or $\defopara^{i - r + n} \delta_{n - m}^{\infty}([y] - [y]') = \defopara^{i + l - b} [w]$ for certain essential differential $d_{r - m + l - a - b + 1}(z) = w$ with $0 < a < r - m - i$, $0 \leq b \leq r - n + l$. In the former case, the Blueprint Theorem \ref{blueprint-GLR} implies that
    \[\defopara^i f[x_{\infty}] = \defopara^{i - r + n} \delta_{n - m}^{\infty} [y] = \defopara^{i + l} [y_{\infty}]\]
    which witnesses the extension ${\color{black} d^{f, E_{\infty}, i}_{l}(x_{\infty}) = y_{\infty}}$. The latter case actually contradicts assumption item 5: from item 5 we know any such crossing satisfies $b > 0$, and the same argument leads to
    \[\defopara^i f[x_{\infty}] = \defopara^{i - r + n} \delta_{n - m}^{\infty} [y] = \defopara^{i + l - b}([w] +  \defopara^b [y_{\infty}]').\]
    Here $\defopara^{i + l - b}[w] \in \IIm(\defopara^{i + l - b}) \setminus \IIm(\defopara^{i + l - b + 1})$ since $w$ receives an essential $d_{r - m + l - a - b + 1}$-differential with $a < r - m - i$. Thus $b \leq l$, and we obtain a nontrivial extension $d_{l - b}^{f, E_{\infty}, i}(x_{\infty}) = w$ which is forbidden by item 5.
\end{proof}

\begin{remark}
    We briefly discuss how to choose parameters when we apply Theorem \ref{GLR-for-SS-2}.
    \begin{itemize}
        \item The parameters $r, m, l$ are fixed by the desideratum. 
        \item The choice of the parameters $n, j$ follows from the same discussion as in Remark \ref{explanation-of-the-parameters-GLR}. 
        \item The parameter $i$ is fine-tuned towards the $i$-complementary crossing condition in assumption item 5. The largest possible choice of $i$ is $i = r - m - 1$, which matches the statement of Theorem \ref{GLR-for-SS}. The smallest possible choice of $i$ is $i = r - n$, which is also the default choice. In practice, the most common situation is $n = r$, $i = 0$.
    \end{itemize}
    Also, the complementary crossing conditions in Theorem \ref{GLR-for-SS-2} can be simplified via Example \ref{complementary-crossing-examples}. 
    \begin{itemize}
        \item Assumption item 4 holds true if the differential in item 1 has no crossing on the $E_{n + 1}$ page or the extension in item 2 has no crossing.
        \item Assumption item 5 holds true if the differential in item 3 has no crossing on the $E_{c + 1}$ page for some $c \in \{r - m - i, \ldots, r - m\}$. In particular, if $i = r - m - 1$, then it holds true automatically.
    \end{itemize}
\end{remark}

\begin{example} \label{GLR-example}
    In this example we will deduce the differential $\color{ForestGreen} d_4(h_0 h_3^2 D_2) = g^2 n$ in the $\mathrm{H}\Fb_2$-Adams SS of $\Sb^0$ through GLR. Our notational convention will be the same as in Example \ref{example-extn-crossing}. \parr 
    \begin{figure}[htbp]
    \centering
    \scalebox{1}{\begin{tikzpicture}[line width=0.1pt]
    \tikzset{
        diff/.style={-{Stealth},line width=0.3pt,shorten <=3pt,shorten >=3pt},
        diff1/.style={-{Stealth},dashed,line width=0.6pt,shorten <=3pt,shorten >=3pt},
        extarr/.style={-{Stealth},line width=0.3pt,shorten <=3pt,shorten >=3pt},
        extarr1/.style={-{Stealth},dashed,line width=0.6pt,shorten <=3pt,shorten >=3pt},
        mapbelow/.style={-{Straight Barb[scale=0.8]},line width=0.7pt,shorten <=.7cm,shorten >=.7cm},
    }
    \draw (70.5,7.5) rectangle (72.5,13.5);
    \node at (71,7) {$71$};
    \node at (72,7) {$72$};
    \node at (70,8) {$8$};
    \node at (70,9) {$9$};
    \node at (70,10) {$10$};
    \node at (70,11) {$11$};
    \node at (70,12) {$12$};
    \node at (70,13) {$13$};
    \begin{scope}
    \clip (70.5,7.5) rectangle (72.5,13.5);
    \draw[black!10] (71.5,7.5) -- (71.5,13.5);
    \draw[black!10] (70.5,8.5) -- (72.5,8.5);
    \draw[black!10] (70.5,9.5) -- (72.5,9.5);
    \draw[black!10] (70.5,10.5) -- (72.5,10.5);
    \draw[black!10] (70.5,11.5) -- (72.5,11.5);
    \draw[black!10] (70.5,12.5) -- (72.5,12.5);
      % \draw (67.887,8.041) -- (71,9);
      % \draw (68.887,8.041) -- (72,9);
      % \draw (68,9) -- (71,10);
      % \draw (69,10) -- (72,11);
      % \draw (68.113,10.959) -- (71.113,11.959);
      \draw (72.113,7.959) -- (73,9);
      \draw (71.887,8.041) -- (72,9);
      \draw (71,9) -- (71,10);
      % \draw (68.887,11.041) -- (72,12);
      \draw (72,10) -- (72,11);
      % \draw (72,10) -- (75,11);
      \draw (71,11) -- (71.113,11.959);
      \draw (71,11) -- (72,12);
      % \draw (71,11) -- (74,12);
      \draw (72,11) -- (72,12);
      \draw (70.887,12.041) -- (72,13);
      \draw (71,13) -- (71.113,13.959);
      % \draw (71,13) -- (74.113,13.959);
      \draw (70.774,13.082) -- (70.887,14.041);
      % \draw (70.774,13.082) -- (73.887,14.041);
      \draw[diff, ForestGreen] (72,9) -- (71.226,12.918);
      \draw[diff, cyan] (70.774,13.082) -- (70,15);
      \draw[diff, cyan] (71.887,8.041) -- (71,10);
      \draw[diff, cyan] (72,10) -- (71.113,11.959);
      \fill (72.113,7.959) circle (0.064);
      \fill (71.887,8.041) circle (0.064) node[left] {\resizelabel{0.667}{$h_3^2D_2$}};
      \fill (71,9) circle (0.064) node[below] {\resizelabel{0.667}{$h_2 \Delta g_2$}};
      \fill (72,9) circle (0.064);
      \fill (71,10) circle (0.064);
      \fill (72,10) circle (0.064) node[below] {\resizelabel{0.667}{$d_0D_2$}};
      \fill (71,11) circle (0.064) node[below] {\resizelabel{0.667}{$d_0 Q_2$}};
      \fill (72,11) circle (0.064);
      \fill (71.113,11.959) circle (0.064);
      \fill (72,12) circle (0.064);
      \fill (71.226,12.918) circle (0.064) node[above right] {\resizelabel{0.667}{$g^2n$}};
      \fill (72,13) circle (0.064);
      \fill (70.887,12.041) circle (0.064) node[below] {};
      \fill (71,13) circle (0.064) node[below] {};
      \fill (70.774,13.082) circle (0.064) node[below] {};
    \end{scope}\end{tikzpicture}}
    \caption{Deducing $\color{ForestGreen} d_4(h_0 h_3^2 D_2) = g^2 n$ through generalized Leibniz rule.}
    \label{figure-GLR-example}
    \end{figure}
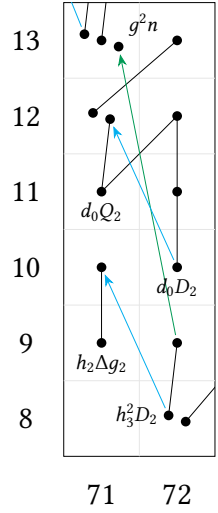

    We choose our inputs as follows: 

    \begin{itemize}
        \item The differential ${\color{cyan} d_2(h_3^2 D_2) = h_0 h_2 \Delta g_2}$ in the Adams SS of $\Sb^0$. 
        \item The distinguished triangle (from \cite[Proposition 3.20]{Lin-Wang-Xu-kervaire})
        \begin{equation*}
            \mathbb{S}^{1,2}_{\mathrm{H}\mathbb{F}_2}\xrightarrow{[h_0]}\mathbb{S}^{0,0}_{\mathrm{H}\mathbb{F}_2}\xrightarrow{i} \cofib(2)_{\mathrm{H}\mathbb{F}_2} \xrightarrow{p}\mathbb{S}^{1,1}_{\mathrm{H}\mathbb{F}_2}
        \end{equation*}
        together with the extension $d_2^{[h_0], E_{4}}(h_0 h_2 \Delta g_2) = g^2 n$. The latter can be deduced via GMT (Theorem \ref{GMT-for-SS}) from the differential ${\color{red} d_3(h_0 h_2 \Delta g_2[1]) = g^2 n[0]}$ in the Adams SS of $\cofib(2)$. 
        \item The differential ${\color{cyan} d_2(d_0 D_2) = h_0 d_0 Q_2}$ in the Adams SS of $\Sb^0$. Applying GMT to the above triangle, this leads to an extension $d_{1}^{p, E_{3}}(d_0 D_2[0]) = d_0 Q_2$. Furthermore, $d_0 D_2[0]$ is a $d_3$-cycle\footnote{Note that $d_0 D_2[0] = h_1^2 D_3'[1]$ and $D_3'[1]$ supports an essential $d_2$. Therefore, if $d_0 D_2[0]$ supports an essential $d_3$, then either $h_1 D_3'[1]$ supports an essential $d_2$ whose target is $h_1$-divisible, or $h_1 D_3'[1]$ supports an essential $d_3$ whose target is not $h_1$-torsion. In fact, neither of these is possible by $E_2$-page data.}, 
        and $d_3(d_0 D_2[0]) = 0$ gives rise to (via GMT) another extension $d_1^{[h_0], E_{3},1}(d_0 Q_2) = 0$. 
    \end{itemize}

    We will apply GLR (Theorem \ref{GLR-for-SS}) to the map $[h_0]\colon \mathbb{S}^{0,0}_{\mathrm{H}\mathbb{F}_2} \to \mathbb{S}^{-1, -2}_{\mathrm{H}\mathbb{F}_2}$ to propagate the differential ${\color{cyan} d_2(h_3^2 D_2) = h_0 h_2 \Delta g_2}$ in its domain. In this case we take the parameters as follows: $r = n = 1, m = i = j = 0, l = 2, U = 3$. As the source extension $d_0^{[h_0], E_2}(h_3^2 D_2) = h_0 h_3^2 d_2$ has no crossing for degree reasons, it suffices to show the target extension $d_2^{[h_0], E_{4}}(h_0 h_2 \Delta g_2) = g^2 n$ has no crossing\footnote{The extension has a crossing in the sense of Lin--Wang--Xu, namely $d_1^{[h_0], E_{3}}(d_0 Q_2) = h_0 d_0 Q_2$. Thus this is a place where our notion of $j$-th layer extensions leads to extra flexibility.} in the sense of Definition \ref{extensions-def}. In fact, a crossing (if it exists) has to be either of the form $d_0^{[h_0], E_{2}, 2}(z) = h_0 z = w$ for some $z$ of stem $71$, filtration $12$, or of the form $d_1^{[h_0], E_{3}, 1}(z) = w$ for some $z$ of stem $71$, filtration $11$. The former case is impossible since the classes of filtration $12$ are $h_0$-torsion, while the latter case is impossible since the only candidate of filtration $11$ is $d_1^{[h_0], E_{3},1}(d_0 Q_2) = 0$. Thus we conclude that $\color{ForestGreen} d_4(h_0 h_3^2 D_2) = g^2 n$ in the Adams SS of $\Sb^0$.
    
    \FloatBarrier
\end{example}

\begin{remark} \label{credits-for-theorem-B}
We end this section with some further historical comments and attributions concerning the generalized Leibniz rule. The relationship with Lin--Wang--Xu's Adams SS version \cite[Theorem 6.1]{Lin-Wang-Xu-kervaire} has already been discussed above.

\begin{itemize}

\item Suppose $R$ is a filtered ring spectrum and $[x] \in \pi_{t - s, t}(R)$ is a lift of $x \in E_{\infty}^{s, t}(R)$. Applying the GLR to the map $[x]\colon \Sigma^{t-s,t}R \to R$ recovers a special case of the classical Leibniz rule, namely $d_r(xy) = (-1)^{t-s}x \cdot d_r(y)$ for any $y \in E_r^{*,*}(R)$. This principle is classical and has been used throughout the study of multiplicative SS.

\item The main innovation underlying the GLR is the naturality of \emph{total} differentials $\delta$ with respect to arbitrary maps $f\colon X \to Y$ in $\Fil\Sp$. In practice, this naturality is used in two ways:
\begin{enumerate}
\item[(a).] If $\delta [x] = \defopara^? [x_{\infty}]$ in $X$, and if there are extensions $f[x] = \defopara^? [y]$ and $f[x_{\infty}] = \defopara^? [y_{\infty}]$, then one obtains $\delta[y] = \defopara^? [y_{\infty}]$ in $Y$.
\item[(b).] If $\delta [x] = \defopara^? [x_{\infty}]$ in $X$, $\delta[y] = \defopara^? [y_{\infty}]$ in $Y$, and there is an extension between the sources, $f[x] = \defopara^? [y]$, then one gets an extension between the targets, $f[x_{\infty}] = \defopara^? [y_{\infty}]$.
\end{enumerate}

\item Perspective (b) appears in several computations:
\begin{itemize}
\item Burklund--Hahn--Senger \cite[Theorem A.20(14)]{BHS1} use this perspective to deduce the extension
\[
2[d_0] = \defopara [h_0][d_0] = \delta_1^{\infty}(h_0 h_4) = [h_0] \delta_1^{\infty}(h_4)  = [h_0]^2 [h_3]^2
\]
in the $\HH \Fb_2$-synthetic stems.
\item Isaksen--Kong--Li--Ruan--Zou \cite{IKLRZ} use this perspective to deduce hidden extensions in the ($\Cb$-motivic = $\MU$-synthetic) ANSS of (the synthetic lift of) $\mathrm{tmf}$.
\item Marek \cite{Marek} uses this perspective to deduce hidden extensions in the ($\HH \Fb_2$-synthetic) Adams SS of (the synthetic lift of) $\mathrm{tmf}$.
\end{itemize}

\item Perspective (a) also occurs in several places:
\begin{itemize}
\item Hill--Hopkins--Ravenel \cite{HHR-BPC4<1>} use a stem-level form of this perspective to deduce exotic transfer differentials in the $C_4$-slice SS of $\BPCfour\langle 1 \rangle$. Their argument also uses hidden extension input of the kind related to the generalized Mahowald trick; see Remark \ref{credits-for-theorem-C} for a more detailed discussion of this point.
\item Burklund--Xu \cite[\S~7]{Burklund-Xu} use this perspective, applied to various $\HH \Fb_2$-synthetic stems lifting $h_j^2$, to deduce the Adams differentials on the family $h_j^3$.
\item Carrick--Davies--van Nigtevecht describe this perspective in \cite[Remark 3.41]{vN25}, and use it systematically in \cite{CDvN2} to deduce differentials in the descent SS of $\mathrm{Tmf}$.
\end{itemize}

\item The applications above operate at the level of bigraded stems. To crystallize the GLR into a statement that takes spectral sequence data as input, one has to translate carefully between bigraded stems and spectral sequence terms, and impose appropriate no-crossing conditions.

An earlier formulation in this direction appears in Chua \cite[Theorem 12.9]{Chua} for the Adams SS. As stated, however, this result cannot be correct without a no-crossing hypothesis; see \cite[Example 6.10]{Lin-Wang-Xu-kervaire} for a counterexample.
\end{itemize}
\end{remark}

%% file: GMT.tex
\label{sec:4}

In this section, we prove the generalized Mahowald trick, which translates between hidden extensions along a map $f\colon X \to Y$ in $\Fil\Sp$ and differentials in the standard SS of $Z = \cofib(f)$. Roughly speaking, this is done in four steps: 
\begin{enumerate}
    \item We construct from the two distinguished triangles 
    \[X \xrightarrow{f} Y \xrightarrow{g} Z \xrightarrow{h} \Sigma^{1, 0} X \ \ \text{ and } \ \ \Sb^{0, -k} / \defopara^{r - k} \xrightarrow{\defopara^k} \Sb^{0, 0} / \defopara^r \xrightarrow{\rho} \Sb^{0, 0} / \defopara^k \xrightarrow{\delta_k} \Sb^{1, -k} / \defopara^{r - k}\] 
    a ``$2$d distinguished triangle'' in $\Fil\Sp$, i.e. a diagram $F\langle -, -\rangle\colon \overrightarrow{\Zb} \times \overrightarrow{\Zb} \to \Fil\Sp$ 
    \[\begin{tikzpicture}[baseline= (a.base)]
        \node[scale=.75] (a) at (0,0){
            \begin{tikzcd}
                & \vdots \ar[d] & \vdots \ar[d] & \vdots \ar[d] & \vdots \ar[d] & \vdots \ar[d] & \\
                \cdots \ar[r] & * \ar[r] \ar[d] & * \ar[r] \ar[d] & * \ar[r] \ar[d] & \Sigma^{0, -k} X / \defopara^{r - k} \ar[r] \ar[d]&  {\color{blue} \Sigma^{0, -k} Y / \defopara^{r - k}} \ar[r] \ar[d, blue, "{\defopara^k}"] & \cdots \\
                \cdots \ar[r] & * \ar[r] \ar[d] & * \ar[r] \ar[d] &  \Sigma^{-1, 0} Z / \defopara^r \ar[r] \ar[d] &  {\color{blue} X / \defopara^r} \ar[r, blue, "{f}"] \ar[d, blue, "{\rho}"]&  {\color{blue} Y / \defopara^r} \ar[r] \ar[d] & \cdots \\
                \cdots \ar[r] & * \ar[r] \ar[d] & \Sigma^{-1, 0} Y / \defopara^{k} \ar[r] \ar[d] & {\color{blue} \Sigma^{-1, 0} Z / \defopara^{k}}  \ar[r, blue, "{h}"] \ar[d, blue, "{\delta_k}"] & {\color{blue} X / \defopara^{k}} \ar[r] \ar[d]& Y / \defopara^k \ar[r] \ar[d] & \cdots \\
                \cdots \ar[r] & {\color{black} \Sigma^{0, -k} X / \defopara^{r - k}} \ar[r] \ar[d] &{\color{blue} \Sigma^{0, -k} Y / \defopara^{r - k}} \ar[r, blue, "{g}"] \ar[d, blue, "{\defopara^k}"] & {\color{blue} \Sigma^{0, -k} Z / \defopara^{r - k}} \ar[r] \ar[d] & \Sigma^{1, -k} X / \defopara^{r - k} \ar[r] \ar[d]& * \ar[r] \ar[d] & \cdots \\
                \cdots \ar[r] &  {\color{blue} X / \defopara^r} \ar[r, blue, "{f}"] \ar[d] & {\color{blue} Y / \defopara^r} \ar[r] \ar[d] & Z / \defopara^r \ar[r] \ar[d] & * \ar[r] \ar[d]& * \ar[r] \ar[d] & \cdots \\
                & \vdots & \vdots & \vdots & \vdots & \vdots & 
            \end{tikzcd}
        };
    \end{tikzpicture}\] 
    together with the coherence data witnessing the fact that every three horizontally (resp. vertically) consecutive arrows form a distinguished triangle. 
    \item We then construct two families of isomorphisms between pullbacks of various {\color{blue} cospans} inside a ``$2$d triangle'' and clarify their relation by Theorem \ref{coherent-Mahowald-trick}, the ``coherent Mahowald trick''. 
    \item Using the previous steps, we extract a pair of isomorphisms from the above ``$2$d triangle''
    \[X / \defopara^r  \times_{f, Y / \defopara^r, \defopara^k}  \Sigma^{0, -k} Y/\defopara^{r - k} \xrightarrow{\sim} 
    \Sigma^{-1, 0} Z/\defopara^{k} \times_{h, X/\defopara^{k}, \rho} X / \defopara^r \xrightarrow{\sim} 
    \Sigma^{0, -k} Y/\defopara^{r - k} \times_{g, \Sigma^{0, -k} Z/\defopara^{r - k}, \delta_k} \Sigma^{-1, 0} Z/\defopara^{k}\]
    compatible with all six projections to $X / \defopara^r$, $\Sigma^{0, -k} Y / \defopara^{r - k}$ and $\Sigma^{-1, 0} Z / \defopara^k$, cf. Theorem \ref{blueprint-GMT}. 
    \item Finally, we apply the Bockstein dictionary in \S~\hyperref[subsec:2.2]{2.2} to deduce the  generalized Mahowald trick in terms of standard spectral sequences.
\end{enumerate}

We organize this section as follows: In \S~\hyperref[subsec:4.1]{4.1} we define the notion of ``$2$d triangles'' and construct the internal isomorphisms. We also state the coherent Mahowald trick and defer its proof to Appendix \ref{app:B}. In \S~\hyperref[subsec:4.2]{4.2} we specialize the discussions to $\Fil\Sp$ and prove the generalized Mahowald trick in Theorems \ref{GMT-for-SS} and \ref{GMT-for-SS-2}.

\subsection{\texorpdfstring{$2$d triangles and the coherent Mahowald trick}{2d triangles and the coherent Mahowald trick}}

\label{subsec:4.1}

To propose a definition for ``$2$d distinguished triangles'', we start from a definition that coherently lifts a classical ($1$d) distinguished triangle, together with all its rotations, in a stable $\infty$-category $\CC$. 

\begin{definition}
    Let $\CP$ be the sub-poset of $(\Zb, \leq) \times (\Zb, \leq)$ consisting of pairs $(m, n)$ such that $-1 \leq m - n \leq 2$. For any stable $\infty$-category $\CC$, a \textbf{coherent unrolled distinguished triangle (cuDT)} in $\CC$ is a functor $F\colon \CP \to \CC$ such that $F(m, n) = 0$ if $m - n = 2$ or $-1$, and the image of each square
    \[\begin{tikzcd}
        (m, n) \ar[r] \ar[d] & (m + 1, n) \ar[d] \\
        (m, n + 1) \ar[r] & (m + 1, n + 1) 
    \end{tikzcd}\]
    is cartesian in $\CC$. We write $\DT^u(\CC)$ for the full subcategory of $\Fun(\CP, \CC)$ spanned by cuDTs. 
\end{definition} 

\begin{remark} \label{drawing-1d-cuDTs}
    One can draw a cuDT as follows: consider the ``diagonal embedding'' map of posets
    \[t\colon (\Zb, \leq) \to \CP, \qquad n \mapsto \left(\bigg\lfloor \frac{n + 1}{2} \bigg\rfloor, \bigg\lceil \frac{n - 1}{2}  \bigg\rceil\right) =  \begin{cases}
        (\frac{n}{2}, \frac{n}{2}), & 2 \mid n. \\
        (\frac{n + 1}{2}, \frac{n - 1}{2}), & 2 \nmid n.
    \end{cases} \]
    and write $F\langle n \rangle \coloneq F(t(n))$. Then  $F \in \DT^u(\CC)$ is a diagram
    \[\begin{tikzpicture}[baseline= (a.base)]
        \node[scale=.8] (a) at (0,0){
            \begin{tikzcd}[row sep = {3.5em, between origins}, column sep = {5em, between origins}]
                & \ddots \ar[r] \ar[d] &  0 \ar[d] \ar[rd, equal] & \\
                \ddots \ar[r] \ar[rd, equal] & F\langle -2 \rangle \ar[r] \ar[d] & F\langle -1 \rangle \ar[r]  \ar[d] & 0 \ar[d] \ar[rd, equal]  & & \\
                & 0 \ar[r] \ar[rd, equal] & F\langle 0 \rangle \ar[r] \ar[d] & F\langle 1 \rangle \ar[r] \ar[d] & 0  \ar[d] \ar[rd, equal] & \\
                && 0 \ar[r] \ar[rd, equal] & F\langle 2 \rangle \ar[r] \ar[d] & F\langle 3 \rangle \ar[r] \ar[d] & \ddots\\
                &&& 0 \ar[r] & \ddots & 
            \end{tikzcd}
        };
    \end{tikzpicture}\]  
    in which each square is cartesian. Often for convenience we suppress null homotopies in the picture, and the resulting diagram of $F$ (more precisely, $t^* F$) becomes 
    \[\begin{tikzcd}
        \cdots \ar[r] & F\langle -2 \rangle \ar[r] & F\langle -1 \rangle \ar[r] & F\langle 0 \rangle \ar[r] & F\langle 1 \rangle \ar[r] & F\langle 2 \rangle \ar[r] & F\langle 3 \rangle \ar[r] & \cdots
    \end{tikzcd}\]
\end{remark}

We treat a cuDT as a coherent refinement of DT due to the following 

\begin{theorem}\label{1d-cuDT-generation}
    Consider the inclusions of posets 
    \[\begin{tikzcd}
        \CP'' \ar[r, "{i_1}"] & \CP' \ar[r, "{i_2}"] & \CP
    \end{tikzcd}\]
    where $\CP' = \{(m, n) \in \CP \;|\; (0, 0) \leq (m, n) \leq (2, 1)\} \cong \Deltatwo \times \Deltaone$ and $\CP'' = \{(0, 0) \leq (1, 0)\} \cong \Deltaone$. Restriction along these functors induces equivalences of $\infty$-categories 
    \[\begin{tikzcd}
        \DT^u(\CC) \ar[r, "{i_2^*}", "{\cong}"'] & \DT(\CC) \ar[r, "{i_1^*}", "{\cong}"'] & \CC^{\Deltaone}.
    \end{tikzcd}\]
    Here $\DT(\CC)$ is the full subcategory of $\Fun(\Deltatwo \times \Deltaone, \CC)$ spanned by diagrams
    \begin{equation*} \label{star-rectangle}
        \begin{tikzcd}
            X_0 \ar[r] \ar[d] & X_1 \ar[r] \ar[d] & 0 \ar[d] \\
            0 \ar[r] & X_2 \ar[r] & X_3
        \end{tikzcd}
        \tag{$\star$}
    \end{equation*}
    in which the upper-right corner and the lower-left corner are $0$, and each square is cartesian in $\CC$. 
\end{theorem}

Intuitively, this says a cuDT is freely generated by a map $X_0 \to X_1$ via taking cofibers and fibers recursively. It also explains the relation between a cuDT and a DT: according to \cite[Definition 1.1.2.11 and Theorem 1.1.2.14]{HA}, if $\CC$ is stable, then a triangle $X_0 \to X_1 \to X_2 \to \Sigma X_0$ in $h\CC$ is distinguished iff there is a diagram (\ref{star-rectangle}) in $\CC$, in which each square is cartesian, so that 
\begin{itemize}
    \item The maps $X_0 \to X_1$ and $X_1 \to X_2$ in the triangle coincide with the two maps in the diagram.
    \item The map $X_2 \to \Sigma X_0$ in the triangle is the composite of $X_2 \to X_3$ in the diagram and the inverse to the comparison map $\Sigma X_0 \to X_3$ induced by the outer rectangle. 
\end{itemize}
Therefore, the theorem also says a cuDT contains the same information as a coherent refinement of a DT (and the latter is also freely generated by the map $X_0 \to X_1$ via taking cofibers twice). 

\begin{proof}
    It suffices to show the composite $i_1^* i_2^*$ is an equivalence, as the same proof would imply $i_1^*$ is an equivalence, while $i_2^*$ is an equivalence by $2$-out-of-$3$. Write $\partial \CP = \CP \setminus \mathrm{Im} (t)$, and write $\CP^{\geq 0} = \{(m, n) \in \CP \;|\; (m, n) \geq (0, 0)\}$. 
    Note that $i_1^* i_2^*$ is the composite of restriction functors
    \begin{align*}
        \DT^{u}(\CC) &= \left\{F \colon \CP \to \CC \;\middle|\;
        \begin{aligned}
            &F(m, n) = 0 \text{ if } (m, n) \in \partial \CP \text{ and the} \\ 
            &\text{image of each square is cartesian in } \CC
        \end{aligned} \right\} \\[6pt]
        &\to \left\{F \colon \CP^{\geq 0} \cup \partial \CP \to \CC \;\middle|\;
        \begin{aligned}
            &F(m, n) = 0 \text{ if } (m, n) \in \partial \CP \text{ and the} \\ 
            &\text{image of each square is cartesian in } \CC
        \end{aligned} \right\} \\[6pt]
        &\to \left\{F \colon \CP''\cup \partial \CP \to \CC \;\middle|\; F(m, n) = 0 \text{ if } (m, n) \in \partial \CP\right\} \\[6pt]
        &\to \left\{F \colon \CP''\cup (\partial \CP \cap \CP^{\geq 0}) \to \CC \;\middle|\; F(m, n) = 0 \text{ if } (m, n) \in \partial \CP\cap \CP^{\geq 0}\right\} \\[6pt]
        &\to \Fun({\CP''}, \CC) = \CC^{\Deltaone}. 
    \end{align*}
    We claim each functor here is an equivalence. Actually, the first functor has an inverse, which sends $F \colon \CP^{\geq 0} \cup \partial \CP \to \CC$ to its right Kan extension along the inclusion $\CP^{\geq 0} \cup \partial \CP \to \CP$. Similarly, the fourth restriction functor also has an inverse given by right Kan extension, while the second functor and the third functor admit inverses given by left Kan extensions. 
\end{proof}

Rotations of DTs correspond to certain autoequivalences of $\DT^u$. 

\begin{fact}
    There are isomorphisms of posets, inverse to each other: 
    \[+1\colon \CP \to \CP, \qquad (m, n) \mapsto (n + 1, m)\]
    \[-1\colon \CP \to \CP, \qquad (m, n) \mapsto (n, m - 1)\]
    for which precomposition functors induce autoequivalences $(+1)^*$ and $(-1)^*$ on $\DT^u(\CC)$. \parr 
    
    For $F \in \DT^u(\CC)$, $((+1)^* F)\langle n \rangle = F\langle n + 1 \rangle$ and $((-1)^* F)\langle n \rangle = F\langle n - 1 \rangle$. Intuitively, $(+1)^*$ rotates the DT once to the right, while $(-1)^*$ rotates it once to the left. \parr 

    In general, for each $k \in \Zb$ there is a poset automorphism $+k\colon \CP \to \CP$ given as $(+1)^{\circ k}$, and it induces an autoequivalence $(+k)^*$ on $\DT^u(\CC)$. This defines a $\Zb$-action on $\DT^u(\CC)$ natural on $\CC \in \Cat_{\st}$. For notational convenience, we will refer to $(+k)^* F$ as $F[+k]$. 
\end{fact}

In a triangulated category, if we rotate a DT $(h, g, f)$ three times, we see $(-\Sigma h, - \Sigma g, - \Sigma f)$ is another DT. This suggests that for $F \in \DT^u(\CC)$, there is an isomorphism between $\Sigma F$ and $F[+3]$. 

\begin{construction}[Alpha isomorphisms] \label{alpha-iso}
    Consider the following map of posets
    \[\Theta\colon \CP \times \Deltaone \times \Deltaone \to \CP, \qquad (m, n, a, b) \mapsto 
    \begin{cases}
        (m, n), & \text{ if } a = b = 0. \\
        (n + 2, n), & \text{ if }  a = 1, b = 0. \\
        (m, m + 1), & \text{ if }  a = 0, b = 1. \\
        (n + 2, m + 1), & \text{ if }  a = b = 1.
    \end{cases}\]
    Precomposition with $\Theta$ induces a functor $\Theta^*\colon \DT^u(\CC) \to \DT^u(\CC)^{\Deltaone \times \Deltaone}$, so that $\Theta^* F |_{a = b = 0} = F$, $\Theta^* F|_{a = b = 1} = F[+3]$, while $\Theta^* F|_{a = 1, b = 0}$ and $\Theta^* F|_{a = 0, b = 1}$ are both constant functors $\CP \to \CC$ taking value $0$. Post-composing $\Theta^*$ with the pushout comparison functor
    \[\cp\colon \DT^{u}(\CC)^{\Deltaone \times \Deltaone} \to \DT^{u}(\CC)^{\Deltaone}, \qquad 
    \begin{tikzcd}
        F \ar[r] \ar[d] \ar[dr, phantom] & G \ar[d] \\
        H \ar[r] & K
    \end{tikzcd}
    \mapsto 
    \begin{tikzcd}
        \colim(H \gets F \to G) \ar[d]\\
        K
    \end{tikzcd}\]
    we get a family of isomorphisms $\alpha = \alpha_F \colon \Sigma F \to F[+3]$ in $\DT^u(\CC)$ natural in $F$. 
\end{construction}

\begin{remark} \label{alpha[+1]-is-minus-alpha}
    Suppose $F \in \DT^u(\CC)$ and $x \in \Zb$. Then the map $\alpha_{F[+1], x}\colon \Sigma F[+1]\langle x \rangle = \Sigma F\langle x + 1 \rangle \to  F[+1]\langle x + 3 \rangle = F\langle x + 4 \rangle$ comes from the transpose of the square defining $\alpha_{F, x + 1}\colon \Sigma F\langle x + 1 \rangle \to F\langle x + 4 \rangle$. Therefore, $\alpha_{F[+1], x} \cong - \alpha_{F, x + 1}$ due to \cite[Lemma 1.1.2.10]{HA}. 
\end{remark}

\begin{remark} \label{DT-with-alpha=id}
    The explicit Kan extensions in the proof of Theorem \ref{1d-cuDT-generation} imply that every distinguished triangle $X \xrightarrow{f} Y \xrightarrow{g} Z \xrightarrow{h} \Sigma X$ in $\CC$ lifts to some $F \in \DT^u(\CC)$ with
    \[t^* F\colon \cdots \xrightarrow{\Sigma^{-1} f} X \xrightarrow{f} Y \xrightarrow{g} Z \xrightarrow{h} \Sigma X \xrightarrow{\Sigma f} \Sigma Y \xrightarrow{\Sigma g} \Sigma Z \xrightarrow{\Sigma h} \Sigma^2 X \xrightarrow{\Sigma^2 f} \cdots\]
    so that $\alpha_{F, x}\colon \Sigma F\langle x\rangle \to F\langle x + 3\rangle$ is $\id$ for each $x \in \Zb$. Therefore, the shape of $\CP$ (and \cite[Lemma 1.1.2.13]{HA}) suggests $F\langle x \rangle \to F\langle x + 1 \rangle \to F\langle x + 2 \rangle \to F\langle x + 3\rangle = \Sigma F\langle x \rangle$ is a DT if $x \equiv 0$ mod $2$ while it is an ADT otherwise. This matches the fact in a triangulated category that 
    \[X \xrightarrow{f} Y \xrightarrow{g} Z \xrightarrow{h} \Sigma X \text{ is a DT} \Longleftrightarrow Y \xrightarrow{g} Z \xrightarrow{h} \Sigma X \xrightarrow{\Sigma f} \Sigma Y \text{ is an ADT}.\] 
\end{remark}

Up to this point, the theory seems to be a faithful translation of what happens in triangulated categories. The merit of this coherent reformulation is that it extends naturally to diagram $\infty$-categories.

\begin{fact} \label{diagram-DT}
    \begin{itemize}
        \item The subcategory $\DT^u(\CC) \subset \Fun(\CP, \CC)$ is closed under limits and colimits. In particular, finite (co)limits exist in $\DT^u(\CC)$ and they can be computed entrywise. 
        \item For any small $\infty$-category $\CK$, the equivalence $\Fun(\CP, \Fun(\CK, \CC)) \cong \Fun(\CK, \Fun(\CP, \CC))$ restricts to an equivalence $\DT^u(\Fun(\CK, \CC)) \cong \Fun(\CK, \DT^u(\CC))$. 
    \end{itemize}
\end{fact}
\begin{proof}
    The first part follows from the stability of $\CC$ and the fact that (co)limits commute with (co)limits. For the second part, after identifying the double functor categories with $\Fun(\CK \times \CP, \CC)$, we see $\DT^u(\Fun(\CK, \CC))$ and $\Fun(\CK, \DT^u(\CC))$ both correspond to the full subcategory spanned by functors $F\colon \CK \times \CP \to \CC$ such that the slices $F(k, -)\colon \CP \to \CC$ are cuDTs for all $k \in \CK$. 
\end{proof}

We can thus formulate a natural $2$d generalization:

\begin{definition} \label{2d-cuDT}
    A \textbf{$2$-dimensional coherent unrolled distinguished triangle ($2$d cuDT)} in $\CC$ is a functor $F\colon \CP^2 \to \CC$ such that for each $(m, n) \in \CP$, the slices $F(-, -, m, n)\colon \CP \to \CC$ and $F(m, n, -, -)\colon \CP \to \CC$ are both cuDTs. We write $2\DT^u(\CC)$ for the full subcategory of $\Fun(\CP^2, \CC)$ spanned by $2$d cuDTs. Under the identification in Fact \ref{diagram-DT}, we have $2\DT^u(\CC) \cong \DT^u(\DT^u(\CC))$. 
\end{definition}

\begin{remark} \label{drawing-2d-cuDTs}
    For convenience, when drawing $F \in 2\DT^u(\CC)$ we usually suppress all nullhomotopies. In other words, write $t \times t \colon (\Zb, \leq) \times (\Zb, \leq) \to \CP \times \CP$ for the square of the diagonal embedding $t$ in Remark \ref{drawing-1d-cuDTs}, and write $F \langle m, n \rangle \coloneq F(t(m), t(n))$, then $F$ takes the form 
    \[\begin{tikzpicture}[baseline= (a.base)]
        \node[scale=.8] (a) at (0,0){
            \begin{tikzcd}
                & \vdots \ar[d] & \vdots \ar[d] & \vdots \ar[d] & \vdots \ar[d] & \vdots \ar[d] & \\
                \cdots \ar[r] & F\langle 0, 0\rangle \ar[r] \ar[d] & F\langle 1, 0\rangle \ar[r] \ar[d] & F\langle 2, 0\rangle \ar[r] \ar[d] & F\langle 3, 0\rangle \ar[r] \ar[d]& F\langle 4, 0\rangle \ar[r] \ar[d] & \cdots \\
                \cdots \ar[r] & F\langle 0, 1\rangle \ar[r] \ar[d] & F\langle 1, 1\rangle \ar[r] \ar[d] & F\langle 2, 1\rangle \ar[r] \ar[d] & F\langle 3, 1\rangle \ar[r] \ar[d]& F\langle 4, 1\rangle \ar[r] \ar[d] & \cdots \\
                \cdots \ar[r] & F\langle 0, 2\rangle \ar[r] \ar[d] & F\langle 1, 2\rangle \ar[r] \ar[d] & F\langle 2, 2\rangle \ar[r] \ar[d] & F\langle 3, 2\rangle \ar[r] \ar[d]& F\langle 4, 2\rangle \ar[r] \ar[d] & \cdots \\
                \cdots \ar[r] & F\langle 0, 3\rangle \ar[r] \ar[d] & F\langle 1, 3\rangle \ar[r] \ar[d] & F\langle 2, 3\rangle \ar[r] \ar[d] & F\langle 3, 3\rangle \ar[r] \ar[d]& F\langle 4, 3\rangle \ar[r] \ar[d] & \cdots \\
                \cdots \ar[r] & F\langle 0, 4\rangle \ar[r] \ar[d] & F\langle 1, 4\rangle \ar[r] \ar[d] & F\langle 2, 4\rangle \ar[r] \ar[d] & F\langle 3, 4\rangle \ar[r] \ar[d]& F\langle 4, 4\rangle \ar[r] \ar[d] & \cdots \\
                & \vdots & \vdots & \vdots & \vdots & \vdots & 
            \end{tikzcd}
        };
    \end{tikzpicture}\]
\end{remark}

\begin{proposition} \label{2DT-from-DTxDT}
    Suppose $\odot\colon \CC \times \CC \xrightarrow{} \CC$ preserves finite colimits separately in each variable. Then for $F_1, F_2 \in \DT^u(\CC)$, the composite $F_1 \odot F_2\colon \CP \times \CP \xrightarrow{F_1 \times F_2} \CC \times \CC \xrightarrow{\odot} \CC$ lies in $2\DT^u(\CC)$. 
\end{proposition}

\begin{proof}
    This follows directly from the definition. 
\end{proof}

\begin{theorem}\label{2d-cuDT-generation}
    Consider the map of posets 
    \(i \times i\colon \Deltaone \times \Deltaone \cong (\CP'')^2 \to \CP^2\),
    where 
    \[i = i_2 i_1\colon \CP'' = \{(0, 0) \leq (1, 0)\} \to \CP\] 
    is the poset inclusion in Theorem \ref{1d-cuDT-generation}. Restriction along $i \times i$ yields an equivalence $2\DT^u(\CC) \cong \CC^{\Deltaone \times \Deltaone}$. 
\end{theorem}
\begin{proof}
    As $2\DT^u(\CC) \cong \DT^u(\DT^u(\CC))$, this follows from applying Theorem \ref{1d-cuDT-generation} twice.
\end{proof}

Intuitively, this says a $2$d cuDT is freely generated by a square 
\[\begin{tikzcd}
    X_{00} \ar[r] \ar[d] & X_{10} \ar[d] \\
    X_{01} \ar[r] & X_{11}
\end{tikzcd}\]
by taking iterated cofibers and fibers in both directions. \parr 

\begin{remark}
    In general, for each $n \in \Nb$, we can define the $\infty$-category of \textbf{$n$-dimensional cuDTs} in $\CC$ to be $n\DT^u(\CC) = (\DT^u)^{\circ n}(\CC)$, and realize it as a full subcategory of $\Fun(\CP^{n}, \CC)$. The same argument then provides an equivalence $n\DT^u(\CC) \cong \Fun(\Deltaone^n, \CC)$. 
\end{remark}

\begin{remark} \label{2d-alpha}
    The $\Zb$-action on $\CP$ leads to a $\Zb^2$-action on $\CP^2$ and thus a $\Zb^2$-action on $2\DT^u(\CC)$, which we denote by $F \mapsto F[+m, +n]$. The comparison map $\alpha$ in Construction \ref{alpha-iso} gives rise to two natural isomorphisms here, the ``horizontal'' comparison $\alpha_h\colon \Sigma F \to F[+3, 0]$ and the ``vertical'' comparison $\alpha_v\colon  \Sigma F \to F[0, +3]$.  We write $\alpha_{vh}$ for the composite $\alpha_v\alpha_h^{-1}\colon F[+3, 0] \to F[0, +3]$. Taking specific values in $\CP^2$, this gives rise to isomorphisms $F[+3, 0] \langle 1, 0 \rangle = F \langle 4, 0 \rangle \to F[0, +3] \langle 1, 0 \rangle = F \langle 1, 3 \rangle$, etc. that are compatible with the functoriality of $F$. Therefore, it induces an isomorphism of pullbacks $\alpha_{vh}\colon F \langle 4, 0 \rangle \times_{F \langle 4, 1 \rangle} F \langle 3, 1 \rangle \to F \langle 1, 3 \rangle \times_{F \langle 1, 4 \rangle} F \langle 0, 4 \rangle$, etc.
\end{remark}

\begin{remark} \label{1dx1d-alpha}
    For instance, if $\odot\colon \CC \times \CC \to \CC$ preserves finite colimits separately in each variable and $F_1, F_2 \in \DT^u(\CC)$, then for $F = F_1 \odot F_2$, $\alpha_h\colon F\langle 0, 0\rangle \to F\langle 3, 0\rangle$ is the composite $\Sigma (F_1\langle 0\rangle \odot F_2\langle 0\rangle) \cong (\Sigma F_1\langle 0\rangle) \odot F_2\langle 0\rangle \xrightarrow{\alpha \odot \id} F_1\langle 3\rangle \odot F_2\langle 0\rangle$. A similar discussion works for $\alpha_v$.
\end{remark}

\begin{construction}[Beta isomorphisms] \label{beta-iso}
    Write $F_\partial\langle x, y\rangle = F\langle x + 1, y \rangle \times_{F\langle x + 1, y + 1\rangle} F\langle x, y + 1 \rangle$ for any $F \in 2\DT^u(\CC)$. Under this notation, the above construction yields $\alpha_{vh}\colon F_\partial \langle x + 3, y\rangle \xrightarrow{\sim} F_\partial \langle x, y + 3\rangle$. Below we will construct another family of natural isomorphisms $\beta\colon F_\partial\langle x + 1, y\rangle \xrightarrow{\sim} F_\partial\langle x, y + 1\rangle$. \parr 
    
    Consider the following posets and maps: Take $\CV = \Lambda^2_2 = \{a_1 < b_1 > c_1\}$ and write $\CW = \{a_2 < b_2 > c_2 < d_2 > e_2\}$. Write $\mu\colon \CV \times \CW \to \CV \times \CV$ for the map of posets given on elements by:

    \begin{table}[htbp]
        \centering
        \begin{tabular}{||c||c|c|c|c|c||}
            \hline
            \diagbox[width=0.28\linewidth]{$\CV$ component}{$\CW$ component} & $a_2$ & $b_2$ & $c_2$ & $d_2$ & $e_2$\\
            \hline 
            \hline 
            $a_1$ & $(a_1, a_1)$ & $(a_1, a_1)$ & $(a_1, a_1)$  & $(a_1, a_1)$  & $(a_1, a_1)$ \\
            \hline
            $b_1$ & $(b_1, a_1)$ & $(b_1, b_1)$ & $(b_1, b_1)$  & $(b_1, b_1)$  & $(a_1, b_1)$ \\
            \hline
            $c_1$ & $(c_1, a_1)$ & $(c_1, b_1)$ & $(c_1, c_1)$  & $(b_1, c_1)$  & $(a_1, c_1)$ \\
            \hline
        \end{tabular}
        % \caption{Caption}
        % \label{tab:placeholder}
    \end{table}

    \FloatBarrier

    Also, for $(m, n, a, b) \in \CP \times \CP$ with $0 \leq m - n \leq 1$, $0 \leq a - b \leq 1$, write $\psi = \psi_{mnab}\colon \CV \times \CV \to \CP \times \CP$ for the map of posets given on elements by:

    \begin{table}[htbp]
        \centering
        \begin{tabular}{||c|c||c||}
            \hline
            $1$st component & $2$nd component & $\psi$-value \\
            \hline 
            \hline
            \multirow{3}{*}{$a_1$} & $a_1$ & $(n + 1, m, b + 1, a)$ \\
            \cline{2-3} 
            & $b_1$ & $(n + 1, m , a + 1, b + 1)$\\
            \cline{2-3}
            & $c_1$ & $(m, n, a + 1, b + 1)$\\
            \hline
            \multirow{3}{*}{$b_1$} & $a_1$ & $(m + 1, n + 1, b + 1, a)$ \\
            \cline{2-3} 
            & $b_1$ & $(m + 1, n + 1, a + 1, b + 1)$\\
            \cline{2-3}
            & $c_1$ & $(2m - n, 2  n - m + 1, a + 1, b + 1)$\\
            \hline
            \multirow{3}{*}{$c_1$} & $a_1$ & $(m + 1, n + 1, a, b)$ \\
            \cline{2-3} 
            & $b_1$ & $(m + 1, n + 1, 2a - b, 2b - a + 1)$\\
            \cline{2-3}
            & $c_1$ & $(2m - n, 2n - m + 1, 2a - b, 2b - a + 1)$\\
            \hline
        \end{tabular}
        % \caption{Caption}
        % \label{tab:placeholder}
    \end{table}
    \FloatBarrier
    Write $m + n = x, a + b = y$. The composite functor 
    \[2\DT^u(\CC) \subset \Fun(\CP^2, \CC) \xrightarrow{\psi^*} \Fun(\CV^2, \CC) \xrightarrow{\mu^*} \Fun(\CV \times \CW, \CC) \cong \Fun(\CW, \CC^{\CV}) \xrightarrow{\lim_{\CV}} \Fun(\CW, \CC)\]
    sends $F \in 2\DT^u(\CC) $ to the following zig-zag diagram in $\CC$:
    \[\begin{tikzcd}
        F_\partial\langle x + 1, y\rangle \ar[r] & \fib(F\langle x + 1, y + 1 \rangle \to F\langle x + 2, y + 2\rangle) \ar[d, equal] \\
        & \fib(F\langle x + 1, y + 1\rangle \to F\langle x + 2, y + 2\rangle) \ar[d, equal] \\ 
        F_\partial\langle x, y + 1\rangle \ar[r] & \fib(F\langle x + 1, y + 1\rangle \to F\langle x + 2, y + 2 \rangle)
    \end{tikzcd}\]
    Here, the two vertical maps are both $\id$. Furthermore, part of the diagram $\psi^* F$ takes the shape 
    \[\begin{tikzcd}
        & F\langle x + 2, y \rangle \ar[r] \ar[d] & 0 = F(m + 1, n + 1, 2a - b, 2b - a + 1)\ar[d] \\
        F\langle x + 1, y + 1 \rangle \ar[r] & F\langle x + 2, y + 1 \rangle \ar[r] & F\langle x + 2, y + 2 \rangle
    \end{tikzcd}\]
    in which the right square is already cartesian, so the upper horizontal map in the zig-zag
    \[F_\partial\langle x + 1, y\rangle = F\langle x + 2, y \rangle \times_{F\langle x + 2, y + 1\rangle} F\langle x + 1, y + 1 \rangle \to 0  \times_{F\langle x + 2, y + 2\rangle} F\langle x + 1, y + 1 \rangle\]
    is an isomorphism, and the lower horizontal map is invertible for the same reason. We write $\beta\colon F_\partial\langle x + 1, y\rangle \to F_\partial\langle x, y + 1\rangle$ for the composite map. By construction, this $\beta$ is invertible, and it is compatible with projections to $F\langle x + 1, y + 1 \rangle$ on both sides.     
\end{construction} 

\begin{theorem}[coherent GMT for general $2$d cuDTs] \label{coherent-Mahowald-trick}
    For any $F \in 2\DT^u(\CC)$ and $(x, y) \in \Zb^2$, the two isomorphisms constructed above $(-1)^{x + y + 1}\alpha_{vh}\colon F_\partial\langle x + 3, y\rangle \to F_\partial\langle x, y + 3\rangle$ and $\beta^3\colon F_\partial\langle x + 3, y\rangle \to F_\partial\langle x + 2, y + 1\rangle \to F_\partial\langle x + 1, y + 2\rangle \to F_\partial\langle x, y + 3\rangle$ are homotopic. 
\end{theorem}

We defer its proof to Appendix~\ref{app:B}. 

\begin{remark}
    The same proof there also shows $(-1)^{x + y}\alpha_{hv} = (-1)^{x + y}\Sigma^{-1}(\alpha_h \alpha_v^{-1})\colon F_\partial \langle x + 3, y\rangle \to \Sigma^{-1} F_\partial \langle x + 3, y + 3\rangle \to F_\partial \langle x, y + 3\rangle$ is homotopic to $\beta^3\colon F_\partial\langle x + 3, y\rangle \to F_\partial\langle x,y + 3\rangle$.
\end{remark}

\begin{remark} \label{coherent-Mahowald-trick-as-a-modification}
    With more effort, for each stable $\infty$-category $\CC$ and each $(x, y) \in \Zb^2$ we can construct an invertible $2$-morphism $\beta^3 \cong (-1)^{x + y + 1} \alpha_{vh}$ in $\Fun(2\DT^u(\CC), \CC)$. The proof sketch for $x = y = 0$ goes as follows: We first reduce to presentable $\CC$ by changing universe and taking $\Ind$. For presentable stable $\CC$ we have the Lurie tensor product, such that $2\DT^u(\CC) \cong 2\DT^u(\Sp) \otimes \CC$ and $\alpha_{vh}, \beta$ both lie in the image of $- \otimes \CC\colon \Fun(2\DT^u(\Sp), \Sp) \to \Fun(2\DT^u(\CC), \CC)$, so it suffices to treat the case $\CC = \Sp$. Under the equivalence $2\DT^u(\Sp) \cong \Sp^{\Deltaone \times \Deltaone}$, the composite $\alpha_{vh}^{-1}\beta^3$ is an automorphism of the functor $U\colon \Sp^{\Deltaone \times \Deltaone} \to \Sp$, $T \mapsto \Sigma( T(1,0) \times_{T(1,1)} T(0, 1))$. Its left adjoint $U^L$ sends $A$ to the square $\widehat{A}$ with $\widehat{A}(0, 0) = 0$ and $\widehat{A}(1, 0) = \widehat{A}(0, 1) = \widehat{A}(1,1) = \Sigma^{-1} A$. By $\Sp$-enriched Yoneda lemma, we have $\Msp_{\Fun(\Sp^{\Deltaone \times \Deltaone}, \Sp)}(U, U) \cong \Msp_{\Fun(\Sp, \Sp)}(\id, U U^L) \cong UU^L(\Sb^0) = \Sb^0$ as an isomorphism of $\Eb_1$-ring spectra, thus the automorphism $\alpha_{vh}^{-1}\beta^3$ is homotopic to either $\id$ or $-\id$. Following the proof of Theorem \ref{coherent-Mahowald-trick} in Appendix~\ref{app:B}, we can construct a specific object $\mathrm{H}\Qb^{\phi} \in 2\DT^u(\Sp)$ on which $\beta^3 \cong - \alpha_{vh}$ and $\alpha_{vh} \not\cong -\alpha_{vh}$ (this is the tricky step). Therefore, $\beta^3 \cong -\alpha_{vh}$ globally in the functor category.
\end{remark}

\subsection{Generalized Mahowald trick for SS}
\label{subsec:4.2}

We start by specializing the above discussion for $2$d cuDTs to $\CC = \Fil\Sp$. This leads to a ``prototype'' of GMT, which identifies hidden extensions along a map $f\colon X \to Y$ in $\Fil\Sp$ with (total) differentials in its cofiber $Z$ on the level of bigraded stems.

\begin{theorem}[Blueprint for GMT] \label{blueprint-GMT}
    Suppose $X \xrightarrow{f} Y \xrightarrow{g} Z \xrightarrow{h} \Sigma^{1, 0} X$ is a distinguished triangle in $\Fil\Sp$. Then for $1 \leq k < r \leq \infty$, there are isomorphisms of pullbacks 
    \[X / \defopara^r  \times_{f, Y / \defopara^r, \defopara^k}  \Sigma^{0, -k} Y/\defopara^{r - k} \xrightarrow[\sim]{\beta} 
    \Sigma^{-1, 0} Z/\defopara^{k} \times_{h, X/\defopara^{k}, \rho} X / \defopara^r \xrightarrow[\sim]{\beta} 
    \Sigma^{0, -k} Y/\defopara^{r - k} \times_{g, \Sigma^{0, -k} Z/\defopara^{r - k}, \delta_k} \Sigma^{-1, 0} Z/\defopara^{k}\]
    that are compatible with projections to $X / \defopara^r, \Sigma^{0, -k} Y/\defopara^{r - k}$ and $\Sigma^{-1, 0} Z/\defopara^{k}$. In particular, 
    \begin{itemize}
        \item If $[a] \in \pi_{n, n + s}(X / \defopara^r)$ and $[b] \in \pi_{n, n + s + k} (Y / \defopara^{r - k})$ satisfy $f[a] = \defopara^k [b]$, then there exists $[c] \in \pi_{n + 1, n + s}(Z / \defopara^k)$ with $\delta_k[c] = g[b]$ and $h[c] = \rho[a]$. 
        \item If $[b] \in \pi_{n, n + s}(Y / \defopara^{r - k})$ and $[c] \in \pi_{n + 1, n + s - k} (Z / \defopara^{k})$ satisfy $g[b] = \delta_k [c]$, then there exists $[a] \in \pi_{n, n + s - k}(X / \defopara^r)$ with $f[a] = \defopara^k[b]$ and $\rho[a] = h[c]$. 
        \item If $[c] \in \pi_{n, n + s} (Z / \defopara^{k})$ and $[a] \in \pi_{n - 1, n + s}(X / \defopara^r)$ satisfy $h[c] = \rho[a]$, then there exists $[b] \in \pi_{n - 1, n + s + k}(Y / \defopara^{r - k})$ with $g[b] = \delta_k[c]$ and $\defopara^k [b] = f[a]$. 
    \end{itemize}
\end{theorem}

\begin{proof}
    For the first part (before ``in particular''), we begin by lifting the distinguished triangles 
    \[X \xrightarrow{f} Y \xrightarrow{g} Z \xrightarrow{h} \Sigma^{1, 0} X \ \ \text{ and } \ \ \Sb^{0, -k} / \defopara^{r - k} \xrightarrow{\defopara^k} \Sb^{0, 0} / \defopara^r \xrightarrow{\rho} \Sb^{0, 0} / \defopara^k \xrightarrow{\delta_k} \Sb^{1, -k} / \defopara^{r - k}\] 
    to $F_1, F_2 \in \DT^u(\Fil\Sp)$ according to Theorem \ref{1d-cuDT-generation} and Remark \ref{DT-with-alpha=id}. Then we form $F = F_1 \otimes F_2 \in 2\DT^u(\Fil\Sp)$ following Proposition \ref{2DT-from-DTxDT}. This $F$ (more precisely, $(t \times t)^* F$) takes the form
    \[\begin{tikzpicture}[baseline= (a.base)]
        \node[scale=.75] (a) at (0,0){
            \begin{tikzcd}
                & \vdots \ar[d] & \vdots \ar[d] & \vdots \ar[d] & \vdots \ar[d] & \vdots \ar[d] & \\
                \cdots \ar[r] & * \ar[r] \ar[d] & * \ar[r] \ar[d] & * \ar[r] \ar[d] & \Sigma^{0, -k} X / \defopara^{r - k} \ar[r] \ar[d]&  {\color{blue} \Sigma^{0, -k} Y / \defopara^{r - k}} \ar[r] \ar[d, blue, "{\defopara^k}"] & \cdots \\
                \cdots \ar[r] & * \ar[r] \ar[d] & * \ar[r] \ar[d] &  \Sigma^{-1, 0} Z / \defopara^r \ar[r] \ar[d] &  {\color{blue} X / \defopara^r} \ar[r, blue, "{f}"] \ar[d, blue, "{\rho}"]&  {\color{blue} Y / \defopara^r} \ar[r] \ar[d] & \cdots \\
                \cdots \ar[r] & * \ar[r] \ar[d] & \Sigma^{-1, 0} Y / \defopara^{k} \ar[r] \ar[d] & {\color{blue} \Sigma^{-1, 0} Z / \defopara^{k}}  \ar[r, blue, "{h}"] \ar[d, blue, "{\delta_k}"] & {\color{blue} X / \defopara^{k}} \ar[r] \ar[d]& Y / \defopara^k \ar[r] \ar[d] & \cdots \\
                \cdots \ar[r] & {\color{red} \Sigma^{0, -k} X / \defopara^{r - k}} \ar[r] \ar[d] &{\color{blue} \Sigma^{0, -k} Y / \defopara^{r - k}} \ar[r, blue, "{g}"] \ar[d, blue, "{\defopara^k}"] & {\color{blue} \Sigma^{0, -k} Z / \defopara^{r - k}} \ar[r] \ar[d] & \Sigma^{1, -k} X / \defopara^{r - k} \ar[r] \ar[d]& * \ar[r] \ar[d] & \cdots \\
                \cdots \ar[r] &  {\color{blue} X / \defopara^r} \ar[r, blue, "{f}"] \ar[d] & {\color{blue} Y / \defopara^r} \ar[r] \ar[d] & Z / \defopara^r \ar[r] \ar[d] & * \ar[r] \ar[d]& * \ar[r] \ar[d] & \cdots \\
                & \vdots & \vdots & \vdots & \vdots & \vdots & 
            \end{tikzcd}
        };
    \end{tikzpicture}\] 
    Here $F\langle 0, 0\rangle = {\color{red} \Sigma^{0, -k} X / \defopara^{r - k}}$ is at the lower-left corner. Construction \ref{beta-iso} supplies isomorphisms 
    \[\beta\colon F_\partial \langle 3, -3 \rangle = X / \defopara^r  \times_{Y / \defopara^r}  \Sigma^{0, -k} Y/\defopara^{r - k} \xrightarrow[]{\sim} F_\partial\langle 2, -2 \rangle = \Sigma^{-1, 0} Z/\defopara^{k} \times_{X/\defopara^{k}} X / \defopara^r\]
    \[\beta\colon F_\partial\langle 2, -2 \rangle = \Sigma^{-1, 0} Z/\defopara^{k} \times_{X/\defopara^{k}} X / \defopara^r \xrightarrow{\sim} F_\partial\langle 1, -1 \rangle = \Sigma^{0, -k} Y/\defopara^{r - k} \times_{\Sigma^{0, -k} Z/\defopara^{r - k}} \Sigma^{-1, 0} Z/\defopara^{k}\]
    compatible with projections to $F\langle 3, -2\rangle = X / \defopara^r$ and $F\langle 2, -1\rangle = \Sigma^{-1, 0} Z / \defopara^k$. It remains to show $\beta^2\colon F_\partial \langle 3, -3 \rangle \to F_\partial \langle 1, -1 \rangle$ is compatible with projections to $\Sigma^{0, -k} Y / \defopara^{r - k}$ on both sides. Note that
    \[\beta\colon F_\partial\langle 1, -1 \rangle = \Sigma^{0, -k} Y/\defopara^{r - k} \times_{\Sigma^{0, -k} Z/\defopara^{r - k}} \Sigma^{-1, 0} Z/\defopara^{k} \xrightarrow[]{\sim} F_\partial \langle 0, 0 \rangle = X / \defopara^r  \times_{Y / \defopara^r} \Sigma^{0, -k} Y/\defopara^{r - k}\]
    is compatible with projections to $F\langle 1, 0\rangle = \Sigma^{0, -k} Y / \defopara^{r - k}$ and the identification $F\langle 4, -3 \rangle = \Sigma^{1, 0} Y \otimes  \Sb^{-1, -k} / \defopara^{r - k} \cong  \Sigma^{0, -k} Y / \defopara^{r - k} \to F\langle 1, 0\rangle =  \Sigma^{0, -k} Y / \defopara^{r - k}$ is precisely $\alpha_{vh}$ by Remark \ref{1dx1d-alpha}. Thus, this result follows from the homotopy $\beta^3 \cong \alpha_{vh}$, which is constructed in Theorem \ref{coherent-Mahowald-trick} with $(x, y) = (0, -3)$. \parr 

    The second part follows directly from the first part together with Lemma \ref{pi*-pullback-to-pullback-pi*}. For instance, say there are $[a] \in \pi_{**}(X / \defopara^r)$ and $[b] \in \pi_{**}(Y / \defopara^{r - k})$ with $f[a] = \defopara^k [b]$. Surjectivity in Lemma \ref{pi*-pullback-to-pullback-pi*} yields a certain $[w] \in \pi_{**}(X / \defopara^r  \times_{Y / \defopara^r}  \Sigma^{0, -k} Y/\defopara^{r - k})$ that projects to both $[a]$ and $[b]$, and the projection $[c]$ of $\beta[w] \in \pi_{**}(\Sigma^{-1, 0} Z/\defopara^{k} \times_{X/\defopara^{k}} X / \defopara^r)$ to $\pi_{**}(Z/\defopara^{k})$ then satisfies $\delta_k[c] = g[b]$, $h[c] = \rho[a]$. 
\end{proof}

\begin{remark}
    \begin{itemize}
        \item We can treat the bigraded stems of $V_{f, r, k} \coloneq X / \defopara^r  \times_{f, Y / \defopara^r, \defopara^k}  \Sigma^{0, -k} Y/\defopara^{r - k}$ as \textit{coherent extensions} on the $E_{r + 1}$-page with filtration jump $\geq k$ along $f\colon X \to Y$. The first part of Theorem \ref{blueprint-GMT} identifies $V_{f, r, k}$ with $\Sigma^{0, -k} Y/\defopara^{r - k} \times_{g, \Sigma^{0, -k} Z/\defopara^{r - k}, \delta_k} \Sigma^{-1, 0} Z/\defopara^{k}$, of which the bigraded stems are \textit{coherent (total) differentials} in $Z = \cofib(f)$ with length $\in [k + 1, r]$ whose targets come from $Y$. Therefore, coherent extensions contain the exact same information as coherent differentials.
        \item On the other hand, \emph{actual} finite page extensions of bigraded stems (e.g. $[a] \in \pi_{**}(X / \defopara^r)$, $[b] \in \pi_{**} (Y / \defopara^{r - k})$ with $f[a] = \defopara^k [b]$) yield total differentials in $Z$ by the second part of Theorem \ref{blueprint-GMT}, and vice versa, but the correspondence is not one-to-one. This is because the kernel in Lemma \ref{pi*-pullback-to-pullback-pi*} is usually nonzero, namely there are multiple ways to enhance an extension (or a total differential) to a coherent one, resulting in the indeterminacy along the correspondence. 
        \item In practice, the input of GMT will come from SS. This suffers from another source of indeterminacy, namely nonuniqueness of the lift from SS to bigraded stems. In some cases, we have ``no-crossing conditions'' at hand to eliminate this indeterminacy, resulting in Theorem \ref{GMT-for-SS} and Theorem \ref{GMT-for-SS-2} below. For delicate situations where ``crossings'' are inevitable, we need to carry out further analysis with a certain explicit choice of lift via Theorem \ref{blueprint-GMT}.
    \end{itemize}
\end{remark}

We also record a variant of Theorem \ref{blueprint-GMT} for the ``limit case'' triangle in Construction \ref{limit-distinguished-triangle}.

\begin{theorem}[Blueprint for the ``abutment type'' GMT]\label{blueprint-GMT-variant}
    Suppose $X \xrightarrow{f} Y \xrightarrow{g} Z \xrightarrow{h} \Sigma^{1, 0} X$ is a distinguished triangle in $\Fil\Sp$. Then there are natural isomorphisms of pullbacks 
    \[X[\defopara^{-1}] \times_{f, Y[\defopara^{-1}], \iota}  Y \xrightarrow[\sim]{\beta} 
    Z^{\defopara\textup{-tors}} \times_{h, \Sigma^{1, 0} X^{\defopara\textup{-tors}}, \rho^{\infty}_{\infty}} X[\defopara^{-1}] \xrightarrow[\sim]{\beta} 
    Y \times_{g, Z, \delta_{\infty}^{\infty}} Z^{\defopara\textup{-tors}}\]
    that are compatible with projections to $X[\defopara^{-1}], Y$ and $Z^{\defopara\textup{-tors}}$. 
\end{theorem}

\begin{proof}
    This follows from the proof of Theorem \ref{blueprint-GMT} upon replacing the second input DT by
    \[\oneb \xrightarrow{\iota \, = \, \langle \defopara^{\infty} \rangle} \oneb [\defopara^{-1}]\xrightarrow{\rho_{\infty}^{\infty}} \Sigma^{1, 0} \oneb ^{\defopara\textup{-tors}} = \colim_{k \to \infty} \Sigma^{0, k} \oneb / \defopara^{k} \xrightarrow{\delta_{\infty}^{\infty}} \Sigma^{1, 0} \oneb\]
    coming from Construction \ref{limit-distinguished-triangle}.
\end{proof}

As in the GLR case, The Blueprint Theorem \ref{blueprint-GMT} for GMT also organizes the information in two directions. On the one hand, given  $f\colon X \to Y$ in $\Fil\Sp$, if there is a differential $d(\overline{x}) = \overline{y}$ in its cofiber $Z$ together with two extensions $\overline{x} \rightsquigarrow x$ and $y \rightsquigarrow \overline{y}$, then it suggests an $f$-extension $x \rightsquigarrow y$.

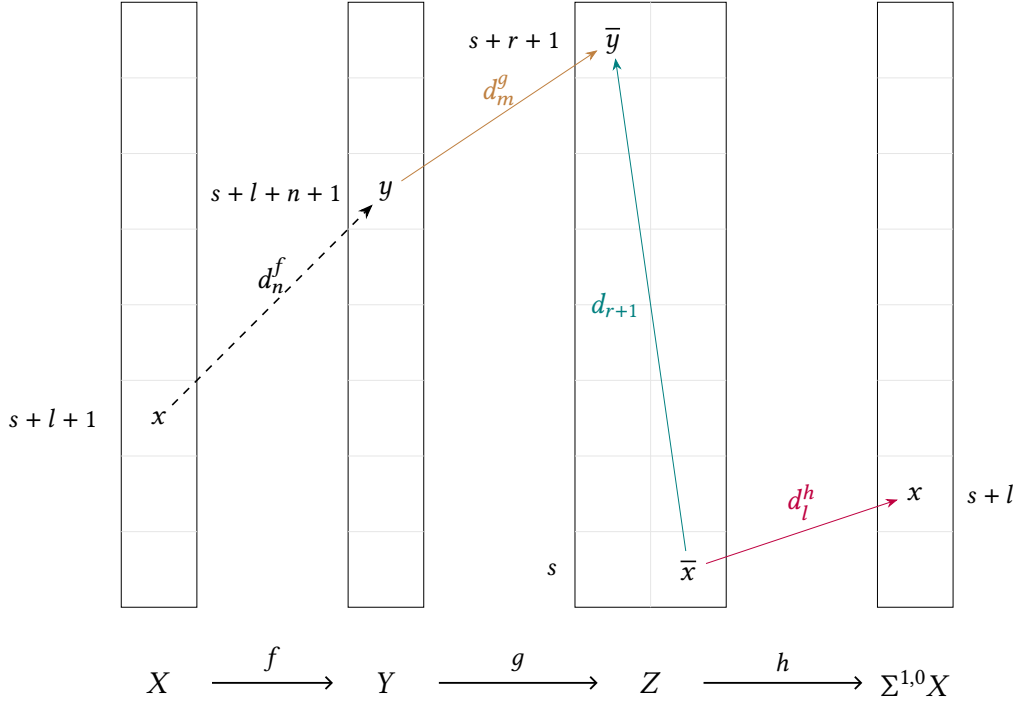
\begin{figure}[htbp]
    \centering
    \scalebox{1}{%
        \begin{tikzpicture}[line width=0.1pt]
            \tikzset{
                diff/.style={-{Stealth},line width=0.3pt,shorten <=7pt,shorten >=7pt},
                diff1/.style={-{Stealth},dashed,line width=0.6pt,shorten <=7pt,shorten >=7pt},
                extarr/.style={-{Stealth},line width=0.3pt,shorten <=7pt,shorten >=7pt},
                extarr1/.style={-{Stealth},dashed,line width=0.6pt,shorten <=7pt,shorten >=7pt},
                mapbelow/.style={-{Straight Barb[scale = .8]},line width=0.7pt,shorten <=.7cm,shorten >=.7cm},
            }
            
            % left copy: 
            
            \begin{scope}
                \draw (49.5,11.5) rectangle (50.5,19.5);
                \node at (48.6,14) {\footnotesize $s + l + 1$};
                % \node at (48.6,16) {$s + r + 1$};
                \begin{scope}
                    % \clip (49.5,11.5) rectangle (50.5,19.5);
                    \draw[black!10] (49.5,12.5) -- (50.5,12.5);
                    \draw[black!10] (49.5,13.5) -- (50.5,13.5);
                    \draw[black!10] (49.5,14.5) -- (50.5,14.5);
                    \draw[black!10] (49.5,15.5) -- (50.5,15.5);
                    \draw[black!10] (49.5,16.5) -- (50.5,16.5);
                    \draw[black!10] (49.5,17.5) -- (50.5,17.5);
                    \draw[black!10] (49.5,18.5) -- (50.5,18.5);

                    \coordinate (X-class) at (50, 14);
                    
                    \node at (X-class) {\small $x$};
                \end{scope}
            \end{scope}
            
            % middle-left copy: 
            
            \begin{scope}[xshift=3cm, yshift=0cm]
                \draw (49.5,11.5) rectangle (50.5,19.5);
                \node at (48.55,17) {\footnotesize $s + l + n + 1$};
                % \node at (48.6,16) {$s + r + 1$};
                \begin{scope}
                    % \clip (49.5,11.5) rectangle (50.5,19.5);
                    \draw[black!10] (49.5,12.5) -- (50.5,12.5);
                    \draw[black!10] (49.5,13.5) -- (50.5,13.5);
                    \draw[black!10] (49.5,14.5) -- (50.5,14.5);
                    \draw[black!10] (49.5,15.5) -- (50.5,15.5);
                    \draw[black!10] (49.5,16.5) -- (50.5,16.5);
                    \draw[black!10] (49.5,17.5) -- (50.5,17.5);
                    \draw[black!10] (49.5,18.5) -- (50.5,18.5);

                    \coordinate (Y-class) at (50, 17);
                    
                    \node at (Y-class) {\small $y$};
                \end{scope}
            \end{scope}

            % middle-right copy: 
            
            \begin{scope}[xshift=6cm, yshift=0cm]
                \draw (49.5,11.5) rectangle (51.5,19.5);
                \node at (49.2,12) {\footnotesize $s$};
                \node at (48.7,19) {\footnotesize $s + r + 1$};
                \begin{scope}
                    % \clip (49.5,11.5) rectangle (51.5,19.5);
                    \draw[black!10] (50.5,11.5) -- (50.5,19.5);
                    \draw[black!10] (49.5,12.5) -- (51.5,12.5);
                    \draw[black!10] (49.5,13.5) -- (51.5,13.5);
                    \draw[black!10] (49.5,14.5) -- (51.5,14.5);
                    \draw[black!10] (49.5,15.5) -- (51.5,15.5);
                    \draw[black!10] (49.5,16.5) -- (51.5,16.5);
                    \draw[black!10] (49.5,17.5) -- (51.5,17.5);
                    \draw[black!10] (49.5,18.5) -- (51.5,18.5);

                    \coordinate (Z-source) at (51, 12);
                    \coordinate (Z-target) at (50, 19);
                    
                    \node at (Z-source) {\small $\overline{x}$};
                    \node at (Z-target) {\small $\overline{y}$};
                    \color{teal} \draw[diff] (Z-source) -- (Z-target) node[midway, left] {\small $d_{r + 1}$};
                \end{scope}
            \end{scope}

            % right copy: 
            
            \begin{scope}[xshift=10cm, yshift=0cm]
                \draw (49.5,11.5) rectangle (50.5,19.5);
                \node at (51,13) {\footnotesize $s + l$};
                % \node at (48.6,16) {$s + r + 1$};
                \begin{scope}
                    % \clip (49.5,11.5) rectangle (50.5,19.5);
                    \draw[black!10] (49.5,12.5) -- (50.5,12.5);
                    \draw[black!10] (49.5,13.5) -- (50.5,13.5);
                    \draw[black!10] (49.5,14.5) -- (50.5,14.5);
                    \draw[black!10] (49.5,15.5) -- (50.5,15.5);
                    \draw[black!10] (49.5,16.5) -- (50.5,16.5);
                    \draw[black!10] (49.5,17.5) -- (50.5,17.5);
                    \draw[black!10] (49.5,18.5) -- (50.5,18.5);

                    \coordinate (X-class-2) at (50, 13);
                    
                    \node at (X-class-2) {\small $x$};
                \end{scope}
            \end{scope}

            % extension arrows

            \draw[extarr1] (X-class) -- (Y-class) node[midway, above] {\small $d_{n}^f$}; 
            \color{brown} \draw[extarr] (Y-class) -- (Z-target) node[midway, above] {\small $d_{m}^g$};
            \color{purple} \draw[extarr] (Z-source) -- (X-class-2) node[midway, above] {\small $d_{l}^h$};

            % underlying maps

            \coordinate (label-X) at (50, 10.5); 
            \coordinate (label-Y) at (53, 10.5);
            \coordinate (label-Z) at (56.5, 10.5);
            \coordinate (label-SX) at (60, 10.5);
            
            \color{black}
            \node at (label-X) {$X$};
            \node at (label-Y) {$Y$};
            \node at (label-Z) {$Z$};
            \node at (label-SX) {$\Sigma^{1, 0} X$};
            \draw[mapbelow] (label-X) -- (label-Y) node[midway,above] {\footnotesize $f$};
            \draw[mapbelow] (label-Y) -- (label-Z) node[midway,above] {\footnotesize $g$};
            \draw[mapbelow] (label-Z) -- (label-SX) node[midway,above] {\footnotesize $h$};
            
        \end{tikzpicture}
    }
    \caption{Generalized Mahowald trick.}
    \label{figure-GMT}
\end{figure}

\begin{theorem}[Generalized Mahowald trick, part 1] \label{GMT-for-SS}
    Suppose $X \xrightarrow{f} Y \xrightarrow{g} Z \xrightarrow{h} \Sigma^{1, 0} X$ is a distinguished triangle in $\Fil\Sp$. Take $n, l, m, i, j \in \Nb$ with $i < n$ and write $r = n + m + l$. Suppose there are classes $x \in Z^{s + l + 1, t + l}_{n}(X)$, $y \in Z^{s + n + l + 1, t + n + l}_{m + 1}(Y)$, $\overline{x} \in Z^{s, t}_{r}(Z)$, $\overline{y} \in Z^{s + r + 1, t + r}_{\infty}(Z)$ subject to the following assumptions:
    \raggedcolumns
    \begin{multicols}{2}
    \begin{enumerate}
        \item ${\color{teal} d_{r + 1}(\overline{x}) = \overline{y}}$. 
        \item ${\color{purple} d^{h, E_{r' + 1}, j}_{l}(\overline{x}) = x}$, where $r' = n + l - i$.
        \item ${\color{brown} d^{g, E_{m + 2}, i}_{m}(y) = \overline{y}}$.
        \item The {\color{teal} differential} in item 1 and the {\color{purple} extension} in item 2 are complementary on the source.
    \end{enumerate}
    \end{multicols}
    \noindent Then $x$ lies in $Z^{s + l + 1, t + l}_{n + m + 1}(X)$, 
    and there is a $(j + l)$-th extension $d_{n}^{f, E_{n + m + 2}, j + l}(x) = y$ (see Figure \ref{figure-GMT}). 
\end{theorem}

\begin{proof}
    Items 1, 2, and 4 provide two bigraded stems $[\overline{x}] \in \pi_{t - s, t}(Z / \defopara^{r'})$ and $[x] \in \pi_{t - s - 1, t + l} (X / \defopara^{n - i})$ lifting $\overline{x}$ and $x$, so that 
    \[ \defopara^j h[\overline{x}] = \defopara^{j + l} [x] \  \text{ in } \  \pi_{t - s - 1, t - j}(X / \defopara^{r' + j}) \quad \text{ and } \quad\delta_{r'}^{m + i + 1} [\overline{x}] = \defopara^{m + i} \overline{y} \ \text{ in } \ \pi_{t - s - 1, t + r'}(Z / \defopara^{m + i + 1}), \]
    cf. Lemma \ref{combining-partial-nocrossings-in-source-refined}. On the other hand, by item 3 we can find $[y] \in \pi_{t - s - 1, t + r'}(Y / \defopara^{m + 1})$ lifting $y$, so that $\defopara^{i} g[y] = \defopara^{m + i} \overline{y}$ in $\pi_{t - s - 1, t + r'}(Z / \defopara^{m + i + 1})$. Apply Theorem \ref{blueprint-GMT} with $1 \leq r' < r + 1$ to the stems
    \[g(\defopara^i [y]) = \defopara^{m + i} \overline{y} = \delta_{r'}^{m + i + 1}[\overline{x}],\]
    we obtain $[a] \in \pi_{t - s - 1, t}(X / \defopara^{r + 1})$ 
    such that $\rho^{r + 1}_{r'}[a] = h[\overline{x}]$ and $f[a] = \defopara^{r'}\!(\defopara^i [y]) = \defopara^{n + l}[y]$. Feeding
    \[\rho^{r + 1 + j}_{r'+ j}(\defopara^{j}[a]) = \defopara^{j}\rho^{r + 1}_{r'}[a] = \defopara^{j} h[\overline{x}] = \defopara^{l + j}[x]\] 
    into the pullback in Lemma \ref{defopara-rho-cartesian} (via Lemma \ref{pi*-pullback-to-pullback-pi*}) produces a refined lift $[x] \in \pi_{t - s - 1, t + l}(X / \defopara^{n + m + 1})$ of $x$ with $\defopara^{l + j} [x] = \defopara^j [a]$. Thus, we have $\defopara^{j + l} f[x] = \defopara^j f[a] = \defopara^{j + l + n} [y]$ which witnesses the desired extension. 
\end{proof}

\begin{remark} \label{explanation-of-the-parameters-GMT-1}
    We briefly discuss how to choose parameters when we apply Theorem \ref{GMT-for-SS}.
    \begin{itemize}
        \item The parameters $n, m, l$ are fixed by the desideratum. 
        \item The parameter $j$ can be as large as possible to guarantee the extension from $\overline{x}$ to $x$ is unobstructed by early boundaries (i.e. $\defopara$-power torsion lifts). In practice, often $j = 0$ suffices. 
        \item The parameter $i$ is fine-tuned towards two goals: the existence of the extension from $\overline{x}$ to $x$, and the ``complementary on the source'' condition in item 4. The largest possible choice of $i$ yields $r' + 1 = n - i  + l + 1 = l + 2$, the earliest page on which extensions of filtration jump $l$ make sense. It can be as small as $0$, and we usually take $i = 0$ in practice.
    \end{itemize}
\end{remark}

\begin{remark} \label{reduction-to-GMT-in-LWX}
    The conclusion of Theorem \ref{GMT-for-SS} can be refined using Theorem \ref{stretching-extensions-across-pages} or Lemma \ref{reducing-refined-extensions-to-extensions}.
    \begin{itemize}
        \item In the setup of Theorem \ref{GMT-for-SS}, suppose furthermore $f(x) = 0$ in $E_2^{s + l + 1, t + l}(Y)$, and there exists $0 \leq v \leq j + l$ so that $B_{1 + j + l + a}^{s + a + l + 1, t + a + l}(Y) / B_{1 + v + a}^{s + a + l + 1, t + a + l}(Y) = 0$ for $1 \leq a \leq n$. Then $d_{n}^{f, E_{n + m + 2}, v}(x) = y$ according to Theorem \ref{stretching-extensions-across-pages}. 
        \item Alternatively, apart from $f(x) = 0$ in $E_2^{s + l + 1, t + l}(Y)$, suppose we only have certain $0 \leq v \leq j + l$ with $B_{1 + j + l + a}^{s + a + l + 1, t + a + l}(Y) / B_{1 + v + a}^{s + a + l + 1, t + a + l}(Y) = 0$ for $1 \leq a \leq n - 1$. Then $d_{n}^{f, E_{n + m + 2}, v}(x) = y$ has a witness $\defopara^v [x] = \defopara^{v + n} [y]$, and the difference $[x] - \defopara^n [y]$ is in $\IIm(\defopara^n) \cap \ker(\defopara^v)$ according to Lemma \ref{reducing-refined-extensions-to-extensions}. In other words, the conclusion becomes $d^{f, E_{n + m + 2}, v}_{n}(x) = y'$ for some $y' \in Z^{s + n + l + 1, t + n + l}_{m + 1}(Y)$ such that $y - y' \in B^{s + n + l + 1, t + n + l}_{1 + v + n}(Y)$. 
    \end{itemize}
    On the other hand, the complementary crossing condition in assumption item 4 there can be simplified via Example \ref{complementary-crossing-examples}. For instance, it holds true if the differential in item 1 has no crossing on the $E_{{r'} + 1}$-page or the extension in item 2 has no crossing. Consequently, under the translation in Remark \ref{synthetic-spectra-as-special-case},  Remark \ref{comparison-with-extension-SS} and Remark \ref{comparison-with-extension-SS-contd}, we recover \cite[Theorem 6.12]{Lin-Wang-Xu-kervaire} for Adams SS. 
\end{remark}

There is also a dual perspective extracted from the Blueprint Theorem \ref{blueprint-GMT}: Given a map $f\colon X \to Y$ in $\Fil\Sp$, if there is an $f$-extension $x \rightsquigarrow y$ together with two other extensions $\overline{x} \rightsquigarrow x$ and $y \rightsquigarrow \overline{y}$ out of/into the cofiber $Z$, then it suggests a differential $d(\overline{x}) = \overline{y}$ in $Z$. 

\begin{theorem}[Generalized Mahowald trick, part 2] \label{GMT-for-SS-2}
    Suppose $X \xrightarrow{f} Y \xrightarrow{g} Z \xrightarrow{h} \Sigma^{1, 0} X$ is a distinguished triangle in $\Fil\Sp$. Take $n, l, m, j \in \Nb$ with $n \geq 1$, and write $r = n + m + l$. Suppose there are classes $x \in Z^{s + l + 1, t + l}_{n + m + 1}(X)$, $y \in Z^{s + n + l + 1, t + n + l}_{m + 1}(Y)$, $\overline{x} \in Z^{s, t}_{l + 1}(Z)$, $\overline{y} \in E^{s + r + 1, t + r}_{2}(Z)$ subject to the following assumptions:
    % \raggedcolumns
    \begin{multicols}{2}
    \begin{enumerate}
        \item $d_{n}^{f, E_{n + m + 2}, j + l}(x) = y$. 
        \item $d^{h, E_{l + 2}, j}_{l}(\overline{x}) = x$. 
        \item $d^{g, E_{m + 2}, n + l - 1}_{m}(y) = \overline{y}$.
        \item Item 3 has no crossing. 
        \item {$B_{a + l + j + 2}^{s + a + l + 2, t + a + l + 1}(Y) / B_{a + l + 1}^{s + a + l + 2, t + a + l + 1}(Y) = 0$} in the standard SS of $Y$ for each $0 \leq a < n$.
    \end{enumerate}
    \end{multicols}
    \noindent Then $\overline{x}, \overline{y}$ are actually $r$-cycles in the standard SS of $Z$, and $d_{r + 1}(\overline{x}) = \overline{y}$. 
\end{theorem}

\begin{proof}
    Item 1 supplies $[x] \in \pi_{t - s - 1, t + l}(X / \defopara^{n + 
    m + 1}), [y] \in \pi_{t - s - 1, t + n + l}(Y / \defopara^{m + 1})$ lifting $x, y$ so that $\defopara^{j + l} f[x] = \defopara^{j + n + l}[y]$, and item 2 yields a lift $[\overline{x}] \in \pi_{t - s, t}(Z / \defopara^{l + 1})$ of $\overline{x}$ so that $\defopara^j h[\overline{x}] = \defopara^{j + l} x$. Apply Theorem \ref{blueprint-GMT} with parameters $1 \leq j + l + 1 < r + j + 1$ to the stems 
    \[h(\defopara^j[\overline{x}]) = \defopara^{j + l} \rho^{n + m + 1}_{1}[x]= \rho^{r + j + 1}_{j + l + 1}(\defopara^{j + l}[x]),\]
    we obtain $[b] \in \pi_{t - s - 1, t + l + 1}(Y / \defopara^{n + m})$ so that $f(\defopara^{j + l}[x]) = \defopara^{j + l + 1}[b], g[b] = \delta_{j + l + 1}^{n + m}(\defopara^j [\overline{x}]) = \delta_{l + 1}^{n + m}[\overline{x}]$. Thus, the class $\defopara^{l}([b] - \defopara^{n - 1}[y]) \in \pi_{t - s - 1, t + 1}(Y / \defopara^{r})$ is $\defopara^{j + 1}$-torsion and $\defopara^{l}$-divisible, and item 5 implies (through Lemma \ref{reducing-refined-extensions-to-extensions}) that all such stems in this bidegree are also $\defopara^{n + l}$-divisible. Therefore, $\defopara^{l} [b] = \defopara^{n + l - 1}[y]'$ where $[y]' = [y] + \defopara \varepsilon$ for some $\varepsilon \in \pi_{t - s - 1, t + n + l + 1}(Y / \defopara^{m})$. Item 3 and item 4 imply that $\defopara^{n + l - 1} g[y]' = \defopara^{r - 1} \overline{y}$ due to Corollary \ref{free-choice-lemma-for-refined-extensions-with-full-no-crossings}, so $\delta_1^{r} \overline{x} = \defopara^l\delta_{l + 1}^{n + m}[\overline{x}] = \defopara^l g[b] = \defopara^{r - 1}\overline{y}$. By Theorem \ref{delta-as-total-diff}, $\overline{x}, \overline{y}$ are actually $r$-cycles in the standard SS of $Z$ and $d_{r + 1}(\overline{x}) = \overline{y}$.
\end{proof}

\begin{remark} \label{explanation-of-the-parameters-GMT-2}
    We briefly discuss how to choose parameters when we apply Theorem \ref{GMT-for-SS-2}.
    \begin{itemize}
        \item The parameters $n, m, l$ are fixed by the desideratum. 
        \item The parameter $j$ can be as large as possible to guarantee the extension from $\overline{x}$ to $x$ is unobstructed by early boundaries (i.e. $\defopara$-power torsion lifts). In practice, often $j = 0$ suffices.
    \end{itemize}
\end{remark}

\begin{example} \label{GMT-example}
    Recall the differential $\color{blue} d_5 (h_5 Pe_0) = d_0 \Delta h_0^2 e_0$ in Example \ref{example-diff-crossing}. This turns out to be the first essential $d_5$ in the $\mathrm{H}\Fb_2$-Adams SS of $\Sb^0$, which is established in \cite[Theorem 3.3.55]{Isaksen-stable-stems}. In this example we provide an alternative argument through GMT. Our notational convention will be the same as in Example \ref{example-extn-crossing}. \parr
    
    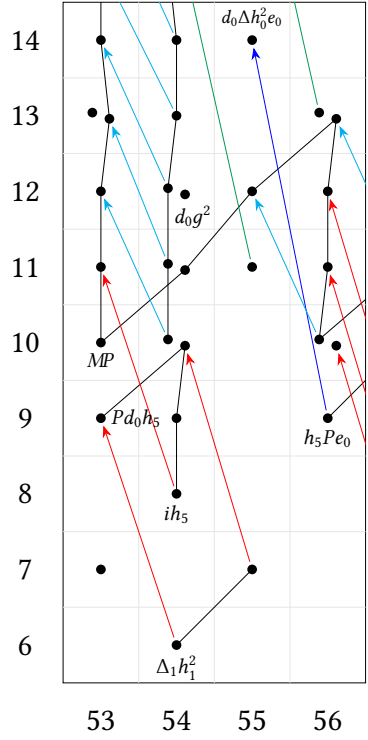
\begin{figure}[htbp]
        \centering
        \scalebox{1}{\begin{tikzpicture}[line width=0.1pt]
        \tikzset{
            diff/.style={-{Stealth},line width=0.3pt,shorten <=3pt,shorten >=3pt},
            diff1/.style={-{Stealth},dashed,line width=0.6pt,shorten <=3pt,shorten >=3pt},
            extarr/.style={-{Stealth},line width=0.3pt,shorten <=3pt,shorten >=3pt},
            extarr1/.style={-{Stealth},dashed,line width=0.6pt,shorten <=3pt,shorten >=3pt},
            mapbelow/.style={-{Stealth},line width=0.4pt,shorten <=.7cm,shorten >=.7cm},
        }
        \draw (52.5,5.5) rectangle (56.5,14.5);
        \node at (53,5) {$53$};
        \node at (54,5) {$54$};
        \node at (55,5) {$55$};
        \node at (56,5) {$56$};
        \node at (52,6) {$6$};
        \node at (52,7) {$7$};
        \node at (52,8) {$8$};
        \node at (52,9) {$9$};
        \node at (52,10) {$10$};
        \node at (52,11) {$11$};
        \node at (52,12) {$12$};
        \node at (52,13) {$13$};
        \node at (52,14) {$14$};
        \begin{scope}
        \clip (52.5,5.5) rectangle (56.5,14.5);
        \draw[black!10] (53.5,5.5) -- (53.5,14.5);
        \draw[black!10] (54.5,5.5) -- (54.5,14.5);
        \draw[black!10] (55.5,5.5) -- (55.5,14.5);
        \draw[black!10] (52.5,6.5) -- (56.5,6.5);
        \draw[black!10] (52.5,7.5) -- (56.5,7.5);
        \draw[black!10] (52.5,8.5) -- (56.5,8.5);
        \draw[black!10] (52.5,9.5) -- (56.5,9.5);
        \draw[black!10] (52.5,10.5) -- (56.5,10.5);
        \draw[black!10] (52.5,11.5) -- (56.5,11.5);
        \draw[black!10] (52.5,12.5) -- (56.5,12.5);
        \draw[black!10] (52.5,13.5) -- (56.5,13.5);
          % \draw (50,6) -- (53,7);
          \draw (54,6) -- (55,7);
          \draw (54,8) -- (54,9);
          % \draw (54,8) -- (56.887,9.041);
          \draw (53,9) -- (54.113,9.959);
          \draw (54,9) -- (54.113,9.959);
          % \draw (54,9) -- (57,10);
          \draw (53,10) -- (53,11);
          \draw (53,10) -- (54.113,10.959);
          % \draw (53,10) -- (56,11);
          \draw (53.887,10.041) -- (53.887,11.041);
          % \draw (53.887,10.041) -- (57,11);
          \draw (53,11) -- (53,12);
          % \draw (53,11) -- (56,12);
          \draw (56,9) -- (57,10);
          \draw (54.113,10.959) -- (55,12);
          \draw (53.887,11.041) -- (53.887,12.041);
          \draw (53,12) -- (53.113,12.959);
          % \draw (53,12) -- (56.113,12.959);
          \draw (55.887,10.041) -- (56,11);
          \draw (55.887,10.041) -- (57,11);
          % \draw (55.887,10.041) -- (59,11);
          \draw (53.887,12.041) -- (54,13);
          \draw (53.113,12.959) -- (53,14);
          \draw (56,11) -- (56,12);
          % \draw (56,11) -- (59,12);
          \draw (55,12) -- (56.113,12.959);
          \draw (54,13) -- (54,14);
          \draw (53,14) -- (53,15);
          \draw (56,12) -- (56.113,12.959);
          \draw (54,14) -- (53.887,15.041);
          \draw[diff, blue] (56,9) -- (55,14);
          \draw[diff, ForestGreen] (55,11) -- (54.113,14.959);
          \draw[diff, ForestGreen] (55.887,13.041) -- (55,17);
          \draw[diff, red] (54,6) -- (53,9);
          \draw[diff, red] (54,8) -- (53,11);
          \draw[diff, red] (55,7) -- (54.113,9.959);
          \draw[diff, red] (57,7) -- (56.113,9.959);
          \draw[diff, red] (56.887,8.041) -- (56,11);
          \draw[diff, red] (56.887,9.041) -- (56,12);
          \draw[diff, cyan] (53.887,10.041) -- (53,12);
          \draw[diff, cyan] (53.887,11.041) -- (53.113,12.959);
          \draw[diff, cyan] (53.887,12.041) -- (53,14);
          \draw[diff, cyan] (54,13) -- (53,15);
          \draw[diff, cyan] (54,14) -- (53.113,15.959);
          \draw[diff, cyan] (55.887,10.041) -- (55,12);
          \draw[diff, cyan] (57,11) -- (56.113,12.959);
          \fill (53,7) circle (0.064);
          \fill (55,7) circle (0.064);
          \fill (54,8) circle (0.064) node[below] {\resizelabel{0.667}{$ih_5$}};
          \fill (53,9) circle (0.064) node[right] {\shortstack{\resizelabel{0.667}{$Pd_0h_5$}}};
          \fill (54,9) circle (0.064);
          \fill (54.113,9.959) circle (0.064);
          \fill (53,11) circle (0.064);
          \fill (56,9) circle (0.064) node[below] {\resizelabel{0.667}{$h_5Pe_0$}};
          \fill (54.113,10.959) circle (0.064);
          \fill (53.887,11.041) circle (0.064);
          \fill (53,12) circle (0.064);
          \fill (56.113,9.959) circle (0.064) node[below] {};
          \fill (55,11) circle (0.064) node[below] {};
          \fill (54.113,11.959) circle (0.064) node[below] {\resizelabel{0.667}{$\ \ \ d_0g^2$}};
          \fill (53.887,12.041) circle (0.064);
          \fill (53.113,12.959) circle (0.064);
          \fill (52.887,13.041) circle (0.064) node[below] {};
          \fill (56,11) circle (0.064);
          \fill (55,12) circle (0.064);
          \fill (54,13) circle (0.064);
          \fill (53,14) circle (0.064);
          \fill (56,12) circle (0.064);
          \fill (54,14) circle (0.064);
          \fill (56.113,12.959) circle (0.064);
          \fill (55.887,13.041) circle (0.064) node[below] {};
          \fill (55,14) circle (0.064) node[above] {\resizelabel{0.667}{$d_0\Delta h_0^2e_0$}};
          \fill (54,6) circle (0.064) node[below] {\resizelabel{0.667}{$\Delta_1h_1^2$}};
          \fill (53,10) circle (0.064) node[below] {\resizelabel{0.667}{$M\!P$}};
          \fill (53.887,10.041) circle (0.064) node[below] {};
          \fill (55.887,10.041) circle (0.064) node[below] {};
        \end{scope}\end{tikzpicture}}
        \caption{Deducing $\color{blue} d_5(h_5 P e_0) = d_0 \Delta h_0^2 e_0$ through generalized Mahowald trick.}
        \label{figure-GMT-example}
    \end{figure}

    We choose our inputs as follows: 
    \begin{itemize}
        \item The distinguished triangle (from \cite[Proposition 3.20]{Lin-Wang-Xu-kervaire})
        \[\Sb^{1, 2}_{\mathrm{H}\Fb_2} \xrightarrow{[h_1]} \Sb^{0, 0}_{\mathrm{H}\Fb_2} \xrightarrow{i} \cofib(\eta)_{\mathrm{H}\Fb_2} \xrightarrow{p} \Sb^{2, 2}_{\mathrm{H}\Fb_2}\]
        together with the differentials ${\color{cyan} d_2(h_1 \Delta_1 h_1^2[2]) = h_5 Pe_0[0]}$, ${\color{cyan} d_2(d_0 g^2[2]) = d_0 \Delta h_0^2 e_0[0]}$ in the $\mathrm{H}\Fb_2$-Adams SS of $\cofib(\eta)$.
        \item Burklund's \cite{Burklund-hidden-extension} hidden extension 
        \[[h_0]\cdot[h_0 i h_5] = \delta_2^{4}[h_1 \Delta_1 h_1^2] + \defopara^2 [d_0 g^2]\] 
        in $\pi_{54, 54 + 10}(\Sb_{\mathrm{H}\Fb_2} / \defopara^4)$.  Here $[h_1 \Delta_1 h_1^2]$ and $[d_0 g^2]$ stand for certain specific lifts of $h_1 \Delta_1 h_1^2$ and $d_0 g^2$ in $\pi_{**}(\Sb_{\mathrm{H}\Fb_2} / \defopara^2)$, while $[h_0]$ and $[h_0 i h_5]$ stand for certain specific lifts in $\pi_{**}(\Sb_{\mathrm{H}\Fb_2} / \defopara^4)$.
    \end{itemize}

    According to GMT (Theorem \ref{GMT-for-SS}), the two differentials in $\cofib(\eta)_{\mathrm{H}\Fb_2}$ lead to two extensions 
    \[d_1^{[h_1], E_3}(h_1 \Delta_1 h_1^2) = h_5 Pe_0 \qquad \text{and} \qquad d_1^{[h_1], E_3}(d_0 g^2) = d_0 \Delta h_0^2 e_0.\]
    Here we use the parameters $n = 1, m = l = i = j = 0$. The two extensions both have no crossing as all classes in $E_{2, \mathrm{ASS}}^{8, 8 + 55}(\Sb^0)$ and $E_{2, \mathrm{ASS}}^{13, 13 + 54}(\Sb^0)$ are $h_1$-torsion. Thus, we have $[h_1] \cdot [h_1 \Delta_1 h_1^2] = \defopara h_5 P e_0$ and $[h_1] \cdot [d_0 g^2] = \defopara d_0 \Delta h_0^2 e_0$ in $\pi_{**}(\Sb_{\mathrm{H}\Fb_2} / \defopara^2)$ due to Corollary \ref{free-choice-lemma-for-refined-extensions-with-full-no-crossings}. As $[h_1] \cdot [h_0] = 0$, it follows that 
    \[0 = [h_1] \delta_2^{4}[h_1 \Delta_1 h_1^2] + \defopara^2 [h_1][d_0 g^2] = \delta_2^4(\defopara h_5 Pe_0) + \defopara^3 d_0 \Delta h_0^2 e_0 = \delta_1^4(h_5 Pe_0) + \defopara^3 d_0 \Delta h_0^2 e_0\]
    in $\pi_{55, 55 + 11}(\Sb_{\mathrm{H}\Fb_2} / \defopara^4)$, in other words $\color{blue} d_5(h_5 P e_0) = d_0 \Delta h_0^2 e_0$ in the Adams SS of $\Sb^0$.

\FloatBarrier
\end{example}

\begin{remark} \label{credits-for-theorem-C}
We end this section with some further historical comments and attributions concerning the generalized Mahowald trick. The relationship with Lin--Wang--Xu's Adams SS version \cite[Theorem 6.12]{Lin-Wang-Xu-kervaire} has already been discussed above.

\begin{itemize}

\item The idea behind the Mahowald trick goes back to Mahowald \cite{MT67}.

In the terminology of the present paper, the rough intuition behind this primitive form of the Mahowald trick can be described as follows. Given a map $f\colon X \to Y$ between filtered spectra and $Z = \cofib(f)$, there are two long exact sequences:
\begin{itemize}
\item applying $\pi_{**}(- / \defopara)$ gives a long exact sequence on the $E_2$-pages of the standard spectral sequences of $X$, $Y$, and $Z$;
\item applying $\pi_{*}(-[\defopara^{-1}])$ gives a long exact sequence on the abutments of these standard spectral sequences.
\end{itemize}
In practice, there can be mismatches between the maps induced by $f$ on the $E_2$-page and on the abutment. Concretely, suppose that $x \in E_2(X)$ detects $\alpha$, and $y \in E_2(Y)$ detects $\beta$, with a relation $f(\alpha) = \beta$ in the abutment which is not visible on the $E_2$-page, for instance with $x \in \ker(f)$ and $y \in \coker(f)$ on $E_2$. Since the kernel and cokernel of $f$ on the $E_2$-page, respectively on the abutment, are controlled by the $E_2$-page, respectively the abutment, of $Z$, one is led to expect that such mismatches should be accounted for by the information governing the passage from $E_2(Z)$ to its abutment, namely the differentials in $E_*(Z)$.

This idea serves as a heuristic tool in many of Mahowald's computations. For example, in \cite[Theorem 4.2.1]{MT67}, Mahowald and Tangora deduce $d_4(e_0 g) = d_0 P d_0$ in the Adams SS of $\Sb^0$ through an argument guided by this philosophy, using what amounts to the triangle
\[
\Sb^{1, 2}_{\mathrm{H}\Fb_2} \xrightarrow{[h_1]} \Sb^{0, 0}_{\mathrm{H}\Fb_2} \xrightarrow{i} \cofib(\eta)_{\mathrm{H}\Fb_2} \xrightarrow{p} \Sb^{2, 2}_{\mathrm{H}\Fb_2}
\]
to relate this differential to the $p$-extension on the abutment $e_0 g[0] \rightsquigarrow h_0^2 m$.

\item Hill--Hopkins--Ravenel prove a version of the generalized Mahowald trick, more precisely a version of the Blueprint Theorem \ref{blueprint-GMT}, in \cite[Theorem 4.5]{HHR-BPC4<1>}. They use it to produce the stem-level extension input for their exotic transfer differentials in the $C_4$-slice SS of $\BPCfour\langle 1 \rangle$.

The technical core of this result is \cite[Lemma 4.4]{HHR-BPC4<1>}, which plays a role analogous to our coherent Mahowald trick, Theorem \ref{coherent-Mahowald-trick}. As pointed out by Meier--Shi--Zeng \cite[Footnote 3]{MSZ-localized-slice-SS}, this lemma requires some additional care as stated. Rephrased in the language of the present paper, the issue is that the input in \cite[Lemma 4.4]{HHR-BPC4<1>} is not a $2$d cuDT in our sense, but only a $\overrightarrow{\Zb} \times \overrightarrow{\Zb}$-diagram in which every three consecutive arrows, either vertically or horizontally, form a DT or an ADT. Since such a diagram does not keep track of the nullhomotopy data, it need not be uniquely determined by its restriction to $\{0 < 1\}^2$.

For instance, let $X$ be a $2$d cuDT with $X\langle 1,1 \rangle = X\langle 3,2 \rangle = X\langle 2,3 \rangle = 0$. After restricting $X$ to $\overrightarrow{\Zb} \times \overrightarrow{\Zb}$, one may choose an automorphism $f$ of $X\langle 2,2 \rangle$ and replace the map $X\langle 2,1 \rangle \to X\langle 2,2 \rangle$ by its postcomposition with $f$. The resulting diagram $X'$ still satisfies the assumptions of \cite[Lemma 4.4]{HHR-BPC4<1>} and has the same restriction to $\{0 < 1\}^2$. However, Theorem \ref{coherent-Mahowald-trick} applied to $X$ shows that the claimed conclusion for $X'$ can only hold after inserting the automorphism $f$ in the appropriate place. Thus the statement does not hold in this level of generality without additional coherence data.

\item The geometric boundary theorem forms a related line of development. In its earlier forms, including Johnson--Miller--Wilson--Zahler \cite[Theorem 1.7]{JMWZ}, Bruner \cite{Bru77}, and Ravenel \cite[Theorem 2.3.4]{GreenBook}, it compares algebraic boundary maps on the $E_2$-pages with geometric maps on abutments. Behrens \cite[Lemma A.4.1]{Behrens-EHP} formulates a version involving filtration jumps, bringing it closer to the intuition described above. This filtration-jump version is refined and corrected by Ma \cite{Ma2}. For a more extensive history, see \cite[Introduction]{Ma2}.

In its most general form in \cite{Ma2}, the geometric boundary theorem describes possible behaviors of the image of a differential $d(x) = x_{\infty}$ in the standard SS of $X$ along a map $f\colon X \to Y$ in $\Fil\Sp$, or more precisely the $f$-extension of the target $x_{\infty}$, in terms of the differentials in $X$, $Y$, and $\cofib(f)$, together with the maps between their $E_2$-pages. In this formulation, bigraded stem lifts are built into the statement, while no-crossing hypotheses are not imposed. 

The technical key of \cite{Ma2} is \cite[Lemma 5.1]{Ma1}, which is essentially equivalent to our coherent Mahowald trick, Theorem \ref{coherent-Mahowald-trick}. The proof there is quite different from ours in Appendix \ref{app:B}: it works directly in the quasicategorical model and verifies the required compatibilities through explicit choices of horn fillers, whereas our proof is model-independent and reduces the coherence check to a single family of concrete examples (or even to a single representative of the family). 

As no-crossing conditions are out of the story, the geometric boundary theorems in \cite{Ma2} are not directly comparable with Theorem \ref{GMT-for-SS} or Theorem \ref{GMT-for-SS-2}. Rather, they can be viewed as arising from iterations of the Blueprint Theorem \ref{blueprint-GMT}, together with the following two facts: for any map $f\colon X \to Y$ in $\Fil\Sp$,
\begin{enumerate}
    \item[(a).] The map $f$ is compatible with total differentials. In other words, $f(\delta_r^{\infty}[x]) = \delta_r^{\infty} f[x]$ in $\pi_{**}(Y)$ for every $[x] \in \pi_{**}(X / \defopara^r)$ and every $r \in \Nb_{\geq 1}$.
    \item[(b).] The map $f$ is compatible with detection phenomena. In other words, $f(\iota[x]) = \iota f[x]$ in $\pi_{**}(Y[\defopara^{-1}])$ for every $[x] \in \pi_{**}(X)$, where $\iota$ denotes the $\defopara$-inversion map.
\end{enumerate}
More precisely, \cite[Theorem 1]{Ma2} is an instance of (a). Theorem 2 of \cite{Ma2} reduces to (a) once the extension $p_* y \rightsquigarrow x$ is obtained using Theorem \ref{blueprint-GMT}; the same argument shows that \cite[Theorem 5]{Ma2} reduces to (b), that \cite[Theorem 3]{Ma2} reduces to \cite[Theorem 2]{Ma2}, and that \cite[Theorem 4]{Ma2} reduces to \cite[Theorem 5]{Ma2}. 

Conversely, starting from a situation of type (a) or type (b), one may recursively replace an extension by an extension-differential-extension zigzag using Theorem \ref{blueprint-GMT}, which suggests an infinite family of possible geometric boundary patterns.

\end{itemize}

\end{remark}

%% file: Slice.tex
\label{sec:5}

In the previous sections, we have developed the generalized Leibniz rule, the generalized Mahowald trick, and the Leibniz rule for total differentials in the realm of $\Fil\Sp$ using BIPWX's cofiber-of-$\tau$ formalism. We summarize these computational tools as the Burklund--Lin--Wang--Xu (BLWX) methods.  In this section we apply BLWX methods to compute the equivariant slice spectral sequences of Hill--Hopkins--Ravenel theories. In \S~\hyperref[subsec:5.1]{5.1} we set up our convention for general $G$-equivariant filtrations and ($\RO(G)$-graded, Mackey-functor-valued) spectral sequences. In \S~\hyperref[subsec:5.2]{5.2} we introduce the equivariant slice SS and the HHR theories $\BPG \langle m \rangle$. In \S~\hyperref[subsec:5.3]{5.3} and \S~\hyperref[subsec:5.4]{5.4}, we record known information for the slice SS of $\BP_{\Rb}$ and $\BPCfour$. In \S~\hyperref[subsec:5.5]{5.5}, we use GLR and GMT to deduce exotic transfer differentials for the $C_4$-slice SS of $\BPCfour \langle m \rangle$ in Theorem \ref{exotic-transfer-paradigm} and the subsequent examples. Finally, in  \S~\hyperref[subsec:5.6]{5.6} we end the paper with a satisfactory summary of generating differentials in the $C_4$-slice SS of $\BPCfour \langle 1 \rangle$.

\subsection{Equivariant filtrations and equivariant SS}
\label{subsec:5.1}

Let $G$ be a finite group. In Appendix \ref{app:C} we discuss Picard-graded SS in general, which in particular leads to a lax symmetric monoidal functor $E_*^{*,\filledstar}$ sending a filtered genuine $G$-spectrum $X \in \Fil\Sp^G$ to its $\RO(G)$-graded Mackey-functor-valued standard SS $E_*^{*, \filledstar}(X)$, cf. Example \ref{bigraded-picard-SS-examples}. In this section, for better compatibility with convergence results in Theorem \ref{Bockstein-dictionary-infinite}, we phrase this construction in a slightly different (but equivalent) manner, which does not leave the realm of abelian groups.

\begin{construction}[parametrized SS] \label{parametrized-SS}
    Suppose $\CC \in \Mod_{\Fil\Sp}(\PPr_{\st})$ and $\CT \subset \CC$ is a small full subcategory. We write $\psi_{\CT}^{\Fil}\colon \CC \to \Fun(\CT^{\op}, \Fil\Sp)$ for the functor $X \mapsto (A \in \CT \mapsto \psi_{A}^{\Fil}(X) = \Msp^{\Fil}(A, X))$. Moreover, we denote the composite 
    \[ \CC \xrightarrow{\psi^{\Fil}} \Fun(\CT^{\op}, \Fil\Sp) \xrightarrow{E_*^{*,*}} \Fun(\CT^{\op}, \SpSeq)\]
    by $X \mapsto \{E_r^{s, t, A}(X)\}_{r \geq 2, s, t \in \Zb, A \in \CT}$, where $E_r^{s, t, A}(X) = E_r^{s, t}(\Msp^{\Fil}(A, X))$. By construction, $\psi^{\Fil}$ controls both this parametrized SS and hidden extensions along the maps induced by each $A \to B$ in $\CT$. Furthermore, if $\CC \in \CAlg(\PPr_{\st})_{\Fil\Sp / }$ and $\CT \subset \CC$ is a symmetric monoidal full inclusion, then $\psi^{\Fil}$ and $E_*^{*,*} \circ \psi^{\Fil}$ are lax symmetric monoidal under Day convolutions. 
\end{construction}

\begin{construction}[equivariant SS] \label{equivariant-SS}
    Let $G$ be a finite group.  
    \begin{itemize}
        \item The $\infty$-category $\Fil\Sp^G$ of filtered genuine $G$-spectra (equipped with the Day convolution symmetric monoidal structure) receives two symmetric monoidal left adjoints
        \[i\colon \Fil\Sp \to \Fil\Sp^G, \qquad j\colon \Sp^G = \Fun(\{0\}, \Sp^G) \xrightarrow{\mathrm{LKE}} \Fil\Sp^G.\]
        We write $\oneb = j(\Sb^0) = i(\Sb^{0, 0}) \in \Fil\Sp^G$ for the tensor unit.  For $V \in \RO(G), w \in \Zb$, we write  $\Sb^{V, w} = j(\Sb^V) \otimes i(\Sb^{0, w})$, and for $Y \in \Fil\Sp^G, P \in \Fin_G$ we write $Y[P] = Y \otimes j(\Sigma_+^{\infty} P)$. 
        \item Write $U\colon \Sp^G \to \Sp, \Sb^V \mapsto \Sb^{|V|}$ for the underlying spectra functor. Take $\CT'_G = \{\Sb^V\}_{V \in \RO(G)} \subset \Sp^G$ and $\CT_G = \{\Sb^{V, w}[P]\}_{V \in \RO(G), w \in \Zb, P \in \Fin_G} \subset \Fil\Sp^G$. Applying Construction \ref{parametrized-SS} to $\CT_G$,  we obtain a lax symmetric monoidal functor $\psi = \psi^{\Fil}\colon \Fil\Sp^G \to \Fun(\CT_G^{\op}, \SpSeq)$. Furthermore, by precomposing the symmetric monoidal  functor 
        \begin{align*}
            \CT'_G \times \Span(\Fin_G) &\xrightarrow{(\id, U) \times \id} \CT'_G \times \Pic(\Sp) \times \Span(\Fin_G) \\
            &\xrightarrow{\id \times \pi_0 \times \id} \CT'_G \times \Zb^{\delta} \times \Span(\Fin_G) \xrightarrow{j \times i \times j\Sigma^{\infty}_+} \CT_G \times \CT_G \times \CT_G \xrightarrow{\otimes} \CT_G.
        \end{align*}
        we obtain a lax symmetric monoidal functor $\Fil\Sp^G \to \Fun({\CT'_G}^{\op}, \Mack_G(\SpSeq))$. Thus, for each $Y \in \Fil\Sp^G$ we extract an $\RO(G)$-graded Mackey-functor-valued SS\footnote{or rather, an $\RO(G)$-family of SS-valued Mackey functors}. 
        
        We introduce some notational conventions: take $Y \in \Fil\Sp^G, V \in \RO(G), P \in \Fin_G$, then
        \begin{itemize}
            \item We write $\psi_{V} = \psi^{\Fil}_{\Sb^{V, |V|}}$ and $E_r^{s, t, V}(Y[P]) = E_r^{s, t, \Sb^{V, |V|}[P]}(Y) = E_r^{s, t}(\psi_{V}(Y[P]))$. 
            \item We further write $E_r^{s, V}(Y[P]) = E_r^{s, 0, V}(Y[P])$. Conversely, under this convention, we have $E_r^{s, t, V}(Y[P]) \cong E_r^{s, t + V}(Y[P])$. The differential now goes as $d_r\colon E_r^{s, V} \to E_r^{s + r, V + r - 1}$. 
        \end{itemize}
        Therefore, it makes sense to say $\{E_r^{s, V}(Y[P])\} \Rightarrow \pi_{V - s}((Y[P])[\defopara^{-1}])$ if $Y$ is $\defopara$-complete and the corresponding (weak, whole-plane) obstructions vanish. 
    \end{itemize}
\end{construction}

\begin{remark}
    For $x \in E_r^{s, V}$, following the convention in \cite{HSWX-BPC4<2>} we will write $|x| = (V - s, s)$ and refer to this pair as the \textbf{bidegree} of $x$. This corresponds to the folklore convention for coordinates in the pictures of equivariant SS, cf. the figures in \cite{HSWX-BPC4<2>}.
\end{remark}

\subsection{Equivariant slice SS and Hill--Hopkins--Ravenel theories}
\label{subsec:5.2}

From this point to the end of this section, we assume each $\Sp^G$ is localized at the prime $2$. 

\begin{recollection} \label{slice-filtrations}
    Recall the notion of \emph{regular slice filtrations} from \cite{Ull13, HHR-book}. 
    \begin{itemize}
        \item The \textbf{regular slice cells} in $\Sp^G$ are genuine $G$-spectra of the form $\Indrm_H^G\Sb^{k \rho_H}$, where $H$ is a subgroup of $G$, $\rho_H$ is the regular representation of $H$ and $k \in \Zb$. 
        \item Write $\Sp^G_{\geq n} \subset \Sp^G$ for the full subcategory of \textbf{regular slice $n$-connective $G$-spectra}. Concretely, this is generated under colimits by the slice cells $\Indrm_H^G\Sb^{k \rho_H}$ with $k|H| \geq n$. By \cite[Theorem A]{Hill-Yarnall}, $X \in \Sp^G$ is regular slice $n$-connective iff for each $H \subset G$, the geometric fixed point $\Phi^H(X) \in \Sp$ is $\lceil n / |H| \rceil$-connective. 
        \item The inclusion $j_n\colon \Sp^G_{\geq n} \hookrightarrow \Sp^G$ preserves all colimits, so it admits a right adjoint $P_n$. For $X \in \Sp^G$, we refer to $P_n X$ as the \textbf{regular slice $n$-connective cover}. Since $\{\Sp^G_{\geq n}\}_{n \in \Zb}$ forms a nested sequence with $\bigcap_n \Sp^G_{\geq n} = \{0\}, \bigcup_n \Sp^G_{\geq n} = \Sp^G$, for each $X \in \Sp^G$ we obtain a filtration 
        \[P_{\bullet} X \colon \cdots \to  P_{n + 1} X \to P_n X \to P_{n - 1} X \to \cdots\]
        so that $\lim_{n \to \infty} P_n X = 0, \colim_{n \to -\infty} P_n X = X$. This  is the \textbf{regular slice filtration} of $X$. Note that if $H \subset G$, then $\Res^G_H (P_\bullet X) \cong P_\bullet (\Res^G_H X)$. In particular, $(P_\bullet X)^e = \tau_{\geq \bullet} (X^e)$. 
        \item The assignment $X \mapsto P_\bullet X$ yields a functor $P_\bullet\colon \Sp^G \to \Fil(\Sp^G)$. Furthermore, according to \cite[Construction C.6 and Example C.7]{BHS2}, $P_\bullet$ is lax symmetric monoidal. 
        \item We refer to $\{E_r^{s, V}(P_\bullet X)\}_{r \geq 2, s \in \Zb, V \in \RO(G)}$ as  the \textbf{regular $G$-slice SS} of $X$. Furthermore, for each $H \subset G$, we refer to $\{E_r^{s, V}(P_\bullet X[G / H])\}_{r \geq 2, s \in \Zb, V \in \RO(G)}$ as the \textbf{regular $H$-slice SS} of $X$. Note that $E_r^{s, V}(P_\bullet X[G / H]) \cong E_r^{s, \Res^G_H V} (P_\bullet (\Res^G_H X))$. In particular, for $H = e$ the SS degenerates at the $E_2$-page.
    \end{itemize}
    In the following, we will suppress all ``regular'', referring e.g. to $P_\bullet X$ just as the \emph{slice filtration} of $X$. 
\end{recollection}

\begin{remark} \label{filtered-HHR-norm}
    In \cite[\S~3]{Car25}, the author equips $\Fil(\Sp^G)$ with a $G$-symmetric monoidal structure using equivariant Day convolution, and proves that $P_\bullet$ is $G$-symmetric monoidal. 
\end{remark}

\begin{recollection}
    We introduce the \textbf{Hill--Hopkins--Ravenel theories} following \cite{Hu-Kriz, HHR-paper}.
    \begin{itemize}
        \item The complex $j$-homomorphism $j_{\Cb}\colon \mathrm{BU} \to \mathrm{BGL}_1(\Sb)$ is $C_2$-equivariant for complex conjugation on the left and trivial action on the right. Taking the Thom construction (and employing its functoriality as in \cite[Definition 2.20]{ABGHR}), we obtain a Borel complete $C_2$-spectrum $\MU_{\Rb}$. Concretely, $\MU_{\Rb}^e \cong \MU, \MU_{\Rb}^{C_2} \cong \MU_{\Rb}^{hC_2}$ and $\Phi^{C_2}(\MU_{\Rb}) \cong \MO$, cf. \cite[\S~2]{Hu-Kriz}.
        \item According to \cite[Proof of Theorem 2.33]{Hu-Kriz}, there is a $C_2$-equivariant map $e_{\Rb}\colon \MU_{\Rb} \to \MU_{\Rb}$ whose underlying map (resp. map on $C_2$-geometric fixed points) is the Quillen idempotent on $\MU$ (resp. on $\MO$), and we write $\BP_{\Rb}$ for its mapping telescope. By construction, $\BP_{\Rb}$ is Borel complete (thus $\BP_{\Rb}^e \cong \BP, \BP_{\Rb}^{C_2} \cong \BP_{\Rb}^{hC_2}$), and $\Phi^{C_2}(\BP_{\Rb}) \cong \mathrm{H}\Fb_2$.
        Note that $e_{\Rb}$ is an $\Eb_1$-ring map due to \cite[Theorem 6.2.1]{QZ25}, so $\BP_{\Rb}$ is an $\Eb_1$-algebra over $\MU_{\Rb}$. 
        \item For any group homomorphism $\rho\colon K \to G$, there is a unique symmetric monoidal sifted colimit preserving functor $\rho_{\otimes}\colon \Sp^K \to \Sp^G$ which sends $\Sigma_+^{\infty} X$ ($X \in \Fin_K$) to\footnote{here $\Map^K(G, X)$ is the set of left $K$-equivariant maps from $G$ to $X$, with the left $G$-action $(g_0 \cdot f)(g) = f(g g_0)$} $\Sigma_+^{\infty} \Map^K(G, X)$ (cf. \cite[\S~9.2]{BH21}). If $\rho$ is surjective with $G = K / A$ then $\rho_{\otimes}$ recovers the geometric fixed point $\Phi^A$. If $\rho$ is injective, then $\rho_{\otimes} = \N_K^G$ is the \textbf{Hill--Hopkins--Ravenel norm} functor in \cite[\S~2.2.3]{HHR-paper}. By construction, $(\N_K^G(X))^e \cong (X^e)^{\otimes|G / K|}$ for each $X \in \Sp^K$, and $\N_K^G(\Sb^V) \cong \Sb^{\Indrm_K^G(V)}$ for each (virtual) $K$-representation $V$. 
        \item For any finite $2$-group $G$ equipped with a fixed inclusion $C_2 \subset G$, we can construct two \textbf{Hill--Hopkins--Ravenel (HHR) theories} $\MUG = \N_{C_2}^{G}\MU_{\Rb}$ and $\BPG = \N_{C_2}^{G}\BP_{\Rb}$. By symmetric monoidality, $\MUG \in \CAlg(\Sp^G)$ and $\BPG$ is an $\Eb_1$-algebra over $\MUG$.
        \item According to \cite[\S~5]{HHR-paper}, for $G = C_{2^n}$ there is a family of classes ${\tbar}_i \in \pi_{i \rho_2}^{C_2}\MUG$ (where $\rho_2$ is the regular representation of $C_2$) lifting classes $t_i \in \pi_{2i}^{C_2}\MUG$ for each $i \geq 1$, so that 
        \[\pi_{*}^{e}\MUG \cong \Zb[t_1, \gamma t_1, \ldots, \gamma^{2^{n - 1} - 1} t_1, t_2 , \gamma t_2, \ldots]\]  
        where the generator $\gamma$ of $G = C_{2^n}$ acts on $\pi_{*}^{e}\MUG$ by the formulas $\gamma \cdot \gamma^k t_i = \gamma^{k + 1} t_i$ for $0 \leq k < 2^{n - 1} - 1$ and $\gamma \cdot \gamma^{2^{n - 1} - 1}t_i = -t_i$. Under this notation, we also have \[\pi_{*}^{e}\BPG \cong \Zb[t_1, \gamma t_1, \ldots, \gamma^{2^{n - 1} - 1} t_1, t_3 , \gamma t_3, \ldots]\]
        where the indices of $t_i$ take values $i = 2^k - 1$ for $k > 0$.
        \item For each $S \subset \Zb_{\geq 1}$, the collection $\big\{{\tbar}_i \colon \Sb^{i \rho_2} \to \Res^G_{C_2} \MUG\big\}_{i \in S}$ yields a map\footnote{here if $S  \subset \Zb_{\geq 1}$ is infinite, we write $\bigotimes_{i \in S} R_i$ for $\colim_{S' \subset S \text{ finite}} \bigotimes_{i \in S'} R_i$} 
        $\oneb[\{{\tbar}_i\}_{i \in S}] = \bigotimes_{i \in S} \mathrm{Free}_{\Eb_1}(\Sb^{i \rho_2}) \to \bigotimes_{i \in S}  \Res^G_{C_2} \MUG \to \Res^G_{C_2} \MUG$
        in $\Alg(\Sp^{C_2})$, which further gives rise to a composite
        \[A_S = \oneb[\{G \cdot {\tbar}_i\}_{i \in S}] = \N_{C_2}^G\oneb[\{{\tbar}_i\}_{i \in S}] \to \N_{C_2}^{G} \Res^G_{C_2}\MUG \to \MUG\]
        in $\Alg(\Sp^G)$, where the first map comes from monoidality of $\N_{C_2}^G$ and the second map comes from the normed $\Eb_{\infty}$-algebra structures on $\MU_{\Rb}$ and $\MUG = \N_{C_2}^G \MU_{\Rb}$ due to \cite[Theorem 2.4]{HM17} and the norm-forgetful adjunction in \cite[Theorem 8.5 and the discussion after Remark 8.6]{BH21}. Also, each $A_S$ admits a canonical $\Eb_1$-algebra map to $\oneb$ which kills all ${\tbar}_i$. Thus, for $1 \leq m \leq \infty$, it makes sense to introduce the following constructions:
        \begin{itemize}
            \item We write $S'_m = \Zb_{> m}, A'_{> m} = A_{S'_m}$ and set $\MUG\langle m \rangle = \MUG \otimes_{A'_{> m}} \oneb$. 
            \item We write $S_m = \{2^k - 1\}_{k > m}, A_{> m} = A_{S_m}$ and set $\BPG\langle m \rangle = \BPG \otimes_{A_{> m}} \oneb$. 
        \end{itemize}
        We refer to these as the \textbf{truncated HHR-theories}.
    \end{itemize}
\end{recollection}

\begin{remark} \label{(total)-Leibniz-rule-for-slice-SS}
    By construction, all HHR theories, truncated or not, are modules over $\MUG$. Also, all HHR theories constructed out of $\BPG$ are modules over $\BPG$. Thus, the slice filtrations of all these objects are modules over $P_\bullet \MUG$ (resp. $P_\bullet \BPG$), so due to Theorem \ref{Burklund's-Leibniz-rule}, we have Leibniz rules for ordinary and total differentials along the maps $P_\bullet\BPG \otimes P_\bullet \BPG \to P_\bullet\BPG$,  $P_\bullet\BPG \otimes P_\bullet(\BPG\langle m \rangle) \to P_\bullet(\BPG\langle m \rangle)$, etc. \parr
    
    As the maps $P_\bullet \MUG \to P_\bullet \BPG$ and $P_\bullet \BPG \to P_\bullet(\BPG\langle m \rangle)$ are all surjective on the slice $E_2$-page, we can treat the Leibniz rule from the module structure as if there is an $\Eb_1$-ring structure on each $P_\bullet(\BPG\langle m \rangle)$ which leads to the Leibniz rule on their equivariant slice SS, cf. \cite[Remark 2.2]{HSWX-BPC4<2>}.
\end{remark}

\subsection{Slice SS of \texorpdfstring{$\BP_{\Rb}$}{BPR}}
\label{subsec:5.3}

\begin{fact}[Homology of a point] \label{HZ-pt-C2}
    The real representation ring $\RO(C_2)$ is generated by the trivial representation $1$ and the sign representation $\sigma = \sigma_2$. The regular representation $\rho_2$ decomposes as $\rho_2 = 1 + \sigma_2$. Furthermore, for the constant Mackey functor $\underline{\Zb}$, the $\RO(C_2)$-graded homology of a point takes the form $\HH\underline{\Zb}_\filledstar = \pi^{C_2}_\filledstar(\HH\underline{\Zb}) = \Zb[a_{\sigma}, u_{2 \sigma}] / (2 a_\sigma) \oplus \CN_2$, where 
    \begin{itemize}
        \item $a_\sigma$ comes from the \textbf{Euler class} with $|a_\sigma| = -\sigma$. 
        \item $u_{2\sigma}$ comes from the \textbf{Thom class} with $|u_{2 \sigma}| = 2 - 2\sigma$. 
        \item $\CN_2$ stands for the \textbf{negative cone}\footnote{which we ignore in this paper} square-zero ideal.
    \end{itemize}
    For reference, see \cite[Theorem 2.8 and Appendix B]{Dug05}, where the earlier computation is traced back to Stong; for a modern account in the notation used here, see also \cite[\S 3]{HHR-BPC4<1>}. Throughout, we use the multiplicative shorthand $a_{k\sigma}=a_\sigma^k$ and $u_{2k\sigma}=u_{2\sigma}^k$; the same convention applies to all Euler classes, Thom classes, and related notation.

\end{fact}

\begin{recollection}[degree filtration] \label{degree-filtration-C2}
    Take $G = C_2$ and $1 \leq m \leq \infty$. 
    \begin{itemize}
        \item We write $\Zb_{[1, m]} = \{k \in \Zb \mid 1 \leq k < m + 1 \}$ and take $\CI_m = \Map^{\mathrm{fs}}(\Zb_{[1, m]}, \Nb)$, the set of maps $\alpha\colon \Zb_{[1, m]} \to \Nb$ so that $\alpha(k) = 0$ for all but finitely many $k$. We equip $\CI_m$ with the alphabetic order and pointwise addition. For $\alpha \in \CI_m$, we write $|\alpha| = \sum_{k \in \Zb_{[1, m]}} \alpha(k) (2^k - 1) \in \Nb$. 
        \item We denote by $A_{C_2, \leq m}$ the $\Eb_1$-ring $\bigotimes_{1 \leq  k < m + 1}\Free_{\Eb_1}(\Sb^{(2^{k} - 1)\rho_2})$ in $\Sp^{C_2}$. Due to \cite[Proposition 4.1.1.18]{HA}, $A_{C_2, \leq m} \cong \bigoplus_{\substack{\alpha \in \CI_m}} \Sb^{|\alpha| \rho_2}$. Therefore, $A_{C_2, \leq m}$ comes with an $\Eb_1$-multiplicative \textbf{degree filtration} $A_{C_2, \leq m, \bullet}$ in which $A_{C_2, \leq m, w} = \bigoplus_{\alpha \in \CI_m,  |\alpha| \geq w / 2} \Sb^{|\alpha| \rho_2}$. 
    \end{itemize}
\end{recollection}

\begin{remark}
    Conceptually, $A_{C_2, \leq m, \bullet} \cong \bigotimes_{1 \leq  k < m + 1} \Free_{\Eb_1}(\Sb^{(2^{k} - 1)\rho_2, 2^{k + 1} - 2})$ in $\Fil\Sp^{C_2}$.
\end{remark}

\begin{theorem}[Slice theorem]
    For each $1 \leq m \leq \infty$, $\BP_{\Rb} \langle m \rangle \otimes_{A_{C_2, \leq m}} A_{C_2, \leq m, \bullet} \cong P_\bullet (\BP_{\Rb} \langle m \rangle)$. Moreover, on the associated graded we have $\gr_\bullet(P_\bullet (\BP_{\Rb} \langle m \rangle)) \cong \HH\underline{\Zb} \otimes \gr_\bullet(A_{C_2, \leq m, \bullet})$.
\end{theorem}

\begin{proof} 
    This is essentially the result of \cite[\S~6]{HHR-paper} for $G = C_2$, first proved in \cite[\S~3]{Hu-Kriz}.
\end{proof}

\begin{corollary}
    The $E_2$-page of the $C_2$-slice SS of $X = \BP_{\Rb}\langle m \rangle$ takes the form
    \[E_2^{s, V}(P_\bullet X) \cong \pi_{V - s} (\gr_{|V|} P_\bullet X) \cong \bigoplus_{\substack{\alpha \in \CI_m \\ |\alpha| = |V| / 2}} \HH \underline{\Zb}_{V - s - |\alpha| \rho_2}\]
    in particular $E_2^{s, V} = 0$ if $|V|$ is odd. Concretely, each slice cell $\Sb^{|\alpha| \rho_2}$ contributes a copy of $\HH\underline{\Zb}_{\filledstar - |\alpha| \rho_2}$ with generator ${\tbar}^{\alpha} = \prod_{k \geq 1} ({\tbar}_{2^k - 1})^{\alpha(k)}$ at bidegree $(|\alpha| \rho_2 , 0) = (|\alpha| + |\alpha| \sigma, 0)$, so the full $E_2$-page is a polynomial ring $\HH\underline{\Zb}_\filledstar [{\tbar}_1, {\tbar}_3, {\tbar}_7, \ldots, {\tbar}_{2^m - 1}]$, with $|a_\sigma| = (-\sigma, 1)$, $|u_{2\sigma}| = (2 - 2\sigma, 0)$ and $|{\tbar}_{2^k - 1}| = (2^k - 1 + (2^k - 1)\sigma, 0)$. 
\end{corollary}

\begin{theorem}[Slice differential theorem] \label{BPR-slice-diffs}
    In the $C_2$-slice SS of $X = \BP_{\Rb}\langle m \rangle$, we have a family of nontrivial differentials $d_{2^{k + 1} - 1}(u_{2^k \sigma}) = {\tbar}_{2^k - 1} a_{(2^{k + 1} - 1) \sigma}$ for $1 \leq k \leq m$. Together with the fact that $1$, $a_\sigma$, and ${\tbar}_{2^k - 1}$ are permanent cycles for $1 \leq k \leq m$, this family determines all differentials in this SS by the Leibniz rule.
\end{theorem}

\begin{proof}
    This follows from \cite[Theorem 9.9]{HHR-paper} for $G = C_2$, also first proved in \cite[\S~3]{Hu-Kriz}.
\end{proof}

\subsection{Slice SS of \texorpdfstring{$\BPCfour$}{BPC4}: known partial data}
\label{subsec:5.4}

From here to the end of this section, we will work with $G = C_4$ unless otherwise specified. We write $\Res$ for $\Res^{C_4}_{C_2}$, and similarly for $\Indrm, \N, \tr, \res$. We start with the slice $E_2$-pages.

\begin{fact}[Homology of a point] \label{HZ-pt-C4}
    The real representation ring $\RO(C_4)$ is generated by the trivial representation $1$, the sign representation $\sigma = \sigma_4$ and the $\pi/2$-rotation representation $\quarterrep$, and the regular representation $\rho_4$ decomposes as $\rho_4 = 1 + \sigma + \quarterrep$. Moreover, for the constant $C_4$-Mackey functor $\underline{\Zb}$:  
    \begin{itemize}
        \item $\HH\underline{\Zb}_\filledstar = \pi^{C_4}_\filledstar(\HH\underline{\Zb}) = \Zb[a_{\quarterrep}, a_{\sigma}, u_{\quarterrep}, u_{2 \sigma}] / (4 a_{\quarterrep}, 2 a_{\sigma}, u_{\quarterrep} a_{2\sigma} - 2 u_{2 \sigma} a_{\quarterrep}) \oplus \CN_4$, where 
        \begin{itemize}
            \item $a_{\quarterrep}, a_{\sigma}$ come from \textbf{Euler classes} with $|a_{\quarterrep}| = - \quarterrep, |a_{\sigma}| = - \sigma$.
            \item $u_{\quarterrep}, u_{2\sigma}$ come from \textbf{Thom classes} with $|u_{\quarterrep}| = 2 - \quarterrep, |u_{2\sigma}| = 2 - 2\sigma$. 
            \item $\CN_4$ stands for the \textbf{negative cone}\footnote{which we ignore in this paper} square-zero ideal.
        \end{itemize}
        \item For each $V \in \RO(C_4)$, $\pi^{C_4}_V (\HH\underline{\Zb}[C_4 / C_2]) \cong \pi^{C_2}_{\Res(V)} (\HH \underline{\Zb})$ as in Fact \ref{HZ-pt-C2}. Therefore,
        \[\HH\underline{\Zb}[C_4 / C_2]_\filledstar = \pi^{C_4}_\filledstar (\Indrm \Res(\HH\underline{\Zb})) \cong \Zb[a_{2\sigma_2}, u_{2 \sigma_2}, u_\sigma^{\pm 1}] / (2 a_{2 \sigma_2}) \oplus \CN'_2.\] 
        Here $u_{\sigma}$ comes from the isomorphism $\Sigma^{1 - \sigma} X[C_4 / C_2] \cong \Indrm(\Res(\Sb^{1 - \sigma}) \otimes \Res(X)) \cong \Indrm\Res X \cong X[C_4 / C_2]$ for each $X \in \Sp^{C_4}$. We take the convention that $|a_{2 \sigma_2}| = - \quarterrep, |u_{2 \sigma_2}| = 2 - \quarterrep$ and $|u_{\sigma}| = 1 - \sigma$. The natural $\HH\underline{\Zb}_\filledstar \!$-action on $\HH\underline{\Zb}[C_4 / C_2]_\filledstar$ is presented by the formulas $a_{\quarterrep} \cdot x = a_{2 \sigma_2} x, u_{\quarterrep}\cdot x = u_{2\sigma_2} x, a_{\sigma} \cdot x = x$ and $u_{2 \sigma} \cdot x = u_{\sigma}^2 x$.  
    \end{itemize}
    For reference, see \cite[\S~3]{HHR-BPC4<1>}.
\end{fact}

\begin{remark} \label{induction-convention-for-RO(C2)-graded-abelian-groups}
    For any $\RO(C_2)$-graded abelian group $A = A_\filledstar$, we write $\Ev(A)$ for its restriction to the subring of $\RO(C_2)$ generated by $1$ and $2\sigma_2$, and we denote by $\Indrm_{\filledstar}(A) = \Ev(A) \otimes_{\Zb} \Zb[u_\sigma^{\pm 1}]$ the $\RO(C_4)$-graded abelian group with $x \in A_{p + 2 q \sigma_2}$ in $\RO(C_4)$-degree $p + q \quarterrep$ and $|u_\sigma| = 1 - \sigma$. Under these notations, $\pi^{C_4}_\filledstar(\Indrm(X)) \cong \Indrm_{\filledstar}(\pi^{C_2}_\filledstar(X))$, in particular $\HH\underline{\Zb}[C_4 / C_2]_\filledstar = \Indrm_{\filledstar}(\Res(\HH\underline{\Zb})_\filledstar)$.
\end{remark}

\begin{lemma} \label{binomial-theorem-for-norms}
    Suppose $H \subset G$ is an index-$2$ subgroup. Then for $X, Y \in \Sp^H$, there is a natural identification $\N_H^G(X \oplus Y) \cong \N_H^G(X) \oplus \N_H^G(Y) \oplus \Indrm_H^G(X \otimes Y)$.
\end{lemma}
\begin{proof}
    This follows from the model category discussion in \cite[Proposition A.37]{HHR-paper}. Alternatively, according to \cite[Proposition 3.4.19]{EH23}, from the distributivity diagram in $\Fin_G$ 
    \[\begin{tikzcd}
        & {G / H \sqcup G / H \sqcup (G / H \sqcup G / H)} \ar[ld, "{(p_1, p_2, p_1, p_2)}"'] \ar[rr, "{i_H^G \sqcup i_H^G \sqcup \nabla}"] \ar[dd, "{\nabla}"] && {G / G \sqcup G / G \sqcup G / H} \ar[dd, "{(\id, \id, i_H^G)}"] \\
        G / H \sqcup G / H \ar[rd, "{\nabla}"] &&& \\
        & G / H \ar[rr, "{i_H^G}"] && G / G 
    \end{tikzcd}\]
    we can extract a commutative pentagon
    \[\begin{tikzcd}
        & {\Sp^H \times \Sp^H \times (\Sp^H \times \Sp^H)} \ar[rr, "{\N_H^G \times \N_H^G \times  \otimes}"] & & {\Sp^G \times \Sp^G \times \Sp^H} \ar[dd, "{\id \oplus \id \oplus \Indrm_H^G}"] \\
        {\Sp^H \times \Sp^H} \ar[rd, "{\oplus}"] \ar[ru, "{(p_1, p_2, p_1, p_2)}"] &&& \\
        & {\Sp^H} \ar[rr, "{\N_H^G}"] && {\Sp^G}. 
    \end{tikzcd}\]
    Thus the identification $\N_H^G(X \oplus Y) \cong \N_H^G(X) \oplus \N_H^G(Y) \oplus \Indrm_H^G(X \otimes Y)$ follows.
\end{proof}

\begin{recollection}[degree filtration] \label{degree-filtration-C4}
    We write $A_{C_4, \leq m} = \N_{C_2}^{C_4} (A_{C_2, \leq m})$. Using Lemma \ref{binomial-theorem-for-norms} and the fact that $\N$ commutes with filtered colimits, we deduce that (under the alphabetic order on $\CI_m$):
    \[A_{C_4, \leq m} \cong \bigoplus_{\alpha \in \CI_m} \N(\Sb^{|\alpha| \rho_2}) \oplus \bigoplus_{\substack{\alpha, \beta \in \CI_m \\ 
    \alpha > \beta}} \Indrm(\Sb^{|\alpha| \rho_2} \otimes \Sb^{|\beta| \rho_2}) \cong \bigoplus_{\alpha \in \CI_m} \Sb^{|\alpha| \rho_4} \oplus \bigoplus_{\substack{\alpha, \beta \in \CI_m \\ 
    \alpha > \beta}} \Indrm(\Sb^{|\alpha + \beta| \rho_2}) \]
    here $\rho_4$ is the regular representation for $C_4$. Therefore, $A_{C_4, \leq m}$ admits a ($\Eb_1$-multiplicative) \textbf{degree filtration} $A_{C_4, \leq m, \bullet}$ with $A_{C_4, \leq m, w} = \bigoplus_{|\alpha| \geq w / 4} \Sb^{|\alpha| \rho_4} \oplus \bigoplus_{\alpha > \beta, |\alpha + \beta| \geq w / 2} \Indrm(\Sb^{|\alpha + \beta| \rho_2})$.
\end{recollection}

\begin{theorem}[slice theorem] \label{slice-theorem-C4}
    Take $1 \leq m \leq \infty$. Treating $\BPCfour\langle m \rangle$ as an $A_{C_4, \leq m}$-module using the $\Eb_1$-ring map $A_{C_4, \leq m} \to \MUCfour$ in the construction of $\BPCfour\langle m \rangle$, we have 
    \[\BPCfour \langle m \rangle \otimes_{A_{C_4, \leq m}} A_{C_4, \leq m, \bullet} \cong P_\bullet (\BPCfour \langle m \rangle), \qquad \gr_\bullet(P_\bullet (\BPCfour \langle m \rangle)) \cong \HH\underline{\Zb} \otimes \gr_\bullet(A_{C_4, \leq m, \bullet}).\]
\end{theorem}
\begin{proof}
    This is essentially the result in \cite[\S~6]{HHR-paper} for $G = C_4$.
\end{proof}

\begin{remark} \label{slice-theorem-in-general}
    Conceptually, under the terminology of filtered HHR normed functor $\N_{C_2}^G\colon \Fil\Sp^{C_2} \to \Fil\Sp^{G}$ coming from the $G$-symmetric monoidal structure in Remark \ref{filtered-HHR-norm}, $A_{C_4, \leq m, \bullet} \cong \N_{C_2}^{C_4} (A_{C_2, \leq m, \bullet})$. In general, for $G = C_{2^k}$, write $A_{G, \leq m, \bullet} = \N_{C_2}^{G} (A_{C_2, \leq m, \bullet})$, then we also have slice theorems 
    \[\BPG \langle m \rangle \otimes_{A_{G, \leq m}} A_{G, \leq m, \bullet} \cong P_\bullet (\BPG \langle m \rangle), \qquad \gr_\bullet(P_\bullet (\BPG \langle m \rangle)) \cong \HH\underline{\Zb} \otimes \gr_\bullet(A_{G, \leq m, \bullet})\] 
    following the same discussion as in \cite[\S~6]{HHR-paper}. The concrete enumeration for slice cells (in particular the induced cells) appears to be involved, cf. \cite[\S~2.4.1]{HHR-paper}. 
\end{remark}

\begin{corollary}
    \begin{itemize}
        \item Note that $\Res(A_{C_4, \leq m, \bullet}) \cong A_{C_2, \leq m, \bullet}^{\otimes 2}$ due to the general formula $\Res^G_H \circ \N^G_H \cong (-)^{\otimes |G / H|}$. Therefore, the $E_2$-page of the $C_2$-slice SS of $X = \BPCfour\langle m \rangle$ takes the form
        \[E_2^{s, V}(P_\bullet X [C_4 / C_2]) \cong \pi_{\Res(V) - s} (\gr_{|V|} \Res(P_\bullet X)) \cong \bigoplus_{\substack{\alpha, \beta \in \CI_m \\ |\alpha + \beta| = |V| / 2}} \Res(\HH \underline{\Zb})_{\Res(V) - s - |\alpha + \beta| \rho_2}.\] 
        In particular, $E_2^{s, V}(P_\bullet X [C_4 / C_2]) = 0$ if $|V|$ is odd. Concretely, $E_2^{*, \filledstar}(\Res (P_\bullet X))$ is the $\RO(C_2)$-graded polynomial ring $\Res(\HH \underline{\Zb})_\filledstar [{\tbar}_1, \gamma {\tbar}_1, {\tbar}_3, \gamma {\tbar}_3, \ldots, {\tbar}_{2^m - 1}, \gamma {\tbar}_{2^m - 1}]$, where $|a_{\sigma_2}| = (-\sigma_2, 1)$, $|u_{2 \sigma_2}| = (2 - 2\sigma_2, 0)$, and $|{\tbar}_{2^k - 1}| = |\gamma {\tbar}_{2^k - 1}| = ((2^k - 1) + (2^k - 1)\sigma_2, 0)$, while $E_2^{s, \filledstar}(P_\bullet X [C_4 / C_2]) = \Indrm_{\filledstar}(E_2^{s, \filledstar}(\Res (P_\bullet X)))$ in terms of Remark \ref{induction-convention-for-RO(C2)-graded-abelian-groups}. 
        \item The $E_2$-page of the $C_4$-slice SS of $X = \BPCfour\langle m \rangle$ takes the form
        \[E_2^{s, V}(P_\bullet X) \cong \pi_{V - s} (\gr_{|V|} P_\bullet X) \cong \bigoplus_{\substack{\alpha \in \CI_m \\ |\alpha| = |V| / 4}} \HH \underline{\Zb}_{V - s - |\alpha| \rho_4} \oplus \bigoplus_{\substack{\alpha, \beta \in \CI_m, \alpha > \beta \\ |\alpha + \beta| = |V| / 2}} \Res( \HH \underline{\Zb})_{\Res(V) - s - |\alpha + \beta| \rho_2}. \] 
        In particular, $E_2^{s, V}(P_\bullet X) = 0$ if $|V|$ is odd. Concretely, there are contributions from \textbf{norm (i.e. noninduced) cells} $\Sb^{|\alpha| \rho_4}$ and \textbf{induced cells} $\Indrm(\Sb^{|\alpha + \beta|\rho_2})$:
        \begin{itemize}
            \item For each $\alpha \in \CI_m$, the norm cell $\Sb^{|\alpha| \rho_4}$ gives rise to a copy of $\HH\underline{\Zb}_{\filledstar - |\alpha| \rho_4}$ generated by the \textbf{norm class}\footnote{Norm classes come from the normed algebra structure map $\N(\Res(\MUCfour)) \to \MUCfour$} $\N({\tbar}^{\alpha}) = \prod_{k \geq 1} \N({\tbar}_{2^k - 1})^{\alpha(k)}$ in bidegree $(|\alpha| \rho_4 , 0) = (|\alpha| + |\alpha| \sigma + |\alpha|\quarterrep, 0)$. In other words, all classes from noninduced cells in the $E_2$-page form a polynomial subring $\HH\underline{\Zb}_\filledstar[\dfbar_1, \dfbar_3, \dfbar_7, \ldots, \dfbar_{2^m - 1}]$, where $|a_{\sigma}| = (-\sigma, 1), |a_{\quarterrep}| = (- \quarterrep, 2), |u_{2\sigma}| = (2 - 2 \sigma, 0), |u_{\quarterrep}| = (2 - \quarterrep, 0)$ and $\dfbar_{2^k - 1} = \N({\tbar}_{2^k - 1})$ is of bidegree $(2^{k} - 1 + (2^k - 1)\sigma + (2^k - 1) \quarterrep, 0)$. 
            \item For each pair $\alpha, \beta \in \CI_m$ such that $\alpha > \beta$, the induced cell $\Indrm(\Sb^{|\alpha + \beta| \rho_2})$ gives rise to a copy of $\Indrm_{\filledstar}(\Res(\HH\underline{\Zb})_{\filledstar - |\alpha + \beta|\rho_2})$ spanned by \textbf{transfer classes}\footnote{Transfer classes come from the transfer map $\BPCfour[C_4 / C_2] \to \BPCfour$}
            \[\tr(\tbar^{\alpha} \gamma \tbar^{\beta} u_{2p\sigma_2} a_{q\sigma_2}) u_{k\sigma}\] 
            of bidegree $(|\alpha| + |\beta| + 2 p + k - k\sigma + \frac{|\alpha| + |\beta| - 2p - q}{2} \quarterrep, q)$, where $\gamma\tbar^{\beta} = \prod_{k \geq 1} (\gamma \tbar_{2^k - 1})^{\beta(k)}$. Furthermore, the classes from the induced cell form an ideal of the $E_2$-page. Concretely, multiplication involving transfer classes is computed by the Frobenius relation and the fact $\res(\dfbar_i) = \tbar_i \gamma \tbar_i$, for instance $\tr(\tbar_1 a_{\sigma_2}) \cdot \dfbar_3 u_{3 \quarterrep} a_{\quarterrep} = \tr(\tbar_1 \tbar_3 \gamma \tbar_3 a_{\sigma_2}) u_{\sigma}^{-3} u_{3 \quarterrep}a_{\quarterrep} = \tr(\tbar_1 \tbar_3 \gamma \tbar_3 u_{6 \sigma_2}a_{3\sigma_2}) u_{-3\sigma}$. In particular, transfer classes are $a_\sigma$-torsion since $\res(a_{\sigma}) = 0$. 
        \end{itemize}
    \end{itemize}
\end{corollary}

\begin{remark}
    All classes from the induced cells are transfer classes, but the converse might not hold true, for instance $2 = \tr(1)$ comes from the norm cell $\Sb^0$.
\end{remark}

We then record two sources of $C_4$-slice differentials, the \emph{transfer differentials} and the \emph{sheared differentials}. Transfer differentials come from a complete description of the  $C_2$-slice SS. 

\begin{fact} \label{norm-forgetful-adjunction-unit-map-for-BP}
    The unit of the norm-forgetful adjunction $\MU_{\Rb} \to \Res(\MUCfour)$ commutes with the real Quillen idempotent $e_{\Rb}$ in the $\infty$-category of $\Eb_1$-algebras in $\Sp^{C_2}$. Therefore, it induces an $\Eb_1$-ring map between the summands $\BP_{\Rb} \to \Res(\BPCfour)$, which further leads to a map of $\Eb_1$-multiplicative filtrations $P_\bullet \BP_{\Rb} \to \Res(P_\bullet \BPCfour)$. 
\end{fact}

\begin{proof}
    The claim above actually holds true in a stronger sense, with $\Eb_1$ replaced by the $C_2$-operad $\Eb_{\rho}$ in the sense of \cite{Hill-Erho}, where $\rho = \rho_2$ is the regular representation of $C_2$. In fact, the second claim follows directly from the first one through a slice connectivity argument. To prove the first claim, we systematically use the discussions in \cite{QZ25}. As $\Res(\MUCfour)$ is a strongly even $\Eb_{\infty}^{C_2}$-ring spectrum in the sense of \cite[Definition 5.1.3]{QZ25}, an application of \cite[Theorem 6.2.1]{QZ25} shows it suffices to construct a homotopy that fills the square 
    \[\begin{tikzcd}
        \MU \ar[r, "u"] \ar[d, "{e}"] & \MU \otimes \MU \ar[d, "{e \, \otimes \, e}"] \\
        \MU \ar[r, "u"] & \MU \otimes \MU
    \end{tikzcd}\]
    in $\Alg_{\Eb_2}(\Sp)$, where $e$ is the (nonequivariant) Quillen idempotent and $u$ is the underlying map of $\MU_{\Rb} \to \Res(\MUCfour)$. On the other hand, after passing to Borel completion, the norm-forgetful adjunction for $C_2 \subset C_4$ can be identified with the LKE-forgetful adjunction between $\CAlg(\Sp)^{\mathrm{B} C_2}$ and $\CAlg(\Sp)^{\mathrm{B} C_4}$. Therefore, the underlying map $u$ of the adjunction unit $\MU_{\Rb} \to \Res(\MUCfour)$ is the unit-insertion map $\MU \to \MU^{\otimes |C_4/C_2|} \cong \MU \otimes \MU$ associated to the finite set inclusion $* = C_2/C_2 \subset C_4/C_2$ singled out by the subgroup inclusion $C_2\subset C_4$, and hence is functorial in $\MU$ as an $\Eb_2$-ring. Thus, the above square has a canonical filler as long as $e$ is a map of $\Eb_2$-ring spectra, which is true by \cite[Corollary 1.3]{ChadwickMandell-EnGenera}.
\end{proof}

\begin{theorem} \label{diffs-in-C2-sliceSS}
    In the $C_2$-slice SS of $\BPCfour\langle m \rangle$, there are differentials for each $1 \leq k \leq 2m$ that take the form
    \[d_{2^{k + 1} - 1}(u_{2^{k} \sigma_2}) = {\vbar}_{k} a_{(2^{k + 1} - 1)\sigma_2}\]
    here ${\vbar}_k$ is the image of the $k$-th \textbf{Araki generator}\footnote{in our convention this is precisely ${\tbar}_{2^k - 1} \in \pi_{(2^{k} - 1) \rho_2}^{C_2} \BP_{\Rb} \langle 2m\rangle$} in $\BP_{\Rb}$ under the map in Fact \ref{norm-forgetful-adjunction-unit-map-for-BP}. Concretely, write $I_{k}$ for the ideal $(2, {\vbar}_1, \ldots, {\vbar}_{k - 1}) \subset \pi_{\filledstar}^{C_2} (\Res \, \BPCfour \langle m \rangle)$, then 
    \[{\vbar}_k = \begin{cases}
        \tbar_{2^{k} - 1} + \gamma \tbar_{2^{k} - 1} + \sum_{j = 1}^{k - 1} \gamma \tbar_{2^{j}-1} \cdot \tbar_{2^{k - j}-1}^{2^{j}} \mod I_{k},  & 1 \leq k \leq m. \\
        \sum_{j = k - m}^{m} \gamma \tbar_{2^{j}-1} \cdot \tbar_{2^{k - j}-1}^{2^{j}} \mod I_{k}, & m + 1 \leq k \leq 2m.
    \end{cases}\]
    Furthermore, the pattern for all differentials in this SS is determined by the Leibniz rule, the fact that $a_{\sigma_2}$ and all $\tbar_{2^{k} - 1}, \gamma \tbar_{2^{k} - 1}$ are permanent cycles, and the above family of differentials.
\end{theorem}
\begin{proof}
    This is \cite[Theorem 2.3]{DLS-vanishing-lines} in the case $G = C_4$. The formulas for $\vbar_k$ come from \cite[Theorem 1.1]{BHSZ} for $G = C_4$ or the discussion in \cite[\S~2.4]{HSWX-BPC4<2>}. 
\end{proof}

\begin{remark} \label{transfer-diffs}
    Under the map $\tr\colon  P_\bullet X [C_4 / C_2] \to P_\bullet X$, every differential $d_r(x) = y$ in the $C_2$-slice SS gives rise to a family of \textbf{transfer differentials} $d_r(\tr(x)u_{k\sigma}) = \tr(y)u_{k\sigma}$ in the $C_4$-slice SS. By Theorem \ref{diffs-in-C2-sliceSS}, each essential transfer differential for $X = \BPCfour \langle m \rangle$ has length $2^{k + 1} - 1$ for some $1 \leq k \leq 2m$. 
\end{remark}

Sheared differentials, on the other hand, come from the transchromatic ``shearing'' isomorphism  between certain regions of two different slice spectral sequences. To describe this we shall introduce the following notion:

\begin{definition} \label{ss-isomorphism-on-or-above-L}
    Suppose $\CL: s = k(t - s) + b$ is a line on the $(t - s, s)$-plane, where $k > -1$. 
    \begin{itemize}
        \item We say $F \in \Fil\Sp$ is \textbf{$\CL$-connective} if $E_2^{s, t}(F) = 0$ whenever $(t - s, s)$ is above  $\CL$, i.e. $s > k(t - s) + b$. 
        \item We say $f\colon X \to Y$ in $\Fil\Sp$ is \textbf{$\CL$-connective} if $F = \fib(f)$ is $\CL$-connective. In other words, $f$ is $\CL$-connective if $E_2^{s, t}(X) \to E_2^{s, t}(Y)$ is an isomorphism for $(t - s, s)$ above $\CL$ and a surjection for $(t - s, s)$ on $\CL$. In this case, by induction (cf. \cite[Proof of Theorem 2.6]{MSZ-stratification}), every essential differential in $\{E_{r}^{s, t}(X)\}$ whose source is on or above the line $\CL$ goes to an essential differential in $\{E_{r}^{s, t}(Y)\}$, and conversely each essential differential in $\{E_{r}^{s, t}(Y)\}$ whose source is on or above the line $\CL$ comes from an essential differential in $\{E_{r}^{s, t}(X)\}$. We refer to this situation as an \textbf{isomorphism of SS on or above the line $\CL$}.
    \end{itemize}
\end{definition}

\begin{remark}
    In general, for each line $\CL$ on the $(t - s, s)$-plane whose slope $\neq -1$, we say $F \in \Fil\Sp$ is \textbf{$\CL$-connective} if $E_2^{s, t}(F) = \pi_{t - s, t}(F / \defopara) = 0$ whenever $(t - s, s)$ and $(k, -k)$ lie in different connected components of $\Rb^2 \setminus \CL$ for $k \gg 0$. 
    \begin{itemize}
        \item The collection of $\CL$-connective objects forms the connective part of a $t$-structure on $\Fil\Sp$, whose $(-1)$-coconnective part consists of complete filtrations $F$ so that $E_2^{s, t}(F) = 0$ whenever $(t - s, s)$ lies on $\CL$ or in the same connected component of $\Rb^2 \setminus \CL$ as $(k, -k)$ for $k \gg 0$. In fact, there is a $t$-structure on $\Fil\Sp_{\hat{\defopara}}$ in which $F$ is (co)connective iff each spectrum $(F / \defopara)(w)$ is $\lceil x(w) \rceil$-(co)connective, where $x(w)$ is the horizontal coordinate of the intersection point in $\Rb^2$ between the line $\CL$ and the unique line of slope $-1$ passing through $(w, 0)$. This is straightforward under the identification between complete filtered spectra and coherent cochain complexes \cite[Theorem 3.19]{Ariotta}, and the desired $t$-structure on $\Fil\Sp$ follows by gluing it with the trivial $t$-structure on $\Sp$ (in which all spectra are connective) along the recollement in Remark \ref{recollements-and-Nakayama-lemma} using \cite[Theorem 5.10]{Ariotta}.
        \item For $\CL\colon s = -2(t - s)$, this recovers the Beilinson $t$-structure in Construction \ref{accelerations}. 
        \item By construction, for every $F \in \Fil\Sp$, $E_2^{s, t}(\tau_{\geq 0}^{\CL} F) = 0$ if $(t - s, s)$ and $(k, -k)$ lie in different connected components of $\Rb^2 \setminus \CL$ for $k \gg 0$, otherwise it is isomorphic to $E_2^{s, t}(F)$. Thus, this is not the same as the linear $t$-structure in \cite[Example 2.8]{Lee-Levy} % (due to \cite[discussions before Lemma 2.7]{Lee-Levy}) 
        or \cite[Definition 3.1]{CD26}. %(due to \cite[Proposition 3.10]{CD26}).
        \item For $\CL$ with slope $< -1$, we can reproduce the discussion in Definition \ref{ss-isomorphism-on-or-above-L}, but here an $\CL$-connective map $f\colon X \to Y$ would induce an \emph{isomorphism of SS on or below the line $\CL$}.
    \end{itemize}
\end{remark}

\begin{theorem}[transchromatic shearing isomorphism \& slice recovery theorem] \label{sheared-diffs}
    Take $1 \leq m < \infty$. Recall the functor $\psi_V = \psi^{\Fil}_V \colon \Fil\Sp^{G} \to \Fil\Sp,  Y \mapsto \Msp^{\Fil}(\Sb^{V, |V|}, Y)$ in Construction \ref{equivariant-SS}. 
    \begin{itemize}
        \item There is an isomorphism $P_\bullet \BPCfour\langle m \rangle[a_{\sigma}^{-1}] \cong \infl(\Sl_2(P_\bullet \BP_{\Rb} \langle m \rangle [a_{\sigma_2}^{-1}]))$ in $\Fil\Sp^{C_4}$, for which $\infl = \infl_{C_4 / C_2}^{C_4}$ is the inflation and $\Sl_2$ is the $2$-fold slowdown in Construction \ref{accelerations}. 
        \item For any $X \in \Sp^{C_2}$, $V \in \RO(C_2)$, the map $\psi_V(P_\bullet X) \to \psi_V(P_\bullet X[a_{\sigma_2}^{-1}])$ is $\CL^{C_2}_{V, 0}$-connective and the map $\psi_V(P_\bullet X[a_{\sigma_2}^{-1}]) \to 0$ is $\CL^{C_2}_{V, 1}$-connective, where 
        \[\CL^{C_2}_{V, 0} \colon s = 0, \qquad \CL^{C_2}_{V, 1}\colon s = (t - s) + 2|V^{C_2}| - |V|. \]
        \item For any $X \in \Sp^{C_4}, V \in \RO(C_4)$, the maps $\psi_V(P_\bullet X) \to \psi_V(P_\bullet X[a_{\quarterrep}^{-1}])$, $\psi_V(P_\bullet X[a_{\quarterrep}^{-1}]) \to \psi_V(P_\bullet X[a_{\sigma}^{-1}])$ and\footnote{the map $X[a_{\quarterrep}^{-1}] \to X[a_{\sigma}^{-1}]$ exists since $\oneb[a_{\quarterrep}^{-1}] = \widetilde{\mathrm{E}\CF\!}_{\mathrm{triv}}$ with $\CF_{\mathrm{triv}} = \{e\}$ and $\oneb[a_{\sigma}^{-1}] = \widetilde{\mathrm{E}\CF\!}_{\mathrm{pr}}$ with $\CF_{\mathrm{pr}} = \{e, C_2\}$} $\psi_V(P_\bullet X[a_{\sigma}^{-1}]) \to 0$ are respectively $\CL^{C_4}_{V, 0}, \CL^{C_4}_{V, 1}, \CL^{C_4}_{V, 2}$-connective, where 
        \[\CL^{C_4}_{V, 0} \colon s = 0, \quad \CL^{C_4}_{V, 1}\colon s = (t - s) + 2|V^{C_2}| - |V|, \quad \CL^{C_4}_{V, 2}\colon s = 3(t - s) + 4|V^{C_4}| - |V|. \]
    \end{itemize}
    Consequently, for each $V \in \RO(C_4), V' \in \RO(C_2)$ and $s, t \in \Zb$ with $t - s \geq 0$,  $\{E_r^{s, t + V}(P_\bullet \BPCfour \langle m \rangle)\}$ vanishes above the line $\CL^{C_4}_{V,2}$ and $\{E_r^{s, t + V'}(P_\bullet \BP_{\Rb} \langle m \rangle)\}$ vanishes above the line $\CL^{C_2}_{V', 1}$. Furthermore, from these comparisons we deduce a \text{``shearing'' isomorphism} between $\{E_{r}^{s, t + V^{C_2}}(P_\bullet \BP_{\Rb} \langle m \rangle)\}$ on or above the line $s = 0$ and $\{E_{r'}^{s', t' + V}(P_\bullet \BPCfour \langle m \rangle)\}$ on or above the line $\CL^{C_4}_{V,1}$, under which 
    \begin{gather*}
        E_{r}^{s, t + V^{C_2}} \leftrightarrow E_{2r - 1}^{s + t + 2|V^{C_2}| - |V|, 2t + V + 2|V^{C_2}| - |V|}, \ \ d_r \leftrightarrow d_{2r - 1}, \\ 
        a_{\sigma_2} \leftrightarrow a_{\sigma}, \ \ u_{2\sigma_2} \leftrightarrow u_{2\sigma}, \ \ \text{ and }  \ \  \tbar_{2^k - 1} \leftrightarrow \dfbar_{2^k - 1} a_{(2^k - 1)\quarterrep}.
    \end{gather*}
    Thus, there is a \textbf{sheared differential} $d_{2^{k + 2} -3}(u_{2^k \sigma}) = \dfbar_{2^k - 1} a_{(2^k - 1) \quarterrep} a_{(2^{k + 1} - 1)\sigma}$ for each $1 \leq k \leq m$ in the $C_4$-slice SS of $\BPCfour \langle m \rangle$, and these generate all differentials on or above $\CL^{C_4}_{\filledstar, 1}$ by the Leibniz rule. 
\end{theorem}

\begin{proof}
    The isomorphism in $\Fil\Sp^{C_4}$ is \cite[Theorem 4.1]{MSZ-transchromatic} for $G = C_4$, and the $\CL$-connectivity results follow from \cite[Theorem 2.1]{MSZ-stratification} for $G = C_2, C_4$. The precise translation along the shearing isomorphism is due to \cite[Theorem 5.5 and Theorem 6.2]{MSZ-transchromatic}. Finally, the sheared differentials were first established in \cite[Theorem 9.9]{HHR-paper}, and they can be deduced from Theorem \ref{BPR-slice-diffs} via the shearing isomorphism (as is recorded in \cite[Theorem 6.7]{MSZ-transchromatic}). The latter approach also proves these generate all differentials on or above the lines $\CL_{\filledstar, 1}^{C_4}$ by the Leibniz rule. 
\end{proof}

\subsection{Exotic transfer differentials}
\label{subsec:5.5}

In this subsection we use Burklund--Lin--Wang--Xu methods to construct new families of differentials in the \(C_4\)-slice spectral sequence of \(\BPCfour\langle m\rangle\). Building on the generalized Leibniz rule and the generalized Mahowald trick from the previous sections, we ``splice'' transfer differentials with sheared differentials to produce \emph{exotic transfer differentials}. This recovers, and then extends to all heights \(m\), the exotic transfer differentials found for \(\BPCfour\langle 1\rangle\) in \cite{HHR-BPC4<1>}.

\begin{definition}
    Take $1 \leq m < \infty$. 
    \begin{itemize}
        \item We write $B_m = \Zb[T_1, T_3, \ldots, T_{2^m - 1}]$ for the graded polynomial ring with $|T_{2^i - 1}| = 2^{i} - 1$. 
        \item We write $R_m = \Zb[\tbar_1, \gamma \tbar_1, \tbar_3, \gamma \tbar_3, \ldots, \tbar_{2^m - 1}, \gamma \tbar_{2^m - 1}]$ for the graded ring with $|\tbar_{2^i - 1}| = |\gamma\tbar_{2^i - 1}| = 2^{i + 1} - 2$. We equip it with the $C_2 = \{e, \gamma\}$ action by $\gamma(\tbar_{2^i - 1}) = \gamma \tbar_{2^i - 1}, \gamma(\gamma \tbar_{2^i - 1}) = -\tbar_{2^i - 1}$.
    \end{itemize}
\end{definition}

\begin{theorem}[exotic transfer paradigm] \label{exotic-transfer-paradigm}
    Take $1 \leq m < \infty$. Suppose for $h \in \{1, \ldots, 2m\}$, there is a homogeneous polynomial $P \in B_m$ with $|P| = p \geq 2^{h - 1}$ and two homogeneous classes $x, c \in R_m$ with $|x| = 4p - 2^{h + 1} + 2, |c| = 4p$, so that in $R_m$ we have
    \[x \cdot \vbar_h = P(\tbar_1 \gamma \tbar_1, \tbar_3 \gamma \tbar_3, \ldots, \tbar_{2^m - 1} \gamma \tbar_{2^m - 1}) + c + \gamma(c) \mod I_h = (2, \vbar_1, \ldots, \vbar_{h - 1}).\]
    Then for each $q \in \{1, \ldots, m\}, j \in \Nb$, there is a differential in the $C_4$-slice SS of $\BPCfour \langle m \rangle$
    \[d_{2^{h + 1} + 2^{q + 2} - 5}(\tr(x a_{\sigma_2})u_{2^{h - 1}\quarterrep} u_{k\sigma}) = P(\dfbar_1, \dfbar_3, \ldots, \dfbar_{2^m - 1}) \dfbar_{2^q - 1} u_{2^{q + 1} \! j \sigma} a_{(2^h + 2^q - 1)\quarterrep}a_{(2^{q + 1} - 2)\sigma}\]
    here $k = k(q, p, j) =  2^{q + 1}\! j + 2^q - p - 1 \in \Zb$ (which can be negative). 
\end{theorem}

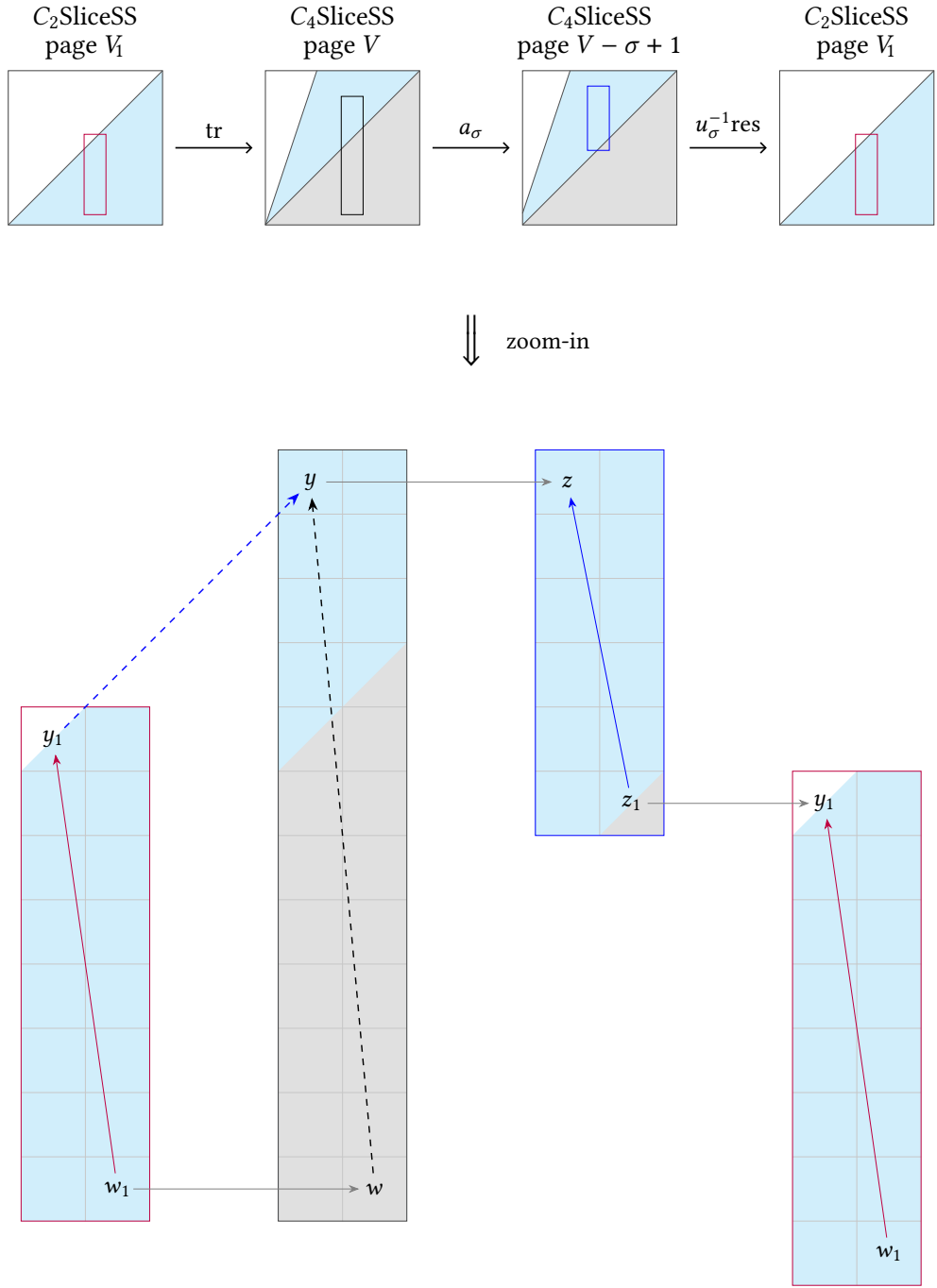
\begin{figure}[htbp]
    \centering
    \scalebox{0.9}{
        \begin{tikzpicture}[line width=0.2pt]
            \tikzset{
                diff/.style={-{Stealth},line width=0.3pt,shorten <=7pt,shorten >=7pt},
                diff1/.style={-{Stealth},dashed,line width=0.6pt,shorten <=7pt,shorten >=7pt},
                extarr/.style={-{Stealth},line width=0.3pt,shorten <=7pt,shorten >=7pt},
                extarr1/.style={-{Stealth},dashed,line width=0.6pt,shorten <=7pt,shorten >=7pt},
                mapabove/.style={-{Straight Barb[scale=0.8]},line width=0.7pt, shorten <=1.4cm,shorten >=1.4cm},
                skel/.style={black!22, line width=0.18pt},
                thumbborder/.style={black!75, line width=0.28pt},
                sepline/.style={black!70, line width=0.22pt},
                selectA/.style={purple, line width=0.42pt},
                selectB/.style={black, line width=0.42pt},
                selectC/.style={blue, line width=0.42pt},
            }

            \colorlet{thumbblue}{cyan!16}
            \colorlet{thumbgray}{black!12}

            %==================================================
            % centers of the four thumbnail squares
            %==================================================
            \coordinate (T1c) at (50.5,28.2);
            \coordinate (T2c) at (54.5,28.2);
            \coordinate (T3c) at (58.5,28.2);
            \coordinate (T4c) at (62.5,28.2);

            %==================================================
            % arrows and labels above the thumbnails
            %==================================================

            \node at ($(T1c) + (0, 1.8)$) {\shortstack{$C_2\mathrm{SliceSS}$ \\ $\mathrm{page}\ V_{\!1}$}};
            \node at ($(T2c) + (0, 1.8)$) {\shortstack{$C_4\mathrm{SliceSS}$ \\ $\mathrm{page}\ V$}}; 
            \node at ($(T3c) + (0, 1.8)$) {\shortstack{$C_4\mathrm{SliceSS}$ \\ $\mathrm{page}\ V - \sigma + 1$}}; 
            \node at ($(T4c) + (0, 1.8)$) {\shortstack{$C_2\mathrm{SliceSS}$ \\ $\mathrm{page}\ V_{\!1}$}}; 
            
            \draw[mapabove]
                ($(T1c)$) -- ($(T2c)$)
                node[midway, above] {\small $\tr$};

            \draw[mapabove]
                ($(T2c)$) -- ($(T3c)$)
                node[midway, above] {\small $a_{\sigma}$};

            \draw[mapabove]
                ($(T3c)$) -- ($(T4c)$)
                node[midway, above] {\small $u_\sigma^{-1}\res$};

            %==================================================
            % upper half: four square thumbnails
            %==================================================

            %----------------------------------------
            % thumbnail 1
            % square: [49.3,51.7] x [27.0,29.4]
            % separator of slope 1 through lower-left corner
            %----------------------------------------
            \begin{scope}
                \clip (49.3,27.0) rectangle (51.7,29.4);
                \fill[thumbblue] (49.3,27.0) -- (51.7,27.0) -- (51.7,29.4) -- cycle;
            \end{scope}
            \draw[thumbborder] (49.3,27.0) rectangle (51.7,29.4);
            \draw[sepline] (49.3,27.0) -- (51.7,29.4);
            \draw[selectA] (50.48,27.16) rectangle (50.82,28.41);

            %----------------------------------------
            % thumbnail 2
            % same slope-1 line through lower-left corner
            % second separator of slope 3 through lower-left corner
            % blue = under slope-3
            % gray = below slope-1
            %----------------------------------------
            \begin{scope}
                \clip (53.3,27.0) rectangle (55.7,29.4);
                \fill[thumbblue] (53.3,27.0) -- (55.7,29.4) -- (54.1,29.4) -- cycle;
                \fill[thumbgray] (53.3,27.0) -- (55.7,27.0) -- (55.7,29.4) -- cycle;
            \end{scope}
            \draw[thumbborder] (53.3,27.0) rectangle (55.7,29.4);
            \draw[sepline] (53.3,27.0) -- (55.7,29.4);
            \draw[sepline] (53.3,27.0) -- (54.1,29.4);
            \draw[selectB] (54.48,27.16) rectangle (54.82,29);

            %----------------------------------------
            % thumbnail 3
            % almost the same background as thumbnail 2
            % apart from the second separator of slope 3, which slightly moves above
            %----------------------------------------
            \begin{scope}
                \clip (57.3,27.0) rectangle (59.7,29.4);
                \fill[thumbblue] (57.3,27.0) -- (59.7,29.4) -- (58.04,29.4) -- (57.3,27.18) -- cycle;
                \fill[thumbgray] (57.3,27.0) -- (59.7,27.0) -- (59.7,29.4) -- cycle;
            \end{scope}
            \draw[thumbborder] (57.3,27.0) rectangle (59.7,29.4);
            \draw[sepline] (57.3,27.0) -- (59.7,29.4);
            \draw[sepline] (57.3,27.18) -- (58.04,29.4);
            \draw[selectC] (58.31,28.16) rectangle (58.65,29.16);

            %----------------------------------------
            % thumbnail 4
            % identical to thumbnail 1
            %----------------------------------------
            \begin{scope}
                \clip (61.3,27.0) rectangle (63.7,29.4);
                \fill[thumbblue] (61.3,27.0) -- (63.7,27.0) -- (63.7,29.4) -- cycle;
            \end{scope}
            \draw[thumbborder] (61.3,27.0) rectangle (63.7,29.4);
            \draw[sepline] (61.3,27.0) -- (63.7,29.4);
            \draw[selectA] (62.48,27.16) rectangle (62.82,28.41);

            %==================================================
            % middle
            %================================================== 

            \node at (56.5,25.2) {\rotatebox{-90}{\scalebox{1.4}[1.2]{$\Longrightarrow$}}};
            \node[right] at (56.9,25.2) {\small zoom-in};

            %==================================================
            % lower half
            %==================================================

            %----------------------------------------
            % left copy
            %----------------------------------------
            \begin{scope}
                % fill first
                \begin{scope}
                    \clip (49.5,11.5) rectangle (51.5,19.5);
                    % below slope-1 line through (50,19)
                    \fill[thumbblue]
                        (49.5,11.5) --
                        (51.5,11.5) --
                        (51.5,19.5) --
                        (50.5,19.5) --
                        (49.5,18.5) -- cycle;
                \end{scope}

                % skeleton after fill
                \draw[skel] (50.5,11.5) -- (50.5,19.5);
                \draw[skel] (49.5,12.5) -- (51.5,12.5);
                \draw[skel] (49.5,13.5) -- (51.5,13.5);
                \draw[skel] (49.5,14.5) -- (51.5,14.5);
                \draw[skel] (49.5,15.5) -- (51.5,15.5);
                \draw[skel] (49.5,16.5) -- (51.5,16.5);
                \draw[skel] (49.5,17.5) -- (51.5,17.5);
                \draw[skel] (49.5,18.5) -- (51.5,18.5);

                \coordinate (X-source) at (51,12);
                \coordinate (X-target) at (50,19);

                \node at (X-source) {\small ${w}_1$};
                \node at (X-target) {\small ${y}_1$};
                \color{purple}
                \draw[diff] (X-source) -- (X-target); %node[midway, left] {\small $d_{r+1}$};

                % border last
                \draw[purple, line width=0.35pt] (49.5,11.5) rectangle (51.5,19.5);
            \end{scope}

            %----------------------------------------
            % middle-left copy
            %----------------------------------------
            \begin{scope}[xshift=4cm]
                % fills first
                \begin{scope}
                    \clip (49.5,11.5) rectangle (51.5,23.5);

                    % whole selected rectangle is under the slope-3 separator in the thumbnail,
                    % so start with blue everywhere
                    \fill[thumbblue] (49.5,11.5) rectangle (51.5,23.5);

                    % then gray below the slope-1 line through the same vertical position
                    % as in the first box, i.e. through (50,19)
                    \fill[thumbgray]
                        (49.5,11.5) --
                        (51.5,11.5) --
                        (51.5,20.5) --
                        (49.5,18.5) -- cycle;
                \end{scope}

                % skeleton after fill
                \draw[skel] (50.5,11.5) -- (50.5,23.5);
                \draw[skel] (49.5,12.5) -- (51.5,12.5);
                \draw[skel] (49.5,13.5) -- (51.5,13.5);
                \draw[skel] (49.5,14.5) -- (51.5,14.5);
                \draw[skel] (49.5,15.5) -- (51.5,15.5);
                \draw[skel] (49.5,16.5) -- (51.5,16.5);
                \draw[skel] (49.5,17.5) -- (51.5,17.5);
                \draw[skel] (49.5,18.5) -- (51.5,18.5);
                \draw[skel] (49.5,19.5) -- (51.5,19.5);
                \draw[skel] (49.5,20.5) -- (51.5,20.5);
                \draw[skel] (49.5,21.5) -- (51.5,21.5);
                \draw[skel] (49.5,22.5) -- (51.5,22.5);

                \coordinate (Y-source) at (51,12);
                \coordinate (Y-target) at (50,23);

                \node at (Y-source) {\small ${w}$};
                \node at (Y-target) {\small ${y}$};
                \color{black}
                \draw[diff1] (Y-source) -- (Y-target); % node[midway, left] {\small $d_{r+1}$};

                % border last
                \draw[black!75, line width=0.35pt] (49.5,11.5) rectangle (51.5,23.5);
            \end{scope}

            %----------------------------------------
            % middle-right copy
            %----------------------------------------
            \begin{scope}[xshift=8cm, yshift=0cm]
                % fills first
                \begin{scope}
                    \clip (49.5,17.5) rectangle (51.5,23.5);

                    % whole selected rectangle is again under the slope-3 separator,
                    % so start with blue everywhere
                    \fill[thumbblue] (49.5,17.5) rectangle (51.5,23.5);

                    % gray below the slope-1 line through \bar x = (51,18)
                    \fill[thumbgray]
                        (50.5,17.5) --
                        (51.5,17.5) --
                        (51.5,18.5) -- cycle;
                \end{scope}

                % skeleton after fill
                \draw[skel] (50.5,17.5) -- (50.5,23.5);
                \draw[skel] (49.5,18.5) -- (51.5,18.5);
                \draw[skel] (49.5,19.5) -- (51.5,19.5);
                \draw[skel] (49.5,20.5) -- (51.5,20.5);
                \draw[skel] (49.5,21.5) -- (51.5,21.5);
                \draw[skel] (49.5,22.5) -- (51.5,22.5);

                \coordinate (Z-source) at (51,18);
                \coordinate (Z-target) at (50,23);

                \node at (Z-source) {\small ${z}_1$};
                \node at (Z-target) {\small ${z}$};
                \color{blue}
                \draw[diff] (Z-source) -- (Z-target); % node[midway, left] {\small $d_{r+1}$};

                % border last
                \draw[blue, line width=0.35pt] (49.5,17.5) rectangle (51.5,23.5);
            \end{scope}

            %----------------------------------------
            % right copy
            %----------------------------------------
            \begin{scope}[xshift=12cm,yshift=-1cm]
                % fill first
                \begin{scope}
                    \clip (49.5,11.5) rectangle (51.5,19.5);
                    \fill[thumbblue]
                        (49.5,11.5) --
                        (51.5,11.5) --
                        (51.5,19.5) --
                        (50.5,19.5) --
                        (49.5,18.5) -- cycle;
                \end{scope}

                % skeleton after fill
                \draw[skel] (50.5,11.5) -- (50.5,19.5);
                \draw[skel] (49.5,12.5) -- (51.5,12.5);
                \draw[skel] (49.5,13.5) -- (51.5,13.5);
                \draw[skel] (49.5,14.5) -- (51.5,14.5);
                \draw[skel] (49.5,15.5) -- (51.5,15.5);
                \draw[skel] (49.5,16.5) -- (51.5,16.5);
                \draw[skel] (49.5,17.5) -- (51.5,17.5);
                \draw[skel] (49.5,18.5) -- (51.5,18.5);

                \coordinate (SX-source) at (51,12);
                \coordinate (SX-target) at (50,19);

                \node at (SX-source) {\small ${w}_1$};
                \node at (SX-target) {\small ${y}_1$};
                \color{purple}
                \draw[diff] (SX-source) -- (SX-target); % node[midway, left] {\small $d_{r+1}$};

                % border last
                \draw[purple, line width=0.35pt] (49.5,11.5) rectangle (51.5,19.5);
            \end{scope}

            %==================================================
            % extension arrows in the lower half
            %==================================================
            \color{blue} \draw[extarr1] (X-target) -- (Y-target);
            \color{gray} \draw[extarr] (X-source) -- (Y-source);
            \color{gray} \draw[extarr] (Y-target) -- (Z-target);
            \color{gray} \draw[extarr] (Z-source) -- (SX-target);

        \end{tikzpicture}
    }
    \caption{A sketch of the proof of Theorem \ref{exotic-transfer-paradigm}.}
    \label{Figure-proof-illustration}
\end{figure}

In fact, Theorem~\ref{exotic-transfer-paradigm} produces families of \textbf{transchromatic differentials} in the \(C_4\)-slice spectral sequence of \(\BPCfour\langle m\rangle\), namely differentials whose source and target lie on opposite sides of one of the lines \(\CL^{C_4}_{\filledstar,1}\). Figure \ref{Figure-proof-illustration} is a schematic illustration of the proof below: Roughly speaking, the generalized Mahowald trick turns a certain {\color{blue} sheared differential $d({z}_1) = {z}$} into a {\color{blue} $\tr$-extension from ${y}_1$ to ${y}$}, and the generalized Leibniz rule blends this with a certain {\color{purple} transfer differential $d({w}_1) = {y}_1$} to produce the desired exotic transfer differential $d({w}) = {y}$.

\begin{proof}
    As $C_4 / C_2$ is the unit sphere of $\sigma \in \RO(C_4)$, the cofiber of $\tr\colon \Sb^0[C_4 / C_2] \to \Sb^0$ in $\Sp^{C_4}$ is the Euler class $a_\sigma\colon \Sb^0 \to \Sb^\sigma$. Taking the Spanier--Whitehead dual, the cofiber of $a_\sigma$ is the composite $\Sb^\sigma \to \Sb^{\sigma}[C_4 / C_2] \cong \Sb^1[C_4 / C_2]$, where the first map is $\res$ and the second isomorphism is $u_\sigma^{-1}$. Thus, 
    \[\Sb^0[C_4 / C_2] \xrightarrow{\tr} \Sb^0 \xrightarrow{a_\sigma} \Sb^{\sigma} \xrightarrow{u_\sigma^{-1} \res} \Sb^1 [C_4 / C_2] \]
    forms a distinguished triangle in $\Sp^{C_4}$. It further induces a distinguished triangle in $\Fil\Sp^{C_4}$
    \begin{equation*} \label{diamond-dt} \tag{$\blacklozenge$}
        Y[C_4 / C_2] \xrightarrow{\tr} Y \xrightarrow{a_\sigma} \Sigma^{\sigma, 0} Y \xrightarrow{u_\sigma^{-1} \res} \Sigma^{1, 0} Y [C_4 / C_2]
    \end{equation*}
    for $Y = P_\bullet \BPCfour\langle m \rangle$, from which the expected differential will be deduced via GLR and GMT. \parr 

    The source of our expected differential sits in bidegree $(V + 2^{h + 1} + 1, 1)$ for $V = (2p + k - 2^{h + 1}) - k \sigma + (p - 2^h)\quarterrep$. A preimage of this class along the map $\tr$ (more precisely, along the map $\tr(-)u_{k \sigma}$) is $w_1 = xu_{2^h \sigma_2} a_{\sigma_2}$, which supports a differential in the $C_2$-slice SS
    \[{\color{purple} d_{2^{h + 1} - 1}(x u_{2^h \sigma_2} a_{\sigma_2}) = x\vbar_h a_{2^{h + 1}\sigma_2} =  (P(\tbar_1 \gamma \tbar_1, \ldots, \tbar_{2^m - 1} \gamma \tbar_{2^m - 1}) + c + \gamma(c))a_{2^{h + 1}\sigma_2}}\] 
    according to Theorem \ref{diffs-in-C2-sliceSS}. Here we get a strict equality on the $E_{2^{h + 1} - 1}$-page since $2 a_{\sigma_2} = 0$ and each $\vbar_i a_{2^{h + 1} \sigma_2}$ ($i < h$) is already wiped out in some earlier page. The target of this {\color{purple} $d_{2^{h + 1} - 1}$} is equal to 
    \[\res(P(\dfbar_1, \ldots, \dfbar_{2^m - 1})u_{(k + p + 1) \sigma} a_{2^h\quarterrep} + \tr(c)u_{(k + 1) \sigma}a_{2^h\quarterrep}).\]
    The classes inside $\res(-)$ sit in bidegree $((2p + k + 1) - (k + 1)\sigma + (p - 2^h) \quarterrep, 2^{h + 1})$, which is on the line $\CL^{C_4}_{V - \sigma + 1, 1}\colon s = (t - s)$. By Theorem \ref{sheared-diffs}, the differential on these classes can be read off in the localized slice SS $E_*^{*, \filledstar}(Y[a_\sigma^{-1}])$. The transfer class is thus a permanent cycle as every transfer class is $a_\sigma$-torsion. On the other hand, as $k + p + 1 = 2^q(2j + 1)$, the norm class supports a differential 
    {\color{blue}\begin{align*} \color{blue}
        d_{2^{q + 2} - 3}(P(\dfbar_1, \ldots, \dfbar_{2^m - 1})u_{(k + p + 1) \sigma} a_{2^h\quarterrep}) &= P(\dfbar_1, \ldots, \dfbar_{2^m - 1}) u_{2^{q + 1}\! j \sigma}  a_{2^h \quarterrep} \cdot d_{2^{q + 2} - 3} (u_{2^q\sigma}) \\ 
        &=  P(\dfbar_1, \ldots, \dfbar_{2^m - 1}) \dfbar_{2^q - 1} u_{2^{q + 1} \! j \sigma} a_{(2^h + 2^q - 1)\quarterrep}a_{(2^{q + 1} - 1)\sigma}
    \end{align*}}
    which is again due to Theorem \ref{sheared-diffs}. Applying the generalized Mahowald trick (i.e. Theorem \ref{GMT-for-SS}) to the $\psi_V$-image of (\ref{diamond-dt}), we can convert this {\color{blue}$ d_{2^{q + 2} - 3}$ differential} into an extension 
    \[{\color{blue} d_{2^{q + 2} - 4}^{\tr(-)u_{k\sigma}, E_{2^{q + 2} - 2}}(\text{target of the above }d_{2^{h + 1} - 1}) = P(\dfbar_1, \ldots, \dfbar_{2^m - 1}) \dfbar_{2^q - 1} u_{2^{q + 1} \! j \sigma} a_{(2^h + 2^q - 1)\quarterrep}a_{(2^{q + 1} - 2)\sigma}}\]
    Here the no-crossing conditions are satisfied as there is no filtration jump along $\res$ and $a_{\sigma}$. Furthermore, the target of the above $d_{2^{h + 1} - 1}$ sits in bidegree $(2p + (2p - 2^{h + 1}\sigma_2), 2^{h + 1})$, which is on the line $\CL^{C_2}_{V_{1}, 1}: s = (t - s)$ for $V_{1} = \Res(V) = (2p - 2^{h + 1}) + (2p - 2^{h + 1})\sigma_2$. On the $E_2$-page there is already no class above this bidegree by Theorem \ref{sheared-diffs}, so we can apply Theorem \ref{stretching-extensions-across-pages} to upgrade the $E_{2^{q + 2} - 2}$-page extension to an extension on the $E_{\infty}$-page. The expected differential now follows from an application of the generalized Leibniz rule (i.e. Theorem \ref{GLR-for-SS}) along the $\psi_V$-image of the map $\tr$ (or more precisely, $\tr(-)u_{k\sigma}$) to the {\color{purple} $d_{2^{h + 1} - 1}$ differential} and the {\color{blue} $E_{\infty}$-page extension with filtration jump $2^{q + 2} - 4$}, in which no-crossing conditions are satisfied since there is no filtration jump on the source and the $E_{\infty}$-page extension on the target has no-crossing for degree reasons. 
\end{proof}

\FloatBarrier

In the following, we list several families arising from Theorem~\ref{exotic-transfer-paradigm} as illustrative examples.

\begin{example} \label{exotic-transfer-diff-omnipresence}
    Take $1 \leq m < \infty$. We first construct a collection of examples that work for each pair $(h, q)$ with $2 \leq h \leq 2m, 1 \leq q \leq m$. Note that for each $z, w \in R_m$, we have $(z + w)\cdot \gamma(z + w) = z \cdot \gamma(z) + w \cdot \gamma(w) + (1 + \gamma)(z \cdot \gamma(w))$ \textup{mod} $2$. Thus for $x = \gamma(\vbar_h)$, there exists $c_h \in R_m$ so that  $\gamma(\vbar_h) \cdot \vbar_h = P_h(\tbar_1 \gamma \tbar_1, \ldots, \tbar_{2^m - 1} \gamma \tbar_{2^m - 1}) + c_h + \gamma(c_h)$ \textup{mod} $I_h$, where
    \[P_h(T_1, \ldots, T_{2^m - 1}) = \begin{cases}
        \sum_{j = 1}^{h - 1} T_{2^{j} - 1} T_{2^{h - j} - 1}^{2^j}, &  2 \leq h \leq m. \\
        \sum_{j = h - m}^{m} T_{2^{j} - 1} T_{2^{h - j} - 1}^{2^j}, & m + 1 \leq h \leq 2m.
    \end{cases}  \]
    In particular, $p = |P_h| = 2^h - 1$. Consequently, for each $j \in \Nb$, 
    \[d_{2^{h + 1} + 2^{q + 2} - 5}(\tr(\vbar_h a_{\sigma_2})u_{2^{h - 1}\quarterrep} u_{(2^{q + 1}\! j + 2^q - 2^h)\sigma}) = P_h(\dfbar_1, \ldots, \dfbar_{2^m - 1}) \dfbar_{2^q - 1} u_{2^{q + 1} \! j \sigma} a_{(2^h + 2^q - 1)\quarterrep}a_{(2^{q + 1} - 2)\sigma}.\]
\end{example}

\begin{example} \label{exotic-transfer-diff-longest}
    There is a refined choice for splicing the longest transfer differential and the longest sheared differential for $\BPCfour \langle m \rangle$: note that in $R_m$
    \[\gamma \tbar_{2^{m} - 1}^{2^m - 1} \cdot \vbar_{2m} = \gamma \tbar_{2^{m} - 1}^{2^m - 1} \cdot \tbar_{2^{m} - 1}^{2^m} \gamma \tbar_{2^m - 1} = \tbar_{2^m - 1}^{2^m} \gamma \tbar_{2^m - 1}^{2^m} \mod I_{2m}.\]
    Thus for $h = 2m$, the tuple $x = \gamma \tbar_{2^m - 1}^{2^m - 1}, P = T_{2^m - 1}^{2^m}, c = 0$ satisfies the hypotheses of Theorem \ref{exotic-transfer-paradigm}, where $p = |P| = 2^{2m} - 2^m$. Taking $q = m$, we deduce that
    \[d_{2^{2m + 1} + 2^{m + 2} - 5} (\tr(\tbar_{2^m - 1}^{2^m - 1} a_{\sigma_2}) u_{2^{2m - 1} \quarterrep} u_{(2^{m + 1}\! j + 2^{m + 1} - 2^{2m} - 1) \sigma}) = \dfbar_{2^m - 1}^{2^m + 1} u_{2^{m + 1} \! j \sigma} a_{(2^{2m} + 2^m - 1)\quarterrep}a_{(2^{m + 1} - 2)\sigma}\]
    for each $j \in \Nb$. Setting $j = 2^{m} - 1$ and multiplying both sides by $\dfbar_{2^m - 1}^{2^{m} + 1} a_{(2^{2m} - 2^m - 1)\quarterrep}$ yields a differential in the integer-graded page of the $C_4$-slice SS
    \begin{align*}
        d_{2^{2m + 1} + 2^{m + 2} - 5} (\dfbar_{2^m - 1}^{2^{m} + 1} \tr(\tbar_{2^m - 1}^{2^m - 1} a_{\sigma_2}) u_{2^{2m - 1} \quarterrep} u_{(2^{2m} - 1) \sigma} & a_{(2^{2m} - 2^m - 1)\quarterrep}) \\ 
        = \ &\dfbar_{2^m - 1}^{2^{m + 1} + 2} u_{(2^{2 m + 1} - 2^{m + 1}) \sigma} a_{(2^{2m + 1} - 2)\quarterrep} a_{(2^{m + 1} - 2)\sigma}
    \end{align*}
    whose source lies in bidegree $(2^{2m + 2} - 2^{m + 1} - 1, 2^{2m + 1} - 2^{m + 1} - 1)$. 
    \begin{itemize}
        \item When $m = 1$, the family of differentials above specializes to
        \[d_{11}(\tr(\tbar_1 a_{\sigma_2})u_{2 \quarterrep} u_{(4j - 1)\sigma}) = \dfbar_{1}^{3} u_{4j \sigma} a_{5 \quarterrep} a_{2\sigma}, \ \ j \in \Nb\]
        which generates all essential $d_{11}$ differentials in the (positive cone of the) $C_4$-slice SS of $\BPCfour \langle 1 \rangle$ under the Leibniz rule, and the integer-graded differential becomes 
        \[d_{11}(\dfbar_{1}^{3} \tr(\tbar_1 a_{\sigma_2}) u_{2 \quarterrep} u_{3\sigma} a_{\quarterrep}) = \dfbar_{1}^{6} u_{4 \sigma} a_{6 \quarterrep} a_{2\sigma}\]
        out of bidegree $(11, 3)$, which is the first essential $d_{11}$ in the (positive cone of the) integer-graded $C_4$-slice SS, cf. \cite[Theorem 3.18 and Figure 9]{HSWX-BPC4<2>}. 
        \item When $m = 2$, the family of differentials above specializes to
        \[d_{43}(\tr(\tbar_3^3 a_{\sigma_2})u_{8 \quarterrep} u_{(8j - 9)\sigma}) = \dfbar_{3}^{5} u_{8j \sigma} a_{19 \quarterrep} a_{6\sigma}, \ \ j \in \Nb\]
        which generates all essential $d_{43}$ differentials in the (positive cone of the) $C_4$-slice SS of $\BPCfour \langle 2 \rangle$ under the Leibniz rule, and the integer-graded differential becomes 
        \[d_{43}(\dfbar_{3}^{5} \tr(\tbar_3^3 a_{\sigma_2}) u_{8 \quarterrep} u_{15\sigma} a_{11\quarterrep}) = \dfbar_{3}^{10} u_{24 \sigma} a_{30 \quarterrep} a_{6\sigma}\]
        out of bidegree $(55, 23)$, which is the first essential $d_{43}$ in the (positive cone of the) integer-graded $C_4$-slice SS, cf. \cite[Theorem 9.5 and Figure 35]{HSWX-BPC4<2>}. 
    \end{itemize}
    In both known cases, this is the family of transchromatic differentials whose targets lie highest above the ``transchromatic line'' $\CL^{C_4}_{\filledstar, 1}$. Indeed, this assertion will be true for all $m \geq 1$ for degree reasons as long as we know such differentials are essential.
\end{example}

\begin{example}
    Not all families of differentials produced by Theorem \ref{exotic-transfer-paradigm} are guaranteed to be essential. Take $2 \leq m < \infty$ and consider the same identity as in Example \ref{exotic-transfer-diff-longest} above
    \[\gamma \tbar_{2^{m} - 1}^{2^m - 1} \cdot \vbar_{2m} = \tbar_{2^m - 1}^{2^m} \gamma \tbar_{2^m - 1}^{2^m} \mod I_{2m}.\]
    Setting $h = 2m$, $x = \gamma \tbar_{2^m - 1}^{2^m - 1}$, $P = T_{2^m - 1}^{2^m}$, $c = 0$, and $q = m - 1$, we also deduce a family of differentials 
    \[d_{2^{2m + 1} + 2^{m + 1} - 5} (\tr(\tbar_{2^m - 1}^{2^m - 1} a_{\sigma_2}) u_{2^{2m - 1} \quarterrep} u_{(2^{m}\! j + 2^{m + 1} - 2^{2m} - 2^{m - 1} - 1) \sigma}) = \dfbar_{2^{m - 1} - 1} \dfbar_{2^m - 1}^{2^m} u_{2^{m} \! j \sigma} a_{(2^{2m} + 2^{m - 1} - 1)\quarterrep}a_{(2^{m} - 2)\sigma}\]
    for $j \in \Nb$. These are usually not essential differentials, for instance:
    \begin{itemize}
        \item When $m = 2$, the family specializes to 
        \[d_{35}(\tr(\tbar_{3}^{3} a_{\sigma_2}) u_{8 \quarterrep} u_{(4j - 11) \sigma}) = \dfbar_{1} \dfbar_{3}^{4} u_{4j \sigma} a_{17\quarterrep}a_{2\sigma}\]
        for $j \in \Nb$. None of these is essential because of the family 
        \[d_{27}(\tr(\tbar_1 a_{\sigma_2})u_{6 \quarterrep} u_{(4j - 5) \sigma}) = \dfbar_1 \dfbar_3^2 u_{4j \sigma} a_{13\quarterrep} a_{2\sigma}\] 
        in the $C_4$-slice SS of $\BPCfour \langle 2 \rangle$, cf. \cite[Proposition 10.5]{HSWX-BPC4<2>}.
    \end{itemize}
\end{example}

\begin{example} \label{exotic-transfer-diff-second-longest}
    For some $m \geq 2$, there exists a homogeneous class $c_0 \in R_m$ such that
    \[\tbar_1 \vbar_{2 m - 1} = \tbar_{2^{m - 1} - 1}^{2^{m - 1} - 1} \gamma \tbar_{2^{m - 1} - 1}^{2^{m - 1} - 1} \cdot \tbar_{2^m - 1} \gamma \tbar_{2^m - 1} + c_0 + \gamma(c_0) \mod I_{2m - 1}.\]
    This is the case for $m = 2$ due to \cite[Theorem 2.11]{HSWX-BPC4<2>}, and also for $m = 3, 4$ by a brute force Gröbner basis computation (cf. the scripts in \cite{WuZenodoC2Poly2026}). We expect such $c_0$ to exist for all $m \geq 2$. \parr 

    In the case such $c_0$ exists, the collection $x = \tbar_1$, $P = T_{2^{m - 1} - 1}^{2^{m - 1} - 1} T_{2^{m} - 1}$, $c = c_0$ satisfies the hypotheses of Theorem \ref{exotic-transfer-paradigm}, here $p = |P| = 2^{2m - 2}$. Taking $q = m - 1$, we deduce that 
    \[d_{2^{2m} + 2^{m + 1} - 5} (\tr(\tbar_1 a_{\sigma_2}) u_{2^{2m - 2} \quarterrep} u_{(2^{m}\! j + 2^{m - 1} - 2^{2m - 2} - 1) \sigma}) = \dfbar_{2^{m - 1} - 1}^{2^{m - 1}} \dfbar_{2^{m} - 1} u_{2^{m} \! j \sigma} a_{(2^{2m - 1} + 2^{m - 1} - 1)\quarterrep}a_{(2^{m} - 2)\sigma}\]
    for each $j \in \Nb$. Setting $j = 2^{m} - 2^{m - 2}$ and multiplying both sides by $\dfbar_{2^{m} - 1}^{2^{m - 1} + 1} a_{(2^{2 m - 2} + 2^{m - 1} - 1)\quarterrep}$ yields a differential in the integer-graded page of the $C_4$-slice SS 
    \begin{align*}
        d_{2^{2m} + 2^{m + 1} - 5} (\dfbar_{2^m - 1}^{2^{m - 1} + 1} \tr(\tbar_{1} a_{\sigma_2}) &u_{2^{2m - 2} \quarterrep} u_{(2^{2m - 1} + 2^{m - 1} - 1)\sigma} a_{(2^{2m - 2} + 2^{m - 1} - 1)\quarterrep}) \\ 
        = \ &\dfbar_{2^{m - 1} - 1}^{2^{m - 1}} \dfbar_{2^{m} - 1}^{2^{m - 1} + 2}  u_{(2^{2 m} - 2^{2m - 2}) \sigma} a_{(2^{2m} - 2^{2 m - 2} + 2^{m} - 2)\quarterrep} a_{(2^{m} - 2)\sigma}
    \end{align*}
    whose source lies in bidegree $(2^{2m + 1} - 2^{2m - 1} + 2^m - 1, 2^{2m - 1} + 2^m - 1)$. 
    \begin{itemize}
        \item  When $m = 2$, the family of differentials above specializes to
        \[d_{19}(\tr(\tbar_1 a_{\sigma_2})u_{4 \quarterrep} u_{(4j - 1)\sigma}) = \dfbar_{1}^{2} \dfbar_{3} u_{4j \sigma} a_{9 \quarterrep} a_{2\sigma}, \ \ j \in \Nb\]
        which generates all essential transchromatic $d_{19}$ differentials in the (positive cone of the) $C_4$-slice SS of $\BPCfour \langle 2 \rangle$ under the Leibniz rule, and the integer-graded differential becomes 
        \[d_{19}(\dfbar_{3}^{3} \tr(\tbar_1 a_{\sigma_2}) u_{4 \quarterrep} u_{9\sigma} a_{5\quarterrep}) = \dfbar_{1}^{2} \dfbar_{3}^{4} u_{12 \sigma} a_{14 \quarterrep} a_{2\sigma}\]
        out of bidegree $(27, 11)$, which is the first essential transchromatic $d_{19}$ in the (positive cone of the) integer-graded $C_4$-slice SS, cf. \cite[Theorem 10.5 and Figure 31]{HSWX-BPC4<2>}. 
    \end{itemize}
\end{example}

\subsection{Epilogue: classifying slice differentials for \texorpdfstring{$\BPCfour\langle 1 \rangle$}{BPC4<1>}}
\label{subsec:5.6}

In general, the slice SS of each $\BPG \langle m \rangle$ collapses after finitely many pages. Also, there are various periodicities in these slice SS, cf. \cite{DHLLSWX-periodicity}. Together, these imply that in each slice SS there are finitely many \textbf{generating differentials} that give rise to all differentials under the Leibniz rule. \parr

\begin{figure}[htbp]
\begin{center}
\makebox[\textwidth]{
  \includegraphics[
    page=1,
    clip,
    trim={13cm 13.5cm 0cm 0cm},
    scale=1.2
  ]{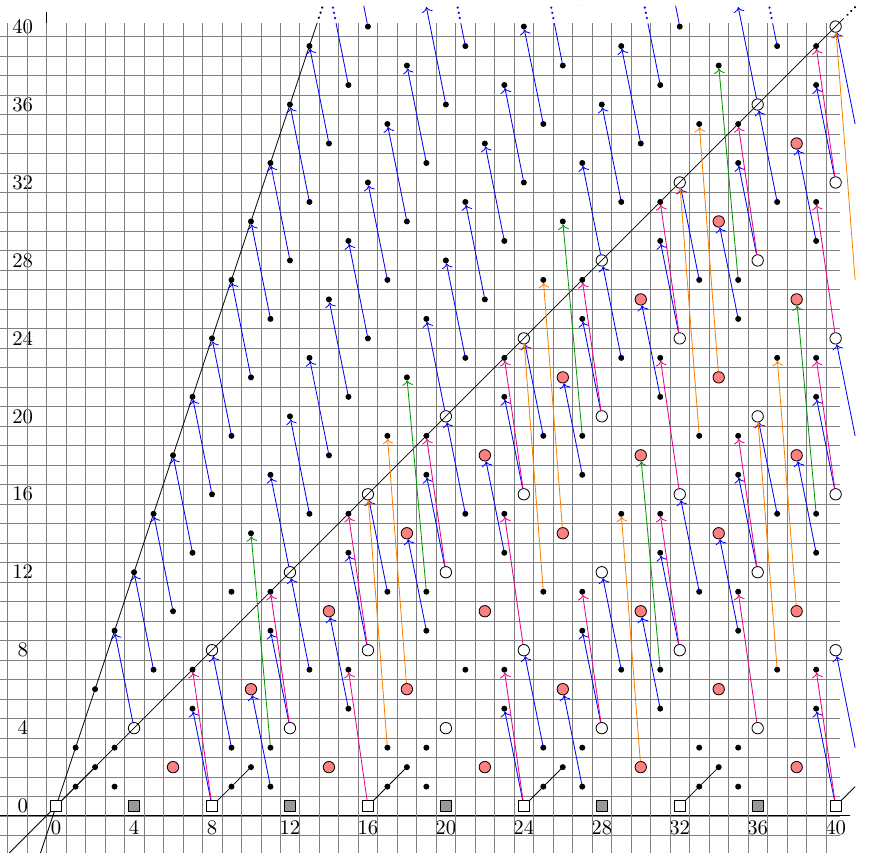}
}
\caption{The $C_4$-slice spectral sequence of $\BPCfour \langle 1 \rangle$, starting from its $E_5$-page.}
\hfill
\label{figure-BPC4<1>}
\end{center}
\end{figure}

Table \ref{table-generating-diffs-in-BPC4<1>} (whose first column is from \cite[\S~3]{HSWX-BPC4<2>}) gives a full list of generating differentials in the $C_4$-slice SS of $Y = \BPCfour \langle 1 \rangle$ with short indications for the proofs of their existence:

    \begin{table}[htbp]
        \centering
        \renewcommand{\arraystretch}{1.35}
        \begin{tabular}{||c||c||}
            \hline 
            Generating differential & Reason \\
            \hline
            \hline
             $d_3(u_{\quarterrep})$ & restriction  \\
            \hline
             $\color{blue} d_5(u_{2\sigma})$ & transchromatic isomorphism \\
            \hline
             {$\color{blue} d_5(u_{2 \quarterrep})$} &  {
             restriction}  \\
            \hline
             $\color{magenta} d_7(2u_{2 \quarterrep})$ & transfer \\
            \hline
             $\color{magenta} d_7(u_{4 \quarterrep})$ & total Leibniz rule\\
            \hline
             {
             $\color{green!60!black} d_{11}(\tr(\tbar_1 a_{\sigma_2}) u_{2\quarterrep} u_{3\sigma})$} &  {
             exotic transfer}  \\
            \hline
             $\color{orange} d_{13}(u_{4\quarterrep} a_\sigma)$ & HHR norm \\ 
            \hline
        \end{tabular}
        \caption{Generating differentials in the $C_4$-slice SS of $\BPCfour\langle 1 \rangle$.}
        \label{table-generating-diffs-in-BPC4<1>}
    \end{table}    

    \FloatBarrier

In Table \ref{table-generating-diffs-in-BPC4<1>}, ``restriction'' stands for a direct comparison with the $C_2$-slice SS along the map $\res\colon Y \to Y[C_4 / C_2]$, ``transfer'' stands for Remark \ref{transfer-diffs} (i.e. a direct comparison with the $C_2$-slice SS along the map $\tr\colon Y[C_4 / C_2] \to Y$), ``transchromatic isomorphism'' stands for Theorem \ref{sheared-diffs}, ``exotic transfer'' stands for Theorem \ref{exotic-transfer-paradigm} (more precisely, here it comes from Example \ref{exotic-transfer-diff-longest}), and ``HHR norm'' stands for the norm differential theorem (cf. \cite[Theorem 4.7]{HHR-BPC4<1>}). The remaining one, $d_7(u_{4\quarterrep})$, can be proved via the Leibniz rule for total differentials (i.e. Theorem \ref{Burklund's-Leibniz-rule}) as follows:

\begin{proposition}
    $d_7(u_{4\quarterrep}) = \dfbar_1 \tr(\tbar_1 a_{\sigma_2}) u_{2\quarterrep} u_\sigma a_{3\quarterrep}$.
\end{proposition}
\begin{proof}
    Write $X = P_\bullet \MUCfour$ and $Y = P_\bullet \BPCfour \langle 1 \rangle$. As $u_{2\quarterrep}$ is a $d_4$-cycle, there exists $[u_{2\quarterrep}]' \in \pi_{4 - 2 \quarterrep, 0}(Y / \defopara^4)$ lifting $u_{2\quarterrep}$, which further lifts\footnote{Strictly speaking, we should also check that $u_{2\quarterrep}$ is a $d_4$-cycle in the $C_4$-slice SS of $\MUCfour$. But this is automatic: For degree reasons the only other possibility is that $u_{2 \quarterrep}$ supports an essential $d_3$ here, which does not happen because the potential target bidegree $(3 - 2 \quarterrep, 3)$ of the $C_4$-slice SS of $\MUCfour$ is isomorphic to that of $\BPCfour \langle 1 \rangle$. In fact, both of them are $\Fb_2$-vector spaces of dimension $1$ generated by the class $\tr(\tbar_1 u_{2 \sigma_2} a_{3 \sigma_2})$.} to a bigraded stem $[u_{2\quarterrep}] \in \pi_{4 - 2 \quarterrep, 0}(X / \defopara^4)$. Applying the total Leibniz rule (i.e. Theorem \ref{Burklund's-Leibniz-rule}) to the map in $\Fil\Sp$
    \[\psi_0(\Sigma^{-V, 0} X) \otimes \psi_0(\Sigma^{-V, 0} X) \to \psi_0(\Sigma^{-2V, 0} X)\]
    induced from the $\Eb_{\infty}$-ring structure on $X$ and the lax symmetric monoidality of $\psi_0 = \Msp^{\Fil}(\Sb^{0, 0}, -)$ (here $V = 4 - 2 \quarterrep$), we compute that 
    \[\delta_4^4([u_{2\quarterrep}]^2) = 2 [u_{2\quarterrep}] \delta_4^4([u_{2\quarterrep}]) = [u_{2\quarterrep}] \delta_4^4(2[u_{2\quarterrep}])\]
    in the slice SS of $X$. Therefore, by projecting along the map $X \to Y$, we find a bigraded stem $[u_{4\quarterrep}] \in \pi_{8 - 4 \quarterrep, 0}(Y / \defopara^4)$ lifting $u_{4 \quarterrep}$ whose $\delta_4^4$ is equal to $[u_{2\quarterrep}] \delta_4^4(2[u_{2\quarterrep}]')$. Note that there is a transfer differential $d_7(2u_{2\quarterrep}) = \dfbar_1 \tr(\tbar_1 a_{\sigma_2})u_\sigma a_{3\quarterrep}$, for which the only possible crossing is a $d_3$-differential into the same target bidegree $(3 - 2 \quarterrep, 7)$. As this transfer $d_7$ has no crossing on the $E_5$-page, Corollary \ref{free-choice-lemma-for-diff-with-full-no-crossing} implies that $\delta_4^{4}(2[u_{2\quarterrep}]') = \defopara^2 [y]$ where $[y] \in \pi_{3 - 2 \quarterrep, 6}(Y / \defopara^4)$ is a lift of $y = \dfbar_1 \tr(\tbar_1 a_{\sigma_2})u_\sigma a_{3\quarterrep}$. Thus, $\delta_4^4([u_{4\quarterrep}]) = \defopara^2[u_{2\quarterrep}][y]$, so $d_7(u_{4\quarterrep}) = u_{2\quarterrep} y $ due to Theorem \ref{delta-as-total-diff}.
\end{proof}

In summary, using Burklund--Lin--Wang--Xu methods, we obtain a satisfactory interpretation for the $C_4$-slice SS computation of $\BPCfour\langle 1 \rangle$: each generating differential in the $C_4$-slice SS can be deduced from a systematic construction that also produces new families of differentials for all $m \geq 1$ (and potentially all groups $G$), which replaces all ad-hoc case-by-case arguments in previous computations. In ongoing work, we will try to develop a similar understanding for the $C_4$-slice SS computation of $\BPCfour\langle 2 \rangle$ using Burklund--Lin--Wang--Xu methods, together with the norm structure on the slice filtration. This will be our second step towards a systematic computation of the slice SS of $\BPG\langle m \rangle$.

%% file: SynASS.tex
Let $E$ be a homotopy ring spectrum of Adams type (cf. \cite[Definition 3.14]{Pst23}). Write $\Syn_{E}$ for the $\infty$-category of $E$-synthetic spectra, and write $\nu\colon \Sp \to \Syn_E$ for the synthetic analogue functor. In this appendix, for any reasonable spectrum $X$, we identify the synthetic $\nu E$-Adams SS of $\nu X$ and $\nu X / \defopara^r$ (as well as their abutments) with trigraded versions of the $\defopara$-Bockstein SS, cf. Theorem \ref{synthetic-Adams-vs-Bockstein}. We thus obtain a complete description of these synthetic $\nu E$-Adams SS in terms of the $E$-Adams SS of $X$ due to ``rigidity'' of the trigraded $\defopara$-Bockstein SS, cf. Theorem \ref{rigidity-for-trigraded-defopara-Bockstein-SS} and Remark \ref{rigidity-for-filtered-SS-by-definition}. 

\begin{remark}
    The results of this appendix are used only in Remark \ref{comparison-with-extension-SS} to compare Definition \ref{extensions-def} with the extension spectral sequence approach in \cite{Lin-Wang-Xu-kervaire}; the remainder of the paper is independent of Appendix \ref{app:A}.
\end{remark}

\begin{definition}
    We write $\biFil\Sp = \Fil(\Fil\Sp)$ for the $\infty$-category of \textbf{bifiltered spectra}. 
    \begin{itemize}
        \item The $\infty$-category $\biFil\Sp$ is presentable and stable. 
        \item We equip $\biFil\Sp = \Fun(\overrightarrow{\Zb}^{\op}, \Fil\Sp)$ with the Day convolution symmetric monoidal structure, which makes it presentably symmetric monoidal. This is the same as the Day convolution symmetric monoidal structure on $\Fun(\overrightarrow{\Zb}^{\op} \times \overrightarrow{\Zb}^{\op}, \Sp)$. We denote by $\oneb$ its tensor unit.
        \item There are three families of invertible objects in $\biFil\Sp$: 
        \begin{itemize}
            \item By construction, the spectral Yoneda embedding functor 
            \[\gamma\colon \overrightarrow{\Zb} \times \overrightarrow{\Zb} \to \biFil\Sp, \quad (a, b) \mapsto \gamma_{(a, b)} = \Sigma^{\infty}_+ \Map(-, (a, b))\] 
            is symmetric monoidal, thus each $\gamma_{(a, b)}$ is invertible in $\biFil\Sp$ with inverse $\gamma_{(-a, -b)}$. We write $\Sb^{0, a, b} = \gamma_{(a, b)}$. Concretely, the image of $(x, y) \in \overrightarrow{\Zb}^{\op} \times \overrightarrow{\Zb}^{\op}$ under $\Sb^{0, a, b}$ is $\Sb^0$ if $x \leq a, y \leq b$ or $0$ otherwise. In particular, $\Sb^{0, 0, 0} = \oneb$. 
            \item On the other hand, each suspension $\Sigma^{n} \oneb$ is also invertible. We write $\Sb^{n, 0, 0} = \Sigma^{n} \oneb$. 
        \end{itemize}
        In general, we write $\Sb^{n, a, b} = \Sb^{n, 0, 0} \otimes \Sb^{0, a, b}$, and we write $\Sigma^{n, a, b} X = X \otimes \Sb^{n, a, b}$. 
        \item We write $\defopara^r_{h} \colon \Sb^{0, -r, 0} \to \Sb^{0, 0, 0}$ (resp. $\defopara^r_v \colon  \Sb^{0, 0, -r} \to \Sb^{0, 0, 0}$) for the image of $(-r, 0) \to (0, 0)$ (resp. $(0, -r) \to (0, 0)$) under $\gamma$, and refer to it as the \textbf{horizontal} (resp. \textbf{vertical}) \textbf{deformation parameter} in $\biFil\Sp$. For general $Y \in \biFil\Sp$, we also write $\defopara^r_h$ for $Y \otimes \defopara^r_h \colon \Sigma^{0, -r, 0} Y \to Y$, and similarly for $\defopara^r_v$.  It follows from the discussion in Construction \ref{accelerations} that both $\oneb / \defopara^r_h$ and $\oneb / \defopara^r_v$ acquire the structure of an $\Eb_{\infty}$-algebra in $\biFil\Sp$. 
        \item Also, $\oneb[\defopara_h^{-1}]$ and $\oneb[\defopara_v^{-1}]$ are both idempotent $\Eb_{\infty}$-algebras in $\biFil\Sp$, so that $\Mod_{\oneb[\defopara_h^{-1}]}\biFil\Sp$ and $\Mod_{\oneb[\defopara_v^{-1}]}\biFil\Sp$ are both equivalent to $\Fil\Sp$. The functor $-[\defopara_h^{-1}]\colon \biFil\Sp \to \Fil\Sp$ sends $\Sb^{n, a, b}$ to $\Sb^{n, b}$, $\defopara_h$ to $\id$ and $\defopara_v$ to $\defopara$. Similar results hold true for the functor $-[\defopara_v^{-1}]$. 
        \item In practice, we write $Y = \{Y_w \}_{w \in \Zb}$ for an arbitrary object in $\biFil\Sp = \Fil(\Fil\Sp)$. Here we treat $w$ as the vertical coordinate, while each $Y_w \in \Fil\Sp$ occupies a horizontal row. Under this convention, $\defopara_h$ is the map $\{\defopara\colon \Sigma^{0, -1}Y_w \to Y_w\}_{w \in \Zb}$, while $\defopara_v$ is the map $\{Y_{w + 1} \to Y_w\}_{w \in \Zb}$. 
    \end{itemize}
\end{definition}

\begin{construction} \label{trigraded-defopara-Bockstein-SS}
    We construct below the trigraded $\defopara$-Bockstein SS. 
    \begin{itemize}
        \item We denote by $\Fil(\SpSeq) = \Fun(\overrightarrow{\Zb}^{\op}, \SpSeq)$ the $1$-category of \textbf{filtered SS}. Concretely, a filtered SS $E = \{E_r^{s, t, w}\}$ consists of the following data:
        \begin{itemize}
            \item A family of abelian groups $\{E_r^{s, t, w}\}_{r \geq 2, s, t, w \in \Zb}$.
            \item A family of maps $d_r\colon E_r^{s, t, w} \to E_r^{s + r, t + r - 1, w}$, together with isomorphisms to witness that $\{(E_r^{s, t, w}, d_r)\}_{r \geq 2, s, t \in \Zb}$ is a SS in the sense of Definition \ref{SS} for each fixed $w$. 
            \item A family of maps $\defopara\colon E_r^{s, t, w} \to  E_r^{s, t, w - 1}$ such that $d_r \defopara = \defopara d_r$.
        \end{itemize}
       By combining Theorem \ref{categorified-total-Leibniz-rule} and \cite[Example 2.2.6.9]{HA}, we equip $\Fil(\SpSeq)$ with the structure of a $1$-truncated $\infty$-operad under Day convolution.
        \item The functor $E_*^{*,*}\colon \Fil\Sp \to \SpSeq$ in Theorem \ref{standard-SS} induces a \textbf{trigraded standard SS} functor 
        \[E_*^{*,*,*}\colon \biFil\Sp = \Fil(\Fil\Sp) \to \Fil(\SpSeq), \quad Y = \{Y_w\}_{w \in \Zb} \mapsto \{E_r^{s, t, w}(Y) = E_r^{s, t}(Y_w)\}_{r \geq 2, s, t, w \in \Zb}.\] 
        Theorem \ref{categorified-total-Leibniz-rule} implies $E_*^{*,*,*}$ is also lax symmetric monoidal under Day convolutions. In particular, a multiplicative bifiltration gives rise to a multiplicative trigraded standard SS.
        \item For $X \in \Fil\Sp$ and $w \in \Zb$, we write $\kappa_{\geq w} X \in \Fil\Sp$ for the filtration 
        \[\cdots \to X(w + 2) \to X(w + 1) \to X(w) \xrightarrow{\id} X(w) \xrightarrow{\id} \cdots \]
        where the first $X(w)$ sits in filtration $w$. The assignments $X \mapsto \{\kappa_{\geq w} X\}_{w \in \Zb}$ assemble into a functor $\kappa \colon \Fil\Sp \to \biFil\Sp$. Formally, write $\Delta\colon \overrightarrow{\Zb} \hookrightarrow \overrightarrow{\Zb} \times \overrightarrow{\Zb}$ for the diagonal inclusion, then $\kappa\colon \Fil\Sp = \Fun(\overrightarrow{\Zb}^{\op}, \Sp) \to \Fun(\overrightarrow{\Zb}^{\op} \times \overrightarrow{\Zb}^{\op}, \Sp) = \biFil\Sp$ computes the left Kan extension along $\Delta$. Since $\Delta$ is a map of commutative monoids, $\kappa$ is symmetric monoidal under Day convolutions. Additionally, for each $X \in \Fil\Sp$, $\kappa X[\defopara_h^{-1}]$ and $\kappa X[\defopara_v^{-1}]$ can both be identified with $X$, while $(\kappa X / \defopara_v^{r})[\defopara_h^{-1}] \cong X / \defopara^r$. We call $\kappa X$ the \textbf{$\defopara$-Bockstein bifiltration} of $X$.
        \item For $X \in \Fil\Sp$, we refer to $E_*^{*,*,*}(\kappa X)$ as the \textbf{trigraded $\defopara$-Bockstein SS} of $X$. Also, for $1 \leq r < \infty$, we refer to $E_*^{*,*,*}(\kappa X / \defopara_v^r)$ as the \textbf{(refined) trigraded $\defopara$-Bockstein SS} of $X / \defopara^r$. The two functors $E_*^{*,*,*}(\kappa(-)), E_*^{*,*,*}(\kappa(-)/ \defopara_v^r)\colon \Fil\Sp \to \Fil(\SpSeq)$ are both lax symmetric monoidal, so multiplicative filtrations lead to multiplicative (refined) trigraded $\defopara$-Bockstein SS. 
    \end{itemize}
\end{construction}

\begin{remark} \label{diagonal-bifiltrations}
    The functor $\kappa\colon \Fil\Sp \to \biFil\Sp$ admits a right adjoint $\kappa^R$ given by restriction along $\Delta\colon \overrightarrow{\Zb} \hookrightarrow \overrightarrow{\Zb} \times \overrightarrow{\Zb}$. The unit $\id \to \kappa^R \kappa$ is an isomorphism, so $\kappa$ is fully faithful, and it induces an equivalence
    $\Fil\Sp \cong \biFil\Sp_{\mathrm{de}}$ where $\biFil\Sp_{\mathrm{de}} \subset \biFil\Sp$ consists of bifiltrations $Y = \{Y_w\}_{w \in \Zb}$ such that for each $w \in \Zb, k \in \Nb$ the maps $Y_w(w) \to Y_w(w - k)$ and $Y_w(w) \to Y_{w - k}(w)$ are isomorphisms\footnote{In other words, such bifiltrations are ``diagonally extended'', i.e. left Kan extended from the diagonal filtration}. 
    According to Remark \ref{recollements-and-Nakayama-lemma} part 2, if $Y = \{Y_w\}_{w \in \Zb}$ is \textbf{horizontally complete} (i.e. each $Y_w$ is complete), then $Y \in \biFil\Sp_{\mathrm{de}}$ iff $E_2^{s, t, w}(Y) = 0$ for $w > t$ and $\defopara\colon E_2^{s, t, w}(Y) \to E_2^{s, t, w - 1}(Y)$ is an isomorphism for $w \leq t$. 
\end{remark}

\begin{remark}
    Some literature reserves the name ``trigraded $\defopara$-Bockstein SS'' for the trigraded SS $E' = \{(E')_r^{s, t, w}\}$ whose weight-$w$ part comes from $\Msp(\Sb^{0, w}, -)$ image of 
    \[\cdots \xrightarrow{\defopara} \Sigma^{0, -3} X \xrightarrow{\defopara} \Sigma^{0, -2} X \xrightarrow{\defopara} \Sigma^{0, -1} X \xrightarrow{\defopara} X \xrightarrow{\id} X \xrightarrow{\id} \cdots  \in \Fil(\Fil\Sp).\]
    In other words, $(E')_r^{s, t, w} = E_r^{s, t}(\Sigma^{0, -w} \kappa_{\geq w} X)$. This is not very different from the construction above, since for each fixed weight $w$ there is an isomorphism of SS $(E')_r^{s, t, w} \cong E_r^{s + w, t + w, w}(\kappa X)$. We prefer our construction for the following reasons: 
    \begin{itemize}
        \item The natural maps $(E')_r^{s, t, w} \to (E')_r^{s, t, w - 1}$ induced by $\Sigma^{0, -w} \kappa_{\geq w} X \to \Sigma^{0, 1-w} \kappa_{\geq w - 1} X$ are zero for each $r$. To recover the $\defopara$-multiplication in $E_r^{s, t, w}(\kappa X)$, one has to set up another operator ``$\defopara$'' on $E'$ with tridegree $(1, 1, -1)$, which leads to unnecessary complications. 
        \item In Theorem \ref{synthetic-Adams-vs-Bockstein} we identify synthetic Adams SS with the trigraded Bockstein SS as in Construction \ref{trigraded-defopara-Bockstein-SS}. If we use the $E'$-convention instead, they are no longer isomorphic on the nose, but only after an appropriate reindexing.
    \end{itemize}
\end{remark}

\begin{theorem}[Rigidity] \label{rigidity-for-trigraded-defopara-Bockstein-SS}
    Take $X \in \Fil\Sp$. The trigraded $\defopara$-Bockstein SS $E_*^{*,*,*}(\kappa X)$ and $E_*^{*,*,*}(\kappa X / \defopara_v^n)$ are completely determined by the standard SS $E_*^{*,*} = E_*^{*,*}(X)$ (together with $Z_*^{*,*}, B_*^{*,*}$ in Remark \ref{cycles-and-boundaries}).
    \begin{enumerate}
        \item $E_2^{*,*,*}(\kappa X) \cong E_2^{*,*} \otimes_{\Zb} \Zb[\defopara]$, where $x \in E_2^{s, t}$ sits in tridegree $(s, t, t)$ and $|\defopara| = (0, 0, -1)$. More generally, for $1 \leq r \leq \infty$, 
        \[E_{r + 1}^{s, t, w}(\kappa X) \cong \begin{cases}
            Z_r^{s, t} / B_{\min\{t - w + 1, r\}}^{s, t}, & t \geq w. \\
            0, & t < w.
        \end{cases}\]
        For differentials, there is a differential $d_{r + 1}(x) = y$ in the standard SS iff there is a differential $d_{r + 1}(x) = \defopara^r y$ in the trigraded $\defopara$-Bockstein SS. Moreover, in this case we also have $d_{r + 1}(\defopara^k x) = \defopara^{r + k} y$ for $k \geq 0$, and all nontrivial differentials in the trigraded $\defopara$-Bockstein SS are of this form. 
        \item If $X$ is complete, then each $\kappa_{\geq w} X$ is complete. If furthermore all weak obstructions $RE_{\infty}^{s, t}(X) = \lim^1_r Z_r^{s, t}$ vanish, then $\{E_r^{s, t, w}(\kappa X)\}$ converges strongly to $\pi_{t - s, w}(X)$ equipped with the $\IIm(\defopara^a)$ filtration in Theorem \ref{Bockstein-dictionary-infinite}, where $\IIm(\defopara^a) \subset \pi_{m, w}(X)$ sits in filtration $w + a$. 
        \item $E_2^{*,*,*}(\kappa X / \defopara_v^n) \cong E_2^{*,*} \otimes_{\Zb} \Zb[\defopara] / \defopara^n$ where $x \in E_2^{s, t}$ sits in tridegree $(s, t, t)$ and $|\defopara| = (0, 0, -1)$. More generally, for $1 \leq r \leq \infty$,
        \[E_{r + 1}^{s, t, w}(\kappa X / \defopara_v^n) \cong \begin{cases}
            0, & t \geq w + n. \\
            Z_{\min\{w + n - t, r\}}^{s, t} / B_{\min\{t - w + 1, r\}}^{s, t}, & w \leq t < w + n. \\
            0, & t < w.
        \end{cases}\]
        For differentials, there is a differential $d_{r + 1}(x) = y$ in the standard SS with $r < n$ iff there is a differential $d_{r + 1}(x) = \defopara^r y$ in the refined trigraded $\defopara$-Bockstein SS. Moreover, in this case we also have $d_{r + 1}(\defopara^k x) = \defopara^{r + k} y$ for $0 \leq k < n - r$, and all nontrivial differentials in the refined trigraded $\defopara$-Bockstein SS are of this form. In particular, there is no nontrivial differential starting from the $E_{n + 1}$-page. Also, $\{E_{r + 1}^{s, t, w}(\kappa X / \defopara_v^n)\}$ always converges strongly to $\pi_{t - s, w}(X / \defopara^n)$ equipped with the $\IIm(\defopara^a)$-filtration in Theorem \ref{Bockstein-dictionary-finite}, where $\IIm(\defopara^a) \subset \pi_{m, w}(X / \defopara^n)$ sits in filtration $w + a$. 
    \end{enumerate}
\end{theorem}

\begin{proof}
    For part 1, note that $(\kappa_{\geq w} X)(t) = X(\max\{t, w\})$. For fixed $w$, $E_{r + 1}^{s, t, w}(\kappa X) \cong E_{r + 1}^{s, t}(\kappa_{\geq w} X) \cong Z_r^{s, t}(\kappa_{\geq w} X) / B_r^{s, t}(\kappa_{\geq w} X)$. Expanding the construction in Theorem \ref{standard-SS}, $Z_r^{s, t}(\kappa_{\geq w} X)$ is the image of $\rho\colon \pi_{t - s}(X(\max\{t, w\}) / X(\max\{t + r, w\})) \to  \pi_{t - s}(X(\max\{t, w\}) / X(\max\{t + 1, w\}))$. Thus, $Z_r^{s, t}(\kappa_{\geq w} X)$ is isomorphic to $Z_r^{s, t}$ if $t \geq w$, and it vanishes otherwise. Furthermore, for $t \geq w$, $B_r^{s, t}(\kappa_{\geq w} X)$ is the image of $\delta_r\colon \pi_{t - s - 1}(X(\max\{t - r + 1, w\}) / X(t)) \to \pi_{t - s}(X(t)) \to \pi_{t - s}(X(t) / X(t + 1))$, which is precisely $B_{t - \max\{w, t - r + 1\} + 1}^{s, t} =  B_{\min\{t - w + 1, r\}}^{s, t}$. We have thus computed $E_{r + 1}^{s, t, w}(\kappa X)$ for each $1 \leq r < \infty$ (the result for $r = \infty$ also follows by taking limits), which in particular recovers the expression for $E_2^{*,*,*}(\kappa X)$. As for differentials,  the natural map $\kappa_{\geq w} X \to \kappa X [\defopara_v^{-1}]$ defines a map of SS $E_*^{**w}(\kappa X) \to E_*^{*,*}(X)$ which is injective on each $E_{r + 1}$ entry that could receive a nontrivial $d_{r + 1}$ and surjective on each $E_{r + 1}$ entry that could support a nontrivial $d_{r + 1}$, so it determines the pattern of $d_{r + 1}$ for each $r \geq 1$ as stated above. \parr 
    
    Part 2 follows from applying the result on strong convergence in Theorem \ref{Bockstein-dictionary-infinite} item 2 to each filtration $\kappa_{\geq w} X$. Here the whole plane obstructions vanish due to Remark \ref{vanishing-of-the-obstructions-to-strong-convergence} and the fact that each $\kappa_{\geq w} X$ gives rise to a half-plane SS.\parr 
    
    For part 3, one computes $E_{r + 1}^{s, t, w} (\kappa X / \defopara_v^n) \cong E_{r + 1}^{s, t}(\kappa_{\geq w} X / \kappa_{\geq w + n} X)$ using the same method as in the proof of part 1. For differentials, the map $\kappa X \to \kappa X / \defopara_v^n$ induces a map of filtered SS $E_*^{*,*,*}(\kappa X) \to E_*^{*,*,*}(\kappa X / \defopara_v^n)$ that is injective on each $E_{r + 1}$ entry that could receive a nontrivial $d_{r + 1}$ and surjective on each $E_{r + 1}$ entry that could support a nontrivial $d_{r + 1}$, so it determines the pattern of $d_{r + 1}$ for each $r \geq 1$ as stated above.  The strong convergence here also follows from Theorem \ref{Bockstein-dictionary-infinite} item 2 applied to $\kappa_{\geq w} X / \kappa_{\geq w + n} X$ (these are complete filtrations with vanishing obstructions for degree reasons). 
\end{proof}

\begin{remark} \label{circularity}
    If one chooses to avoid $\defopara, \rho, \delta$ when setting up the theory of standard SS and strong convergence while recovering the relation to these operators afterwards, then one can establish Theorem \ref{rigidity-for-trigraded-defopara-Bockstein-SS} first (by combining the argument above and an appropriate paraphrase of \cite[\S~7]{Boa99}) and use it to prove Theorem \ref{Bockstein-dictionary-infinite} and Theorem \ref{Bockstein-dictionary-finite}, cf. \cite[\S 3]{vN25}.
\end{remark}

\begin{remark} \label{rigidity-for-filtered-SS-by-definition}
    In this remark we go through the story of \emph{rigidity} for spectral sequences. 
    \begin{enumerate}
        \item Since filtered colimits are exact in $\Ab$, there is a functor $-[\defopara^{-1}]\colon \Fil(\SpSeq) \to \SpSeq$ sending $E = \{E_r^{s, t, w}\}_{r \geq 2, s, t, w \in \Zb}$ to $E[\defopara^{-1}] = \{\colim_{w \to -\infty} E_r^{s, t, w}\}_{r \geq 2, s, t\in \Zb}$. We say a filtered SS $E$ is \textbf{rigid} (over $E[\defopara^{-1}]$) if $E$ is determined by $E[\defopara^{-1}]$ in the sense of Theorem \ref{rigidity-for-trigraded-defopara-Bockstein-SS} part 1. 
        \item Suppose $E = \{E_r^{s, t, w}\}$ is a filtered SS, $E' = \{(E')_r^{s, t}\}$ is a SS, $\{f_w\colon E_*^{**w} \to (E')_*^{*,*}\}_{w \in \Zb}$ is a family of maps of SS so that $f_{w} \defopara = f_{w + 1}$ for each $w \in \Zb$, and these maps induce an isomorphism $E_2^{*,*,*} \cong (E'_2)^{*,*} \otimes_{\Zb} \Zb[\defopara]$. Then by induction on $1 \leq r < \infty$, the map $\{f_w\colon E_{r + 1}^{**w} \to (E')_{r + 1}^{*,*}\}_{w \in \Zb}$ is injective on each $E_{r + 1}$ entry that could receive a nontrivial $d_{r + 1}$ and surjective on each $E_{r + 1}$ entry that could support a nontrivial $d_{r + 1}$, so the $d_{r + 1}$-patterns in $E$ and $E'$ determine each other. Thus, a page-by-page computation shows $E$ is rigid over $E[\defopara^{-1}] \cong E'$. 
        \item According to Theorem \ref{rigidity-for-trigraded-defopara-Bockstein-SS} part 1, if $Y \in \biFil\Sp$ is isomorphic to $\kappa X$ for some $X \in \Fil\Sp$, then $E_*^{*,*,*}(Y)$ is rigid over $E_*^{*,*,*}(Y)[\defopara^{-1}] \cong E_*^{*,*}(X)$. 
        \item Conversely, suppose $Y = \{Y_w\}_{w\in \Zb} \in \biFil\Sp$ is horizontally complete (cf. Remark \ref{diagonal-bifiltrations}), $X \in \Fil\Sp$ is complete, and there is a family of compatible maps $\{Y_w \to X\}_{w \in \Zb}$ which induce an isomorphism $E_2^{*,*,*}(Y) \cong E_2^{*,*}(X) \otimes_{\Zb} \Zb[\defopara]$. Then $Y[\defopara_v^{-1}] \cong X$ due to Remark \ref{recollements-and-Nakayama-lemma} part 2. Moreover, Remark \ref{diagonal-bifiltrations} shows $Y \in \biFil\Sp_{\mathrm{de}}$, and $-[\defopara_v^{-1}] \colon \biFil\Sp_{\mathrm{de}} \to \Fil\Sp$ is an equivalence since it is a left inverse to $\kappa \colon \Fil\Sp \to \biFil\Sp_{\mathrm{de}}$. Therefore, $Y \cong \kappa X$.
    \end{enumerate}
    Consequently, trigraded spectral sequences with the rigidity property are precisely the shadows of $\defopara$-Bockstein bifiltrations. Examples of this rigidity property include the motivic Adams--Novikov SS (cf. \cite[\S~6.1 and \S~6.2]{Isaksen-stable-stems}) which is a special case of Remark \ref{rigidity-of-synthetic-Adams-SS}, and the motivic Cartan--Eilenberg SS (cf. \cite[Theorem 2.8]{Burklund-Xu}) which follows directly from item 2 (and item 4) above.
\end{remark}

We now compare these with the synthetic Adams SS. Let $E$ be a homotopy ring spectrum of Adams type. Recall the forgetful functor $U = \Msp^{\Fil}(\Sb^{0, 0}_E, -)\colon \Syn_E \to \Fil\Sp$ as in Remark \ref{synthetic-spectra-as-special-case}.

\begin{theorem} \label{synthetic-Adams-vs-Bockstein}
    Suppose $X \in \Sp$ is $E$-nilpotent complete. 
    \begin{enumerate}
        \item The $E$-Adams SS of $X$, starting from its $E_2$-page, coincides with the standard SS of $U(\nu X)$.
        \item The synthetic $\nu E$-Adams SS of $\nu X$, starting from its $E_2$-page, coincides with the trigraded $\defopara$-Bockstein SS of $U(\nu X)$.
        \item For each $1 \leq r < \infty$, the synthetic $\nu E$-Adams SS of $\nu X / \defopara^r$, starting from its $E_2$-page, coincides with the refined trigraded $\defopara$-Bockstein SS of $U(\nu X) / \defopara^r$. 
    \end{enumerate}
    The identifications above are also compatible with their abutments, so the  $\IIm(\defopara^a)$-filtrations on $\pi_{**}(\nu X)$ and $\pi_{**}(\nu X / \defopara^r)$, as in Theorems \ref{Bockstein-dictionary-infinite} and \ref{Bockstein-dictionary-finite}, coincide with their synthetic $\nu E$-Adams filtrations. 
\end{theorem}

\begin{proof}
    To identify the SS and the abutment altogether, it suffices to identify the filtered spectra behind them. The result of part 1 is standard: we have
    \begin{align*}
        U(\nu X)(w) &= \Msp(\Sb^{0, w}, \nu X) \cong \lim_{[k] \in \Delta} \Msp(\Sb^{0, w}, \nu X \otimes \nu E^{\otimes k}) \\
        &\cong \lim_{[k] \in \Delta} \Msp(\Sb^{0, w}, \nu(X \otimes E^{\otimes k})) \cong \lim_{[k] \in \Delta} \tau_{\geq w}(X \otimes E^{\otimes k}).
    \end{align*}
    Here the first $\cong$ follows from \cite[Proposition A.13]{BHS1}, the second $\cong$ follows from \cite[Lemma 4.24]{Pst23}, and the third $\cong$ follows from the fact that $\Msp(\Sb^{0, w}, \nu M) \cong \tau_{\geq w} M$ for each $M \in \Sp$ with a homotopy $E$-module structure (which boils down to \cite[Proposition 4.60]{Pst23}). On the other hand, the $E$-Adams SS of $X$ is constructed using the $E$-Adams tower coming from the fiber of $\Sb^0 \to E$. According to \cite[Proposition 2.14]{MNN17} and \cite[Proposition 6.3]{Lev15}, the filtered spectrum for this SS, starting from its $E_2$-page, is also given by $\lim_{[k] \in \Delta} \tau_{\geq \bullet} (X \otimes E^{\otimes k})$. We thus conclude part 1.\parr 
    
    For part 2, the weight $= w$ part of the trigraded $\defopara$-Bockstein SS of $U(\nu X)$ comes from the filtration $\kappa_{\geq w} U(\nu X)$.
    On the other hand, the synthetic $\nu E$-Adams SS of $\nu X$ is constructed using the $\nu E$-Adams tower coming from the fiber of $\nu\Sb^0 \to \nu E$. Thus, the same argument shows the filtered spectrum for the weight $= w$ part of this SS, starting from its $E_2$-page, 
    is given by $\lim_{[k] \in \Delta} \tau_{\geq \bullet}(\Msp(\Sb^{0, w}, \nu X \otimes \nu E^{\otimes k}))$. Note that $\Msp(\Sb^{0, w}, \nu X  \otimes \nu E^{\otimes k}) \cong \Msp(\Sb^{0, w}, \nu(X \otimes E^{\otimes k})) \cong \tau_{\geq w} (X \otimes E^{\otimes k})$, and $\tau_{\geq \bullet}(\tau_{\geq w} (X \otimes E^{\otimes k})) \cong \kappa_{\geq w} (\tau_{\geq \bullet} (X \otimes E^{\otimes k}))$. 
    Since $\kappa_{\geq w}$ preserves limits, the totalization yields $\kappa_{\geq w} U(\nu X)$ as well. \parr 
    
    For part 3, the weight $= w$ part of the refined trigraded $\defopara$-Bockstein SS for $U(\nu X) / \defopara^r$ comes from the filtration $\kappa_{\geq w} U(\nu X) / \kappa_{\geq w + r} U(\nu X)$.
    Also, the identical argument as before shows the filtered spectrum for the weight $= w$ part of the synthetic $\nu E$-Adams SS of $\nu X / \defopara^r$, starting from its $E_2$-page, 
    is given by $\lim_{[k] \in \Delta} \tau_{\geq \bullet}(\Msp(\Sb^{0, w}, (\nu X / \defopara^r) \otimes \nu E^{\otimes k}))$. Note that 
    \begin{align*}
        \Msp(\Sb^{0, w}, (\nu X / \defopara^r)  \otimes \nu E^{\otimes k}) &\cong \Msp(\Sb^{0, w}, \nu(X \otimes E^{\otimes k}) / \defopara^r) \\ 
        &\cong \tau_{\geq w} (X \otimes E^{\otimes k}) / \tau_{\geq w + r} (X \otimes E^{\otimes k}) 
    \end{align*}
    and $\tau_{\geq \bullet}(\tau_{\geq w} (X \otimes E^{\otimes k}) / \tau_{\geq w + r} (X \otimes E^{\otimes k}) ) \cong \kappa_{\geq w}(\tau_{\geq \bullet}(X \otimes E^{\otimes k})) / \kappa_{\geq w + r}(\tau_{\geq \bullet}(X \otimes E^{\otimes k}))$. 
    Thus, the totalization also gives rise to $\kappa_{\geq w}(U(\nu X)) / \kappa_{\geq w + r}(U(\nu X))$.
\end{proof}

\begin{remark} \label{rigidity-of-synthetic-Adams-SS}
    In terms of Remark \ref{rigidity-for-filtered-SS-by-definition}, the above theorem implies that the synthetic $\nu E$-Adams SS of $\nu X$ is rigid (over the $E$-Adams SS of $X$). 
\end{remark}

%% file: cMT-Proof.tex
In this appendix we prove Theorem \ref{coherent-Mahowald-trick}. Our strategy goes as follows:
For any stable $\infty$-category $\CC$ and any $X \in \CC$, we construct 
an ``auxiliary'' $2$d cuDT of the form 
    \[\begin{tikzpicture}[baseline= (a.base)]
        \node[scale=1.0] (a) at (0,0){
            \begin{tikzcd}
                & \vdots \ar[d] & \vdots \ar[d] & \vdots \ar[d, equal] & \vdots \ar[d] & \vdots \ar[d] & \\
                \cdots \ar[r] & 0 \ar[r] \ar[d] & X \ar[r, equal] \ar[d, equal] & X \ar[r] \ar[d] & 0 \ar[r] \ar[d]& \Sigma X \ar[r, equal] \ar[d, equal] & \cdots \\
                \cdots \ar[r] & X \ar[r, equal] \ar[d, equal] & X \ar[r] \ar[d] & 0 \ar[r] \ar[d] & \Sigma X \ar[r, equal] \ar[d, equal]& \Sigma X \ar[r] \ar[d] & \cdots \\
                \cdots \ar[r, equal] & X \ar[r] \ar[d] & 0 \ar[r] \ar[d] & \Sigma X \ar[r, equal] \ar[d, equal] & \Sigma X \ar[r] \ar[d]& 0 \ar[r] \ar[d] & \cdots \\
                \cdots \ar[r] & 0 \ar[r] \ar[d] & \Sigma X \ar[r, equal] \ar[d, equal] & \Sigma X \ar[r] \ar[d] & 0 \ar[r] \ar[d]& \Sigma^2 X \ar[r, equal] \ar[d, equal] & \cdots \\
                \cdots \ar[r] & \Sigma X \ar[r, equal] \ar[d, equal] & \Sigma X \ar[r] \ar[d] & 0 \ar[r] \ar[d] & \Sigma^2 X \ar[r, equal] \ar[d, equal]& \Sigma^2 X \ar[r] \ar[d] & \cdots \\
                & \vdots & \vdots & \vdots & \vdots & \vdots & 
            \end{tikzcd}
        };
    \end{tikzpicture}\]
for which we can pin down all data (objects, maps, squares, $\alpha$-isomorphisms, $\beta$-isomorphisms, etc.) explicitly. We then reduce Theorem \ref{coherent-Mahowald-trick} to the full subcategory $\IIm(\phi^*\colon \Aux(\CC) \hookrightarrow 2\DT^u(\CC))$ spanned by such $2$d cuDTs, and prove it in these cases through a simple diagram chase. 
\begin{remark}
    Actually, according to Remark \ref{coherent-Mahowald-trick-as-a-modification}, we can further reduce to just one concrete example, namely the auxiliary $2$d cuDT for $X = \mathrm{H}\Qb \in \CC = \Sp$ (or $\CC = \D(\Qb)$). 
\end{remark}
We will continue to use the notations from \S~\hyperref[subsec:4.1]{4.1}. The crucial step for our strategy is the following construction.

\begin{construction} \label{phi-map}
    Write $\CQ$ for the sub-poset $(\Zb,\leq) \times (\Zb, \leq)$ consisting of pairs $(i, j)$ such that $-1 \leq i - j \leq 1$. There is a well-defined map of sets 
    \[\phi \colon \CP \times \CP \to \CQ, \qquad (m, n, a, b) \mapsto (k, k) + \epsilon(m, n, a, b),\]
    where $s = m + n + a + b$, $k = \lfloor s/3 \rfloor$ and the correction term $\epsilon = \epsilon(m, n, a, b)$ takes values in the five-element set $\{(0, 0), (1, 0), (0, 1), (-1, 0), (0, -1)\}$. More precisely, the value of $\epsilon(m, n, a, b)$ is determined by the table below:
    \begin{table}[htbp]
        \centering
        \begin{tabular}{||c|c|c|c||c||}
            \hline
            $s \mod 3$  & $a + b - k \mod 2$ & $0 \leq m - n \leq 1$ & $0 \leq a - b \leq 1$ &  value of $\epsilon$\\
            \hline \hline
           
            \multirow{8}{*}{$0$} & \multirow{4}{*}{$0$} & \multirow{2}{*}{T} & T & $(0, -1)$\\
            \cline{4-5}
            && & F & \multirow{3}{*}{$(-1, 0)$} \\
            \cline{3-4}
            && \multirow{2}{*}{F} & T & \\
            \cline{4-4}
            && & F & \\
            \cline{2-5}
            & \multirow{4}{*}{$1$} & \multirow{2}{*}{T} & T & $(-1, 0)$ \\
            \cline{4-5}
            && & F & \multirow{3}{*}{$(0, -1)$} \\
            \cline{3-4}
            && \multirow{2}{*}{F} & T & \\
            \cline{4-4}
            && & F & \\
            \hline

            \multirow{8}{*}{$1$} & \multirow{4}{*}{$0$} & \multirow{2}{*}{T} & T & $(0, 0)$\\
            \cline{4-5}
            && & F & $(-1, 0)$ \\
            \cline{3-5}
            && \multirow{2}{*}{F} & T & $(0, -1)$ \\
            \cline{4-5}
            && & F & $(1, 0)$ \\
            \cline{2-5}
            & \multirow{4}{*}{$1$} & \multirow{2}{*}{T} & T & $(0, 0)$ \\
            \cline{4-5}
            && & F & $(0, -1)$ \\
            \cline{3-5}
            && \multirow{2}{*}{F} & T & $(-1, 0)$ \\
            \cline{4-5}
            && & F & $(0, 1)$ \\
            \hline

            \multirow{8}{*}{$2$} & \multirow{4}{*}{$0$} & \multirow{2}{*}{T} & T & $(0, 0)$\\
            \cline{4-5}
            && & F & $(1, 0)$ \\
            \cline{3-5}
            && \multirow{2}{*}{F} & T & $(0, 1)$ \\
            \cline{4-5}
            && & F & $(-1, 0)$ \\
            \cline{2-5}
            & \multirow{4}{*}{$1$} & \multirow{2}{*}{T} & T & $(0, 0)$ \\
            \cline{4-5}
            && & F & $(0, 1)$ \\
            \cline{3-5}
            && \multirow{2}{*}{F} & T & $(1, 0)$ \\
            \cline{4-5}
            && & F & $(0, -1)$ \\
            \hline
            
        \end{tabular}
        \caption{Values of the correction term $\epsilon = \epsilon(m, n, a, b)$.}
        \label{table-epsilon}
    \end{table}
\end{construction}

\FloatBarrier

The specific formula for \(\phi\) is not essential in itself; what matters is that it gives an explicit map \(\phi\colon \CP \times \CP \to \CQ\) with the following properties.

\begin{lemma} \label{is-map-of-posets}
    \begin{enumerate}
        \item The map $\phi\colon \CP \times \CP \to \CQ$ is an order-preserving surjection. Also, if $(m, n, a, b)$ lies in the $\phi$-preimage of some $(i, i)$, then $0 \leq m - n \leq 1$ and $0 \leq a - b \leq 1$. 
        \item Furthermore, if a square in $\CP \times \CP$ is of the form 
        \[\begin{tikzcd}
            (m, n, a, b) \ar[r] \ar[d] & (m + 1, n, a, b) \ar[d] \\
            (m, n + 1, a, b) \ar[r] & (m + 1, n + 1, a, b) 
        \end{tikzcd} \quad \text{ or } \quad
        \begin{tikzcd}
            (m, n, a, b) \ar[r] \ar[d] & (m, n, a + 1, b) \ar[d] \\
            (m, n, a, b + 1) \ar[r] & (m, n, a + 1, b + 1) 
        \end{tikzcd}\]
        with $0 \leq m - n \leq 1$, $0 \leq a - b \leq 1$, then its $\phi$-image in $\CQ$ has one of the following shapes (up to transpose):
        \[\begin{tikzcd}
            (i, j) \ar[r] \ar[d] & (i + 1, j) \ar[d] \\
            (i, j + 1) \ar[r] & (i + 1, j + 1) 
        \end{tikzcd} \text{ or }
        \begin{tikzcd}
            (i, j) \ar[r] \ar[d] & (i + 1, j) \ar[d] \\
            (i, j) \ar[r] & (i + 1, j) 
        \end{tikzcd} \text{ or }
        \begin{tikzcd}
            (i, j) \ar[r] \ar[d] & (i, j) \ar[d] \\
            (i, j + 1) \ar[r] & (i, j + 1) 
        \end{tikzcd}\]
    \end{enumerate}
\end{lemma}

\begin{proof}

    We start by introducing some notation: 
    
    \begin{itemize}
        \item We say $(m, n) \in \CP$ is \emph{of type T} if $0 \leq m - n \leq 1$, and it is \emph{of type F} otherwise. 
        
        \item For $(m, n, a, b) \in \CP \times \CP$, we say it is of type TT if $(m, n)$ and $(a, b)$ are both of type T. Likewise, in $\CP \times \CP$ we can make sense of type TF, FT and FF. 
        
        \item We say $(m, n) \leq (m', n')$ in $\CP$ is an \emph{increment} if $m' = m + 1, n' = n$ or $m' = m, n' = n + 1$. The order relation on $\CP$ is generated by increments, i.e. each nontrivial $(m, n) \leq (m', n')$ in $\CP$ factors as a finite chain of increments. Note that there are only three types of increments in $\CP$: T to T, T to F, F to T. For degree reasons, there is no increment of type F to F. 
        
        \item The order relation on $\CP \times \CP$ is generated by increments in two directions. We say increments $(m, n, a, b) \leq (m', n', a, b)$ are \emph{horizontal} and increments $(m, n, a, b) \leq (m, n, a', b')$ are \emph{vertical}. 
        \item The poset $\CQ$ admits a $\Zb / 2$ action generated by the \emph{swap} map $\rho\colon \CQ \to \CQ, (i, j) \mapsto (j, i)$. 
    \end{itemize}

    We now prove part 1. The map $\phi$ is surjective since each element of $\CQ$ takes one of the forms $(i, i)$, $(i, i) + (0, -1)$, or $(i, i) + (-1, 0)$, which are hit respectively by $(i + 1, i, \lfloor (i + 1) / 2\rfloor, \lceil (i - 1)/2 \rceil)$, $(i, i, \lfloor (i + 1) / 2\rfloor, \lceil (i - 1)/2 \rceil)$, or $(i + 1, i - 1, \lfloor (i + 1) / 2\rfloor, \lceil (i - 1)/2 \rceil)$. Also, if $\phi(m, n, a, b) = (i, i)$, then $\epsilon(m, n, a, b) = 0$, which forces $(m,n,a,b)$ to be of type TT by definition. \parr 
    
    To show $\phi(m, n, a, b) \leq \phi(m', n', a', b')$ whenever $(m, n, a, b) \leq (m', n', a', b')$, it suffices to consider increments. We first fix $(a, b)$ and treat the horizontal increments $(m, n) \leq (m', n')$. It is safe to assume $a + b - k \equiv 0$ mod $2$, as $\epsilon$-images for the case $a + b - k \equiv 1$ mod $2$ would only differ by a swap. If $(a, b)$ is of type T, then the $\epsilon$-image of this horizontal increment is given by 
    \begin{table}[htbp]
        \centering
        \begin{tabular}{||c||c|c|c||}
            \hline
            \diagbox[width=0.28\linewidth]{$s \mod 3$}{increment type} & TT to TT & TT to FT & FT to TT\\
            \hline 
            \hline 
            $0$ & $(0, -1)$ to $(0, 0)$ & $(0, -1)$ to $(0, -1)$ & $(-1, 0)$ to $(0, 0)$ \\
            \hline
            $1$ & $(0, 0)$ to $(0, 0)$ & $(0, 0)$ to $(0, 1)$ & $(0, -1)$ to $(0, 0)$ \\
            \hline
            $2$ & $(0, 0)$ to $(-1, 0)$ & $(0, 0)$ to $(0, -1)$ & $(0, 1)$ to $(-1, 0)$\\
            \hline
        \end{tabular}
        % \caption{Caption}
        % \label{tab:placeholder}
    \end{table}

    \FloatBarrier

    We access the table as follows: the top left slot reads ``$(0, -1)$ to $(0, 0)$'', which means if $(m, n, a, b)$ satisfies $s = m + n + a + b \equiv 0$ mod $3$ and $(m, n), (m', n')$ are both of type T, then 
    \[\phi(m, n, a, b) = (k, k) + (0, -1) = (k, k - 1) \leq \phi(m', n', a, b) = (k, k) + (0, 0) = (k, k)\]   
    so $\phi$ is order-preserving in this case. Using the first two rows above, this argument shows that $\phi$ is order preserving for $s \equiv 0, 1$ mod $3$. If $s \equiv 2$ mod $3$, then $k' = \lfloor (m' + n' + a + b) / 3 \rfloor = k + 1$. Thus, for instance the bottom right slot reads ``$(0, 1)$ to $(-1, 0)$'', which means if $(m, n, a, b)$ satisfies $s = m + n + a + b \equiv 2$ mod $3$, $(m, n)$ is of type F and $(m', n')$ is of type T, then 
    \[\phi(m, n, a, b) = (k, k) + (0, 1) = (k, k + 1) \leq \phi(m', n', a, b) = (k + 1, k + 1) + (-1, 0) = (k, k + 1).\] Using the third row, this argument shows that $\phi$ is also order preserving for $s \equiv 2$ mod $3$. \parr 

    If $(a, b)$ is of type F, then the $\epsilon$-image of this increment is given by another table:
    
    \begin{table}[htbp]
        \centering
        \begin{tabular}{||c||c|c|c||}
            \hline
            \diagbox[width=0.28\linewidth]{$s \mod 3$}{increment type} & TF to TF & TF to FF & FF to TF\\
            \hline 
            \hline 
            $0$ & $(-1, 0)$ to $(-1, 0)$ & $(-1, 0)$ to $(1, 0)$ & $(-1, 0)$ to $(-1, 0)$ \\
            \hline
            $1$ & $(-1, 0)$ to $(1, 0)$ & $(-1, 0)$ to $(-1, 0)$ & $(1, 0)$ to $(1, 0)$ \\
            \hline
            $2$ & $(1, 0)$ to $(0, -1)$ & $(1, 0)$ to $(0, -1)$ & $(-1, 0)$ to $(0, -1)$\\
            \hline
        \end{tabular}
        % \caption{Caption}
        % \label{tab:placeholder}
    \end{table}

    \FloatBarrier

    A similar argument applied to this table shows that $\phi$ is order preserving in all listed instances, so this concludes the case of horizontal increments. Next we fix $(m, n)$ and treat vertical increments $(a, b) \leq (a', b')$. Again we can assume $(m, n, a, b)$ satisfies $a + b - k \equiv 0$ mod $2$. We can summarize the $\epsilon$-image of this vertical increment in two tables, for $(m, n)$ of type T or F respectively:
    
    \begin{table}[htbp]
        \centering
        \begin{tabular}{||c||c|c|c||}
            \hline
            \diagbox[width=0.28\linewidth]{$s \mod 3$}{increment type} & TT to TT & TT to TF & TF to TT\\
            \hline 
            \hline 
            $0$ & $(0, -1)$ to $(0, 0)$ & $(0, -1)$ to $(0, -1)$ & $(-1, 0)$ to $(0, 0)$ \\
            \hline
            $1$ & $(0, 0)$ to $(0, 0)$ & $(0, 0)$ to $(0, 1)$ & $(-1, 0)$ to $(0, 0)$ \\
            \hline
            $2$ & $(0, 0)$ to $(0, -1)$ & $(0, 0)$ to $(-1, 0)$ & $(1, 0)$ to $(0, -1)$\\
            \hline

            \multicolumn{4}{@{}c@{}}{\rule{0pt}{12pt}}\\[-12pt]

            \hline
            \diagbox[width=0.28\linewidth]{$s \mod 3$}{increment type} & FT to FT & FT to FF & FF to FT\\
            \hline 
            \hline 
            $0$ & $(-1, 0)$ to $(-1, 0)$ & $(-1, 0)$ to $(0, 1)$ & $(-1, 0)$ to $(-1, 0)$ \\
            \hline
            $1$ & $(0, -1)$ to $(1, 0)$ & $(0, -1)$ to $(0, -1)$ & $(1, 0)$ to $(1, 0)$ \\
            \hline
            $2$ & $(0, 1)$ to $(-1, 0)$ & $(0, 1)$ to $(-1, 0)$ & $(-1, 0)$ to $(-1, 0)$\\
            \hline
        \end{tabular} 
        % \caption{Caption}
        % \label{tab:placeholder}
    \end{table}

    \FloatBarrier

    The argument for the first table applies equally well here, which explains that $\phi$ is order-preserving on all vertical increments. Thus $\phi$ is order-preserving, which finishes part 1. \parr 

    For part 2, take $(m, n, a, b) \in \CP \times \CP$ of type TT. Here we can also assume $a + b - k \equiv 0$ mod $2$, since $\phi$-images for the case $a + b - k \equiv 1$ mod $2$ would only differ by a swap. We first consider the horizontal bi-increment square ``$(m, n, a, b)$ to $(m + 1, n, a, b)$ \& $(m, n + 1, a, b)$ to $(m + 1, n + 1, a, b)$'' in $\CP \times \CP$. The type of this square is always ``TT to TT \& FT to TT'' up to transpose, and the possible $\phi$-values are listed in the following table: 
 
    \begin{table}[htbp]
        \centering
        \begin{tabular}{||c||c||}
            \hline
            $s \mod 3$ & $\phi$-value of the square (up to transpose)\\
            \hline
            \hline
            $0$ & $(k, k - 1)$ to $(k, k)$ \& $(k, k - 1)$ to $(k, k)$ \\
            \hline
            $1$ & $(k, k)$ to $(k, k)$ \& $(k, k + 1)$ to $(k, k + 1)$ \\
            \hline
            $2$ & $(k, k)$ to $(k, k + 1)$ \& $(k + 1, k)$ to $(k + 1, k + 1)$ \\
            \hline
        \end{tabular}
        % \caption{Caption}
        % \label{tab:placeholder}
    \end{table}

    \FloatBarrier

    On the other hand, for the vertical bi-increment square ``$(m, n, a, b)$ to $(m, n, a + 1, b)$ \& $(m, n, a, b + 1)$ to $(m, n, a + 1, b + 1)$'', its type should be ``TT to TT \& TF to TT'' up to transpose, and the possible $\phi$-values are as follows: 

    \begin{table}[htbp]
        \centering
        \begin{tabular}{||c||c||}
            \hline
            $s \mod 3$ & $\phi$-value of the square (up to transpose)\\
            \hline
            \hline
            $0$ & $(k, k - 1)$ to $(k, k)$ \& $(k, k - 1)$ to $(k, k)$ \\
            \hline
            $1$ & $(k, k)$ to $(k, k)$ \& $(k, k + 1)$ to $(k, k + 1)$ \\
            \hline
            $2$ & $(k, k)$ to $(k + 1, k)$ \& $(k, k + 1)$ to $(k + 1, k + 1)$ \\
            \hline
        \end{tabular}
        % \caption{Caption}
        % \label{tab:placeholder}
    \end{table}

    \FloatBarrier

    From these we can conclude that part 2 holds true in full generality. 
\end{proof}

\begin{proposition}\label{aux-generation}
    Write $\Aux(\CC) \subset \Fun(\CQ, \CC)$ for the full subcategory spanned by functors $F\colon \CQ \to \CC$ so that $F(i, j) = 0$ if $i \neq j$, and the image of each square 
    \[\begin{tikzcd}
        (i, j) \ar[r] \ar[d] & (i + 1, j) \ar[d] \\
        (i, j + 1) \ar[r] & (i + 1, j + 1) 
    \end{tikzcd}\]
    is cartesian in $\CC$. Then restrictions along the maps 
    \(\CQ'' = \{(0, 0)\} \to \CQ' = \{(0, -1) \leq (0, 0)\} \to \CQ\) induce equivalences $\Aux(\CC) \xrightarrow{\sim} \Aux'(\CC) = \{F\colon \CQ' \to \CC \mid F(0, -1) = 0\} \xrightarrow{\sim} \CC = \Fun(\CQ'', \CC)$. 
\end{proposition}

\begin{proof}
    Apply the same argument as in Theorem \ref{1d-cuDT-generation}. 
\end{proof}

This $\Aux(\CC)$ describes an interesting portion of $2\DT^u(\CC)$.  

\begin{proposition} \label{aux-to-2DT} 
    The precomposition functor $\phi^*\colon \Fun(\CQ, \CC) \to \Fun(\CP \times \CP, \CC)$ restricts to a fully faithful functor $\phi^*\colon \Aux(\CC) \to 2\DT^u(\CC)$, and its essential image consists of $F \in 2\DT^u(\CC)$ such that
    \begin{enumerate}
        \item If $m + n + a + b \equiv 0$ \textup{mod} $3$, then $F(m, n, a, b)= 0$. 
        \item If $0 \leq m - n \leq 1$, $0 \leq a - b \leq 1$ and $m + n + a + b \equiv  1$ \textup{mod} $3$, then $F(m, n, b + 1, a) \gets F(m, n, a, b) \to F(n + 1, m, a, b)$ are isomorphisms.
    \end{enumerate}
    Also, $\phi^*\colon \Aux(\CC) \to 2\DT^u(\CC)$ admits a left adjoint $\phi_!$ and a right adjoint $\phi_*$.
\end{proposition}

\begin{proof}
    The image of $\Aux(\CC) \subset \Fun(\CQ, \CC) \xrightarrow{\phi^*} \Fun(\CP \times \CP, \CC)$ lies in $2\DT^u(
    \CC)$ due to Lemma \ref{is-map-of-posets}. Furthermore, for $\CP'' = \{(0, 0), (1, 0)\} \subset \CP$, the restriction of $\phi\colon \CP \times \CP \to \CQ$ on $\CP'' \times \CP''$ factors through $\CQ' = \{(0, -1), (0, 0)\} \subset \CQ$. This induces a commutative square
    \[\begin{tikzcd}
        \Aux(\CC) \ar[r, "{\phi^*}"] \ar[d] & 2\DT^u(\CC) \ar[d] \\
        \Aux'(\CC) \ar[r, "{\phi^*}"] & \CC^{\Deltaone \times \Deltaone}
    \end{tikzcd}\]
    in which the two vertical functors are equivalences by Theorem \ref{2d-cuDT-generation} and Proposition \ref{aux-generation}. The lower horizontal functor is explicitly given as
    \begin{align*}
        \phi^*\colon \Aux'(\CC) \quad \ &\to \ \quad \CC^{\Deltaone \times \Deltaone} \\
        \begin{tikzcd}
        0 \ar[d] \\
        X
    \end{tikzcd} \quad \quad &\mapsto 
    \begin{tikzcd}[ampersand replacement = \&]
        0 \ar[d] \ar[r] \& X \ar[d, equal] \\
        X \ar[r, equal] \& X
    \end{tikzcd}
    \end{align*}
    so it admits a right adjoint $\phi_*$ sending $T \in \CC^{\Deltaone \times \Deltaone}$ to $0 \to T(1, 0) \times_{T(1, 1)} T(0, 1)$ and a left adjoint $\phi_!$ sending $T$ to $0 \to \cofib(T(0, 0) \to T(1,1))$. The unit of the adjunction $\phi^* \dashv \phi_*$ is an isomorphism, so $\phi^*$ is fully faithful, and its essential image in $\CC^{\Deltaone \times \Deltaone}$ consists of the squares $T$ whose arrows $T(1, 0) \to T(1, 1) \gets T(0, 1)$ are isomorphisms. This essential image corresponds to the full subcategory of $2\DT^u(\CC)$ described in the statement, since $2\DT^u(\CC) \cong \CC^{\Deltaone \times \Deltaone}$ is constructed by ``restriction along $(\CP'')^2 \to \CP^2$'' and ``taking iterated (co)fibers in both directions''. 
\end{proof}

On the other hand, we can construct objects in $\Aux(\CC)$ with explicit coherence data.

\begin{proposition} \label{explicit-objs-in-aux}
    The poset $\CQ$ is the colimit in $\Cat$ of the diagram 
    \[\begin{tikzcd}[column sep = tiny]
        & {[0]} \ar[ld, "{(1, 1)}"] \ar[rd, "{(0, 0)}"] & & {[0]} \ar[ld, "{(1, 1)}"] \ar[rd, "{(0, 0)}"] & & {[0]} \ar[ld, "{(1, 1)}"] \ar[rd, "{(0, 0)}"] & & {[0]} \ar[ld, "{(1, 1)}"] \ar[rd, "{(0, 0)}"] & \\
        \cdots & & {\Deltaone \times \Deltaone} & & {\Deltaone \times \Deltaone} & & {\Deltaone \times \Deltaone} & & \cdots
    \end{tikzcd}\]
    Therefore, $\Fun(\CQ, \CC)$ is the limit in $\Cat_{\st}$ of the diagram 
    \[\begin{tikzcd}[column sep = tiny]
        \cdots \ar[rd, "{\ev_{(1, 1)}}"] & & \CC^{\Deltaone \times \Deltaone} \ar[rd, "{\ev_{(1, 1)}}"] \ar[ld, "{\ev_{(0, 0)}}"] & & \CC^{\Deltaone \times \Deltaone} \ar[rd, "{\ev_{(1, 1)}}"] \ar[ld, "{\ev_{(0, 0)}}"] & & \CC^{\Deltaone \times \Deltaone} \ar[rd, "{\ev_{(1, 1)}}"]  \ar[ld, "{\ev_{(0, 0)}}"] & & \cdots \ar[ld, "{\ev_{(0, 0)}}"] \\
        & \CC & & \CC & & \CC & & \CC & 
    \end{tikzcd}\]
    Thus, specifying an object $T \in \Fun(\CQ, \CC)$ is the same as providing a family of squares $\{T_n \in \CC^{\Deltaone \times \Deltaone}\}_{n \in \Zb}$ together with a family of isomorphisms $\{\gamma_n\colon T_n(1, 1) \cong T_{n + 1}(0, 0)\}_{n \in \Zb}$. 
\end{proposition}

\begin{proof}
    The second part follows from the first since $\Fun(-, \CC)\colon \Cat^{\op} \to \Cat_{\st}$ preserves limits. For the first part, since $\CQ = \bigcup_{(m, n) \in \Zb \times \Zb, m < n} \{(x, y) \in \CQ \mid (m, m) \leq (x, y) \leq (n, n)\}$ is a filtered union (i.e. a filtered colimit in $\Cat$ by Lemma \ref{filtered-union}), it suffices by induction to show each $\CQ_{[m, n]} = \{(x, y) \in \CQ \mid (m, m) \leq (x, y) \leq (n, n)\}$ is the pushout in $\Cat$ of $\CQ_{[m, n - 1]} \gets [0] \to \CQ_{[n - 1, n]} \cong \Deltaone \times \Deltaone$, where the left (resp. right) map hits the terminal (resp. initial) object. This is due to Lemma \ref{gluing-posets}.
\end{proof}

\begin{remark}
    If we replace $\Deltaone \times \Deltaone$ by $\Deltaone$ in the diagram above, the colimit in $\Cat$ would be $(\Zb, \leq)$. This essentially recovers \cite[Proposition 3.3]{Ariotta}.  
\end{remark}

\begin{lemma}\label{filtered-union}
    Suppose $\CK$ is a small $\infty$-category, and $\{\CK_i\}_{i \in I}$ is a filtered family of full subcategories of $\CK$, so that\footnote{concretely, this means for each $x \in \CK$ there exist $i \in I$ and $y \in \CK_i$ so that $x \simeq y$ in $\CK$} $\CK = \bigcup_{i \in I} \CK_i$. Then $\CK \cong \colim_{i \in I} \CK_i$ in $\Cat$. 
\end{lemma}

\begin{proof}
    The functor $\colim_{i \in I} \CK_i \to \CK$ induced by the inclusions is essentially surjective since  
    $\pi_0(\CK^\simeq) = \bigcup_i \pi_0(\CK_i^\simeq)$. For full faithfulness, it suffices to show each $\CK_{j} \to \colim_i \CK_i$ is fully faithful. Actually, for $X, Y \in \CK_j$,
    \begin{align*}
        \Map_{\colim_{i \in I} \CK_i}(X, Y) &\cong \Map_{\colim\nolimits_{i \geq j} \CK_i}(X, Y) \\
        &= \lim\left({[0]} \xrightarrow{(X, Y)} (\colim\nolimits_{i \geq j} \CK_i)^2 \xleftarrow{(s, t)} \Fun(\Deltaone, \colim\nolimits_{i \geq j} \CK_i)\right) \\
        &\cong \colim\nolimits_{i \geq j} \lim\left({[0]} \xrightarrow{(X, Y)} \CK_i^2 \xleftarrow{(s, t)} \Fun(\Deltaone, \CK_i)\right) \\
        &\cong \colim\nolimits_{i \geq j} \Map_{\CK_i}(X, Y) \cong \Map_{\CK_i}(X, Y)
    \end{align*}
    as filtered colimits commute with finite limits in $\Cat$ and $\Deltaone \in \Cat$ is compact. 
\end{proof}

\begin{lemma}\label{gluing-posets}
    Suppose $\CA, \CB$ are posets, $x \in \CA$ is terminal (i.e. largest) and $z \in \CB$ is initial (i.e. smallest). Then the pushout $\CO$ of $\CA \xleftarrow{x} [0] \xrightarrow{z} \CB$ in $\Cat$ is also a poset. The underlying set of $\CO$ is $\CA \sqcup \CB / \{x \sim z\}$, and the order relation in $\CO$ is as follows: the two induced maps $i\colon \CA \to \CO$, $j\colon \CB \to \CO$ are order embeddings, $i(a) \leq j(b)$ for all $(a, b) \in \CA \times \CB$, while $j(b) \leq i(a)$ only if $a = x, b = z$.
\end{lemma}

\begin{proof}
    Here $\CA \sqcup \CB \to \CO$ is essentially surjective, so it suffices to compute mapping spaces. By \cite[Theorem 0.1]{HRS25}, the two functors $i, j$ are both fully faithful. Moreover, according to the same result, $\Map_{\CO}(i(a), j(b)) = |\CA_{a/} \times_{\CA} * \times_{\CB} \CB_{/b}|$ where 
    \[\CA_{a/} \times_{\CA} [0] \times_{\CB} \CB_{/b} =( \CA_{a/} \times_{\CA} [0]) \times ([0] \times_{\CB} \CB_{/b}) = \Map_{\CA}(a, x) \times \Map_{\CB}(z, b) = [0]\]
    while $\Map_{\CO}(j(b), i(a)) = |{\CB}_{b/} \times_{\CB} [0] \times_{\CA} \CA_{/a}|$ where
    \[\CB_{b/} \times_{\CB} [0] \times_{\CA} \CA_{/a} = \Map_{\CB}(b, z) \times \Map_{\CA}(x, a) = 
    \begin{cases}
        [0], & (a, b) = (x, z). \\
        \emptyset, & \text{otherwise}.
    \end{cases}\]
    Therefore, $\CO$ is a poset. The descriptions of $\CO$ also follow from these computations.
\end{proof}

\begin{proof}[Proof of Theorem \ref{coherent-Mahowald-trick}]
    It suffices to treat the case $x = y = 0$, as the general case then follows from this specific case of $F[+x, +y]$ combined with Remark \ref{alpha[+1]-is-minus-alpha}. Consider the adjunction $\phi^* \dashv \phi_*$ between $\Aux(\CC)$ and $2\DT^u(\CC)$ as in Proposition \ref{aux-to-2DT}. For $F \in 2\DT^u(\CC)$, the counit of this adjunction $\phi^*\phi_* F \to F$ induces an isomorphism $(\phi^*\phi_* F)_\partial \langle 0,0 \rangle \cong F_\partial \langle 0,0 \rangle$ by inspecting the explicit construction. Due to the naturality of $\alpha_v$, $\alpha_h$, and $\beta$, this counit map also induces isomorphisms $(\phi^*\phi_* F)_\partial \langle x,y \rangle \cong F_\partial \langle x,y \rangle$ for all $(x, y) \in \Zb^2$ with $x + y \equiv 0$ \textup{mod} $3$. Thus, it suffices to establish $\beta^3 \cong - \alpha_{vh}$ for each $F$ in the essential image of $\phi^*\colon \Aux(\CC) \to 2\DT^u(\CC)$. According to Proposition \ref{explicit-objs-in-aux}, for each $X \in \CC$ there is an object $\widetilde{X}\colon \CQ \to \CC$ in $\Aux(\CC)$ which takes the following form:
    \[\begin{tikzpicture}[baseline= (a.base)]
        \node[scale=.8] (a) at (0,0){
            \begin{tikzcd}
                % & \ddots \ar[rd, equal] \ar[d]&&& \\
                & & \ddots \ar[d] \ar[rd, equal] &&\\
                & \ddots \ar[r] \ar[rd, equal] & X \ar[r] \ar[d] & 0 \ar[d] \ar[rd, equal] &\\
                && 0 \ar[r] \ar[rd, equal] & \Sigma X \ar[r] \ar[d] & \ddots \\
                &&& \ddots  &  
            \end{tikzcd}
        };
    \end{tikzpicture}\]
    Concretely, $\widetilde{X}$ is determined by a family of squares and a family of isomorphisms, so that the $n$-th square $\widetilde{X}_n$ corresponds to the map $\id\colon \Sigma^{n + 1} X \to \Sigma^{n + 1} X$ under the pushout comparison functor $\cp\colon \CC^{\Deltaone \times \Deltaone} \to \CC^{\Deltaone}, T \mapsto (T(1, 0) \sqcup_{T(0,0)} T(0,1) \to T(1,1))$, and the $n$-th transition isomorphism $\widetilde{X}_n(1,1) = \Sigma^{n + 1} X \to \widetilde{X}_{n + 1}(0, 0)= \Sigma^{n + 1} X $ is $\id$. Furthermore, each object in $\Aux(\CC)$ is isomorphic to one such $\widetilde{X}$ due to Proposition \ref{aux-generation}. So it suffices to show for each $X^\phi = \phi^*\widetilde{X}\colon \CP \times \CP \to \CC$
    \[\begin{tikzpicture}[baseline= (a.base)]
        \node[scale=.8] (a) at (0,0){
            \begin{tikzcd}
                & \vdots \ar[d] & \vdots \ar[d] & \vdots \ar[d, equal] & \vdots \ar[d] & \vdots \ar[d] & \\
                \cdots \ar[r] & 0 \ar[r] \ar[d] & X \ar[r, equal] \ar[d, equal] & X \ar[r] \ar[d] & 0 \ar[r] \ar[d]& \Sigma X \ar[r, equal] \ar[d, equal] & \cdots \\
                \cdots \ar[r] & X \ar[r, equal] \ar[d, equal] & X \ar[r] \ar[d] & 0 \ar[r] \ar[d] & \Sigma X \ar[r, equal] \ar[d, equal]& \Sigma X \ar[r] \ar[d] & \cdots \\
                \cdots \ar[r, equal] & X \ar[r] \ar[d] & 0 \ar[r] \ar[d] & \Sigma X \ar[r, equal] \ar[d, equal] & \Sigma X \ar[r] \ar[d]& 0 \ar[r] \ar[d] & \cdots \\
                \cdots \ar[r] & 0 \ar[r] \ar[d] & \Sigma X \ar[r, equal] \ar[d, equal] & \Sigma X \ar[r] \ar[d] & 0 \ar[r] \ar[d]& \Sigma^2 X \ar[r, equal] \ar[d, equal] & \cdots \\
                \cdots \ar[r] & \Sigma X \ar[r, equal] \ar[d, equal] & \Sigma X \ar[r] \ar[d] & 0 \ar[r] \ar[d] & \Sigma^2 X \ar[r, equal] \ar[d, equal]& \Sigma^2 X \ar[r] \ar[d] & \cdots \\
                & \vdots & \vdots & \vdots & \vdots & \vdots & 
            \end{tikzcd}
        };
    \end{tikzpicture}\]
    we have $\beta^3 \cong -\alpha_{vh}$. Concretely, $X^\phi\langle 4, 0\rangle \to X^\phi\langle 4, 1\rangle \gets X^\phi\langle 3, 1\rangle$ is $\Sigma X \xrightarrow{\id} \Sigma  X \xleftarrow{\id} \Sigma  X$, so we take the pullback $X^\phi_\partial \langle 3, 0\rangle = X^\phi\langle 4, 0\rangle  \times_{X^\phi\langle 4, 1\rangle} X^\phi\langle 3, 1\rangle$ to be $\Sigma X$ and both projections $X^\phi_\partial \langle 3, 0\rangle \to X^\phi\langle 4, 0\rangle$, $X^\phi_\partial \langle 3, 0\rangle \to X^\phi\langle 3, 1\rangle$ to be $\id \colon \Sigma X \to \Sigma X$. We adopt the same convention for the other pullbacks $X^\phi_\partial \langle 2, 1\rangle, X^\phi_\partial \langle 1, 2\rangle, X^\phi_\partial \langle 0, 3\rangle$ and $X^\phi_\partial \langle 0, 0\rangle$ (replacing $\Sigma X$ by $X$ in the last case). By Construction \ref{beta-iso}, $\beta\colon X^\phi_\partial \langle 3, 0\rangle = \Sigma X \to X^\phi_\partial \langle 2, 1\rangle = \Sigma X$ is compatible with projections to $X^\phi \langle 3, 1\rangle = \Sigma X$ on both sides, so $\beta \cong \id_{\Sigma X}$. The same argument shows $\beta^3 \cong \id_{\Sigma X}$. On the other hand, $\alpha_v \colon \Sigma X^\phi_\partial \langle 0, 0\rangle = \Sigma X \to X^\phi_\partial \langle 0, 3\rangle = \Sigma X$ is homotopic to $\alpha_v\colon \Sigma X^\phi \langle 1, 0\rangle = \Sigma X \to X^\phi \langle 1, 3\rangle = \Sigma X$. By the definitions of $\alpha$ and $\phi$ in Constructions \ref{alpha-iso} and \ref{phi-map}, it comes from the square
    \[\begin{tikzcd}
        X^\partial (1, 0, 0, 0) = \widetilde{X} (0, 0) = X \ar[r] \ar[d] & X^\partial (1, 0, 2, 0) = \widetilde{X} (1, 0) = 0 \ar[d] \\
        X^\partial (1, 0, 0, 1) = \widetilde{X} (0, 1) = 0 \ar[r] & X^\partial (1, 0, 2, 1) = \widetilde{X} (1, 1) = \Sigma X        
    \end{tikzcd}\] 
    under the pushout comparison functor, so it is homotopic to $\id_{\Sigma X}$ due to our convention on $\widetilde{X}$. Also, $\alpha_h \colon \Sigma X^\phi_\partial \langle 0, 0\rangle = \Sigma X \to X^\phi_\partial \langle 3, 0\rangle = \Sigma X$ is homotopic to $\alpha_h\colon \Sigma X^\phi \langle 1, 0\rangle = \Sigma X \to X^\phi \langle 4, 0\rangle = \Sigma X$, which corresponds under the pushout comparison functor to
    \[\begin{tikzcd}
        X^\partial (1, 0, 0, 0) = \widetilde{X} (0, 0) = X \ar[r] \ar[d] & X^\partial (2, 0, 0, 0) = \widetilde{X} (0, 1) = 0 \ar[d] \\
        X^\partial (1, 2, 0, 0) = \widetilde{X} (1, 0) = 0 \ar[r] & X^\partial (2, 2, 0, 0) = \widetilde{X} (1, 1) = \Sigma X
    \end{tikzcd}\] 
    which is the transpose of the previous square. Thus, $\alpha_h$ is homotopic to $-\id_{\Sigma X}$ by \cite[Lemma 1.1.2.10]{HA}. Putting all these together, we get $\beta^3 \cong - \alpha_{v} \alpha_{h}^{-1}$ for $X^\phi$. 
\end{proof}

%% file: Picard.tex
In this appendix we study the lax \(\Eb_n\)-monoidality of the standard SS functor and its variants with general coefficient categories and additional grading data, with Picard-graded SS arising as the main examples. We begin with the case of spectra: in Theorem \ref{categorified-total-Leibniz-rule}, we prove the lax symmetric monoidality of \(E_*^{*,*}\colon \Fil\Sp \to \SpSeq\), together with the corresponding statement for total differentials. We then record the analogous statement for presentably \(\Eb_n\)-monoidal \(t\)-\(\infty\)-categories in Corollary \ref{categorified-total-Leibniz-rule-with-general-coefficient}. The rest of the appendix incorporates grading data. First, we introduce a coarse form of grading data and prove lax \(\Eb_n\)-monoidality for the resulting trigraded standard SS in Theorem \ref{categorified-total-Leibniz-rule-with-picard-trigrading}. We then refine the grading data and explain how the Koszul sign rule can be incorporated through Day convolution, allowing the trigraded statement to be compressed into lax \(\Eb_n\)-monoidality for the resulting bigraded standard SS in Theorem \ref{categorified-total-Leibniz-rule-with-picard-bigrading}.

\begin{construction}\label{SS-as-an-operad}
    Take $r \in \Nb$, $r \geq 1$.
    \begin{itemize}
        \item We write $\Gr\Ch(\Ab)^{E_{r}}$ for the abelian $1$-category of pairs $(M^{*, *}, d_{r})$ such that $\{M^{s, t}\}_{s, t \in \Zb}$ is a bigraded abelian group and $d_{r}\colon M^{s, t} \to M^{s + r, t + r - 1}$ is a family of abelian group homomorphisms with $d_r d_r = 0$. For $r = 1$, we have a canonical identification 
        \[\Gr\Ch(\Ab)^{E_{1}} \to \Gr(\Ch(\Ab)), \qquad (\{M^{s, t}\}_{s, t \in \Zb}, d_1) \mapsto \{(M^{*, w}, d_1)\}_{w \in \Zb}\]
        and we equip $\Gr\Ch(\Ab)^{E_{1}}$ with the symmetric monoidal structure from $\Gr(\Ch(\Ab)^\Kos)^\Kos$. In general, there is an isomorphism of abelian groups $v_r\colon \Zb \times \Zb \to \Zb \times \Zb, (s, t) \mapsto (s + (r - 1)(t - s), t + (r - 1)(t - s))$, such that reindexing (i.e. precomposition) along $v_r$ leads to a symmetric monoidal equivalence $v_r^*\colon \Gr(\Gr(\Ab)^\Kos)^\Kos \cong \Gr(\Gr(\Ab)^\Kos)^\Kos$. We thus define $\Gr\Ch(\Ab)^{E_{r}}$ to be the pullback of the cospan diagram of symmetric monoidal functors
        \[\Gr(\Ch(\Ab)^{\Kos})^{\Kos} \xrightarrow{U} \Gr(\Gr(\Ab)^{\Kos})^{\Kos} \xleftarrow{v_r^*} \Gr(\Gr(\Ab)^{\Kos})^{\Kos}\]
        where $U\colon \{(M^{*, w}, d_r)\}_{w \in \Zb} \mapsto \{M^{*, w}\}_{w \in \Zb}$ is the forgetful functor. By construction, we see:
        \begin{itemize}
            \item The functor $E_2 = \pi^B_*\colon \Fil\Sp \to \Gr(\Ch(\Ab)^{\Kos})^{\Kos}$ in Example \ref{lax-symmetric-monoidality-of-pi-star} can be identified with the composite of two lax symmetric monoidal functors $\Fil\Sp \to \Gr\Ch(\Ab)^{E_{2}}, X \mapsto (\{E_2^{s, t}(X)\}_{s, t \in \Zb}, d_2)$ and $v_2^*\colon \Gr\Ch(\Ab)^{E_{2}} \to \Gr(\Ch(\Ab)^{\Kos})^{\Kos}$. 
            \item The forgetful functor 
            \[U\colon \Gr\Ch(\Ab)^{E_{r}} \to \Gr(\Gr(\Ab)^{\Kos})^{\Kos}, \quad (\{M^{s, t}\}_{s, t\in \Zb}, d_{r}) \mapsto \{M^{s, t}\}_{s, t\in \Zb}\]
            is faithful and symmetric monoidal. Also, taking cycles and taking homology yield lax symmetric monoidal functors $Z, H\colon \Gr\Ch(\Ab)^{E_{r}} \to \Gr(\Gr(\Ab)^{\Kos})^{\Kos}$. 
            \item A map $F\colon X \otimes Y \to T$ in $\Gr\Ch(\Ab)^{E_{r}}$ amounts to a family of bilinear maps $\{X^{s_1, t_1} \times Y^{s_2,t_2} \to T^{s_1 + s_2, t_1 + t_2}\}$ satisfying the Leibniz rule $d_{r} F(x, y) = F(d_{r}x, y) + (-1)^{t_1 - s_1} F(x, d_{r}y)$. Similarly, a map of the form $X_1 \otimes \cdots \otimes X_n \to T$ amounts to a family of multilinear maps satisfying the Leibniz rule.
        \end{itemize}
        \item We write $\SpSeq$ for the limit of the diagram of lax symmetric monoidal functors
        \[\begin{tikzcd}[column sep = tiny, row sep = small]
            \Gr\Ch(\Ab)^{E_{2}} \ar[rd, "{H}"] & & \Gr\Ch(\Ab)^{E_{3}} \ar[rd, "{H}"] \ar[ld, "{U}"'] & & \cdots \ar[ld, "{U}"'] \\
            & \Gr(\Gr(\Ab)^{\Kos})^{\Kos} & & \Gr(\Gr(\Ab)^{\Kos})^{\Kos} &   
        \end{tikzcd}\]
        in the $\infty$-category of $\infty$-operads. Concretely,
        \begin{itemize}
            \item The $\infty$-operad $\SpSeq$ is $1$-truncated, i.e. it is a colored operad in the classical sense.
            \item Its underlying category is the $1$-category of spectral sequences in Definition \ref{SS}. 
            \item A multimorphism $F \in \mathrm{Mul}_{\SpSeq}(E, E'; E'')$ corresponds to a \textbf{bilinear map of SS}, which is a family of bilinear pairings 
            \[F_r\colon E_r^{s, t} \times {(E')}_r^{s', t'} \to {(E'')}_r^{s + s', t + t'}\] 
            satisfying the Leibniz rule for each $d_r$, so that the induced pairing of $F_r$ on $d_r$-homology coincides with $F_{r + 1}$. Similarly, a multimorphism $F\in \mathrm{Mul}_{\SpSeq}(E_1,\ldots, E_n; E'')$ corresponds to a \textbf{multilinear map of SS}.
        \end{itemize}
    \end{itemize}
\end{construction}

\begin{lemma} \label{faitful-inclusion-universal-property-absolute}
    Suppose $G\colon \mathcal{D} \to \mathcal{E}$ is a faithful functor between $\infty$-categories (i.e. for every $x_1, x_2 \in \CD$, $\Map_{\CD}(x_1, x_2) \to \Map_{\CE}(Gx_1, Gx_2)$ is a monomorphism\footnote{i.e. inclusion of connected components} between spaces). Then for each $\infty$-category $\CC$, the functor
    \[\Fun(\CC, \CD)^{\simeq} \to \Fun(\CC, \CE)^{\simeq} \times_{\Fun(\CC^{\simeq}, \CE^\simeq)}\Fun(\CC^{\simeq}, \CD^\simeq)\] 
    is a monomorphism between spaces, and $(H\colon \CC \to \CE, F\colon \CC^\simeq \to \CD^\simeq)$ lies in the image iff for every $z_1, z_2 \in \CC$, $\Map_{\CC}(z_1, z_2) \to \Map_{\CE}(Hz_1, Hz_2)$ factors through the inclusion $\Map_{\CD}(F(z_1), F(z_2)) \xrightarrow{G} \Map_{\CE}(GFz_1, GFz_2) \simeq \Map_{\CE}(Hz_1, Hz_2)$.
\end{lemma}

\begin{proof}
    To prove this we utilize the language of \emph{flagged $(\infty, 1)$-categories} in \cite{AF18}. 
    \begin{itemize}
        \item A \textbf{flagged $\infty$-category} is a pair $(\CC, A \to \CC)$ where $\CC$ is an $\infty$-category, $A$ is an $\infty$-groupoid and $A \to \CC$ is an effective epimorphism\footnote{i.e. surjection on $\pi_0$} into the core $\infty$-groupoid $\CC^\simeq$. Such a flagged $\infty$-category is \textbf{univalent} if $A \to \CC^\simeq$ is an equivalence. We write $\Cat^{\mathrm{f}}$ for the $\infty$-category of flagged $\infty$-categories, in which the full subcategory of univalent ones is canonically identified with the $\infty$-category of $\infty$-categories $\Cat$.
        \item A \textbf{simplicial space} is a functor $X\colon \Delta^{\op} \to \mathsf{Spc}$ into the $\infty$-category of spaces. Such $X$ is a \textbf{Segal space} if the natural comparison maps $X_n \to X_1^{\times_{X_0} n}$ are equivalences, and it is furthermore \textbf{complete} if the natural map $X_0 \to X_3 \times_{X_{1} \times X_1} (X_0 \times X_0)$ is an equivalence. We write $\mathsf{Seg}(\mathsf{Spc})$ for the $\infty$-category of Segal spaces and write $\mathsf{CSeg}(\mathsf{Spc})$ for complete Segal spaces.
        \item The Rezk nerve functor 
        \[\mathrm{N}^{\mathrm{R}}\colon \Cat^{\mathrm{f}} \to \mathsf{Spc}^{\Delta^{\op}}, (\CC, A \to \CC) \mapsto ([n] \mapsto \Map_{\Cat^{\mathrm{f}}}(([n], [n]^\simeq \to [n]), (\CC, A \to \CC)))\]
        induces an equivalence between $\Cat^{\mathrm{f}}$ and $\mathsf{Seg}(\mathsf{Spc})$, which further restricts to an equivalence between $\Cat$ and $\mathsf{CSeg}(\mathsf{Spc})$, cf. \cite[Theorem 0.26]{AF18}.
    \end{itemize}
    In this language, the LHS space of the functor
    \[\Fun(\CC, \CD)^{\simeq} \to \Fun(\CC, \CE)^{\simeq} \times_{\Fun(\CC^{\simeq}, \CE^\simeq)}\Fun(\CC^{\simeq}, \CD^\simeq)\] 
    is the mapping space in $\Cat^{\mathrm{f}}$ from $(\CC, \CC^\simeq \to \CC)$ to $(\CD, \CD^\simeq \to \CD)$, while the RHS space is the mapping space from $(\CC, \CC^\simeq \to \CC)$ to $(\CE, \CD^\simeq \to \CD \xrightarrow{G} \CE)$. In the language of Segal spaces, $\mathrm{N}^{\mathrm{R}}(\CE, \CD^\simeq \to \CE)$ is given by the pullback $\mathrm{N}^{\mathrm{R}}(\CE, \CE^\simeq \to \CE) \times_{\mathrm{cosk}_0 \CE^\simeq} \mathrm{cosk}_0 \CD^\simeq$, where $\mathrm{cosk}_0\colon \mathsf{Spc} \to \Fun(\Delta^{\op}, \mathsf{Spc})$ is right adjoint to the ``evaluation at $[0]$'' functor. Consider the map 
    \[X = \mathrm{N}^{\mathrm{R}}(\CD, \CD^\simeq \to \CD) \to Y = \mathrm{N}^{\mathrm{R}}(\CE, \CD^\simeq \to \CE).\] 
    We have $X_0 = \CD^{\simeq} = Y_0$, and for every pair $x_1, x_2 \in \CD$, faithfulness of $G$ implies the induced map on fibers $X_1 \times_{X_0 \times X_0} \{(x_1, x_2)\} = \Map_{\CD}(x_1, x_2) \to \Map_{\CE}(Gx_1, Gx_2) = Y_1 \times_{X_0 \times X_0} \{(x_1, x_2)\}$ is a monomorphism. By Segalness, we see $X_n \to Y_n$ is a monomorphism for all $n$, i.e. $X \to Y$ is a monomorphism in $\mathsf{Seg}(\mathsf{Spc})$. Thus, the above functor between mapping spaces is a monomorphism. Also, a map $\mathrm{N}^{\mathrm{R}}(\CC, \CC^\simeq \to \CC) \to 
    Y = \mathrm{N}^{\mathrm{R}}(\CE, \CD^\simeq \to \CE)$ factors through $X = \mathrm{N}^{\mathrm{R}}(\CD, \CD^\simeq \to \CD)$ iff for each $[n] \in \Delta$ the map $\mathrm{N}^{\mathrm{R}}(\CC, \CC^\simeq \to \CC)_n \to Y_n$ factors through $X_n$. By Segalness this boils down to the case $n = 1$, which is equivalent to the factorization property in the statement.
\end{proof}

\begin{corollary} \label{faitful-inclusion-universal-property-relative}
    Suppose $\CF$ is an $\infty$-category and $G\colon \CD \to \CE$ is a faithful functor over $\CF$. Then for any $\infty$-category $\CC$ over $\CF$, the functor
    \[\Fun_{/ \CF}(\CC, \CD)^{\simeq} \to \Fun_{/ \CF}(\CC, \CE)^{\simeq} \times_{\Fun_{/ \CF^\simeq}(\CC^{\simeq}, \CE^\simeq)}\Fun_{/ \CF^\simeq}(\CC^{\simeq}, \CD^\simeq)\]
    is a monomorphism between spaces. A pair of functors $(H\colon \CC \to \CE, F\colon \CC^\simeq \to \CD^\simeq)$, respectively over $\CF$ and over $\CF^{\simeq}$, lies in the image iff it satisfies the description in Lemma \ref{faitful-inclusion-universal-property-absolute}.
\end{corollary}
\begin{proof}
    It follows by the same argument as in Lemma \ref{faitful-inclusion-universal-property-absolute} applied to the equivalence between slice categories $\mathrm{N}^{\mathrm{R}}\colon\Cat^{\mathrm{f}}_{/ (\CF, \CF^{\simeq} \to \CF)} \cong \mathsf{Seg}(\mathsf{Spc})_{/ \mathrm{N}^{\mathrm{R}}(\CF, \CF^{\simeq} \to \CF) }$.
\end{proof}

\begin{corollary} \label{Lifting-lax-sym-mon-functors-along-a-faithful-map} 
    Suppose $\CD, \CE$ are $\infty$-operads\footnote{By definition, an $\infty$-operad $\CC$ is an inert-cocartesian fibration $\CC^{\otimes} \to \Fin_*$ satisfying certain properties. We will also write $\CC$ for the \emph{underlying category} $\CC^{\otimes}_{\langle 1 \rangle}$, while we refer to $\CC^{\otimes}$ as the \emph{category of operators}.} and $G\colon \CD \to \CE$ is a lax symmetric monoidal functor (i.e. a map of $\infty$-operads), so that for any $x_1,\ldots, x_n, y \in \CD$, the map between multimapping spaces $\mathrm{Mul}_{\CD}(x_1,\ldots, x_n; y) \to \mathrm{Mul}_{\CE}(Gx_1,\ldots, Gx_n; Gy)$ is a monomorphism. Then for any $\infty$-operad $\CC$, the functor \footnote{Here $\Fun(\CC, \CD)$ consists of functors between the underlying categories, while $\Fun^{\otimes\lax}(\CC, \CD)$ consists of the maps between operads (namely, the lax symmetric monoidal functors).}
    \[\Fun^{\otimes\lax}(\CC, \CD)^{\simeq} \to \Fun^{\otimes\lax}(\CC, \CE)^{\simeq} \times_{\Fun(\CC, \CE)^\simeq}\Fun(\CC, \CD)^\simeq\] 
    is a monomorphism of spaces, and $(H\colon \CC \to \CE, F\colon \CC \to \CD)$ lies in the image iff for every tuple $z_1,\ldots, z_n, w \in \CC$, the map $\mathrm{Mul}_{\CC}(z_1,\ldots, z_n; w) \to \mathrm{Mul}_{\CE}(Hz_1,\ldots, Hz_n; Hw)$ factors through the inclusion $\mathrm{Mul}_{\CD}(F(z_1),\ldots, F(z_n); F(w)) \xrightarrow{G} \mathrm{Mul}_{\CE}(H(z_1),\ldots, H(z_n); H(w))$. 
\end{corollary}

\begin{proof}
    We first prove the functor above is a monomorphism between spaces. The assumption implies that the functor between categories of operators $G\colon \CD^{\otimes} \to \CE^{\otimes}$ is faithful, and an arrow in $\CD^{\otimes}$ is inert-cocartesian iff it maps to an inert-cocartesian arrow in $\CE^{\otimes}$. Furthermore, Lemma \ref{faitful-inclusion-universal-property-absolute} applies to the underlying categories, so the functor 
    \[\Fun(\CC, \CD)^{\simeq} \to \Fun(\CC, \CE)^{\simeq} \times_{\Fun(\CC^{\simeq}, \CE^\simeq)}\Fun(\CC^{\simeq}, \CD^\simeq)\] 
    is a monomorphism between spaces. Thus, it suffices to show
    \[\Fun^{\otimes\lax}(\CC, \CE)^{\simeq} \to \Fun^{\otimes\lax}(\CC, \CD)^{\simeq} \times_{\Fun(\CC^{\simeq}, \CD^\simeq)}\Fun(\CC^{\simeq}, \CE^\simeq)\] 
    is a monomorphism. For convenience, we utilize the language of \emph{algebraic patterns} in \cite{CH21}. Under this formalism, an $\infty$-operad $\CC^\otimes \to \Fin_*$ is a weak Segal fibration over the pattern $\CO = \Fin_*^{\flat}$. Write $\Fin_*^{\mathrm{int}}$ for the wide subcategory of the $1$-category of pointed finite sets $\Fin_*$ spanned by inert morphisms, then the pullback $\CC^{\mathrm{int}} = \CC^\otimes \times_{\Fin_*}\Fin_*^{\mathrm{int}}$ is a Segal fibration over $\CO^{\mathrm{int}} = \Fin_*^{\mathrm{int}} \cap \Fin_*^{\flat}$. Using the adjoint equivalence $\Fun(\CO^{\mathrm{el}}, \Cat) \cong \mathsf{Seg}_{\CO^{\mathrm{int}}}(\Cat)$ in \cite[the second paragraph in \S~8]{CH21}, we deduce that for every pair of $\infty$-operads $\CC, \CD$, the restriction functor 
    \[\Fun_{/ \Fin_*^{\mathrm{int}}}^{\mathrm{cc}}(\CC^{\mathrm{int}}, \CD^{\mathrm{int}})^\simeq \to \Fun(\CC^\simeq, \CD^\simeq)\] 
    is an equivalence, where the LHS is the sub-$\infty$-groupoid of $\Fun_{/ \Fin_*^{\mathrm{int}}}(\CC^{\mathrm{int}}, \CD^{\mathrm{int}})^\simeq$ spanned by functors over $\Fin_*^{\mathrm{int}}$ which preserve cocartesian arrows. Furthermore, recall that  $\Fun^{\otimes \lax}(\CC, \CD) = \Fun_{ / \Fin_*}^{\text{int.cc}}(\CC^\otimes, \CD^\otimes)$ is the full sub-$\infty$-category spanned by functors over $\Fin_*$ which preserve cocartesian arrows over inert maps in $\Fin_*$. Now consider the consecutive functors
    \begin{align*}
        \Fun_{ / \Fin_*}^{\text{int.cc}}(\CC^\otimes, \CD^\otimes)^\simeq &\to \Fun_{ / \Fin_*}^{\text{int.cc}}(\CC^\otimes, \CE^\otimes)^\simeq \times_{\Fun_{/ \Fin_*^{\mathrm{int}}}^{\mathrm{cc}}(\CC^{\mathrm{int}}, \CE^{\mathrm{int}})^\simeq}\Fun_{/ \Fin_*^{\mathrm{int}}}^{\mathrm{cc}}(\CC^{\mathrm{int}}, \CD^{\mathrm{int}})^\simeq \\
        &\to \Fun_{ / \Fin_*}^{\text{int.cc}}(\CC^\otimes, \CE^\otimes)^\simeq \times_{\Fun_{/\Fin_*^\simeq}(\CC^{\otimes, \simeq}, \CE^{\otimes,\simeq})}\Fun_{/\Fin_*^\simeq}(\CC^{\otimes, \simeq}, \CD^{\otimes,\simeq})\\
        &\to \Fun_{ / \Fin_*}^{\text{int.cc}}(\CC^\otimes, \CE^\otimes)^\simeq \times_{\Fun(\CC^{\simeq}, \CE^\simeq)}\Fun(\CC^{\simeq}, \CD^\simeq).
    \end{align*}
    The composite of the last two functors is an equivalence by the discussion above, while the composite of the first two functors is a monomorphism between spaces due to Corollary \ref{faitful-inclusion-universal-property-relative}. Therefore, the total composite is a monomorphism. \parr 

    It remains to identify the image of 
    \[\Fun^{\otimes\lax}(\CC, \CD)^{\simeq} \to \Fun^{\otimes\lax}(\CC, \CE)^{\simeq} \times_{\Fun(\CC, \CE)^\simeq}\Fun(\CC, \CD)^\simeq.\] 
    It is clear that a pair $(H\colon \CC \to \CE, F\colon \CC \to \CD)$ in the image would satisfy the condition in the statement. Conversely, suppose $(H, F)$ satisfies this condition. We can extend $F\colon \CC^\simeq \to \CD^{\simeq}$ to a functor $\CC^{\mathrm{int}} \to \CD^{\mathrm{int}}$ over $\Fin_*^{\mathrm{int}}$ due to the discussion on algebraic patterns above, and then restrict this to a functor $\CC^{\otimes, \simeq} \to \CD^{\otimes, \simeq}$ over $\Fin_*^\simeq$. Feeding this functor, together with the functor between the category of operators $H\colon \CC^{\otimes} \to \CE^{\otimes}$ over $\Fin_*$, into Corollary \ref{faitful-inclusion-universal-property-relative}, we obtain a functor $\widehat{F}\colon \CC^{\otimes} \to \CD^{\otimes}$ over $\Fin_*$ whose restrictions recover the pair $(H, F)$. For it to be a lax symmetric monoidal functor, it remains to show $\widehat{F}$ preserves inert-cocartesian arrows, which holds true since an arrow is inert-cocartesian in $\CD^\otimes$ iff its image is so in $\CE^\otimes$, while $H = G\widehat{F}$ has this property.
\end{proof}

\begin{lemma} \label{epimorphsim-domination-lemma}
    Suppose $F\colon \CC \to \CD$ is a functor between $\infty$-categories whose target $\CD$ is a $1$-category, $A$ is an $\infty$-groupoid\footnote{We usually take such $A$ to be the set $\pi_0(\CC^\simeq)$.}, and $i\colon A \to \CC$ is an essential surjection onto $\CC^\simeq$. Take $n \in \Nb_{\geq 1} \cup \{\infty\}$.
    \begin{itemize}
        \item Write $\Fun(\CC, \CD)^{\mathrm{epi}}_{F / }$ for the full subcategory of the slice category $\Fun(\CC, \CD)_{F / }$ spanned by natural transformations $\alpha\colon F \Rightarrow G$ so that $\alpha_x\colon Fx \to Gx$ is an epimorphism\footnote{Here an \emph{epimorphism} in $\CD$ is an arrow $z_1 \to z_2$ so that for every $w \in \CD$, the induced arrow $\Map_{\CD}(z_2, w) \to \Map_{\CD}(z_1, w)$ is a monomorphism between spaces (actually, sets).} for each $x \in \CC$. Then the restriction functor 
        \[\Fun(\CC, \CD)^{\mathrm{epi}}_{F / } \to \Fun(A, \CD)^{\mathrm{epi}}_{F i / }\]
        is fully faithful, and $(G\colon A \to \CD, \beta\colon Fi \Rightarrow G)$ lies in the essential image iff for any $a, b \in A$ and any arrow $i(a) \to i(b)$ in $\CC$, the composite $Fi(a) \to Fi(b) \to Gb$ factors through $Ga$. 
        \item Suppose furthermore $\CC, \CD$ are $\Eb_n$-monoidal and the functor $F$ is oplax $\Eb_n$-monoidal. Write $\Fun^{\otimes \oplax}(\CC, \CD)^{\mathrm{epi}}_{F / }$ for the pullback $\Fun^{\otimes \oplax}(\CC, \CD)_{F / } \times_{\Fun(\CC, \CD)_{F / }}\Fun(\CC, \CD)^{\mathrm{epi}}_{F / }$. Then
        \[\Fun^{\otimes \oplax}(\CC, \CD)^{\mathrm{epi}}_{F / } \to \Fun(\CC, \CD)^{\mathrm{epi}}_{F / }\]
        is fully faithful, and a pair $(G\colon \CC \to \CD, \beta\colon F \Rightarrow G)$ lies in the essential image iff $F(\oneb_{\CC}) \to \oneb_{\CD}$ factors through $G(\oneb_{\CC})$, and for each pair $x, y \in \CC$ the composite $F(x \otimes y) \to Fx \otimes Fy \to Gx \otimes Gy$ factors through $G(x\otimes y)$.
        \item Alternatively, suppose furthermore $\CC, \CD$ are $\Eb_n$-monoidal, the tensor product on $\CD$ preserves epimorphisms in each variable, and the functor $F$ is lax $\Eb_n$-monoidal. Write $\Fun^{\otimes \lax}(\CC, \CD)^{\mathrm{epi}}_{F / }$ for the pullback $\Fun^{\otimes \lax}(\CC, \CD)_{F / } \times_{\Fun(\CC, \CD)_{F / }}\Fun(\CC, \CD)^{\mathrm{epi}}_{F / }$. Then
        \[\Fun^{\otimes \lax}(\CC, \CD)^{\mathrm{epi}}_{F / } \to \Fun(\CC, \CD)^{\mathrm{epi}}_{F / }\]
        is fully faithful, and a pair $(G\colon \CC \to \CD, \beta\colon F \Rightarrow G)$ lies in the essential image iff for each $x, y \in \CC$, the composite $Fx \otimes Fy \to F(x \otimes y) \to G(x\otimes y)$ factors through $Gx \otimes Gy$.
    \end{itemize} 
\end{lemma}

\begin{proof}
    As $\CD$ is a $1$-category, we can replace $\CC$ and $A$ with the homotopy $1$-categories $h\CC$ and $h A$. Now that everything is $1$-categorical, the proof can be done through a direct verification. 
\end{proof}

\begin{theorem} \label{categorified-total-Leibniz-rule}
    Take $r \in \Nb$, $r \geq 1$.
    \begin{itemize}
        \item The assignment $A^{*,*}_{r + 1}\colon \Fil\Sp \to \Gr\Ch(\Ab)^{E_{r + 1}}, X \mapsto \{(A^{s, t}_{r + 1}(X) = \pi_{t - s, t}(X / \defopara^r), d_{r + 1} = \delta_r^r)\}$ determines a lax symmetric monoidal functor. 
        \item The assignment $E^{*,*}_{r + 1}\colon \Fil\Sp \to \Gr\Ch(\Ab)^{E_{r + 1}}, X \mapsto \{(E_{r + 1}^{s, t}(X), d_{r + 1})\}$ determines (uniquely) a lax symmetric monoidal functor, such that the natural epimorphisms $\pi_{t- s, t}(X / \defopara^r) \to E_{r + 1}^{s, t}(X)$ assemble into a symmetric monoidal natural transformation $\alpha\colon A_{r + 1}^{*,*} \Rightarrow E_{r + 1}^{*,*}$.  
        \item The functors $\{E^{*,*}_{r + 1}\}_{r \geq 1}$ assemble into a lax symmetric monoidal functor $E_*^{*,*}\colon \Fil\Sp \to \SpSeq$.
    \end{itemize}
\end{theorem}

\begin{proof}
    Putting the Beilinson homotopy objects $\{\pi_n^B(\Ac_r(\Sigma^{0, w}X / \defopara^r))\}_{n \in \Zb, 0 \leq w < r}$ together, we obtain a functor $A^{*,*}_{r + 1}\colon \Fil\Sp \to \Gr\Ch(\Ab)^{E_{r + 1}}$ which matches the above concrete description. The composite $U A^{*,*}_{r + 1}\colon \Fil\Sp \to \Gr(\Gr(\Ab))$ is also the composite of the functor $X \mapsto \pi_{**}(X / \defopara^r)$ with the reindexing equivalence $(v_{r + 1}^*)^{-1}v^*\colon \Gr(\Gr(\Ab)^{\Day})^{\Kos} \cong \Gr(\Gr(\Ab)^{\Kos})^{\Kos}$, in which $v^*$ comes from Lemma \ref{products-on-the-E2-page} and $v_{r + 1}^*$ comes from Construction \ref{SS-as-an-operad}. Thus, $UA^{*,*}_{r + 1}$ admits a lax symmetric monoidal structure due to Example \ref{lax-symmetric-monoidality-of-pi-star}. Furthermore, for every map $X_1\otimes \cdots\otimes X_n \to T$ in $\Fil\Sp$, the induced multilinear pairing $A^{*,*}_{r + 1}(X_1) \times \cdots \times A^{*,*}_{r + 1}(X_n) \to A^{*,*}_{r + 1}(T)$ satisfies the Leibniz rule (which follows from the first half of Theorem \ref{Burklund's-Leibniz-rule} and Remark \ref{Burklund's-Leibniz-rule-multilinear}). We can thus lift the lax symmetric monoidal structure on $U A_{r + 1}^{*,*}$ uniquely to $A^{*,*}_{r + 1}$ via Corollary \ref{faitful-inclusion-universal-property-relative}. \parr 
    
    We can then apply Lemma \ref{epimorphsim-domination-lemma} items 1 and 3 with $A = \pi_0(\Fil\Sp^\simeq)$ (whose assumptions are satisfied according to the proof of the second half of Theorem \ref{Burklund's-Leibniz-rule}) to the family of epimorphisms $A^{s, t}_{r + 1}(X) \to E^{s, t}_{r + 1}(X)$, which leads to a unique lax symmetric monoidal functor $E_{r + 1}^{*,*}\colon \Fil\Sp \to \Gr\Ch(\Ab)^{E_{r + 1}}$ equipped with a symmetric monoidal natural transformation $A_{r + 1}^{*,*} \Rightarrow E_{r + 1}^{*,*}$. \parr 

    To assemble the family $\{E^{*,*}_{r + 1}\}_{r \geq 1}$ into one lax symmetric monoidal functor $E_*^{*,*}\colon \Fil\Sp \to \SpSeq$, it remains to construct a commutative square of lax symmetric monoidal functors
    \[\begin{tikzcd}
        \Fil\Sp \ar[r, "{E_{r + 2}^{*,*}}"] \ar[d, "{E_{r + 1}^{*,*}}"] & \Gr\Ch(\Ab)^{E_{r + 2}} \ar[d, "{U}"] \\
        \Gr\Ch(\Ab)^{E_{r + 1}} \ar[r, "{H}"] & \Gr(\Gr(\Ab)^{\Kos})^{\Kos}
    \end{tikzcd}\]
    for each $r \geq 1$. By Lemma \ref{epimorphsim-domination-lemma} item 1 and item 3 with $A = \pi_0(\Fil\Sp^\simeq)$, it suffices to find a lax symmetric monoidal functor $F\colon \Fil\Sp \to \Gr(\Gr(\Ab)^{\Kos})^{\Kos}$ with two entrywise epic symmetric monoidal natural transformations $F \Rightarrow H E_{r + 1}^{*,*}$, $F \Rightarrow U E_{r + 2}^{*,*}$, so that for each $X \in \Fil\Sp$ the two arrows $F^{s, t}(X) \to H(E_{r + 1}^{*, *})^{s,t}(X), F^{s, t}(X) \to E_{r + 2}^{s,t}(X)$ have the same kernel. We take our $F$ to be $U A_{r + 2}^{*,*}$, so $F \Rightarrow U E_{r + 2}^{*,*}$ comes for free. On the other hand, according to the proof of Theorem \ref{standard-SS}, the natural transformation $U A_{r + 2}^{*,*} \Rightarrow U A_{r + 1}^{*,*} \Rightarrow U E_{r + 1}^{*,*}$ factors through $ZE_{r + 1}^{*,*}$. The induced transformation $U A_{r + 2}^{*,*} \Rightarrow ZE_{r + 1}^{*,*}$ is entrywise epic, and it can be enhanced uniquely into a symmetric monoidal natural transformation via Lemma \ref{epimorphsim-domination-lemma} item 2 applied to $UE_{r + 1}^{*,*}\colon \Fil\Sp^{\op} \to (\Gr(\Gr(\Ab)^{\Kos})^{\Kos})^{\op}$. We thus obtain the other entrywise epic symmetric monoidal natural transformation $F = U A_{r + 2}^{*,*} \Rightarrow ZE_{r + 1}^{*,*} \Rightarrow HE_{r + 1}^{*,*}$. It is then clear from the construction in Theorem \ref{standard-SS} that for every $X \in \Fil\Sp, s, t \in \Zb$, the corresponding kernels for the two transformations are both given by the image of 
    \[\defopara \oplus \delta_{r}^{r + 1} \colon \pi_{t - s, t}(X / \defopara^{r}) \oplus \pi_{t - s + 1, t - r}(X / \defopara^{r}) \to \pi_{t - s, t}(X / \defopara^{r + 1}) = F^{s, t}.\]
    Therefore, we can conclude with the lax symmetric monoidal functor $E_*^{*,*}\colon \Fil\Sp \to \SpSeq$.
\end{proof}

We wrap up this appendix by generalizing Theorem \ref{categorified-total-Leibniz-rule} to the case where $\Sp$ is replaced by some suitable $t$-$\infty$-category $\CE$, continuing the style in Corollary \ref{standard-SS-with-general-coefficients} and Corollary \ref{Bockstein-dictionary-finite-page-with-general-coefficients}. 

\begin{definition}
    Take $n \in \Nb_{\geq 1} \cup \{\infty\}$. We say $\CE$ is a \textbf{presentably $\Eb_n$-monoidal $t$-$\infty$-category} if $\CE$ is a stable $\Eb_n$-monoidal category with a $t$-structure $(\CE_{\geq 0}, \CE_{\leq 0})$, so that the following holds true: 
    \begin{itemize}
        \item The underlying $\infty$-category $\CE$ is presentable, and the $t$-structure is accessible in the sense of \cite[Definition 1.4.4.12]{HA}, namely the sub-$\infty$-category $\CE_{\geq 0}$ is also presentable. 
        \item The tensor product functor $\otimes\colon \CE \times \CE \to \CE$ preserves colimits separately in each variable. 
        \item The $t$-structure is compatible with the tensor product, i.e. the unit in $\CE$ is connective, and the tensor product of two connective objects is also connective. 
    \end{itemize}
\end{definition}

\begin{corollary} \label{categorified-total-Leibniz-rule-with-general-coefficient}
    Take $n \in \Nb_{\geq 1} \cup \{\infty\}$. Suppose $\CE$ is a presentably $\Eb_n$-monoidal $t$-$\infty$-category\footnote{in fact, here we do not need the $t$-structure to be accessible}. We can construct the $\Eb_n$-monoidal $\infty$-categories $\{\Gr\Ch(\CE^\heartsuit)^{E_{r + 1}}\}_{r \geq 1}$ and assemble them into an $\infty$-operad $\SpSeq(\CE^\heartsuit)$ over $\Eb_n$ by the same discussion as in Construction \ref{SS-as-an-operad}. Furthermore, 
    \begin{itemize}
        \item For each $r \in \Nb_{\geq 1}$, the assignment
        \[A^{*,*}_{r + 1}\colon \Fil(\CE) \to \Gr\Ch(\CE^\heartsuit)^{E_{r + 1}}, X \mapsto \{(A^{s, t}_{r + 1}(X) = \pi_{t - s, t}(X / \defopara^r), d_{r + 1} = \delta_r^r)\}\] 
        determines a lax $\Eb_n$-monoidal functor. 
        \item The assignment $E^{*,*}_{r + 1}\colon  \Fil(\CE) \to \Gr\Ch(\CE^\heartsuit)^{E_{r + 1}}, X \mapsto \{(E_{r + 1}^{s, t}(X), d_{r + 1})\}$ determines (uniquely) a lax $\Eb_n$-monoidal functor, such that the natural epimorphisms $\pi_{t- s, t}(X / \defopara^r) \to E_{r + 1}^{s, t}(X)$ assemble into an $\Eb_n$-monoidal natural transformation $\alpha\colon A^{*,*}_{r + 1} \Rightarrow E_{r + 1}^{*,*}$.  
        \item The functors $\{E^{*,*}_{r + 1}\}_{r \geq 1}$ assemble into a lax $\Eb_n$-monoidal functor\footnote{i.e. a map of $\infty$-operads over $\Eb_n$} $E_*^{*,*}\colon \Fil(\CE) \to \SpSeq(\CE^\heartsuit)$.
    \end{itemize}
    These statements follow by the same argument as in Theorem \ref{Burklund's-Leibniz-rule} and Theorem \ref{categorified-total-Leibniz-rule}.
\end{corollary}

In practice, when one extracts graded homotopy groups and spectral sequences from a filtration $X \in \Fil(\CE)$, the $\Zb$-grading alone is often too coarse. For example, much of the nontrivial structure in equivariant spectral sequences is carried by the $\RO(G)$-degree, and more generally one is led to consider gradings by the Picard groupoid, or even finer grading data. This makes it necessary to establish a lax $\Eb_n$-monoidality result for SS with such general gradings, as well as for their total-differential refinements. In the rest of this appendix, we give two versions of this result.

\begin{setup}[Coarse grading data] \label{picard-grading-coarse-setup}
    Take $n \in \Nb_{\geq 1} \cup \{\infty\}$. We fix the following data:
    \begin{itemize} 
        \item Let $\CE$ be a stable presentably $\Eb_{n + 1}$-monoidal $\infty$-category. 
        \item Let $\CF$ be a presentably $\Eb_{n}$-monoidal $t$-$\infty$-category, and let $T\colon \CE \to \CF$ be an exact lax $\Eb_{n}$-monoidal functor. 
        \item Let $\CI$ be an $\Eb_n$-monoidal $\infty$-category, equipped with two $\Eb_n$-monoidal functors 
        \[P\colon \CI \to \CE \qquad \text{and} \qquad Q\colon \CI \to \Zb^\delta.\]
        \item Let $\CK$ be an $\Eb_n$-monoidal $1$-category equipped with an $\Eb_n$-monoidal functor $R\colon \CK \to h\CI$. We refer to the two composite functors $PR\colon \CK \to h\CI \to h \CE$, $QR\colon \CK \to h\CI \to \Zb^\delta$ respectively as $V \mapsto \Sb^V$ and $V \mapsto |V|$. 
    \end{itemize}
    In practice, usually we have $\CE = \CF$ and $\CK = h\CI$, and we will refer to such a full setup by saying $\CK \to h\CE$ is part of a \textbf{coarse grading datum}.
\end{setup}

\begin{example} \label{coarse-picard-grading-examples}
    Take $n \in \Nb_{\geq 1} \cup \{\infty\}$.
    \begin{itemize}
        \item Let $\CE$ be a presentably $\Eb_{n + 1}$-monoidal $t$-$\infty$-category, equipped with a map\footnote{One can always take $Q = 0$; however, the practical choice of $Q$ is usually the ``virtual dimension''.} of abelian groups $Q \colon \pi_0 \Pic(\CE) \to \Zb^\delta$. Take $\CF = \CE$, $\CI = \Pic(\CE)$, and $\CK = h\CI$. Then the collection satisfies the requirements in Setup \ref{picard-grading-coarse-setup}, with $P$ being the inclusion and $R, T$ being $\id$. 
        \item Let $G$ be a finite group. Take $\CE = \Sp^G = \CF$, the symmetric monoidal $\infty$-category of genuine $G$-spectra with its Mackey $t$-structure. Furthermore, take $\CI = \Pic(\Sp^G)$ and take $\CK = \PG$ to be the pullback of the cospan $\RO(G) \to \pi_0 \Pic(\Sp^G) \gets h \Pic(\Sp^G)$. Writing $P\colon \CI \to \CE$ for the inclusion, $Q\colon \CI \to \Zb^\delta$ to be the $0$-truncation of the forgetful functor $\Pic(\Sp^G) \to \Pic(\Sp)$, $R\colon \CK \to h\CI$ to be the natural comparison map, and $T = \id$, we obtain a collection that satisfies the requirements in Setup \ref{picard-grading-coarse-setup}. In this case, for $V = V_+ - V_- \in \RO(G)$, $|V| = |V_+| - |V_-| \in \Zb^\delta$ corresponds to the virtual dimension. 
    \end{itemize}
\end{example}

\begin{construction} \label{coarse-picard-grading-constructions-1}
    Under Setup \ref{picard-grading-coarse-setup}, we have the following constructions in $\CE$.
    \begin{itemize}
        \item Write $\hom_{\CE}\colon \CE \times \CE^{\op} \to \CE$ for the internal hom functor of $\CE$ (i.e. the two-variable right adjoint of $\otimes_{\CE}$), which is lax $\Eb_n$-monoidal due to \cite[Corollary 3.4.10]{HHLN1}.
        \item For $X \in \CE, V \in \CK, n \in \Zb$, we write $\Sigma^{V + n} X$ for the object $\Sigma^n(X \otimes \Sb^V) \in \CE$. Also, we write $\pi_{V + n}(X) = \pi_0T(\hom_{\CE}(\Sigma^n\Sb^V, X)) \in \CF^\heartsuit$ and refer to these as the \textbf{$\CK$-graded homotopy objects} of $X$. Note that if $V \in \CK$ has an inverse $-V$, then $\pi_{V + n}X = \pi_0 T(\Sigma^{-V-n} X)$.
        \item We denote by $\Gr_{\CI}(\CE)$ (resp. $\Gr_{\CK}(\CF^\heartsuit)$) the $\infty$-category $\Fun(\CI^{\op}, \CE)$ (resp. $\Fun(\CK^{\op}, \CF^\heartsuit)$), and equip it with the Day convolution $\Eb_n$-monoidal structure via \cite[Proposition 2.2.6.16]{HA}. Using the $\Eb_{n}$-monoidal composite $PR\colon V \mapsto \Sb^V$ we can construct a lax $\Eb_n$-monoidal functor $\pi_{\filledstar}\colon \CE \to \Gr_{\CK}(\CF^\heartsuit), X \mapsto (V \mapsto \pi_V(X))$ by transposing the lax $\Eb_n$-monoidal composite 
        \[\CE \times \CI^{\op} \xrightarrow{\CE \times P} \CE \times \CE^{\op} \xrightarrow{\hom_{\CE}(-, -)} \CE \xrightarrow{T} \CF \xrightarrow{\pi_0} \CF^\heartsuit\]
        and postcompose it with $R^*\colon \Gr_{\CI}(\CF^\heartsuit) = \Gr_{h\CI}(\CF^\heartsuit) \to \Gr_{\CK}(\CF^\heartsuit)$. 
    \end{itemize}
\end{construction}

\begin{construction} \label{coarse-picard-grading-constructions-2}
    Under Setup \ref{picard-grading-coarse-setup}, we have the following constructions in $\Fil(\CE)$. 
    \begin{itemize}
        \item The $\infty$-category $\Fil(\CE)$, equipped with the Day convolution $\Eb_{n + 1}$-monoidal structure, receives two $\Eb_{n + 1}$-monoidal left adjoints
        \[i\colon \Fil\Sp \to \Fil(\CE), \qquad j\colon \CE = \Fun(\{0\}, \CE) \xrightarrow{\mathrm{LKE}} \Fil(\CE).\]
        We write $\oneb = j(\oneb_{\CE}) = i(\Sb^{0, 0}) \in \Fil(\CE)$ for the tensor unit.  For $V \in \CK, n, w \in \Zb$, we write  $\Sb^{V + n, w} = j(\Sigma^n\Sb^V) \otimes i(\Sb^{0, w})$; furthermore, for $Y \in \Fil(\CE)$ we write $\Sigma^{V + n, w}Y$ for $Y \otimes \Sb^{V + n, w}$, and we write $\pi_{V + n, w}(Y) \in \CF^\heartsuit$ for $\pi_{V + n} Y(w)$. 
        \item We write $\hom_{\Fil(\CE)}$ for the internal hom functor on $\Fil(\CE)$, which is also lax $\Eb_n$-monoidal. By construction, $\pi_{V + n, w} Y = \pi_{0, 0}\hom_{\Fil(\CE)}(\Sb^{V + n, w}, Y)$.
        \item According to Dunn's additivity \cite[Theorem 5.1.2.2]{HA}, the functor $\otimes\colon \Fil(\CE) \times \Fil(\CE) \to \Fil(\CE)$ is $\Eb_n$-monoidal. From this, we get an $\Eb_n$-monoidal composite 
        \[\CI \xrightarrow{(P, Q)} \CE \times \Zb^\delta \to \CE \times \overrightarrow{\Zb} \to \CE \times \Fil\Sp \xrightarrow{i \times j} \Fil(\CE) \times \Fil(\CE) \xrightarrow{\otimes} \Fil(\CE)\]
        which sends each $R(V)$ for $V \in \CK$ to $\Sb^{V, |V|}$. Therefore, we can construct a lax $\Eb_n$-monoidal functor $\psi\colon \Fil(\CE) \to \Gr_{\CI}(\Fil(\CE))^{\Day}$ by transposing the lax $\Eb_n$-monoidal composite
        \[\Fil(\CE) \times \CI^{\op} \to \Fil(\CE) \times \Fil(\CE)^{\op} \xrightarrow{\hom_{\Fil(\CE)}(-, -)} \Fil(\CE).\]
        \item We write $\Gr_{\CK}(\Gr\Ch(\CF^\heartsuit))^{E_{r + 1}}$ for $\Fun(\CK^{\op}, \Gr\Ch(\CF^\heartsuit)^{E_{r + 1}})$ equipped with the Day convolution $\Eb_n$-monoidal structure. Also, we write $\Gr_{\CK} \SpSeq(\CF^\heartsuit)$ for $\Fun(\CK^{\op}, \SpSeq(\CF^\heartsuit))$ considered as an $\infty$-operad over $\Eb_n$ through Day convolution \cite[Construction 2.2.6.7]{HA}.
    \end{itemize}
\end{construction}

\begin{theorem} \label{categorified-total-Leibniz-rule-with-picard-trigrading}
    Under Setup \ref{picard-grading-coarse-setup}, the following statements hold true: 
    \begin{itemize}
        \item For each $r \in \Nb_{\geq 1}$, the assignment 
        \[A^{*,*, \filledstar}_{r + 1}\colon \Fil(\CE) \to \Gr_{\CK}(\Gr\Ch(\CF^\heartsuit))^{E_{r + 1}}, \quad X \mapsto \{(A^{s, t, V}_{r + 1}(X) = \pi_{V + t - s, |V| + t}(X / \defopara^r), d_{r + 1} = \delta_r^r)\}\]
        determines a lax $\Eb_n$-monoidal functor. 
        \item Write $E_{r + 1}^{s, t, V}(X) = E_{r + 1}^{s, t}(T\hom_{\Fil(\CE)}(\Sb^{V, |V|}, X)) \in \CF^\heartsuit$. For each $r \in \Nb_{\geq 1}$,  the assignment 
        \[E^{*,*, \filledstar}_{r + 1}\colon  \Fil(\CE) \to \Gr_{\CK}(\Gr\Ch(\CF^\heartsuit))^{E_{r + 1}}, \quad X \mapsto \{(E_{r + 1}^{s, t, V}(X), d_{r + 1})\}\] 
        determines (uniquely) a lax $\Eb_n$-monoidal functor, so that the epimorphisms 
        \[A^{s, t, V}_{r + 1}(X) = \pi_{t - s, t}(T\hom_{\Fil(\CE)}(\Sb^{V, |V|}, X) / \defopara^r) \to E_{r + 1}^{s, t}(T\hom_{\Fil(\CE)}(\Sb^{V, |V|}, X)) = E_{r + 1}^{s, t, V}(X)\] 
        assemble into an $\Eb_n$-monoidal natural transformation $\alpha\colon A_{r + 1}^{*,*, \filledstar} \Rightarrow E_{r + 1}^{*,*,\filledstar}$.  
        \item The functors $\{E_{r + 1}^{*,*, \filledstar}\}_{r \geq 1}$ induce a lax $\Eb_n$-monoidal functor $E_*^{*,*, \filledstar}\colon \Fil(\CE) \to \Gr_{\CK}\SpSeq(\CF^\heartsuit)$.
    \end{itemize}
\end{theorem}

\begin{proof}
    For the first assertion we use the lax $\Eb_n$-monoidal composite 
    \begin{align*}
        A_{r + 1}^{*,*, \filledstar}\colon \Fil(\CE) &\xrightarrow{\psi} \Gr_{\CI}(\Fil(\CE)) \xrightarrow{\Gr_{\CI}(T)} \Gr_{\CI}(\Fil(\CF))  \\
        &\xrightarrow{\Gr_{\CI}(A_{r + 1}^{*,*})} \Gr_{\CI}(\Gr\Ch(\CF^\heartsuit))^{E_{r + 1}} = \Gr_{h\CI}(\Gr\Ch(\CF^\heartsuit))^{E_{r + 1}} \to \Gr_{\CK}(\Gr\Ch(\CF^\heartsuit))^{E_{r + 1}}.
    \end{align*}
    Here $A_{r + 1}^{*,*}\colon \Fil(\CF) \to \Gr\Ch(\CF^\heartsuit)^{E_{r + 1}}$ is the lax $\Eb_n$-monoidal functor from Corollary \ref{categorified-total-Leibniz-rule-with-general-coefficient}. Replacing this $A_{r + 1}^{*,*}$ with $E_{r + 1}^{*,*}\colon \Fil(\CF) \to \Gr\Ch(\CF^\heartsuit)^{E_{r + 1}}$ and $E_*^{*,*} \colon \Fil(\CF) \to \SpSeq(\CF^\heartsuit)$, we obtain the lax $\Eb_n$-monoidal functors $E_{r + 1}^{*,*, \filledstar} \colon \Fil(\CE) \to \Gr_{\CK}(\Gr\Ch(\CF^\heartsuit))^{E_{r + 1}}$ and $E_*^{*,*, \filledstar}\colon \Fil(\CE) \to \Gr_{\CK}\SpSeq(\CF^\heartsuit)$. 
\end{proof}

% Redundancy.

\begin{example} \label{trigraded-picard-SS-examples}
    Take $n \in \Nb_{\geq 1} \cup \{\infty\}$.
    \begin{itemize}
        \item Let $\CE$ be a presentably $\Eb_{n + 1}$-monoidal $t$-$\infty$-category, equipped with a map of abelian groups $\pi_0 \Pic(\CE) \to \Zb^\delta$. Applying Theorem \ref{categorified-total-Leibniz-rule-with-picard-trigrading} to the coarse picard grading recipe in Example \ref{coarse-picard-grading-examples} item 1,  we obtain %the assignments $X \mapsto (A^{s, t, V}_{r}(X) = \pi_{V + t - s, |V| + t}(X / \defopara^{r - 1}), \delta_{r - 1}^{r - 1})$ and $X \mapsto (E^{s, t, V}_{r}(X) = E^{s, t}_r(\hom_{\Fil(\CE)}(\Sb^{V, |V|}), X)), d_r)$ assemble into 
        a lax $\Eb_{n}$-monoidal functor 
        \[E_*^{*,*,\filledstar}\colon \Fil(\CE) \to \Gr_{\Pic(\CE)}\SpSeq(\CE^\heartsuit), \quad X \mapsto \{E_r^{s, t, V}(X) = E_r^{s, t}(\hom_{\Fil(\CE)}(\Sb^{V, |V|}, X))
        \}%_{r \geq 2, s, t \in \Zb, V \in \Pic(\CE)}
        \] together with its (lax $\Eb_{n}$-monoidal) total differential refinements $A_{r}^{*,*,\filledstar}$. 
        \item Let $G$ be a finite group. Applying Theorem \ref{categorified-total-Leibniz-rule-with-picard-trigrading} to the coarse $\RO(G)$-grading recipe in Example \ref{coarse-picard-grading-examples} item 2,  we obtain %the assignments $X \mapsto (A^{s, t, V}_{r}(X) = \pi_{V + t - s, |V| + t}(X / \defopara^{r - 1}), \delta_{r - 1}^{r - 1})$ and $X \mapsto (E^{s, t, V}_{r}(X) = E^{s, t}_r(\hom_{\Fil(\CE)}(\Sb^{V, |V|}), X)), d_r)$ assemble into 
        a lax symmetric monoidal functor 
        \[E_*^{*,*,\filledstar}\colon \Fil(\Sp^G) \to \Gr_{\PG}\SpSeq(\Mack_G(\Ab)), \quad X \mapsto \{E_r^{s, t, V}(X) = E_r^{s, t}(\hom_{\Fil(\Sp^G)}(\Sb^{V, |V|}, X))
        \}%_{r \geq 2, s, t \in \Zb, V \in \RO(G)}
        \] together with its (lax symmetric monoidal) total differential refinements $A_{r}^{*,*,\filledstar}$.  
        % \item $\Mod_{\mathrm{KU}}(\Sp)\colon E_2^{s, t, V} = \pi_{V + t - s, t}(X / \defopara) = E_2^{s + 1, t, V + 1}$.
    \end{itemize}
\end{example}

The examples above have a common feature: there exists an object $1 \in \CK$ (necessarily different from the tensor unit, which we denote by $0$), such that $|1| = 1 \in \Zb^\delta$ and $\Sb^1 = \Sigma \oneb \in \CE$. If such an object exists, then for each $X \in \Fil(\CE)$ there are repetitions in the trigraded SS
\[E_r^{s, t, V}(X) = E_r^{s, t}(\hom_{\Fil(\CE)}(\Sb^{V, |V|}, X)) \cong E_r^{s, 0}(\hom_{\Fil(\CE)}(\Sb^{V + t, |V| + t}, X)) = E_r^{s, 0, t + V}(X).\] 
This suggests that the $\Zb \times \Zb \times \CK$-graded spectral sequence actually comes from a $\Zb \times \CK$-graded spectral sequence whose components are $E_r^{s, V}(X) = E_r^{s, 0, V}(X)$, in which there is no redundancy. One would then ask if this leads to a lax symmetric monoidal functor $X \mapsto (\{E_r^{s, V}(X)\}, d_r)_{r \geq 2}$. This is not as straightforward as it seems, due to the following two reasons: 
\begin{itemize}
    \item The $E_r$-page of $\Zb \times \Zb$-graded spectral sequences lies in the $1$-category $\Gr\Ch(\CE^\heartsuit)^{E_r}$, which is a reindexed variant of $\Gr(\Ch(\CE^\heartsuit)^{\Kos})^{\Kos}$. We expect the $E_r$-page of these $\Zb \times \CK$-graded spectral sequences to lie in a version of $\Gr_{\CK}(\Ch(\CE^\heartsuit)^{\Kos})$. A priori, there is only one reasonable symmetric monoidal structure on this category, defined through Day convolution, on which it is not clear how to describe the Koszul sign rule for the $\CK$-parameter. 
    \item On the other hand, to define the $\infty$-operad $\SpSeq$ one needs a reparametrization automorphism $v_r^*$ on $\Gr(\Gr(\CE^\heartsuit)^{\Kos})^{\Kos}$, which is induced by the symmetric monoidal functor $v_r\colon \Zb \times \Zb \to \Zb \times \Zb, (s, t) \mapsto (s + (r - 1)(t - s), t + (r - 1)(t - s))$. For the $\Zb \times \CK$-graded variant, we would also need a reparametrization automorphism $v_r^*$ on $\Gr_{\CK}(\Gr(\CE^\heartsuit)^{\Kos})$. One might expect this to come from a symmetric monoidal endofunctor on $\CK \times \Zb^\delta$, but this is not possible: usually there is no $\Eb_2$-monoidal functor $S\colon \Zb^\delta \to \CK$ so that $|S(1)| = 1$.
\end{itemize}
To resolve these issues, we replace $\Zb^\delta$ by a certain symmetric monoidal $1$-groupoid. The following construction is inspired by Quinn \cite{Qu25}.

\begin{construction}[Koszul grading as a Day convolution idempotent] \label{Koszul-grading-as-a-Day-idempotent}
    Take $n \in \Nb_{\geq 1} \cup \{\infty\}$. There is a symmetric monoidal $1$-groupoid $\Pe$ which unifies both Day and Koszul $\Eb_n$-monoidal structures. 
    \begin{itemize}
        \item The maps of connective spectra (here $J\colon \mathrm{ko} \to \pic(\Sb)$ is the $J$-homomorphism) 
        \[\Sb = \Sb^0 \xrightarrow{\mathrm{unit}} \mathrm{ko} \xrightarrow{J} \pic(\Sb) \to \pic(\mathrm{H}\Zb) \to \pic(\mathrm{H}\Fb_3) \]
        become isomorphisms upon applying $\tau_{\leq 1}$. We write 
        \[\Pe = \Omega^{\infty} \tau_{\leq 1}\pic(\Sb) = h \Pic(\Sb)\] 
        for the corresponding symmetric monoidal $1$-groupoid. Explicitly, $\Pe$ is (canonically) identified with $\Zb^\delta \times \mathrm{B}C_2$ as a monoidal $1$-category, while its symmetric monoidal structure is determined by the self braiding $1 + 1 \to 1 + 1$, which is $f \in \{e, f \mid f^2 = e\} = C_2 = \Aut_{\Pe}\!(2)$. 
        \item We write $\Gr_{\Pe}(\Ab)$ for the $1$-category $\Fun(\Pe^{\op}, \Ab)$ equipped with the Day convolution symmetric monoidal structure. As a monoidal $1$-category, this is (canonically) equivalent to $\Gr(\Ab^{\mathrm{B}C_2})$, so an object of it is a family $M = \{M^s\}_{s \in \Zb}$, here each $M^s$ is an abelian group with $C_2 = \{e, f\}$-action, and a morphism $M \to N$ is a family of maps compatible with the $C_2$-actions on both sides. For $M, N \in \Gr_{\Pe}(\Ab)$, the tensor product $M \otimes N$ is computed as follows:
        \[(M \otimes N)^s = \bigoplus_{p + q = s} (M^p \otimes N^q)_{C_2} = \bigoplus_{p + q = s} M^p \otimes N^q / (x \otimes y - fx \otimes fy)\]
        whose $C_2$-action is determined by $f \cdot (x \otimes y) = (fx) \otimes y = x \otimes (fy)$. The tensor unit in $\Gr_{\Pe}(\Ab)$ is $\Zb[C_2][0]$, the group ring of $C_2$ concentrated in degree $0$. The symmetric monoidal structure is determined by the self-braiding on $\Zb[C_2][1]$ (the group ring of $C_2$ concentrated in degree $1$), which is given by $f \in \Aut_{\Gr_{\Pe}(\Ab)}(\Zb[C_2][2])$. 
        \item The tensor unit $\Zb[C_2][0] \in \Gr_{\Pe}(\Ab)$ has two quotients $\Zb^{\triv}[0]$ and $\Zb^{\sgn}[0]$. Here $\Zb^{\triv}[0]$ is the coequalizer of two maps $f, \id\colon \Zb[C_2][0] \rightrightarrows \Zb[C_2][0]$, which is explicitly a copy of $\Zb$ in degree $0$ with trivial $C_2$-action; $\Zb^{\sgn}[0]$ is the coequalizer of $f, -\id\colon \Zb[C_2][0] \rightrightarrows \Zb[C_2][0]$, which is explicitly a copy of $\Zb$ in degree $0$ whose $C_2$-action is the sign action (i.e. $f x = - x$ for each $x$). These two are both idempotent algebras\footnote{Actually, in an abelian $1$-category with right-exact tensor product, any epimorphism $\oneb \to X$ becomes an isomorphism after tensoring $X$, so it equips $X$ with the structure of an idempotent algebra.} in $\Gr_{\Pe}(\Ab)$. Moreover, $\Mod_{\Zb^{\triv}[0]}$ is the full subcategory of $\Gr_{\Pe}(\Ab)$ consisting of objects $\{M^s\}_{s \in \Zb}$ whose $C_2$-action is trivial, and it identifies canonically with $\Gr(\Ab)^{\Day}$; $\Mod_{\Zb^{\sgn}[0]}$ is the full subcategory consisting of objects $\{M^s\}_{s \in \Zb}$ whose $C_2$-action is the sign action, and it identifies canonically with $\Gr(\Ab)^{\Kos}$. 
        \item In $\Gr_{\Pe}(\Gr_{\Pe}(\Ab)) = \Gr_{\Pe \times \Pe}(\Ab)$ the tensor unit is $\Zb[C_2 \times C_2][0, 0]$, which has four quotients $\Zb^{\epsilon_1, \epsilon_2}[0, 0]$ with $\epsilon_1, \epsilon_2 \in \{\triv, \sgn\}$, corresponding to the two generators $f_1 = (f, e), f_2 = (e, f) \in C_2 \times C_2$ acting by $\id$ or $-\id$. We also have $\Mod_{\Zb^{\epsilon_1, \epsilon_2}[0, 0]} = \Gr(\Gr(\Ab)^{\epsilon_2})^{\epsilon_1}$, where $\triv \leftrightarrow \Day$ and $\sgn \leftrightarrow \Kos$. Furthermore: 
        \begin{itemize}
            \item Consider the functor $v\colon \Pe \times \Pe \to \Pe \times \Pe, (n, a) \mapsto (n - a, a)$. It is symmetric monoidal since each of its two components is a composition of $\id$ with the two functors $\Pe \times \Pe \to \Pe, (x, y) \mapsto x + y$ and $\Pe \to \Pe, x \mapsto -x$. Therefore, it induces a symmetric monoidal equivalence $v^*\colon \Gr_{\Pe \times \Pe}(\Ab) \to \Gr_{\Pe \times \Pe}(\Ab)$. Furthermore, under this functor the maps $f_1, f_2 \in \Aut(\Zb[C_2 \times C_2][0, 0])$ are sent to $f_1$ and $f_1 f_2$, so the idempotent algebra $\Zb^{\sgn, \triv}[0, 0]$ goes to $\Zb^{\sgn, \sgn}[0, 0]$. Therefore, it induces a symmetric monoidal equivalence $v^*\colon \Gr(\Gr(\Ab)^{\Day})^{\Kos} \cong \Gr(\Gr(\Ab)^{\Kos})^{\Kos}$, which is used in Lemma \ref{products-on-the-E2-page}.
            \item For $r \in \Nb_{\geq 1}$, consider the functor $v_r\colon \Pe \times \Pe \to \Pe \times \Pe, (s, t) \mapsto (s + (r - 1)(t - s), t + (r - 1)(t - s))$, which is also symmetric monoidal. It induces another symmetric monoidal equivalence $v_r^*\colon \Gr_{\Pe \times \Pe}(\Ab) \to \Gr_{\Pe \times \Pe}(\Ab)$. Moreover, write $a \equiv r \text{ mod } 2$, then under this functor the maps $f_1, f_2 \in \Aut(\Zb[C_2 \times C_2][0, 0])$ are sent to $f_1^{a}f_2^{1 - a}$ and $f_1^{1 - a}f_2^{a}$, so $\Zb^{\sgn, \sgn}[0, 0]$ goes to $\Zb^{\sgn, \sgn}[0, 0]$. Therefore, it induces a symmetric monoidal equivalence $v_r^*\colon \Gr(\Gr(\Ab)^{\Kos})^{\Kos} \cong \Gr(\Gr(\Ab)^{\Kos})^{\Kos}$, which is used in Construction \ref{SS-as-an-operad}.
        \end{itemize}
        \item The discussion above works ($\Eb_n$-monoidally) with $\Ab$ replaced by any $\Eb_n$-monoidal abelian $1$-category $\CA$ whose tensor product is right exact. More generally, 
        \begin{itemize}
            \item If $\CA$ is cocomplete and $\otimes_{\CA}\colon \CA \times \CA \to \CA$ preserves colimits separately in each variable, then for any $\Eb_{n}$-monoidal $1$-category $\CK$ equipped with an $\Eb_{n}$-monoidal functor $\Pe \to \CK$, the abelian $1$-category $\Gr_{\CK}(\CA) = \Fun(\CK^{\op}, \CA)$ admits an $\Eb_n$-monoidal structure by Day convolution. 
            \item The tensor unit $\oneb \in \Gr_{\CK}(\CA)$ admits two quotients $\oneb_{\CA}^{\triv}[0]$ and $\oneb_{\CA}^{\sgn}[0]$, respectively given by the coequalizer of $f, \id\colon \oneb  \rightrightarrows \oneb$ and $f, -\id\colon \oneb  \rightrightarrows \oneb$, where $f$ comes from the $C_2$-action on $\oneb$ induced by the functor $\Pe \to \CK$. Both of them are idempotent algebras, so their module categories $\Gr_{\CK}^{\triv}(\CA)$, $\Gr_{\CK}^{\sgn}(\CA)$ are also $\Eb_n$-monoidal; they are respectively the full subcategory of $\Gr_{\CK}(\CA)$ of objects $\{M^V\}_{V \in \CK}$ so that the $C_2$-action on each $M^V$ is trivial/sign. In the case $\CK = \Pe$, we have $\Gr_{\Pe}^{\triv}(\CA) = \Gr(\CA)^{\Day}$ and $\Gr_{\Pe}^{\sgn}(\CA) = \Gr(\CA)^{\Kos}$. 
            \item Furthermore, the tensor unit $\oneb$ in $\Gr_{\CK}(\Gr_{\Pe}(\CA)) = \Gr_{\CK \times \Pe}(\CA)$ has four quotients $\oneb_{\CA}^{\epsilon_1, \epsilon_2}[0]$ with $\epsilon_1, \epsilon_2 \in \{\triv, \sgn\}$, which correspond to $f_1 = (f, e), f_2 = (e, f) \in C_2 \times C_2$ acting by $\id$ or $-\id$. These are also idempotent algebras, so their module categories are $\Eb_n$-monoidal; more precisely, $\Mod_{\oneb_{\CA}^{\epsilon_1, \epsilon_2}[0]}$ is identified with $\Gr_{\CK}^{\epsilon_1}(\Gr_{\Pe}^{\epsilon_2}(\CA))$. 
            \item Suppose the category $\CK$ and the functor $\Pe \to \CK$ are $\Eb_{n + 1}$-monoidal. Then the functor $+\colon \CK \times \CK \to \CK$ is $\Eb_n$-monoidal, and 
            \[\begin{tikzcd}
                \Pe \times \Pe \ar[r, "{+}"] \ar[d] & \Pe \ar[d] \\
                \CK \times \CK \ar[r, "{+}"] & \CK
            \end{tikzcd}\]
            is a commutative square in $\Alg_{\Eb_n}(\Cat)$. Thus, the reindex map 
            \[v\colon \CK \times \Pe \to \CK \times \Pe, \qquad (V, a) \mapsto (V - a, a)\] 
            admits the structure of an $\Eb_n$-monoidal functor, which (by precomposition) leads to an $\Eb_n$-monoidal equivalence $v^*\colon \Gr_{\CK \times \Pe}(\CA) \cong \Gr_{\CK \times \Pe}(\CA)$. Note that under $v^*$ the maps $f_1, f_2 \in \Aut(\oneb)$ are sent to $f_1$ and $f_1 f_2$, so the idempotent algebra $\oneb_{\CA}^{\sgn, \triv}[0]$ goes to $\oneb_{\CA}^{\sgn, \sgn}[0]$. Therefore, reparametrization along $v$ gives rise to an $\Eb_n$-monoidal equivalence $v^*\colon \Gr_{\CK}^{\sgn}(\Gr(\CA)^{\Day}) \cong \Gr_{\CK}^{\sgn}(\Gr(\CA)^{\Kos})$.
            \item Suppose moreover there is an $\Eb_{n + 1}$-monoidal functor $|\cdot|\colon \CK \to \Pe$ so that the composite $\Pe \to \CK \to \Pe$ is equivalent to $\id$ as an $\Eb_{n + 1}$-monoidal functor. Then for each $r \in \Nb_{\geq 1}$, the reindex map 
            \[v_r\colon \Pe \times \CK \to \Pe \times \CK, \qquad (s, V) \mapsto (s + (r - 1)(|V| - s), V + (r - 1)(|V| - s))\] 
            admits the structure of an $\Eb_n$-monoidal functor,  which (by precomposition) leads to an $\Eb_n$-monoidal equivalence $v_r^*\colon \Gr_{\CK \times \Pe}(\CA) \cong \Gr_{\CK \times \Pe}(\CA)$. Moreover, write $a \equiv r \text{ mod } 2$, then under this functor the maps $f_1, f_2 \in \Aut(\oneb)$ are sent to $f_1^{a}f_2^{1 - a}$ and $f_1^{1 - a}f_2^{a}$, so the idempotent algebra $\oneb_{\CA}^{\sgn, \sgn}[0]$ goes to $\oneb_{\CA}^{\sgn, \sgn}[0]$. Therefore, reparametrization along $v_r$ gives rise to an $\Eb_n$-monoidal equivalence $v_r^*\colon \Gr_{\CK}^{\sgn}(\Gr(\CA)^{\Kos}) \cong \Gr_{\CK}^{\sgn}(\Gr(\CA)^{\Kos})$. 
        \end{itemize}
    \end{itemize}
\end{construction}

\begin{example} \label{Antieau-lemma-with-picard-grading}
    Take $n \in \Nb_{\geq 1} \cup \{\infty\}$ and let $\CE$ be a presentably $\Eb_{n + 1}$-monoidal $t$-$\infty$-category. Write $\CI = \Pic(\CE)$ and $\CK = h \CI$, then $\CK$ receives a canonical $\Eb_{n + 1}$-monoidal functor $S\colon \Pe \to \CK$. Recall from Construction \ref{coarse-picard-grading-constructions-1} the lax $\Eb_n$-monoidal functor $\pi_\filledstar\colon \CE \to \Gr_{\CK}(\CE^\heartsuit)$ taking $X$ to its $\CK$-graded homotopy objects $\{\pi_V(X)\}_{V \in \CK}$, which is derived from transposing the lax $\Eb_n$-monoidal composite 
    \[\CE \times \CI^{\op} \to \CE \times \CE^{\op} \xrightarrow{\hom_{\CE}(-, -)} \CE \xrightarrow{\pi_0} \CE^\heartsuit.\]
    For each $X \in \CE, V \in \CK$, the $C_2$-action on $\pi_V(X) = \pi_0(\hom_{\CE}(\oneb_{\CE} \otimes \Sb^V, X))$ comes from the action of $C_2 = \pi_0(\Sb)^\times = \Aut_{\Sp}(\Sb^0)$ on $\oneb_{\CE} = S(\Sb^0)$, which is multiplication by $-1$. Therefore, $\pi_\filledstar$ factors canonically through $\Gr_{\CK}^{\sgn}(\CE^\heartsuit)$. In particular, after restriction along $S$, we obtain a lax $\Eb_n$-monoidal functor $\pi_*\colon \CE \to \Gr_{\Pe}^{\sgn}(\CE^\heartsuit) = \Gr(\CE^\heartsuit)^{\Kos}$. The last sentence is slightly weaker than Proposition \ref{Antieau-lemma} which says $\pi_*$ is lax $\Eb_{n + 1}$-monoidal, but this perspective is more feasible for generalizations.
\end{example}

% \begin{remark}
%     In the context of Example \ref{Antieau-lemma-with-picard-grading}, there is a group homomorphism from $\pi_1 \Pic(\CE) = \Aut(\oneb_{\CE})$ to $\Aut(\oneb_{\CE^\heartsuit})$. Therefore, we can construct another quotient $\oneb^{\red}_{\CE^\heartsuit}[0]$ of the tensor unit $\oneb \in \Gr_{\CK}(\CE^\heartsuit)$ by taking the coequalizer of two maps $f_1, f_2\colon \bigoplus_{\pi_1 \Pic(\CE)} \oneb \rightrightarrows \oneb$ induced by $\pi_1 \Pic(\CE)$-actions on $\oneb$ respectively from the indexing category $\CK$ and from the coefficient category $\CE^\heartsuit$. This $\oneb^{\red}_{\CE^\heartsuit}[0]$ is still an idempotent algebra, and its module category $\Gr_{\CK}^{\red}(\CE^\heartsuit)$ is the full subcategory of $\Gr_{\CK}(\CE^\heartsuit)$ spanned by $\CK$-graded objects $\{M^V\}_{V \in \CK}$ such that on each $M^V \in \CE^\heartsuit$ the two $\Aut(V) \cong \Aut(\oneb_{\CE})$ actions (from the index and the coefficient) agree. In fact, the discussion in Example \ref{Antieau-lemma-with-picard-grading} shows that $\pi_\filledstar\colon \CE \to \Gr_{\CK}^{\sgn}(\CE^\heartsuit)$ factors through the more refined target $\Gr_{\CK}^{\red}(\CE^\heartsuit)$, which has the advantage that as a monoidal category $\Gr_{\CK}^{\red}(\CE^\heartsuit)$ is equivalent to $\Gr_{\pi_0\Pic(\CE)}(\CE^\heartsuit)$ with Day convolution. 
% \end{remark}

\begin{setup}[Refined grading data] \label{picard-grading-refined-setup}
    Take $n \in \Nb_{\geq 1} \cup \{\infty\}$. We fix the following data:
    \begin{itemize} 
        \item Let $\CE$ be a stable presentably $\Eb_{n + 1}$-monoidal $\infty$-category. 
        \item Let $\CF$ be a presentably $\Eb_{n + 1}$-monoidal $t$-$\infty$-category, and let $T\colon \CE \to \CF$ be an exact lax $\Eb_{n + 1}$-monoidal functor. 
        % \item Let $\CE^\heartsuit$ be an $\Eb_{n}$-monoidal abelian $1$-category, with a lax $\Eb_{n}$-monoidal functor $U\colon \CF^\heartsuit \to \CE^\heartsuit$.
        \item Let $\CI$ be an $\Eb_{n + 1}$-monoidal $\infty$-category, equipped with three $\Eb_{n + 1}$-monoidal functors 
        \[O\colon \Pic(\Sb) \to \CI, \qquad P\colon \CI \to \CE \qquad \text{and} \qquad Q\colon \CI \to \Pe\]
        so that $PO\colon \Pic(\Sb) \to \CE$ is equivalent to the canonical $\Eb_{n + 1}$-monoidal functor $\Pic(\Sb) \hookrightarrow \Sp \to \CE$, and $QO\colon \Pic(\Sb) \to \Pe$ is equivalent to the truncation functor $\Pic(\Sb) \to h\Pic(\Sb)$. We also assume for each $V \in \CI$, the image $P(V) \in \CE$ is dualizable.
        \item Let $\CK$ be an $\Eb_{n + 1}$-monoidal $1$-category equipped with two $\Eb_{n + 1}$-monoidal functors 
        \[R\colon \CK \to h\CI \qquad \text{and} \qquad S\colon \Pe \to \CK\]
        so that $RS\colon \Pe \to h\CI$ is equivalent to the $1$-truncation of $O\colon \Pic(\Sb) \to \CI$. We refer to the two composites $PR\colon \CK \to h \CE$, $QR\colon \CK \to \Pe$ respectively as $V \mapsto \Sb^V$ and $V \mapsto |V|$.
    \end{itemize}
    In practice, usually we have $\CE = \CF$ and $\CK = h\CI$, and we will refer to such a full setup by saying $\CK \to h\CE$ is part of a \textbf{(refined) grading datum}.
\end{setup}

\begin{example} \label{refined-picard-grading-examples}
    Take $n \in \Nb_{\geq 1} \cup \{\infty\}$.
    \begin{itemize}
        \item Let $\CE$ be a presentably $\Eb_{n + 1}$-monoidal $t$-$\infty$-category, equipped with an $\Eb_{n + 1}$-monoidal retract $Q$ of the canonical map $\Pic(\Sb) \to \Pic(\CE)$. Take $\CF = \CE$, $\CI = \Pic(\CE)$ and $\CK = h\CI$. Set $O\colon \Pic(\Sb) \to \CI$ and $S\colon \Pe \to \CK$ to come from the unique $\Eb_{n + 1}$-monoidal left adjoint $\Sp \to \CE$, $P\colon \CI\to \CE$ to be the inclusion, and $T = \id$, then the collection is legitimate for Setup \ref{picard-grading-refined-setup}. This recipe refines Example \ref{coarse-picard-grading-examples} item 1.
        \item Let $G$ be a finite group. Take $\CE = \Sp^G = \CF$, the symmetric monoidal $\infty$-category of genuine $G$-spectra with its Mackey $t$-structure. Furthermore, take $\CI = \Pic(\Sp^G)$ and take $\CK = \PG$ to be the pullback of the cospan $\RO(G) \to \pi_0 \Pic(\Sp^G) \gets h \Pic(\Sp^G)$. Setting up $O,P,Q, R, T$ as above (the section $Q$ comes from the forgetful functor $U\colon \Sp^G \to \Sp$), and taking $S\colon \Pe \to \CK$ through the standard refinement $\Zb \to \RO(G)$ of $\Zb \to \pi_0\Pic(\Sp^G)$, we obtain a collection legitimate for Setup \ref{picard-grading-refined-setup}. This recipe refines Example \ref{coarse-picard-grading-examples} item 2.
    \end{itemize}
\end{example}

\begin{construction}
    Under Setup \ref{picard-grading-refined-setup}, we have the following constructions:
    \begin{itemize}
        \item As refined grading data restricts to coarse grading data (in the sense of Setup \ref{picard-grading-coarse-setup}), all discussions in Construction \ref{coarse-picard-grading-constructions-1} and Construction \ref{coarse-picard-grading-constructions-2} make sense in the current setup. In particular, we have a lax $\Eb_n$-monoidal functor $\pi_\filledstar\colon \CE \to \Gr_{\CK}(\CF^\heartsuit)$ sending $X$ to the collection of $\CK$-graded homotopy objects $\{\pi_V(X)\}_{V \in \CK}$. Furthermore, the discussion in Example \ref{Antieau-lemma-with-picard-grading} shows that such a functor factors canonically through $\Gr^{\sgn}_{\CK}(\CF^\heartsuit)$.
        \item As with $\Fil(\Sp)$ and $\Gr(\Sp)$, the stable $\infty$-categories $\Fil(\CF)$ and $\Gr(\CF)$ also admit a pointwise $t$-structure and a Beilinson $t$-structure. On the other hand, there is a lax $\Eb_n$-monoidal functor $\psi_0\colon \Fil(\CE) \to \Gr_{\CI}\Fil(\CE)$ extracted from transposing the composite 
        \[\Fil(\CE) \times \CI^{\op} \xrightarrow{\id \times P} \Fil(\CE) \times \CE^{\op} \xrightarrow{\id \times j} \Fil(\CE) \times \Fil(\CE)^{\op} \xrightarrow{\hom_{\Fil(\CE)}(-,-)} \Fil(\CE).\]
        Combining these, we get two lax $\Eb_n$-monoidal functors: the functor
        \begin{align*}
            \pi_{\filledstar,*}\colon \Fil(\CE) &\xrightarrow{\psi_0} \Gr_{\CI}\Fil(\CE) \xrightarrow{\Gr_{\CI}(T)} \Gr_{\CI}\Fil(\CF) \\
            &\xrightarrow{\mathrm{forget}} \Gr_{\CI}\Gr(\CF) \xrightarrow{\Gr_{\CI}(\pi^P_0)} \Gr_{h\CI}(\Gr(\CF^\heartsuit)^{\Day}) \to \Gr_{\CK}(\Gr(\CF^\heartsuit)^{\Day})
        \end{align*}
        whose target lies in $\Gr_{\CK}^{\sgn}(\Gr(\CF^\heartsuit)^{\Day})$ via the discussion in Example \ref{Antieau-lemma-with-picard-grading}, and the functor 
        \begin{align*}
            E_2 \colon \Fil(\CE) \xrightarrow{\psi_0} \Gr_{\CI}\Fil(\CE) &\xrightarrow{\Gr_{\CI}(T)} \Gr_{\CI}\Fil(\CF) \\ 
            &\xrightarrow{\Gr_{\CI}(\pi^B_0)} \Gr_{h\CI}(\Ch(\CF^\heartsuit)^{\Kos}) \to \Gr_{\CK}(\Ch(\CF^\heartsuit)^{\Kos})
        \end{align*}
        whose target lies in $\Gr_{\CK}^{\sgn}(\Ch(\CF^\heartsuit)^{\Kos})$. 
        \item For $r \in \Nb_{\geq 1}$, we denote by  $\Gr_{\CK}^{E_r}\Ch(\CF^\heartsuit)$ the $1$-category of pairs $(\{M^{s, V}\}_{s \in \Zb, V \in \CK}, d_r)$, where 
        \begin{itemize}
            \item the first component $\{M^{s, V}\}_{s \in \Zb, V \in \CK}$ is a $\Zb \times \CK$ graded object in $\CF$, on which the sign $C_2$-action on each $V \in \CK$ (coming from $S\colon \Pe \to \CK$) induces the sign action on $M^{s, V}$. 
            \item the second component is a family of maps $d_r\colon M^{s, V} \to M^{s + r, V + r - 1}$ in $\CF^\heartsuit$, so that $d_r d_r = 0$, and for each $f\colon W \to V$ in $\CK$ the induced square 
            \[\begin{tikzcd}
                M^{s,V} \ar[r, "{d_r}"] \ar[d,"{f^*}"] & M^{s + r, V + r - 1} \ar[d,"{f^*}"] \\
                M^{s, W} \ar[r, "{d_r}"] & M^{s + r, W + r - 1}
            \end{tikzcd}\]
            commutes in $\CF^\heartsuit$.
        \end{itemize}
        Concretely, for $r = 1$ we have a canonical identification
        \[\Gr_{\CK}^{E_1}\Ch(\CF^\heartsuit) \to \Gr_{\CK}^{\sgn}(\Ch(\CF^\heartsuit)^{\Kos}), \quad (\{M^{s, V}\}_{s \in \Zb, V \in \CK}, d_1) \mapsto \{(M^{*, V}, d_1)\}_{V \in \CK}\]
        and we equip $\Gr_{\CK}^{E_1}\Ch(\CF^\heartsuit)$ with the $\Eb_{n}$-monoidal structure pulled back from the target. In general, we define $\Gr_{\CK}^{E_r}\Ch(\CF^\heartsuit)$ as the pullback of the cospan of $\Eb_n$-monoidal functors
        \[\Gr_{\CK}^{\sgn}(\Ch(\CF^\heartsuit)^{\Kos}) \xrightarrow{U} \Gr_{\CK}^{\sgn}(\Gr(\CF^\heartsuit)^{\Kos}) \xleftarrow{v_r^*} \Gr_{\CK}^{\sgn}(\Gr(\CF^\heartsuit)^{\Kos}),\]
        where $v_r^*\colon \Gr_{\CK}^{\sgn}(\Gr(\CF^\heartsuit)^{\Kos}) \cong \Gr_{\CK}^{\sgn}(\Gr(\CF^\heartsuit)^{\Kos})$ is the reparametrization equivalence from Construction \ref{Koszul-grading-as-a-Day-idempotent}. By construction, the following holds true: 
        \begin{itemize}
            \item The functor $E_2 = \pi^B_\filledstar\colon \Fil(\CE) \to \Gr_{\CK}^{\sgn}(\Ch(\CF^\heartsuit)^{\Kos})$ is isomorphic to the composite of two lax $\Eb_n$-monoidal functors $\Fil(\CE) \to \Gr_{\CK}^{E_2}\Ch(\CF^\heartsuit), X \mapsto (\{E_2^{s, V}(X)\}_{s \in \Zb, V \in \CK}, d_2)$ and $v_2^*\colon \Gr_{\CK}^{E_2}\Ch(\CF^\heartsuit) \to \Gr_{\CK}^{\sgn}(\Ch(\CF^\heartsuit)^{\Kos})$. 
            \item The forgetful functor 
            \[U\colon \Gr_{\CK}^{E_r}\Ch(\CF^\heartsuit) \to \Gr_{\CK}^{\sgn}(\Gr(\CF^\heartsuit)^{\Kos}), \quad (\{M^{s, V}\}_{s\in \Zb, V \in \CK}, d_{r}) \mapsto \{M^{s, V}\}_{s\in \Zb, V \in \CK}\]
            is faithful and $\Eb_n$-monoidal. Also, taking cycles and taking homology yield lax $\Eb_n$-monoidal functors $Z, H\colon \Gr_{\CK}^{E_r}\Ch(\CF^\heartsuit) \to  \Gr_{\CK}^{\sgn}(\Gr(\CF^\heartsuit)^{\Kos})$. 
        \end{itemize}
        \item We write $\SpSeq_{\CK}(\CF^\heartsuit)$ for the limit of the diagram of lax $\Eb_n$-monoidal functors
        \[\begin{tikzcd}[column sep = tiny, row sep = small]
            \Gr_{\CK}^{E_2}\Ch(\CF^\heartsuit) \ar[rd, "{H}"] & & \Gr_{\CK}^{E_3}\Ch(\CF^\heartsuit) \ar[rd, "{H}"] \ar[ld, "{U}"'] & & \cdots \ar[ld, "{U}"'] \\
            & \Gr_{\CK}^{\sgn}(\Gr(\CF^\heartsuit)^{\Kos}) & & \Gr_{\CK}^{\sgn}(\Gr(\CF^\heartsuit)^{\Kos}) &   
        \end{tikzcd}\]
        in the $\infty$-category of $\infty$-operads over $\Eb_n$. The underlying $1$-category of $\SpSeq_{\CK}(\CF^\heartsuit)$ serves as the $1$-category of $\CK$-graded spectral sequences valued in $\CF^\heartsuit$, and the multimorphisms in this $\infty$-operad correspond to multilinear maps between $\CK$-graded spectral sequences.
    \end{itemize}
\end{construction}

There is an analog of Lemma \ref{products-on-the-E2-page} comparing two multiplications on the $\CK$-graded $E_2$-page.

\begin{lemma} \label{products-on-the-E2-page-with-picard-grading}
    Under Setup \ref{picard-grading-refined-setup}, there is a commutative square of lax $\Eb_n$-monoidal functors
    \[\begin{tikzcd}
        \Fil(\CE) \ar[r, "{-\otimes \oneb / \defopara}"] \ar[d, "{E_2 = \pi^B_\filledstar}"]& \Fil(\CE) \ar[r, "{\pi_{\filledstar, *}}"] &  \Gr_{\CK}^{\sgn}(\Gr(\CF^\heartsuit)^{\Day}) \ar[d, "{v^*}"]\\
        \Gr_{\CK}^{\sgn}(\Ch(\CF^\heartsuit)^{\Kos}) \ar[rr, "{\mathrm{forget}}"] && \Gr_{\CK}^{\sgn}(\Gr(\CF^\heartsuit)^{\Kos})
    \end{tikzcd}\]
    where $v^*$ is the reparametrization equivalence along $v\colon \CK \times \Pe \to \CK \times \Pe, (V, a) \mapsto (V - a, a)$,  the $\Eb_n$-monoidal functor in Construction \ref{Koszul-grading-as-a-Day-idempotent}.
\end{lemma}
\begin{proof}
    It suffices to prove this in the case $\CK = h\CI$. The functor $-\otimes \oneb / \defopara\colon \Fil(\CE) \to \Gr(\CE)$ is $\Eb_n$-monoidal and $t$-exact for both pointwise and Beilinson $t$-structure. Thus, we can rewrite the lower-left composite as 
    \[\Fil(\CE) \xrightarrow{-\otimes \oneb / \defopara} \Gr(\CE) \xrightarrow{\psi_0} \Gr_{\CI} \Gr(\CE) \xrightarrow{T} \Gr_{\CI} \Gr(\CF) \xrightarrow{\pi_\filledstar^B = \Gr_{\CI}(\pi^B_0)}  \Gr_{h\CI} (\Gr(\CF^\heartsuit)^{\Kos}) \]
    where $\psi_0\colon \Gr(\CE) \to \Gr_{\CI}\Gr_{\CE}$ is the transpose of the functor 
    \[\Gr(\CE) \times \CI^{\op} \xrightarrow{\id \times P} \Gr(\CE) \times \CE^{\op} \xrightarrow{\id \times j} \Gr(\CE) \times \Gr(\CE)^{\op} \xrightarrow{\hom_{\Gr(\CE)}(-,-)} \Gr(\CE)\]
    in which $j\colon \CE \to \Gr(\CE)$ is given by left Kan extension along $* = \{0\} \to \CI$. On the other hand, the upper-right composite can also be rewritten as 
    \[\Fil(\CE) \xrightarrow{-\otimes \oneb / \defopara} \Gr(\CE) \xrightarrow{\psi_0} \Gr_{\CI} \Gr(\CE) \xrightarrow{T} \Gr_{\CI} \Gr(\CF) \xrightarrow{\pi_\filledstar^P = \Gr_{\CI}(\pi^P_0)}  \Gr_{h\CI} (\Gr(\CF^\heartsuit)^{\Day}). \]
    Also, recall that $v^*\colon \Gr_{h\CI}^{\sgn}(\Gr(\CF^\heartsuit)^{\Day}) \cong \Gr_{h\CI}^{\sgn}(\Gr(\CF^\heartsuit)^{\Kos})$ is an identification of full subcategories restricted from the equivalence $v^*\colon \Gr_{h\CI}\Gr_{\Pe}(\CF^\heartsuit) \cong \Gr_{h\CI}\Gr_{\Pe}(\CF^\heartsuit)$. Thus, it remains to construct a commutative square of lax $\Eb_n$-monoidal functors 
    \begin{equation*} \label{dagger-square} \tag{$\dagger$}
        \begin{tikzcd}[column sep=small]
            \Gr(\CE) \ar[d, "{\id}"]  \ar[r, "{\psi_0}"] & \Gr_{\CI}\Gr(\CE) \ar[r, "{T}"] & \Gr_{\CI}\Gr(\CF) \ar[r, "{\pi^P_{\filledstar}}"] & \Gr_{h\CI}({\Gr(\CF^\heartsuit)}^{\Day}) \ar[r] & \Gr_{h\CI}{\Gr_{\Pe}(\CF^\heartsuit)} \ar[d, "{v^*}"] \\
            \Gr(\CE) \ar[r, "{\psi_0}"] & \Gr_{\CI}\Gr(\CE) \ar[r, "{T}"] &  \Gr_{\CI}\Gr(\CF) \ar[r, "{\pi^B_{\filledstar}}"] & \Gr_{h\CI}({\Gr(\CF^\heartsuit)}^{\Kos}) \ar[r] & \Gr_{h\CI}{\Gr_{\Pe}(\CF^\heartsuit)}
        \end{tikzcd}
    \end{equation*}
    To further simplify, we introduce the following notion: for any stable presentably $\Eb_k$-monoidal category $\CC$ we write $\Gr_{\Pic(\Sb)}(\CC)$ for the functor category $\Fun(\Pic(\Sb)^{\op}, \CC)$ equipped with the Day convolution $\Eb_{k}$-monoidal structure. Note that $\Gr_{\Pic(\Sb)}(\CF)$ also admits two $t$-structures: 
    \begin{itemize}
        \item The \textbf{pointwise $t$-structure}, whose connective part $(\Gr_{\Pic(\Sb)}\CF)^{P}_{\geq 0}$ (resp. whose coconnective part $(\Gr_{\Pic(\Sb)}\CF)^{P}_{\leq 0}$) consists of $X \in \Gr_{\Pic(\Sb)}\CF$ such that for each $w \in \Zb^\delta = \pi_0(\Pic(\Sb))$, $X(w)$ is $0$-connective (resp. $0$-coconnective). As for the heart, we have $(\Gr_{\Pic(\Sb)}\CF)^{P,\heartsuit} \cong \Gr_{\Pe}(\CF^\heartsuit)$ as an $\Eb_{n + 1}$-monoidal $1$-category. 
        \item The \textbf{Beilinson $t$-structure}, whose connective part $(\Gr_{\Pic(\Sb)}\CF)^{B}_{\geq 0}$ (resp. whose coconnective part $(\Gr_{\Pic(\Sb)}\CF)^{B}_{\leq 0}$) consists of $X \in \Gr_{\Pic(\Sb)}\CF$ such that for each $w \in \Zb^\delta = \pi_0(\Pic(\Sb))$, $X(-w)$ is $w$-connective (resp. $w$-coconnective). As for the heart, we have $(\Gr_{\Pic(\Sb)}\CF)^{B,\heartsuit} \cong \Gr_{\Pe}(\CF^\heartsuit)$ as an $\Eb_{n + 1}$-monoidal $1$-category. 
    \end{itemize}
    The functor $p^*\colon \Gr(\CF) \to \Gr_{\Sb}(\CF)$ induced from precomposing the map $p\colon \Pic(\Sb) \to \pi_0\Pic(\Sb) = \Zb^\delta$ is both lax $\Eb_{n + 1}$-monoidal and $t$-exact for both $t$-structures. Furthermore, write $\psi_0\colon \Gr_{\Pic(\Sb)} \CE \to \Gr_{\CI}\Gr_{\Pic(\Sb)} \CE $ for the transpose of the functor 
    \[\Gr_{\Pic(\Sb)} (\CE) \times \CI^{\op} \xrightarrow{\id \times P} \Gr_{\Pic(\Sb)} (\CE) \times \CE^{\op} \xrightarrow{\id \times j} \Gr_{\Pic(\Sb)} (\CE) \times \Gr_{\Pic(\Sb)} (\CE)^{\op} \xrightarrow{\hom_{\Gr_{\Pic(\Sb)} (\CE)}} \Gr_{\Pic(\Sb)} (\CE)\]
    then there exists a commutative square of lax $\Eb_{n}$-monoidal functors
    \[\begin{tikzcd}
        \Gr(\CE) \ar[r, "{\psi_0}"] \ar[d, "{p^*}"] & \Gr_{\CI} \Gr(\CE) \ar[d, "{p^*}"] \\
        \Gr_{\Pic(\Sb)} (\CE) \ar[r, "{\psi_0}"] & \Gr_{\CI}\Gr_{\Pic(\Sb)} \CE
    \end{tikzcd}\]
    For this, the only thing to check is $\hom(\Sb^{V, 0}, p^*(X)) = p^*\hom(\Sb^{V, 0}, X)$ for each $V \in h \CI = \CK$ and $X \in \Gr(\CE)$, which holds true by projection formula \cite[Proposition A.4(1)]{BHS2} since the left adjoint of $p^*$ is $\Eb_{n + 1}$-monoidal and each $\Sb^V \in \CE$ is dualizable. Under this notation, we can rewrite the upper horizontal arrow in (\ref{dagger-square}) as the composite 
    \[\Gr(\CE) \xrightarrow{p^*} \Gr_{\Pic(\Sb)}(\CE) \xrightarrow{\psi_0} \Gr_{\CI} \Gr_{\Pic(\Sb)}(\CE) \xrightarrow{T} \Gr_{\CI} \Gr_{\Pic(\Sb)}(\CF) \xrightarrow{\pi^P_\filledstar = \Gr_{\CI}(\pi^P_0)} \Gr_{\CI} \Gr_{\Pic(\Sb)}(\CF^\heartsuit)  \]
    and rewrite the lower horizontal arrow in (\ref{dagger-square}) as the composite 
    \[\Gr(\CE) \xrightarrow{p^*} \Gr_{\Pic(\Sb)}(\CE) \xrightarrow{\psi_0} \Gr_{\CI} \Gr_{\Pic(\Sb)}(\CE) \xrightarrow{T} \Gr_{\CI} \Gr_{\Pic(\Sb)}(\CF) \xrightarrow{\pi^B_\filledstar = \Gr_{\CI}(\pi^B_0)} \Gr_{\CI} \Gr_{\Pic(\Sb)}(\CF^\heartsuit).  \]
    Thus, it remains to produce a commutative square of lax $\Eb_n$-monoidal functors
    \begin{equation*} \label{double-dagger-square} \tag{$\ddagger$}
        \begin{tikzcd}
            \Gr_{\Pic(\Sb)}(\CE) \ar[r, "{\psi_0}"] \ar[d, "{\id}"] &  \Gr_{\CI} \Gr_{\Pic(\Sb)}(\CE) \ar[r, "{T}"] & \Gr_{\CI} \Gr_{\Pic(\Sb)}(\CF) \ar[r, "{\pi^P_\filledstar}"]  & \Gr_{\CI} \Gr_{\Pic(\Sb)}(\CF^\heartsuit)  \ar[d, "{v^*}"]\\
            \Gr_{\Pic(\Sb)}(\CE) \ar[r, "{\psi_0}"] & \Gr_{\CI} \Gr_{\Pic(\Sb)}(\CE) \ar[r, "{T}"] & \Gr_{\CI} \Gr_{\Pic(\Sb)}(\CF) \ar[r, "{\pi^B_\filledstar}"] & \Gr_{\CI} \Gr_{\Pic(\Sb)}(\CF^\heartsuit)
        \end{tikzcd}
    \end{equation*}
    to achieve which we use a twisting automorphism. Recall that $\id\colon \Gr_{\Pic(\Sb)}(\Sp) \to \Gr_{\Pic(\Sb)}(\Sp)$ is the transpose of the lax symmetric monoidal functor 
    \[\ev\colon \Gr_{\Pic(\Sb)}(\Sp) \times \Pic(\Sb) \to \Sp, \quad (X, a) \mapsto X(a).\]
    We can construct another symmetric monoidal left adjoint $V\colon \Gr_{\Pic(\Sb)}(\Sp) \to \Gr_{\Pic(\Sb)}(\Sp), X \mapsto \{\Sigma^{-w}X(w)\}_{w \in \Zb}$ by transposing the lax symmetric monoidal functor 
    \[\Gr_{\Pic(\Sb)}(\Sp) \times \Pic(\Sb) \xrightarrow{\id \times \Delta}  \Gr_{\Pic(\Sb)}(\Sp) \times \Pic(\Sb) \times \Pic(\Sb) \xrightarrow{\ev \times \id} \Sp \times \Pic(\Sb) \to \Sp \times \Sp^{\op} \xrightarrow{\Msp(-,-)}  \Sp.\]
    Its inverse $V^{-1}\colon \Gr_{\Pic(\Sb)}(\Sp) \to \Gr_{\Pic(\Sb)}(\Sp), X \mapsto \{\Sigma^{w}X(w)\}_{w \in \Zb}$ can be constructed in a similar fashion. It follows that for any stable presentably $\Eb_{k}$-monoidal $\infty$-category $\CC$, we get a pair of $\Eb_{k}$-monoidal left adjoints $V, V^{-1}\colon \Gr_{\Pic(\Sb)}(\CC) \to \Gr_{\Pic(\Sb)}(\CC)$ by tensoring $V, V^{-1}$ with $\CC$ over $\Sp$. We will subsequently construct two commutative squares of lax $\Eb_n$-monoidal functors
    \begin{equation*} \label{diamond-squares} \tag{$\lozenge$}
        \begin{tikzcd}
            \Gr_{\Pic(\Sb)}(\CE) \ar[d, "{V}"] \ar[rr, "{\pi_\filledstar^P = \pi_\filledstar^P T\psi_0}"] &&  \Gr_{\CI} \Gr_{\Pic(\Sb)}(\CF^\heartsuit) \ar[d, "{\id}"]\\
            \Gr_{\Pic(\Sb)}(\CE) \ar[rr, "{\pi_\filledstar^B = \pi_\filledstar^B T\psi_0}"] && \Gr_{\CI} \Gr_{\Pic(\Sb)}(\CF^\heartsuit)
        \end{tikzcd} \quad \text{and} \quad 
        \begin{tikzcd}
            \Gr_{\Pic(\Sb)}(\CE) \ar[d, "{V^{-1}}"] \ar[r, "{\pi_\filledstar^P}"] &  \Gr_{\CI} \Gr_{\Pic(\Sb)}(\CF^\heartsuit)  \ar[d, "{v^*}"]\\
            \Gr_{\Pic(\Sb)}(\CE) \ar[r, "{\pi_\filledstar^P}"] & \Gr_{\CI} \Gr_{\Pic(\Sb)}(\CF^\heartsuit)
        \end{tikzcd}
    \end{equation*}
    which paste together to form the desired commutative square (\ref{double-dagger-square}). The left square in (\ref{diamond-squares}) is 
    \[\begin{tikzcd}
        \Gr_{\Pic(\Sb)}(\CE) \ar[r, "{\psi_0}"] \ar[d, "{V}"] &  \Gr_{\CI} \Gr_{\Pic(\Sb)}(\CE) \ar[r, "{T}"] \ar[d, "{V}"] & \Gr_{\CI} \Gr_{\Pic(\Sb)}(\CF) \ar[r, "{\pi^P_\filledstar}"] \ar[d, "{V}"] & \Gr_{\CI} \Gr_{\Pic(\Sb)}(\CF^\heartsuit)  \ar[d, "{\id}"]\\
        \Gr_{\Pic(\Sb)}(\CE) \ar[r, "{\psi_0}"] & \Gr_{\CI} \Gr_{\Pic(\Sb)}(\CE) \ar[r, "{T}"] & \Gr_{\CI} \Gr_{\Pic(\Sb)}(\CF) \ar[r, "{\pi^B_\filledstar}"] & \Gr_{\CI} \Gr_{\Pic(\Sb)}(\CF^\heartsuit)
    \end{tikzcd}\]
    in which the right square commutes since $V\colon \Gr_{\CI} \Gr_{\Pic(\Sb)}(\CF) \to \Gr_{\CI} \Gr_{\Pic(\Sb)}(\CF)$ is $\Eb_{n}$-monoidal and $t$-exact for pointwise $t$-structure in the source and Beilinson $t$-structure in the target, while the left square commutes since both $V \psi_0$ and $\psi_0 V$ transpose to the composite
    \begin{align*}
        \Gr_{\Pic(\Sb)}(\CE) \times \Pic(\Sb) \times \CI^{\op} 
        &\xrightarrow{\id \times \Delta \times \id} \Gr_{\Pic(\Sb)}(\CE) \times \Pic(\Sb) \times \Pic(\Sb) \times \CI^{\op}  \\
        &\xrightarrow{\ev \times \id \times \id} \CE \times \Pic(\Sb) \times \CI^{\op} \to \CE \times \CE^{\op} \times \CE^{\op} \xrightarrow{\hom_{\CE}(-\otimes -, -)}  \CE.
    \end{align*}
    For the middle square, it suffices to establish the commutativity of
    \[\begin{tikzcd}
        \Gr_{\Pic(\Sb)}(\CE) \ar[r, "{T}"] \ar[d, "{V}"] & \Gr_{\Pic(\Sb)}(\CF) \ar[d, "{V}"] \\
        \Gr_{\Pic(\Sb)}(\CE) \ar[r, "{T}"] & \Gr_{\Pic(\Sb)}(\CF)
    \end{tikzcd}\]
    which comes from the commutativity of the transpose rectangle 
    \[\begin{tikzcd}
        \Gr_{\Pic(\Sb)}(\CE) \times \Pic(\Sb) \ar[r, "{\id \times \Delta}"] \ar[d, "T \times \id"] & \Gr_{\Pic(\Sb)}(\CE) \times \Pic(\Sb) \times \Pic(\Sb) \ar[r, "{\ev \times (-\id)}"] \ar[d, "T \times \id \times \id"] & \CE \times \Pic(\Sb) \ar[r, "{-\otimes-}"]  \ar[d, "T \times \id"] & \CE \ar[d, "T"] \\
        \Gr_{\Pic(\Sb)}(\CF) \times \Pic(\Sb) \ar[r, "{\id \times \Delta}"] & \Gr_{\Pic(\Sb)}(\CF) \times \Pic(\Sb) \times \Pic(\Sb) \ar[r, "{\ev \times (-\id)}"] & \CF \times \Pic(\Sb) \ar[r, "{-\otimes-}"] & \CF
    \end{tikzcd}\]
    Here the crucial step is the commutativity of the rightmost square above in $\Alg_{\Eb_n}(\Cat)^{\lax}$, which follows from that $T\colon \CE \to \CF$ is lax $\Eb_{n + 1}$-monoidal. On the other hand, the right square in (\ref{diamond-squares}) is 
    \[\begin{tikzcd}
        \Gr_{\Pic(\Sb)}(\CE) \ar[r, "{\psi_0}"] \ar[d, "{V^{-1}}"] &  \Gr_{\CI} \Gr_{\Pic(\Sb)}(\CE) \ar[r, "{T}"] \ar[d, "{v^*}"] & \Gr_{\CI} \Gr_{\Pic(\Sb)}(\CF) \ar[r, "{\pi^P_\filledstar}"] \ar[d, "{v^*}"] & \Gr_{\CI} \Gr_{\Pic(\Sb)}(\CF^\heartsuit)  \ar[d, "{v^*}"]\\
        \Gr_{\Pic(\Sb)}(\CE) \ar[r, "{\psi_0}"] & \Gr_{\CI} \Gr_{\Pic(\Sb)}(\CE) \ar[r, "{T}"] & \Gr_{\CI} \Gr_{\Pic(\Sb)}(\CF) \ar[r, "{\pi^P_\filledstar}"] & \Gr_{\CI} \Gr_{\Pic(\Sb)}(\CF^\heartsuit)
    \end{tikzcd}\]
    in which the vertical functors $v^*$ are induced by reparametrization along the $\Eb_{n + 1}$-monoidal functor $v\colon \CI \times \Pic(\Sb) \to \CI \times \Pic(\Sb), (V, a) \mapsto (V - a, a)$. Here the left square commutes since both $v^* \psi_0$ and $\psi_0 V^{-1}$ transpose to the composite 
    \begin{align*}
        \Gr_{\Pic(\Sb)}(\CE) \times \Pic(\Sb) \times \CI^{\op} 
        &\xrightarrow{\id \times \Delta \times \id} \Gr_{\Pic(\Sb)}(\CE) \times \Pic(\Sb) \times \Pic(\Sb) \times \CI^{\op}  \\
        &\xrightarrow{\ev \times \id \times \id} \CE \times \Pic(\Sb) \times \CI^{\op} \to \CE \times \CE \times \CE^{\op} \xrightarrow{\hom_{\CE}(-,- \otimes -)}  \CE
    \end{align*}
    while for the other two squares the commutativity is automatic.
\end{proof}

\begin{theorem} \label{categorified-total-Leibniz-rule-with-picard-bigrading}
    Under Setup \ref{picard-grading-refined-setup}, the following statements hold true: 
    \begin{itemize}
        \item For each $r \in \Nb_{\geq 1}$, the assignment 
        \[A^{*, \filledstar}_{r + 1}\colon \Fil(\CE) \to \Gr_{\CK}^{E_{r + 1}}\Ch(\CF^\heartsuit), \quad X \mapsto \{(A^{s, V}_{r + 1}(X) = \pi_{V - s, |V|}(X / \defopara^r), d_{r + 1} = \delta_r^r)\}\]
        determines a lax $\Eb_n$-monoidal functor. 
        \item Write $E_{r + 1}^{s, V}(X) = E_{r + 1}^{s, 0}(T\hom_{\Fil(\CE)}(\Sb^{V, |V|}, X)) \in \CF^\heartsuit$. For each $r \in \Nb_{\geq 1}$,  the assignment 
        \[E^{*, \filledstar}_{r + 1}\colon  \Fil(\CE) \to \Gr_{\CK}^{E_{r + 1}}\Ch(\CF^\heartsuit), \quad X \mapsto \{(E_{r + 1}^{s, V}(X), d_{r + 1})\}\] 
        determines (uniquely) a lax $\Eb_n$-monoidal functor, so that the epimorphisms 
        \[A^{s, V}_{r + 1}(X) = \pi_{- s, 0}(T\hom_{\Fil(\CE)}(\Sb^{V, |V|}, X) / \defopara^r) \to E_{r + 1}^{s, 0}(T\hom_{\Fil(\CE)}(\Sb^{V, |V|}, X)) = E_{r + 1}^{s, V}(X)\] 
        assemble into an $\Eb_n$-monoidal natural transformation $\alpha\colon A_{r + 1}^{*,\filledstar} \Rightarrow E_{r + 1}^{*,\filledstar}$.  
        \item The functors $\{E_{r + 1}^{*,\filledstar}\}_{r \geq 1}$ induce a lax $\Eb_n$-monoidal functor $E_*^{*,\filledstar}\colon \Fil(\CE) \to \SpSeq_{\CK}(\CF^\heartsuit)$.
    \end{itemize}
\end{theorem}

\begin{proof}
    We first prove the analog of Theorem \ref{Burklund's-Leibniz-rule}, which says for any $r \in \Nb_{\geq 1}$ and any map $F\colon X \otimes Y \to N$ in $\Fil(\CE)$, the induced bilinear pairing (extracted from the lax $\Eb_n$-monoidality of $\pi_{\filledstar, *}$)
    \[F = F / \defopara^r\colon \pi_{V_1, w_1}(X / \defopara^r) \times \pi_{V_2, w_2}(Y/\defopara^r) \to \pi_{V_1 + V_2, w_1 + w_2}(N/\defopara^r)\] 
    satisfies the Leibniz rule for total differentials $\delta_r^r F(\alpha, \beta) = F(\delta_r^r(\alpha), \beta) + (-1)^{|V_1|} F(\alpha, \delta_r^r(\beta))$. For $r = 1$ this follows from Lemma \ref{products-on-the-E2-page-with-picard-grading}. In general, we can define the ($\Eb_{n + 1}$-monoidal left adjoint) slowdown functor $\Sl_r\colon \Fil(\CE) \to \Fil(\CE)$ by tensoring $\Sl_r\colon \Fil\Sp \to \Fil\Sp$ in Construction \ref{accelerations} with $\CE$, and take its (lax $\Eb_{n + 1}$-monoidal) right adjoint, the acceleration $\Ac_r\colon \Fil(\CE) \to \Fil(\CE)$. It follows that the statement for general $r$ follows from the case $r = 1$ through the same discussion on $\{\Ac_r(\Sigma^{0, w}X / \defopara^r)\}_{0 \leq w < r}$ as in the proof of Theorem \ref{Burklund's-Leibniz-rule}. Note that the same proof also shows that, for every map $X_1 \otimes \cdots \otimes X_k \to N$, the induced multilinear pairing on $\pi_{\filledstar, *}(- / \defopara^r)$ also satisfies the Leibniz rule for total differentials. \parr 

    Putting together $\{\pi_\filledstar^B(\Sigma^{0, w}\Ac_r(- / \defopara^r))\}_{0 \leq w < r}$, we obtain a functor that matches the description
    \[A^{*, \filledstar}_{r + 1}\colon \Fil(\CE) \to \Gr_{\CK}^{E_{r + 1}}\Ch(\CF^\heartsuit), \quad X \mapsto \{(A^{s, V}_{r + 1}(X) = \pi_{V - s, |V|}(X / \defopara^r), d_{r + 1} = \delta_r^r)\}.\]
    The further composite $UA^{*, \filledstar}_{r + 1}\colon \Fil(\CE) \to \Gr_{\CK}^{\sgn}(\Gr(\CF^\heartsuit)^{\Kos})$ can be identified with the composite of $\pi_{\filledstar, *}(- / \defopara^r)\colon \Fil(\CE) \to \Gr_{\CK}^{\sgn}(\Gr(\CF^\heartsuit)^{\Day})$, $(v_{r + 1}^*)^{-1}v^*\colon\Gr_{\CK}^{\sgn}(\Gr(\CF^\heartsuit)^{\Day}) \to \Gr_{\CK}^{\sgn}(\Gr(\CF^\heartsuit)^{\Kos})$; hence it admits a lax $\Eb_n$-monoidal structure. As 
    $U\colon \Gr_{\CK}^{E_{r + 1}}\Ch(\CF^\heartsuit) \to \Gr_{\CK}^{\sgn}(\Gr(\CF^\heartsuit)^{\Kos})$, the forgetful functor,  
    is faithful and $\Eb_n$-monoidal, and the image of multimorphism sets in the source is precisely the multilinear maps in the target that satisfy the Leibniz rule, it follows from Corollary \ref{Lifting-lax-sym-mon-functors-along-a-faithful-map} that $A^{*, \filledstar}_{r + 1}$ has a unique $\Eb_n$-monoidal structure lifting that of $UA^{*, \filledstar}_{r + 1}$. We then apply Lemma \ref{epimorphsim-domination-lemma} to the epimorphisms
    \[A^{s, V}_{r + 1}(X) = \pi_{- s, 0}(T\hom_{\Fil(\CE)}(\Sb^{V, |V|}, X) / \defopara^r) \to E_{r + 1}^{s, 0}(T\hom_{\Fil(\CE)}(\Sb^{V, |V|}, X)) = E_{r + 1}^{s, V}(X)\] 
    (whose required hypotheses can be deduced by repeating the second half of the proof of Theorem \ref{Burklund's-Leibniz-rule}) to extract the lax $\Eb_n$-monoidal functor 
    \[E^{*, \filledstar}_{r + 1}\colon  \Fil(\CE) \to \Gr_{\CK}^{E_{r + 1}}\Ch(\CF^\heartsuit), \quad X \mapsto \{(E_{r + 1}^{s, V}(X), d_{r + 1})\}.\]
    The last part in the proof of Theorem \ref{categorified-total-Leibniz-rule} also works without essential change in the $\CK$-graded setup; thus the collection $\{E_{r + 1}^{*,\filledstar}\}_{r \geq 1}$ forms a lax $\Eb_n$-monoidal functor $E_*^{*,\filledstar}\colon \Fil(\CE) \to \SpSeq_{\CK}(\CF^\heartsuit)$.
\end{proof}

\begin{example} \label{bigraded-picard-SS-examples}
    Take $n \in \Nb_{\geq 1} \cup \{\infty\}$.
    \begin{itemize}
        \item Let $\CE$ be a presentably $\Eb_{n + 1}$-monoidal $t$-$\infty$-category, equipped with an $\Eb_{n + 1}$-monoidal retraction of the canonical map $\Pic(\Sb) \to \Pic(\CE)$. Applying Theorem \ref{categorified-total-Leibniz-rule-with-picard-bigrading} to the picard grading recipe in Example \ref{refined-picard-grading-examples} item 1,  we obtain %the assignments $X \mapsto (A^{s, t, V}_{r}(X) = \pi_{V + t - s, |V| + t}(X / \defopara^{r - 1}), \delta_{r - 1}^{r - 1})$ and $X \mapsto (E^{s, t, V}_{r}(X) = E^{s, t}_r(\hom_{\Fil(\CE)}(\Sb^{V, |V|}), X)), d_r)$ assemble into 
        a lax $\Eb_{n}$-monoidal functor 
        \[E_*^{*,\filledstar}\colon \Fil(\CE) \to \SpSeq_{\Pic(\CE)}(\CE^\heartsuit), \quad X \mapsto \{E_r^{s, V}(X) = E_r^{s, 0}(\hom_{\Fil(\CE)}(\Sb^{V, |V|}, X))
        \}%_{r \geq 2, s, t \in \Zb, V \in \Pic(\CE)}
        \] together with its (lax $\Eb_{n}$-monoidal) total differential refinements $A_{r}^{*,\filledstar}$. 
        \item Let $G$ be a finite group. Applying Theorem \ref{categorified-total-Leibniz-rule-with-picard-bigrading} to the $\RO(G)$-grading recipe in Example \ref{refined-picard-grading-examples} item 2,  we obtain %the assignments $X \mapsto (A^{s, t, V}_{r}(X) = \pi_{V + t - s, |V| + t}(X / \defopara^{r - 1}), \delta_{r - 1}^{r - 1})$ and $X \mapsto (E^{s, t, V}_{r}(X) = E^{s, t}_r(\hom_{\Fil(\CE)}(\Sb^{V, |V|}), X)), d_r)$ assemble into 
        a lax symmetric monoidal functor 
        \[E_*^{*,\filledstar}\colon \Fil(\Sp^G) \to \SpSeq_{\PG}(\Mack_G(\Ab)), \quad X \mapsto \{E_r^{s, V}(X) = E_r^{s, 0}(\hom_{\Fil(\Sp^G)}(\Sb^{V, |V|}, X))
        \}%_{r \geq 2, s, t \in \Zb, V \in \RO(G)}
        \] together with its (lax symmetric monoidal) total differential refinements $A_{r}^{*,\filledstar}$.  
        % \item $\Mod_{\mathrm{KU}}(\Sp)\colon E_2^{s, t, V} = \pi_{V + t - s, t}(X / \defopara) = E_2^{s + 1, t, V + 1}$.
    \end{itemize}
\end{example}